\documentclass[10pt]{article}

\input{src/project-preamble}

\usepackage[display-counters=section]{declare-style-notes}
\usepackage{declare-maths}
\usepackage{declare-text}
\usepackage{declare-tables-figures}
\usepackage[numbering-within=section]{declare-theorems-notes}
\usepackage{declare-bibliography}

\hypersetup{
	pdfauthor = {Nikita Nikolaev},
	pdftitle = {Spectral Structures on Vector Bundles},
	pdfsubject = {Vector bundles on curves and their moduli, Complex-analytic moduli problems, Stacks and moduli problems},
	pdfkeywords = {spectral structures, holomorphic vector bundles, spectral curves, spectral data, finite flat covers, spectral correspondence, Higgs bundles, analytic moduli stacks, relative Picard stacks, finite sheaf-module correspondence},
}

\newcommand{\MSCSubjectCode}[1]{\href{https://zbmath.org/classification/?q=cc\%3A#1}{#1}}
\usepackage{tocloft}

\title{\Large\bfseries Spectral Structures on Vector Bundles}
\author{Nikita Nikolaev}
\affil{School of Mathematics, University of Edinburgh, United Kingdom}
\date{2 November 2026}

\begin{document}

\maketitle

\begin{abstract}\noindent
We develop an intrinsic theory of spectral structures on holomorphic vector bundles over compact Riemann surfaces, generalising eigenspace decompositions without eigenvalues.
Their geometric realisations are presentations of vector bundles as pushforwards of line bundles along finite branched covers.
We construct their analytic moduli stack locally of finite presentation and identify it with the relative Picard stack of the universal family of finite covers.
We prove that the forgetful morphism to the stack of vector bundles is representable by analytic spaces and is smooth and surjective on a distinguished smooth open substack.
Consequently, every holomorphic family of vector bundles locally admits spectral presentations along smooth connected simply branched covers, and these presentations can be deformed alongside the underlying bundles.
The relative Picard description yields explicit local dimension formulas and recovers the stack of regular Higgs bundles as base change along the Hitchin base.
We also characterise the additional data supplied to a spectral structure by a regular Higgs field as a global twisted generator, and prove that every vector bundle of rank at least three admits smooth connected simply branched spectral structures which cannot arise from any Higgs field.
\end{abstract}

{\footnotesize

\textbf{2020 MSC:}
\MSCSubjectCode{14H60} (primary),	
\MSCSubjectCode{32G13},				
\MSCSubjectCode{14D23}, 			
\MSCSubjectCode{14H30}, 			
\MSCSubjectCode{32L05}, 			
\MSCSubjectCode{32C25}  			
}

\vspace{-10pt}
\enlargethispage{10pt}

\tableofcontents

\setcounter{section}{-1}
\section{Introduction}

A \textit{spectral structure} generalises the operation of forgetting the eigenvalues of an endomorphism while retaining its eigenspaces.
To see the basic idea, imagine $T$ is an endomorphism of an $n$-dimensional complex vector space $\E$ with distinct eigenvalues.
It decomposes $\E$ into one-dimensional eigenspaces:
\begin{eqn}
	\E = \L_1 \oplus \cdots \oplus \L_n
\fullstop
\end{eqn}
The point is that this decomposition makes sense independently of the endomorphism $T$, and many different endomorphisms give rise to the same decomposition.
Indeed, the algebra $\A \coleq \CC[T] \subset \End(\E)$ is precisely the algebra of all endomorphisms that act diagonally with respect to this decomposition.
Changing $T$ by replacing its eigenvalues with any other collection of distinct scalars yields a new endomorphism $T'$ but leaves this algebra unchanged: $\CC[T] = \CC[T']$.
Thus, the algebra $\A$ and its action on $\E$ retain the decomposition without fixing any eigenvalues for its summands.
Spectral structures make this observation intrinsic and extend it to holomorphic vector bundles, allowing monodromy, ramification, singularities, and nonreduced parameter spaces.

Geometrically, a spectral structure on a vector bundle $\cal{E}$ of rank $n$ over a Riemann surface $\X$ describes it as the pushforward $\cal{E} \cong \pi_\ast \cal{L}$ of a line bundle $\cal{L}$ over an $n$-fold branched cover $\pi : \Y \to \X$.
The sheets of $\Y$ take the place of the summands in local decompositions of $\cal{E}$ into line bundles, and the global topology of $\Y$ records how these lines are permuted as we move around the base $\X$.
A global decomposition of $\cal{E}$ into line bundles is the simplest example that corresponds to a disjoint union of $n$ copies of $\X$.
Allowing connected, branched, or more general covers makes this picture much more flexible.

The need to study such structures first arose for the author in his work on the abelianisation of meromorphic connections \cite{MR3809826,Nikolaev2021Abelianisation}.
In that setting, one encounters local decompositions of a vector bundle into line subbundles which play the role of eigenspaces, but no actual endomorphism defining these decompositions is specified.
This naturally suggested isolating and studying such spectral structures in their own right.

The spectral correspondence for Higgs bundles \cite{Hitchin1987, MR998478, Markman1994} supplies an important source of such structures, and suitable twists yield such presentations for arbitrary vector bundles \cite{HirschowitzNarasimhan1994}.
However, a spectral structure remains meaningful even when it cannot be induced by any global Higgs field.
The theory of \textit{abstract Higgs bundles} \cite{BenZviFrenkel2001, DonagiGaitsgory2002} already separates the underlying spectral structure from the particular Higgs field that induces it.
It captures the proper subclass of spectral structures locally generated by one element.
The present paper develops the intrinsic holomorphic geometry and moduli theory of spectral structures.

\clearpage
\enlargethispage{30pt}
\mbox{}
\vspace{-50pt}

\subsubsection*{Results}

The principal result of this paper is an intrinsic analytic moduli theory of spectral structures, valid for arbitrary finite flat covers and arbitrary holomorphic families over complex analytic parameter spaces including nonreduced ones.
More specifically, we construct the analytic moduli stack of spectral structures and prove that its forgetful morphism to the stack of vector bundles is representable by analytic spaces (\autoref{1792494933366}).
On a distinguished smooth open substack of clean spectral structures, this forgetful morphism is smooth and surjective (\autoref{1792494388516}).
Consequently, every holomorphic family of vector bundles admits, locally on its parameter space, a presentation as the pushforward of a family of line bundles along a smooth connected simply branched cover, and these presentations can be deformed alongside the vector bundles they describe (\autoref{1791544458956}).

There is also a natural forgetful morphism to the moduli stack of all finite covers.
Viewed over this stack of covers, the moduli stack of spectral structures is canonically equivalent to the relative Picard stack of the universal covering curve (\autoref{1792502946092}).
This description recovers the stack of regular twisted Higgs bundles by base change along the Hitchin base (\autoref{1792546482295}) and yields explicit local dimension formulas (\autoref{1792503824011}).

We describe precisely the additional data supplied to a spectral structure by a regular Higgs field (\autoref{1792167464001}), and prove that every vector bundle of rank at least three admits smooth, connected, simply branched spectral structures which cannot be induced by any twisted Higgs field (\autoref{1793629098096}). 
We also exhibit clean spectral structures on the trivial rank-three bundle over $\PP^1$ which have this property (\autoref{1792434634960}).
Further examples reveal the geometry beyond the distinguished smooth open substacks: the moduli stack can have obstructed deformations (\autoref{1792516769027}), singular spectral curves can represent smooth moduli points (\autoref{1793011572262}), and the distinguished open substacks need not be dense (\autoref{1793123354876}).

\subsubsection*{Narrative Outline of the Paper}

We have taken particular care with the exposition, aiming to make the constructions transparent and their geometric meaning accessible.
The main arguments are accompanied by detailed explanations, concrete local models, and 28 numbered explicit examples, alongside further examples and calculations embedded in remarks and proofs.
These range from elementary linear algebra and cyclic covers to singular spectral curves and the geometry of the moduli stack, and are intended to build intuition, explain the role of the hypotheses, and clarify the scope and limitations of the results.

In \ref{1792146089353}, we introduce spectral structures and establish their equivalence with spectral realisations.
We begin in \ref{1790804961989} with the simpler setting of \dfn{linear spectral structures}, in which the vector bundle $\cal{E}$ is replaced by a single complex vector space $\E$ (\autoref{1790804983939}).
Already the algebra $\A = \CC[x,y] / \inner{x,y}^2$, acting on itself, gives a rank-three linear spectral structure which cannot be induced by any endomorphism (\autoref{1790846065686}).
Then in \ref{1790715973901}, we give the central definition of this paper: a \dfn{spectral structure} on a holomorphic vector bundle $\cal{E}$ over a compact Riemann surface $\X$ is a pair $(\cal{A}, \rho)$ consisting of a locally free commutative algebra $\cal{A}$ of the same rank with a homomorphism $\rho : \cal{A} \to \cal{End} (\cal{E})$ that makes $\cal{E}$ a locally free rank-one $\cal{A}$-module (\autoref{1789147824525}).

The geometric meaning of this definition comes from the Finite Sheaf-Module Correspondence (\autoref{1793046131649}) which we explain in \ref{1792143161338}: the relative analytic spectrum of the algebra $\cal{A}$ defines a finite flat cover $\pi : \Y \coleq \Spec_\X (\cal{A}) \to \X$, and the rank-one $\cal{A}$-module $\cal{E}$ defines a line bundle $\cal{L}$ on $\Y$ together with an isomorphism $\iota : \pi_\ast \cal{L} \iso \cal{E}$.
In \ref{1790944507299}, we call this equivalent data a \dfn{spectral realisation} of $(\cal{A}, \rho)$ (\autoref{1793121565784}).
Through the resulting equivalence of groupoids (\autoref{1789153885363}), spectral structures interpolate between the geometry of finite flat covers of $\X$ and the geometry of vector bundles over $\X$, summarised in the following diagram of groupoids:
\vspace{-7.5pt}
\begin{eqn}
\begin{tikzcd}
&	\mathbf{Pic} (\Y)
		\ar[r, hook]
&	\mathbf{SpStr}_n (\X)
		\ar[dl, "G"']
		\ar[dr, "F"]
&	\mathbf{SpStr} (\cal{E})
		\ar[l, hook']
\\
	(\Y \tinyoverset{\pi}{\to} \X)
		\ar[r, hook]
&	\mathbf{Cov}_n (\X)
&&	\mathbf{Bun}_n (\X)
&	\cal{E}
		\ar[l, hook']
\fullstop
\end{tikzcd}
\tag{\ref{1790940189747}}
\end{eqn}

\clearpage
\enlargethispage{10pt}
\mbox{}
\vspace{-30pt}

In a subsequent paper, we will fix a smooth connected cover $\pi : \Y \to \X$ and study the pushforward functor $\pi_\ast : \mathbf{Pic} (\Y) \to \mathbf{Bun}_n (\X)$, including the geometry of its image.

In \ref{1791480675169}, we study \dfn{smooth} and \dfn{connected} spectral structures, meaning their spectral curves $\Y$ are smooth and connected, respectively (\autoref{1789157644037}).
Every vector bundle admits a smooth connected spectral structure whose cover is simply branched (\autoref{1790777298686}). We give two distinct proofs of this claim.
The first proof uses a projective-geometric construction suggested by Anand Patel, in which the fibres of the covering curve are hyperplane sections of a rational normal curve.
The second proof appears in \ref{1791465731919} and uses Higgs bundle methods.

In \ref{1792146117234}, we characterise smoothness by normality of the algebra $\cal{A}$ (\autoref{1789157786307}) and give criteria for recognising spectral structures from actions of the algebra of a smooth cover (\autoref{1791408088098}).
We then show in \ref{1792146177359} that every smooth spectral structure is locally generated by a regular endomorphism, including near branch points (\autoref{1791453151472}), and explain why these local generators need not extend globally.
We conclude that section with several explicit nontrivial examples in \ref{1792146206910} arising from cyclic covers, examining their smoothness, connectedness, and ramification, and how the choice of a line bundle on the cover affects its pushforward.

In \ref{1792146312894}, we introduce the off-diagonal bundle, clean spectral structures, and trace duality.
This begins with the observation that the homomorphism $\rho : \cal{A} \to \cal{End} (\cal{E})$ is automatically injective, giving rise to the short exact sequence
\begin{eqn}
	0 \too \cal{A} \overset{\rho}{\too} \cal{End} (\cal{E}) \too \cal{Q} \too 0
\tag{\ref{1789348440662}}
\fullstop
\end{eqn}
The vector bundle $\cal{Q}$ defined by this sequence is called the \dfn{off-diagonal bundle} (\autoref{1790711090032}).
For smooth spectral structures, it has rank $n(n-1)$ and degree $r/2$ where $r$ is the ramification degree of $\pi$ (\autoref{1790716932644}).
In \ref{1791375997197}, we show explicitly that, away from branch points, $\cal{E}$ is locally decomposed into line bundles, so $\cal{A}$ embeds as the diagonal endomorphisms and $\cal{Q}$ is identified precisely with the bundle of off-diagonal endomorphisms, hence the choice of terminology.
At a branch point $x \in \X$, the fibre algebra $\cal{A}(x)$ becomes nonreduced and its action records the corresponding nilpotent structure (\autoref{1790706598412}).

The cohomology of $\cal{Q}$ has a useful geometric meaning: $\H^1 (\X, \cal{Q}) = 0$ precisely when every first-order deformation of $\cal{E}$ is induced by a first-order deformation of the spectral line bundle $\cal{L}$ on the same cover.
This motivates the definition of \dfn{clean spectral structures} in \ref{1792426526441}, whose covers are smooth, connected, and simply branched, and which satisfy this off-diagonal cohomology vanishing condition (\autoref{1792404625757}).
Every vector bundle admits a clean spectral structure (\autoref{1792167000673}).

In \ref{1791466758188}, we also observe that the canonical trace pairing on the endomorphism bundle $\cal{End} (\cal{E})$ of any vector bundle $\cal{E}$ provides a useful embedding of the dual $\cal{Q}^\vee \inj \cal{End} (\cal{E})$ as the trace-annihilator of the subbundle $\cal{A} \inj \cal{End} (\cal{E})$ (\autoref{1791453955298}).

In \ref{1792146501847}, we describe the precise relationship between spectral structures and Higgs bundles.
Namely, in \ref{1792146546398} we explain that a regular $\Lambda$-twisted Higgs field $\Phi$ induces a spectral structure $\cal{A} \inj \cal{End} (\cal{E})$ together with a global $\Lambda$-twisted generator of $\cal{A}$, or equivalently an embedding of its spectral cover into the total space $\tot (\Lambda)$ (\autoref{1792167000670}).
Combined with the diagram \eqref{1790940189747} above, this relationship is summarised in \ref{1793135088636} by the following diagram:
\begin{eqn}
\begin{tikzcd}
&
&	\mathbf{Bun}_n (\X)
\\
	\mathbf{Higgs}^\rm{reg}_{n, \Lambda} (\X)
		\ar[r, "H"]
		\ar[rru, bend left]
		\ar[rrd, bend right]
&	\mathbf{SpStr}_n (\X)
		\ar[ur, "F"]
		\ar[dr, "G"']
\\
&
&	\mathbf{Cov}_n (\X)
\fullstop
\end{tikzcd}
\tag{\ref{1792401591304}}
\end{eqn}

\clearpage
\enlargethispage{20pt}
\mbox{}
\vspace{-30pt}

An interpretation of this diagram is that, although regular Higgs bundles also interpolate between finite covers and vector bundles, spectral structures provide a more intrinsic intermediary.
Indeed, the groupoid of Higgs bundles requires a choice of a twisting line bundle $\Lambda$, but the groupoid of spectral structures and the two forgetful functors to covers and vector bundles are canonical.

We also explain that the fibre of the forgetful functor $H$ above any spectral structure $(\cal{A}, \rho)$ is canonically identified with the space of all global $\Lambda$-twisted generators of $\cal{A}$ (\autoref{1792167464001}).
However, this fibre may be empty: such a generator need not exist for \textit{any} twist $\Lambda$.
This happens even in the simple setting of clean spectral structures on the trivial rank-three bundle over $\PP^1$ (\autoref{1792434634960}).
Nevertheless, suitable positively twisted Higgs fields yield smooth connected spectral structures (\autoref{1791390497007}) and clean spectral structures (\autoref{1792430412943}) on every vector bundle.
These constructions, carried out in \ref{1791465731919}, provide a distinct, alternative proof of the existence of smooth, connected, simply branched spectral structures (\autoref{1790777298686}) and establish the stronger existence of clean spectral structures (\autoref{1792167000673}).

In \ref{1793634279019}, we return to the projective-geometric construction of \ref{1791480675169} and show that, in rank at least three, it produces spectral structures which cannot be induced by any twisted Higgs field (\autoref{1793629098096}).
We also give a determinant obstruction to the existence of a global twisted generator (\autoref{1792167000671}) and illustrate it with explicit examples, including a clean spectral structure on the trivial rank-three bundle over $\PP^1$ (\autoref{1792434634960}).

Finally, in \ref{1792146586506}, we show that a regular Higgs field supplies a natural $\Lambda$-valued perfect alternating form $\eta : \cal{Q} \otimes \cal{Q} \to \Lambda$ on the off-diagonal bundle (\autoref{1791458921327}).
For classical Higgs fields, this gives balanced off-diagonal cohomology (\autoref{1791459820103}).

In \ref{1791914002217}, we extend these constructions to holomorphic families, allowing the base curve itself to vary and the parameter space to be nonreduced.
We begin in \ref{1793124287763} with establishing the relative properties of the Finite Sheaf-Module Correspondence (\autoref{1793046154727}).
In \ref{1792146709610}, we explain how it ensures that spectral structures, their realisations, and off-diagonal bundles are compatible with arbitrary base change (\autoref{1791551111381}).
In \ref{1792146750001} and \ref{1793045933301}, we explain the fibrewise meaning of smoothness, connectedness, cleanliness, regularity of Higgs bundles, and extend the recognition principle to families (\autoref{1792143470690}).
Finally in \ref{1792146750003}, the main result of this section is that every holomorphic family of vector bundles admits a fibrewise clean spectral structure locally on its parameter space (\autoref{1791544458956}).

\autoref{1792401370132} is the heart of the article.
We begin in \ref{1793143183111} by constructing the relative analytic stack of vector bundles and proving that it is smooth over the parameter space (\autoref{1793142766507}).
The proof uses a general representability result for relative spaces of sections (\autoref{1792526939305}), which also underlies the subsequent moduli constructions.
In \ref{1793143227123}, we organise families of spectral structures into the moduli stack $\frak{SpStr}_n (\X / \S)$ and prove that it is an analytic stack locally of finite presentation (\autoref{1792494933366}).
The two forgetful functors $F$ and $G$ introduced above globalise to morphisms of analytic stacks:
\begin{eqn}
\begin{tikzcd}
&	\frak{SpStr}_n (\X / \S)
		\ar[dl, "G"']
		\ar[dr, "F"]
\\
	\frak{Cov}_n (\X / \S)
&&	\frak{Bun}_n (\X / \S)
\fullstop
\end{tikzcd}
\tag{\ref{1792401198207}}
\end{eqn}

The forgetful morphism $F$ to the analytic stack of vector bundles is representable by analytic spaces and locally of finite presentation (\autoref{1792494933366}), which means spectral structures on a prescribed family of vector bundles are parametrised by an analytic space.
Then in \ref{1793136202900} we prove that smooth, smooth connected, and clean spectral structures form open substacks smooth over $\S$, and that the forgetful morphism $F$ restricted to clean spectral structures is smooth and surjective (\autoref{1792494388516}).
Thus, clean spectral presentations exist throughout the moduli of vector bundles and can be deformed alongside the bundles they present.
Examples show that the full moduli stack can have obstructed deformations (\autoref{1792516769027}) and that the distinguished smooth open substacks need not be dense (\autoref{1793123354876}).
They also show that singular spectral curves can represent smooth points of the moduli stack (\autoref{1793011572262}).

In \ref{1793136253999}, we investigate the other forgetful morphism $G$ to the moduli stack of finite flat covers.
We prove that it is smooth and surjective (\autoref{1792511507140}), and that its fibre over a given cover is precisely the Picard stack of that cover.
Globally, it identifies the moduli stack of spectral structures with the relative Picard stack of the universal covering curve over the stack of finite covers (\autoref{1792502946092}).
The utility of this description of $\frak{SpStr}_n (\X / \S)$ is twofold.

First, in \ref{1793045910309}, it crystallises the concrete relationship with the stack of regular $\Lambda$-twisted Higgs bundles.
Namely, the latter is the fibre product of the stack of spectral structures and the relative Hitchin base $\frak{B}_{n, \Lambda} (\X / \S)$ over the stack of finite covers (\autoref{1792546482295}):
\begin{eqn}
\begin{tikzcd}
	\frak{Higgs}_{n,\Lambda}^\rm{reg} (\X / \S)
		\ar[r, "H"]
		\ar[d, "h"']
		\ar[dr, phantom, "\ulcorner" description, near start]
&
	\frak{SpStr}_n (\X / \S)
		\ar[d, "G"]
\\
	\frak{B}_{n, \Lambda} (\X / \S)
		\ar[r, "i"']
&
	\frak{Cov}_n (\X / \S)
\fullstop
\end{tikzcd}
\tag{\ref{1792546191429}}
\end{eqn}
This yields in particular an alternative proof that $\frak{Higgs}_{n,\Lambda}^\rm{reg} (\X / \S)$ is an analytic stack locally of finite presentation and that the Hitchin morphism $h$ is smooth and surjective.
The forgetful morphism $H$ is representable by analytic spaces and locally of finite presentation.
Although every nonempty fibre of $H$ is smooth, the morphism $H$ itself need not be smooth, even over the locus of clean spectral structures (\autoref{1793046252953}).

Second, in \ref{1793136463825} it leads to explicit local dimension formulas for the moduli stack of spectral structures (\autoref{1792503824011}).
In particular, for a fixed curve $\X$, at a smooth connected spectral realisation $(\pi : \Y \to \X, \cal{L})$ with total ramification degree $r$, we obtain (\autoref{1793102347934}):
\begin{eqn}
    \dim_{(\pi,\cal{L})} \frak{SpStr}_n (\X)
    =
    r + g_\Y - 1
    =
    \tfrac{3}{2} r + n (g_\X - 1)
\fullstop
\tag{\ref{1792503826975}}
\end{eqn}
Here, $r$ is the local dimension of the stack of finite covers, whilst $g_\Y - 1$ is the dimension of the Picard stack of the spectral curve.

Finally, on the locus of clean spectral structures, we also compute the relative dimension of $F$, which is $r + h^0 (\X, \cal{Q})$ (\autoref{1793103207826}).
Equivalently, this is the local dimension of the analytic space of spectral structures on the fixed vector bundle $\cal{E} = \pi_\ast \cal{L}$.

\paragraph*{AI Declaration}
OpenAI Codex (GPT-6.1, model identifier \texttt{gpt-6.1-sol}) was used during the preparation of this article.
Its use was limited to editorial feedback, suggestions for wording, bibliographic assistance, and proofreading.
The conceptual framework and mathematical proofs are almost entirely the author's own work, with the author retaining full accountability for the final content.

\paragraph*{Acknowledgements}
The author wishes to thank Tom Bridgeland and Marco Gualtieri for helpful discussions during various stages of this project.
The author is particularly grateful to Anand Patel, who, following the appearance of the first version of this paper, generously suggested the beautiful projective-geometric construction underlying the first proof of \autoref{1790777298686}.
This research was partially funded by the European Union's Horizon 2020 Research and Innovation Programme under the Marie Skłodowska-Curie Grant Agreement No.\,101026083 (\href{https://cordis.europa.eu/project/id/101026083}{\texttt{AbQuantumSpec}}), as well as the Leverhulme Trust through the Leverhulme Research Project grant \textit{Extended Riemann-Hilbert Correspondence, Quantum Curves and Mirror Symmetry}.

\clearpage
\enlargethispage{10pt}
\mbox{}
\vspace{-40pt}

\section{Spectral Structures and Their Realisation}
\label{1792146089353}

We introduce spectral structures on vector bundles and develop their geometric interpretation in terms of line bundles on finite covers.
The bridge between these two descriptions is the Finite Sheaf-Module Correspondence (\autoref{1793046131649}).

\subsection{Linear Spectral Structures}
\label{1790804961989}

We begin with the linear algebra underlying spectral structures.
Diagonal actions and Jordan blocks provide the basic examples, whilst more general commutative algebras show how the notion extends beyond structures induced by a single endomorphism.

\begin{definition}[linear spectral structure]
\label{1790804983939}
A \dfn{linear spectral structure} on an $n$-dimensional complex vector space $\E$ is a pair $(\A, \rho)$ consisting of a commutative unital algebra $\A$ of dimension $n$ and an algebra homomorphism $\rho : \A \to \End (\E)$ which turns $\E$ into a free rank-one $\A$-module; i.e., $\E \simeq \A$ as $\A$-modules.
\end{definition}

In this context, we shall sometimes refer to $\A$ as \dfn{spectral algebra} and $\rho$ as \dfn{spectral action}.

\begin{definition}
\label{1791380118188}
An isomorphism of two linear spectral structures $(\A, \rho) \iso (\A', \rho')$ on the same vector space $\E$ is an algebra isomorphism $\alpha : \A \iso \A'$ such that $\rho = \rho' \circ \alpha$.
More generally, an isomorphism $(\E, \A, \rho) \iso (\E', \A', \rho')$ is a pair $(f, \alpha)$ consisting of a linear isomorphism $f : \E \iso \E'$ and an algebra isomorphism $\alpha : \A \iso \A'$ such that $\Ad_f \circ \rho = \rho' \circ \alpha$.
\end{definition}

Here, $\Ad_f : \End(\E) \to \End(\E')$ is the map $g \mapsto f \circ g \circ f^\inv$.

If $\E$ carries a linear spectral structure $(\A, \rho)$, then an element $a \in \A$ acts on a vector $e \in \E$ by ordinary endomorphism action: $a \cdot e = \rho(a) e$.
An identification $\E \simeq \A$ is noncanonical: a choice of $\A$-module generator $e_0 \in \E$ determines an isomorphism $\epsilon : \A \iso \E$ by sending $a \mapsto e = a \cdot e_0$.
Under this identification, an element $a \in \A$ acts on a vector $e = b \cdot e_0 \in \E$ by ordinary multiplication in $\A$; i.e., $a \cdot e = a \cdot (b \cdot e_0) = (ab) \cdot e_0$.
In particular, it follows automatically from this definition that the homomorphism $\rho$ is injective, for if $a \in \A$ is such that $\rho (a) = 0$ then $a = \epsilon^{-1} (a \cdot e_0) = \epsilon^{-1} (\rho(a) e_0) = 0$.

There is a short exact sequence of vector spaces
\begin{eqntag}
\label{1790842918333}
	0
	\too
	\A
	\overset{\rho}{\too}
	\End(\E)
	\too
	\Q
	\too
	0
\qtext{where}
	\Q \coleq \End(\E) / \rho(\A)
\fullstop
\end{eqntag}
The quotient vector space $\Q$ has dimension $n^2 - n = n(n-1)$, and we call it the \dfn{off-diagonal space}.
The image $\rho(\A)$ is equal to its own centraliser in $\End (\E)$; i.e., $C_{\End (\E)} \big( \rho (\A) \big) = \rho (\A)$.
This is because the centraliser of $\rho (\A)$ in $\End (\E)$ is the same as the set of $\A$-linear endomorphisms of $\E$, but this is isomorphic to $\A$ because $\E$ is free of rank one over $\A$.
This means that $\Q$ can be interpreted as the space of endomorphisms modulo all endomorphisms commuting with the entire spectral algebra.

\begin{example}[diagonal spectral structure]
\label{1790843208023}
Take any vector space $\E = \L_1 \oplus \cdots \oplus \L_n$ decomposed into a direct sum of lines.
Consider the algebra $\A = \CC^n$ with coordinatewise multiplication; i.e., $(a_1, \ldots, a_n) (a'_1, \ldots, a'_n) = (a_1 a'_1, \ldots, a_n a'_n)$.
Define the homomorphism $\rho : \A \to \End(\E)$ by $(a_1, \ldots, a_n) \mapsto \diag(a_1, \ldots, a_n)$.
If we choose a generator $v_i \in \L_i$ for each $i$, then the vector $e_0 = (v_1, \ldots, v_n) \in \E$ is an $\A$-module generator of $\E$, inducing an isomorphism $\A \iso \E$ by $a \mapsto \diag(a) e_0$.
In this case, $\rho(\A)$ is the subspace of diagonal endomorphisms, and the off-diagonal space $\Q$ is canonically isomorphic to the subspace of $\End (\E)$ consisting of strictly off-diagonal endomorphisms:
\begin{eqntag}
\label{1790853303053}
	\Q \cong \Oplus_{i \neq j} \Hom (\L_j, \L_i)
\fullstop
\end{eqntag}
If $\E = \CC^n$, then $\End (\E) = \sfop{Mat}_n (\CC)$ is the vector space of all $n \times n$ complex matrices, the image $\rho (\A) = \sfop{Diag}_n (\CC)$ is the subspace of diagonal matrices, and $\Q$ is the subspace of all strictly off-diagonal matrices.
Notice that $\Q$ is not an algebra.
This is the primary motivation for our terminology.
The spectrum $\Spec (\A)$ consists of $n$ distinct points.
\end{example}

\begin{example}[spectral structure from a single Jordan block]
\label{1790843494586}
Take $\E = \CC^n$, and take $\A = \CC [t] / \inner{t^n}$ with usual multiplication mod $t^n$.
Let $N$ be the standard $n \times n$ nilpotent Jordan block (with $1$s above the diagonal and $0$s everywhere else).
Define $\rho : \A \to \End (\E)$ by sending $t \mapsto N$; i.e., $f(t) \mapsto f(N)$ for any $f(t) \in \A$.
Let $e_0 = (0, \ldots, 0, 1) \in \E$ be the cyclic vector; i.e., $e_0, Ne_0, N^2e_0, \ldots, N^{n-1}e_0$ is a basis of $\E$.
Then the map $\A \to \E$ sending $f(t) \mapsto f(N) e_0$ is an $\A$-module isomorphism.
In this case, $\End (\E) = \sfop{Mat}_n(\CC)$ and the image $\rho (\A) = \CC [N] = \Span_\CC \set{I, N, N^2, \ldots, N^{n-1}}$ consists of upper-triangular Toeplitz matrices.
But the $n(n-1)$-dimensional quotient vector space 
\begin{eqntag}
\label{1790854169563}
	\Q = \sfop{Mat}_n (\CC) / \CC [N] \cong \CC^{n^2} / \CC [N]
\end{eqntag}
is not naturally a subspace of $\sfop{Mat}_n (\CC)$.

The spectrum $\Spec (\A)$ is the so-called \textit{fat point}; i.e., a single nonreduced point with multiplicity $n$, understood as the $(n-1)$-st order infinitesimal neighbourhood of the origin in $\CC$.
Thus, the algebra $\A$ is the algebra of functions on this infinitesimal neighbourhood.
The representation $\rho : \A \to \End (\E)$ realises these functions as multiplication operators on $\E$.
Accordingly, $\Q$ is the quotient of all endomorphisms by these multiplication operators.
Very roughly but more intuitively, $\A$ captures the intrinsic infinitesimal directions at a point, and $\Q$ captures the complementary transverse infinitesimal directions inside $\End (\E)$ (although this complementary nature is certainly not canonical).

A choice of the generator $N$ gives a further interpretation: the commutator identifies $\Q$ with the tangent space at $N$ to its $\GL(\E)$-conjugacy orbit, as formalised below in \autoref{1791381265048}.
This describes the first-order variation of the chosen regular endomorphism with its Jordan form fixed.
\end{example}

\enlargethispage{10pt}

\subsubsection*{Regular endomorphisms and monogenic algebras}

We now examine when a linear spectral structure can be induced by an endomorphism.
The answer is governed by regularity of the endomorphism and the ability to generate the algebra by a single element.

\begin{example}[spectral structures from regular endomorphisms]
\label{1790844612906}
Given an endomorphism $T \in \End (\E)$ of an $n$-dimensional vector space, consider the algebra $\A = \CC [T] \subset \End (\E)$ and let $\rho$ be the inclusion map.
If $m_T \in \CC[t]$ is the minimal polynomial of $T$, then there is a canonical isomorphism $\CC [T] \cong \CC[t] / \inner{ m_T }$ determined by the evaluation map $\CC[t] \to \End (\E)$ sending $f(t) \mapsto f(T)$.
Therefore, $\dim \CC [T] = \deg m_T$ which means the algebra $\CC [T]$ is $n$-dimensional if and only if the minimal and the characteristic polynomials of $T$ coincide.
Equivalently, $T$ admits a cyclic vector: there exists some $e_0 \in \E$ such that $e_0, T e_0, T^2 e_0, \ldots, T^{n-1} e_0$ span the vector space $\E$.
An endomorphism $T$ with this property is called \textit{regular} or \textit{cyclic}.
Thus, \textit{an endomorphism $T \in \End (\E)$ induces a linear spectral structure $(\A,\rho)$ on $\E$, in the sense that $\rho (\A) = \CC [T]$, if and only if it is regular.}

In addition, if $T$ has distinct eigenvalues $\lambda_1, \ldots, \lambda_n$, the vector space $\E$ acquires the eigendecomposition $\E \cong \L_1 \oplus \cdots \oplus \L_n$ into lines.
In this case, the vector space $\Q$ is naturally isomorphic to the subspace of all off-diagonal endomorphisms of $\E$ as in \eqref{1790853303053}.
Thus, we may interpret the vector space $\Q$ as capturing the variations of the eigenlines of $T$ whilst keeping the eigenvalues fixed.

Observe that the diagonal spectral structure on the split vector space $\E \cong \L_1 \oplus \cdots \oplus \L_n$ can be induced by \textit{any} diagonal endomorphism $T = \diag(\lambda_1, \ldots, \lambda_n)$ with distinct eigenvalues.
This example suggests an intuitive interpretation of spectral structures on a vector space as eigendecompositions without a distinguished choice of eigenvalues.

More generally, $\E$ decomposes into the generalised eigenspaces $\E_1 \oplus \cdots \oplus \E_m$ of $T$.
Regularity of $T$ means that each generalised eigenspace corresponds to a single Jordan block.
Therefore, we find that the vector space $\Q$ has two distinct contributions, one coming from the off-diagonal endomorphism pieces like in \autoref{1790843208023} and the other one coming from transverse infinitesimal directions as in \autoref{1790843494586}:

\clearpage
\enlargethispage{100pt}
\mbox{}
\vspace{-55pt}

\begin{eqntag}
\label{1790855272405}
	\Q 
	\cong
	\left( \Oplus_{i \neq j} \Hom (\E_j, \E_i) \right)
	\oplus
	\left( \Oplus_{i} \End (\E_i) / \CC [N_i] \right)
\fullstop{,}
\end{eqntag}
where $T |_{\E_i} = \lambda_i I + N_i$.
Thus, we may interpret the space $\Q$ more generally as capturing the variations of the eigendecomposition of $T$ whilst keeping its conjugacy invariants fixed.
\end{example}

Notice that the algebra $\A = \CC [T]$ has a special property that it is generated by a single element $T$.
The above example suggests a simple criterion to detect when a given linear spectral structure can be induced by a regular endomorphism.

\begin{lemma}[monogenicity criterion]
\label{1790846785676}
A linear spectral structure on a vector space $\E$ is induced by a regular endomorphism of $\E$ if and only if the spectral algebra is monogenic.
In other words, there exists a regular endomorphism $T \in \End (\E)$ such that $\rho (\A) = \CC [T] \subset \End (\E)$ if and only if the algebra $\A$ is generated as a $\CC$-algebra by a single element.
\end{lemma}

\begin{proof}
If $\rho (\A) = \CC [T]$, then $\rho$ defines the algebra isomorphism $\A \iso \CC [T]$ by virtue of being injective.
Since $\CC [T]$ is monogenic with generator $T$, the algebra $\A$ is monogenic with generator $\rho^\inv (T)$.
Conversely, if $a \in \A$ is a generator, let $T \coleq \rho (a)$.
Then $\rho$ defines the algebra isomorphism $\A \iso \CC [T]$, and the fact that $\dim \A = \dim \E$ forces $T$ to be regular.	
\end{proof}

We also record the following observation that, if a linear spectral structure is induced by an endomorphism, this endomorphism is far from unique.

\begin{lemma}[changing the generator]
\label{1791379641510}
Two regular endomorphisms $T,T' \in \End(\E)$ induce the same linear spectral structure on $\E$ if and only if $\CC[T] = \CC[T']$; i.e., each one is a polynomial in the other.
Moreover, the set $\set{ T \in \End (\E) ~\big|~ \rho (\A) = \CC [T]}$ of all endomorphisms that induce the same monogenic spectral structure $(\A, \rho)$ is a Zariski-open and therefore dense subset of $\rho (\A)$.
\end{lemma}

\begin{proof}
The first assertion follows from the definition of the induced spectral structure.
For the second assertion, note that the set of endomorphisms in question is the image under the injective homomorphism $\rho$ of the subset $\set{a \in \A ~\big|~ \A = \CC[a]}$ of all generators of $\A$.
An element $a \in \A$ is a generator if and only if the vectors $1, a, a^2, \ldots, a^{n-1}$ form a basis of $\A$ as a complex vector space.
If we fix one generator $a \in \A$, then for any element $b \in \A$, let us write the vectors $1, b, b^2, \ldots, b^{n-1}$ in terms of the basis $1, a, a^2, \ldots, a^{n-1}$.
Let $B \in \frak{gl}_n (\CC)$ be the matrix of their coordinate vectors.
Then the condition that $b$ is a generator of $\A$ is equivalent to $\det B \neq 0$, which is a Zariski-open condition.
\end{proof}

\begin{lemma}[commutators and trace orthogonality]
\label{1791381265048}
Suppose a linear spectral structure $(\A, \rho)$ on $\E$ is induced by a regular endomorphism.
Then, once such a regular endomorphism $T$ has been specified, the commutator linear map $\ad_T : \End(\E) \to \End(\E)$ induces an isomorphism
\begin{eqn}
	\Q = \End (\E) / \CC[T] \iso \image (\ad_T)
\qqtext{sending}
	[S] \mapstoo \ad_T (S) = [T, S]
\fullstop
\end{eqn}
Moreover, under this identification, $\Q$ and $\A \cong \rho (\A)$ are orthogonal subspaces with respect to the natural trace pairing $(f,g) \mapsto \tr(fg)$ on $\End(\E)$:
\begin{eqn}
	\Q \cong \image (\ad_T) = (\CC[T])^\perp \cong \A^\perp
\fullstop
\end{eqn}
\end{lemma}

Consequently, we have the interpretation that, although many different regular endomorphisms $T$ can induce the same spectral structure, the subspace $\A^\perp$ is independent of the chosen generator $T$.
Only the isomorphism $\Q \cong \A^\perp$ depends on $T$.

\begin{proof}
Since $T$ is regular, the kernel $\ker (\ad_T)$ is the centraliser of $T$ in $\End (\E)$.
But this is equal to $\CC [T]$, so $\Q \cong \image (\ad_T)$ by the First Isomorphism Theorem.
Indeed, for any endomorphisms $S, U \in \End(\E)$, cyclicity of the trace gives $\tr \big([T,S] U \big) = \tr \big(S[U,T]\big)$.
Consequently, every element of $\image(\ad_T)$ is orthogonal to $\ker(\ad_T) = \CC[T]$.
The trace pairing on $\End(\E)$ is nondegenerate, so $\CC[T]^\perp$ has dimension $n^2-n$.
The First Isomorphism Theorem gives the same dimension for $\image(\ad_T)$, and therefore $\image(\ad_T) = \CC[T]^\perp$.
\end{proof}

\clearpage

\begin{remark}[orthogonality $\not\Rightarrow$ decomposition]
\label{1793100364071}
The above orthogonality between $\Q$ and $\A$ with respect to the trace pairing does \textit{not} imply a decomposition $\End(\E) \cong \A \oplus \Q$ in a canonical way.
In particular, despite being orthogonal subspaces, $\rho(\A)$ and $\rho(\A)^\perp$ of $\End(\E)$ can intersect nontrivially, and therefore fail to form an internal direct sum $\End(\E) = \rho(\A) \oplus \rho(\A)^\perp$ in general.
Specifically, the natural addition map
\begin{eqn}
	\rho(\A) \oplus \rho(\A)^\perp
	\too
	\End(\E)
\qqtext{given by}
	(a,b) \mapstoo a+b
\end{eqn}
need not be an isomorphism: its kernel is canonically isomorphic to $\rho(\A) \cap \rho(\A)^\perp$, which may be nonzero.
Consequently, the identification $\Q \cong \A^\perp$ determined by the chosen generator $T$ in \autoref{1791381265048} may not yield a splitting of the exact sequence \eqref{1790842918333}.
For example, if $N$ is a nilpotent Jordan block of size $n \geq 2$ like in \autoref{1790843494586}, then $N$ is a nonzero element of both $\CC[N]$ and $\CC[N]^\perp$: for every polynomial $f$, the matrix $Nf(N)$ is nilpotent and traceless.
\end{remark}

The algebra $\A$ in our \autoref{1790804983939} is not assumed to be generated by a single element, which means there exist spectral structures that are not induced by endomorphisms.
Here is the simplest example of one.

\begin{example}[non-monogenic spectral structure]
\label{1790846065686}
Take $\A = \CC[x,y] / \inner{x,y}^2$.
As a vector space, $\A$ is $3$-dimensional with basis $1,x,y$.
Let $\E \simeq \CC^3$ be this vector space.
Then $\E$ is a free rank-one $\A$-module and $\rho$ is simply the $\A$-multiplication, so $(\A,\rho)$ is a linear spectral structure on $\E$.
However, $\A$ cannot be generated by one element.
Indeed, every element $a = a_0 + a_1 x + a_2 y \in \A$ satisfies $(a - a_0)^2 = 0$, so $a^k = a_0^k + k a_0^{k-1} (a_1 x + a_2 y)$ for every power $k$.
Consequently, $\CC [a] = \Span_\CC \set{1, a_1 x + a_2 y}$ which has dimension at most $2$ hence it cannot possibly be the full three-dimensional algebra $\A$.
So there is \textit{no} endomorphism $T \in \End (\E)$ such that $\A = \CC [T]$.

With respect to the basis $1,x,y$, multiplication by $a = a_0 + a_1 x + a_2 y$ has matrix
\vspace{-5pt}
\begin{eqn}
	\rho (a) = \mtx{ a_0 & 0 & 0 \\ a_1 & a_0 & 0 \\ a_2 & 0 & a_0 }
\fullstop
\vspace{-5pt}
\end{eqn}
The six-dimensional quotient vector space $\Q$ can then be described, with respect to the chosen basis $1,x,y$, as the vector subspace of $3 \times 3$ matrices whose first column is zero:
\begin{eqn}
	\Q = \sfop{Mat}_3 (\CC) / \rho (\A) 
	\simeq
	\set{ \mtx{ 0 & \ast & \ast \\ 0 & \ast & \ast \\ 0 & \ast & \ast } }
	\subset \sfop{Mat}_3 (\CC)
\fullstop
\end{eqn}
The spectrum $\Spec (\A)$ is a fat point of multiplicity $3$.
\end{example}

\enlargethispage{30pt}

\subsection{Spectral Structures on Curves}
\label{1790715973901}

Let $\X$ be a compact Riemann surface (i.e., a smooth connected compact complex curve) and let $\cal{E}$ be a holomorphic vector bundle on $\X$ of rank $n$.
A \textit{spectral structure} on $\cal{E}$ is, roughly speaking, a rank-$n$ algebra acting on $\cal{E}$ in a way that $\cal{E}$ becomes a line bundle over this algebra.
More precisely:

\begin{definition}[spectral structure]
\label{1789147824525}
A \dfn{spectral structure} on $\cal{E}$ is a pair $(\cal{A}, \rho)$ consisting of a commutative unital $\cal{O}_\X$-algebra $\cal{A}$, locally free of rank $n$, and an $\cal{O}_\X$-algebra morphism $\rho : \cal{A} \to \cal{End} (\cal{E})$ such that $\cal{E}$, with the induced $\cal{A}$-module structure, is locally free of rank $1$.
We will also sometimes refer to $\cal{A}$ as a \dfn{spectral algebra}.
\end{definition}

Unpacking parts of this definition, a locally free rank-$n$ commutative unital $\cal{O}_\X$-algebra $\cal{A}$ is a vector bundle on $\X$ of rank $n$ equipped with two bundle maps, called \textit{unit} and \textit{multiplication},
\begin{eqn}
	e : \cal{O}_{\X} \too \cal{A}
\qqtext{and}
	\mu : \cal{A} \otimes \cal{A} \too \cal{A}
\end{eqn}
satisfying the unit, commutativity, and associativity identities:
\begin{eqn}
	\mu (e \otimes a) = a = \mu (a \otimes e)
~,
~~
	\mu (a \otimes b) = \mu (b \otimes a)
~,
~~
	\mu \big( \mu (a \otimes b) \otimes c \big) = \mu \big( a \otimes \mu (b \otimes c) \big)
\fullstop{,}
\end{eqn}
\clearpage
for any local sections $a, b, c \in \cal{A}$; i.e., bundle maps $a, b, c : \cal{O}_\U \to \cal{A}$ for some open $\U \subset \X$.
The homomorphism $\rho$ carries the unit $e$ to the identity endomorphism $\id$ of $\cal{E}$ and it respects the multiplication on $\cal{A}$ and the ordinary composition $\circ$ of endomorphisms of $\cal{E}$:
\begin{eqn}
	\rho (e) = \id
\qqtext{and}
	\rho \big( \mu (a \otimes b) \big) = \rho(a) \circ \rho(b)
\fullstop
\end{eqn}
An element $a \in \cal{A}$ then acts on a local section $v \in \cal{E}$ through $\rho$ by ordinary endomorphism action: $v \mapsto av \coleq \rho(a) v$. 
This is the induced $\cal{A}$-module structure on $\cal{E}$.
That $\cal{E}$ is a locally free rank-one $\cal{A}$-module means $\cal{E} \simeq \cal{A}$ locally as $\cal{A}$-modules; i.e., every point on $\X$ has a neighbourhood $\U$ such that there is an $\cal{A}$-linear isomorphism $\cal{A}_\U \iso \cal{E}_\U$ between their restrictions to $\U$.
Such an identification is induced by choosing a local $\cal{A}$-module generator $e_0 \in \cal{E}_\U$, giving us a map $\epsilon : \cal{A}_\U \iso \cal{E}_\U$ that sends $a \mapsto \rho(a) e_0$.
Under this identification, an element $a \in \cal{A}_\U$ acts on a section $v = b e_0 \in \cal{E}_\U$ for some $b \in \cal{A}_\U$ by ordinary $\cal{A}$-multiplication: $v = b e_0 \mapsto (ab) e_0 \eqcol av$.
In particular, it follows automatically from this definition that $\rho$ is injective, for if $a \in \cal{A}$ is such that $\rho (a) = 0$ then it follows that $a = \epsilon^{-1} (\rho(a) e_0) = 0$.

\begin{definition}[groupoid of spectral structures]
\label{1789150206981}
For a fixed vector bundle $\cal{E}$, an isomorphism $(\cal{A}, \rho) \iso (\cal{A}', \rho')$ of two spectral structures on $\cal{E}$ is an $\cal{O}_\X$-algebra isomorphism $\alpha : \cal{A} \iso \cal{A}'$ such that $\rho = \rho' \circ \alpha$.
This determines the groupoid of spectral structures on the vector bundle $\cal{E}$, denoted by
\begin{eqn}
	\mathbf{SpStr} (\cal{E}) 
	\quad:\quad
	\set{ (\cal{A}, \rho) ~\Big|~  \text{spectral structures on $\cal{E}$} }
\fullstop
\end{eqn}
More generally, unfixing the bundle $\cal{E}$, an isomorphism $(\cal{E}, \cal{A}, \rho) \iso (\cal{E}', \cal{A}', \rho')$ is a pair $(\phi, \alpha)$ consisting of a vector bundle isomorphism $\phi : \cal{E} \iso \cal{E}'$ and an $\cal{O}_\X$-algebra isomorphism $\alpha : \cal{A} \iso \cal{A}'$ such that $\Ad_\phi \circ \rho = \rho' \circ \alpha$.
This determines the groupoid of all vector bundles of rank $n$ with spectral structure over $\X$, denoted by
\begin{eqn}
	\mathbf{SpStr}_n (\X)
	\quad:\quad
	\set{ (\cal{E}, \cal{A}, \rho) ~\Big|~ \text{rank-$n$ vector bundles over $\X$ with spectral structure} }
\fullstop
\end{eqn}
\end{definition}

\subsubsection*{Fibres of the forgetful functor}

There is a forgetful functor to the groupoid $\mathbf{Bun}_n (\X)$ of all rank-$n$ vector bundles on $\X$:
\begin{eqntag}
\label{1791407291197}
	F : \mathbf{SpStr}_n (\X) \too \mathbf{Bun}_n (\X)
\qtext{sending}
	(\cal{E}, \cal{A}, \rho) \mapsto \cal{E}
\fullstop
\end{eqntag}
This gives us a good opportunity to remind the reader of two different notions of the fibre of a functor over an object.
For any $\cal{E}_0 \in \mathbf{Bun}_n (\X)$, the \textit{strict fibre} $F^\inv (\cal{E}_0)$ of this functor above $\cal{E}_0$ is the groupoid that parametrises all triples $(\cal{E}, \cal{A}, \rho)$ where $\cal{E} = \cal{E}_0$.
Morphisms are required to induce the identity morphism of $\cal{E}_0$ under $F$.
Thus, it is precisely the groupoid of spectral structures on $\cal{E}_0$; in symbols, $F^\inv (\cal{E}_0) = \mathbf{SpStr} (\cal{E}_0)$.

On the other hand, its \textit{essential fibre} (or \textit{homotopy fibre}) $F^{[\inv]} (\cal{E}_0)$ keeps track of such triples where $\cal{E}$ is merely isomorphic to $\cal{E}_0$ with isomorphism as part of the data.
Thus, $F^{[\inv]} (\cal{E}_0)$ is the groupoid of tuples $(\cal{E}, \cal{A}, \rho; \iota)$ consisting of an object $(\cal{E}, \cal{A}, \rho) \in \mathbf{SpStr}_n (\X)$ together with an identification $\iota : \cal{E} \iso \cal{E}_0$.
An isomorphism $(\cal{E}, \cal{A}, \rho; \iota) \iso (\cal{E}', \cal{A}', \rho'; \iota')$ is an isomorphism $(\phi, \alpha) : (\cal{E}, \cal{A}, \rho) \iso (\cal{E}', \cal{A}', \rho')$ such that $\phi : \cal{E} \to \cal{E}'$ intertwines the identifications $\iota, \iota'$; i.e., $\iota' \circ \phi = \iota$.
This groupoid is commonly denoted using the $(F \downarrow \cal{E}_0)$ and referred to as the \textit{comma category}.

In general, the strict and essential fibres are not equivalent as groupoids, but in this elementary situation this happens to be so.
Indeed, there is a canonical equivalence of groupoids $F^\inv (\cal{E}_0) \iso F^{[\inv]} (\cal{E}_0)$ given by sending $(\cal{A}, \rho) \mapsto (\cal{E}_0, \cal{A}, \rho; \id_{\cal{E}_0})$.
The inverse equivalence is given by transporting the spectral structure $\rho$ on $\cal{E}$ to $\cal{E}_0$ through the identification $\iota$; i.e., $(\cal{E}, \cal{A}, \rho, \iota) \mapsto (\cal{A}, \Ad_\iota \circ \rho)$.
In summary, for any vector bundle $\cal{E} \in \mathbf{Bun}_n (\X)$,
\begin{eqntag}
\label{1791390052260}
	F^\inv (\cal{E}) = \mathbf{SpStr} (\cal{E})
	\cong
	F^{[\inv]} (\cal{E}) = (F \downarrow \cal{E})
\fullstop
\end{eqntag}

\begin{proposition}[spectral structures as subalgebras]
\label{1791405911510}
For a fixed vector bundle $\cal{E}$, an isomorphism between two spectral structures $(\cal{A}, \rho)$ and $(\cal{A}', \rho')$, if it exists, is unique.
This isomorphism exists precisely when $\rho (\cal{A}) = \rho' (\cal{A}')$ as subsheaves of algebras of $\cal{End} (\cal{E})$.
In particular, every object of $\mathbf{SpStr} (\cal{E})$ has a trivial automorphism group.
Consequently, $\mathbf{SpStr} (\cal{E})$ is canonically equivalent to the discrete groupoid of commutative unital subalgebras $\cal{B} \subset \cal{End} (\cal{E})$ that are locally free of rank $n$ over $\cal{O}_\X$ and make $\cal{E}$ into a locally free rank-one $\cal{B}$-module.
\end{proposition}

Here, the term `discrete groupoid' refers to a groupoid in which the only morphisms that exist are the identity morphisms.

\begin{proof}
Both action maps $\rho$ and $\rho'$ are injective.
If $\rho = \rho' \circ \alpha$ for some $\alpha : \cal{A} \iso \cal{A}'$, then their images in $\cal{End} (\cal{E})$ must coincide.
This equality of images forces $\alpha = \rho'^{-1} \circ \rho$ where $\rho'^{-1}$ is the map $\rho' (\cal{A}') \iso \cal{A}'$ inverting $\rho'$.

Conversely, if $\rho (\cal{A}) = \rho' (\cal{A}')$, then $\alpha \coleq \rho'^{-1} \circ \rho$ defines an isomorphism of $\cal{O}_\X$-algebras $\cal{A} \iso \cal{A}'$ which is clearly compatible with the actions $\rho, \rho'$.
Hence an isomorphism of spectral algebras on $\cal{E}$.
Taking the two spectral structures to be equal shows that their only automorphism is the identity.

Finally, sending a spectral structure $(\cal{A}, \rho)$ to its image $\cal{B} \coleq \rho (\cal{A})$ defines the stated equivalence of groupoids.
The inverse equivalence is given by simply taking $\cal{A} = \cal{B}$ and $\rho$ the inclusion map $\cal{B} \inj \cal{End} (\cal{E})$.
\end{proof}

Notice that the vector bundle $\cal{E}$ in the groupoid $\mathbf{SpStr} (\cal{E})$ is fixed: the isomorphisms between spectral structures are required to induce the identity on $\cal{E}$.
If we allow nontrivial bundle automorphisms of $\cal{E}$, the automorphism groups of spectral structures may be nontrivial.
More generally, the forgetful functor \eqref{1791407291197} is faithful: a given bundle isomorphism $\phi$ uniquely determines any compatible algebra isomorphism through $\rho' \circ \alpha = \Ad_\phi \circ \rho$.

\subsubsection*{Examples}

Let us now illustrate the definition of spectral structures in a number of examples.
First, we look at a direct sum of line bundles (\autoref{1791434155987}) and then the pushforward of a line bundle along a finite cover (\autoref{1789283249898}).
In both cases, the spectral action is given concretely by multiplication of sections by functions.
Then we study a special class of explicit examples called cyclic spectral structures (\autoref{1791437029550} and \autoref{1791450524665}): they are built from a line bundle and a section of one of its tensor powers.
We write out their algebra multiplication and spectral action explicitly in ranks two and three.

\begin{example}[split spectral structures]
\label{1791434155987}
Suppose the vector bundle $\cal{E}$ is a direct sum of line bundles: $\cal{E} = \cal{L}_1 \oplus \cdots \oplus \cal{L}_n$.
Take $\cal{A} = \cal{O}_\X^{\oplus n}$ with its componentwise product algebra structure: $(f_1, \ldots, f_n) \cdot (g_1, \ldots, g_n) = (f_1 g_1, \ldots, f_n g_n)$.
Let it act on $\cal{E}$ diagonally: $(f_1, \ldots, f_n) \cdot (e_1, \ldots, e_n) = (f_1 e_1, \ldots, f_n e_n)$.
Then $\cal{E}$ is a locally free $\cal{A}$-module of rank $1$: local generators $e_i$ of $\cal{L}_i$ form a local section $e = (e_1, \ldots, e_n)$ of $\cal{E}$ which can serve as an $\cal{A}$-module generator of $\cal{E}$ because any local section $s$ of $\cal{E}$ can be written as $s = (f_1 e_1, \ldots, f_n e_n) = (f_1, \ldots, f_n) \cdot e$ for some $(f_1, \ldots, f_n) \in \cal{A}$.
The endomorphism bundle $\cal{End} (\cal{E})$ decomposes accordingly:
\begin{eqn}
	\cal{End} (\cal{E})
	= \Oplus_{i,j = 1}^n \cal{Hom} (\cal{L}_j, \cal{L}_i)
	= \Oplus_{i=1}^n \cal{End} (\cal{L}_i) 
		\oplus 
		\Oplus_{\substack{i,j=1 \\ i \neq j}}^n \cal{Hom} (\cal{L}_j, \cal{L}_i)
\fullstop
\end{eqn}
Since each $\cal{L}_i$ is a line bundle, there is a canonical isomorphism $\cal{End} (\cal{L}_i) \cong \cal{O}_\X$, and so the homomorphism $\rho : \cal{O}_\X^{\oplus n} \to \cal{End} (\cal{E})$ is just inclusion as the diagonal endomorphisms:
\begin{eqn}
	\rho (\cal{A}) = \Oplus_{i=1}^n \cal{End} (\cal{L}_i)
\fullstop
\qedhere
\end{eqn}
In particular, every vector bundle on $\PP^1$ admits a split spectral structure.
This is because, by the Birkhoff-Grothendieck Theorem \cite[Théorème~2.1]{Grothendieck1957}, every vector bundle $\cal{E}$ over $\PP^1$ is isomorphic to a direct sum of line bundles: $\cal{E} \simeq \cal{O} (k_1) \oplus \cdots \cal{O} (k_n)$.
\end{example}

\begin{example}[finite covers]
\label{1789283249898}
Let $\pi : \Y \to \X$ be a smooth cover of degree $n$, and choose a line bundle $\cal{L}$ on $\Y$.
Then the pushforward $\cal{E} = \pi_\ast \cal{L}$ is a vector bundle on $\X$ of rank $n$.
It carries a canonical spectral structure $\rho : \cal{A} \to \cal{End} (\cal{E})$ with $\cal{A} = \pi_\ast \cal{O}_\Y$.
It is induced by the canonical $\cal{O}_\Y$-action on $\cal{L}$ by multiplication of sections by functions.

Concretely, take an open set $\U \subset \X$ and put $\V \coleq \pi^\inv (\U)$.
Then there are canonical isomorphisms 
\begin{eqn}
	\cal{A} (\U) \cong \cal{O}_\Y (\V)
\qtext{and}
	\cal{E}(\U) \cong \cal{L} (\V)
\fullstop{,}
\end{eqn}
so an algebra element $f \in \cal{A} (\U)$ is nothing but a holomorphic function over $\V$, and a section $e \in \cal{E}(\U)$ is nothing but a section of $\cal{L}$ over $\V$.
Consequently, the action $\rho$ is given by ordinary multiplication: 
\begin{eqn}
	\rho(f) (e) \coleq fe
\fullstop
\end{eqn}
Explicitly, the procedure is this: given $f \in \cal{A} (\U)$ and $e \in \cal{E}(\U)$, to calculate $\rho(f) (e)$ we first regard $f$ and $e$ as elements of $\cal{O}_\Y (\V)$ and $\cal{L} (\V)$ via the canonical isomorphisms, then we multiply $f$ and $e$ in the ordinary sense, thus getting a new element $fe \in \cal{L} (\V)$ which we regard as an element of $\cal{E} (\U)$ via the canonical isomorphism.

Even more concretely, if $\U$ is simply connected and avoids branch points, then $\V$ is a disjoint union of isomorphic copies of $\U$; i.e., $\V = \U_1 \sqcup \cdots \sqcup \U_n$ with $\U_i \cong \U$ via $\pi$.
Then
\begin{eqn}
	\cal{A} (\U)
	\cong \cal{O}_\Y (\U_1) \oplus \cdots \oplus \cal{O}_\Y (\U_n)
\qtext{and}
	\cal{E} (\U)
	\cong \cal{L} (\U_1) \oplus \cdots \oplus \cal{L} (\U_n)
\fullstop
\end{eqn}
So an algebra element $f \in \cal{A} (\U)$ is a tuple $(f_1, \ldots, f_n)$ of functions $f_i \in \cal{O}_\Y (\U_i)$
Similarly, a vector bundle section $e \in \cal{E} (\U)$ is a tuple $(e_1, \ldots, e_n)$ of line bundle sections $e_i \in \cal{L} (\U_i)$.
Then the action of $\cal{A} (\U)$ on $\cal{E} (\U)$ is diagonal: 
\begin{eqn}
	\rho(f) (e) = (f_1, \ldots, f_n) \cdot (e_1, \ldots, e_n) = (f_1 e_1, \ldots, f_n e_n)
\fullstop
\end{eqn}
Thus, locally away from the branch locus, the spectral structure $(\cal{A}, \rho)$ behaves like a split spectral structure of \autoref{1791434155987}.
\end{example}

\begin{example}[cyclic spectral structures of rank 2]
\label{1791437029550}
Let $\cal{M}$ be any line bundle of degree $k \geq 0$ on $\X$ and pick a section $s \in \H^0 (\X, \cal{M}^{\otimes 2} )$.
Take $\cal{A}$ and $\cal{E}$ to be the same rank-$2$ vector bundle
\begin{eqn}
	\cal{A} = \cal{O}_\X \oplus \cal{M}^{-1} = \cal{E}
\fullstop
\end{eqn}
We can define an $\cal{O}_\X$-algebra structure on $\cal{A}$ using the map $s : \cal{M}^{-2} \to \cal{O}_\X$ as follows.
Local sections of $\cal{A}$ are of the form $f + \alpha$ where $f \in \cal{O}_\X$ and $\alpha \in \cal{M}^{-1}$, and the product between two such sections of $\cal{A}$ is given by
\begin{eqn}
	(f + \alpha) \cdot (f' + \alpha') \coleq \big(ff' + s (\alpha \otimes \alpha') \big) + (f \alpha' + f' \alpha)
\fullstop
\end{eqn}
Notice that $s (\alpha \otimes \alpha')$ is a function because $\alpha \otimes \alpha'$ is a section of $\cal{M}^{-2}$.
The $\cal{A}$-module structure on $\cal{E}$ is described in terms of endomorphism matrices of $\cal{O}_\X \oplus \cal{M}^{-1}$ like so:
\begin{eqn}
	\rho (f + \alpha) \coleq \mtx{ f & s (\alpha \otimes \dummy) \\ \alpha & f }
\qqtext{i.e.}
	\rho (f + \alpha) \left( \mtx{f' \\ \alpha'} \right)
	= \mtx{ ff' + s(\alpha \otimes \alpha') \\ f \alpha' + f' \alpha }
\fullstop
\end{eqn}
Even more explicitly, let $\mu$ be a local generator of $\cal{M}^{-1}$ over a sufficiently small open set $\U \subset \X$.
Then $\cal{M}^{-1}_\U \cong \cal{O}_\U \mu$,~ $\cal{M}^{-2}_\U \cong \cal{O}_\U \mu^2$, and $\cal{A}_\U \cong \cal{O}_\U \oplus \cal{O}_\U \mu$.
Put $c \coleq s (\mu^2) \in \cal{O}_\U$.
Then the multiplication in $\cal{A}$ is given by 
\begin{eqn}
	(a + b\mu) \cdot (a' + b'\mu) = aa' + (ab' + a'b) \mu + bb' s(\mu^2) = (aa' + bb'c) + (ab' + a'b) \mu
\fullstop
\end{eqn}
In other words, there is an isomorphism of $\cal{O}_\U$-algebras: $\cal{A}_\U \cong \cal{O}_\U [\mu] \big/ \inner{\mu^2 - c}$.
\end{example}

\begin{example}[cyclic spectral structures of higher rank]
\label{1791450524665}
The same construction as in \autoref{1791437029550} extends to higher rank.
For example, let $\cal{M}$ be a line bundle of degree $k \geq 0$ on $\X$ and pick a section $s \in \H^0 (\X, \cal{M}^{\otimes 3})$.
Take
\begin{eqn}
	\cal{A} = \cal{O}_\X \oplus \cal{M}^{-1} \oplus \cal{M}^{-2} = \cal{E}
\end{eqn}
with multiplication determined by the map $s : \cal{M}^{-3} \to \cal{O}_\X$.
Namely, local sections of $\cal{A}$ are of the form $f + \alpha + \beta$ where $f \in \cal{O}_\X$,~ $\alpha \in \cal{M}^{-1}$, and $\beta \in \cal{M}^{-2}$, and the product between two such sections $f + \alpha + \beta$ and $f' + \alpha' + \beta'$ is given by
\begin{eqn}
	\Big(ff' + s(\beta \otimes \alpha') + s(\alpha \otimes \beta') \Big)
	+ \Big( f \alpha' + f' \alpha + s (\beta \otimes \beta') \Big)
	+ \Big( f \beta' + f' \beta + \alpha \otimes \alpha' \Big)
\fullstop
\end{eqn}
Here, the first bracket belongs to $\cal{O}_\X$, second to $\cal{M}^{-1}$, and third to $\cal{M}^{-2}$.
In particular, $s(\beta \otimes \beta') \in \cal{M}^{-1}$ is the image of $\beta \otimes \beta' \in \cal{M}^{-4}$ under the induced map $s \otimes 1 : \cal{M}^{-4} \to \cal{M}^{-1}$.
Similarly, the $\cal{A}$-module structure on $\cal{E}$ is given by
\begin{eqn}
	\rho (f + \alpha + \beta)
	= \mtx{ f & s (\beta \otimes \dummy) & s (\alpha \otimes \dummy) 
		\\ \alpha & f & s (\beta \otimes \dummy)
		\\ \beta & \alpha & f
		} : \mtx{f' \\ \alpha' \\ \beta'}
	\mapsto \mtx{ ff' + s(\beta \otimes \alpha') + s(\alpha \otimes \beta') 
	\\ f \alpha' + f' \alpha + s (\beta \otimes \beta')
	\\ f \beta' + f' \beta + \alpha \otimes \alpha'
	}
\fullstop
\end{eqn}
Finally, if $\mu$ is a local generator of $\cal{M}^{-1}$ over $\U$, let $c \coleq s(\mu^3)$.
Then it induces an $\cal{O}_\U$-algebra isomorphism $\cal{A}_\U \cong \cal{O}_\U [\mu] \big/ \inner{\mu^3 - c}$.
\end{example}

\subsection{Finite Sheaf-Module Correspondence}
\label{1792143161338}

Given a spectral algebra $\cal{A}$ on $\X$ --- i.e., a commutative unital $\cal{O}_\X$-algebra, locally free of rank $n$ --- we can form its relative analytic spectrum $\Y$ over $\X$.
This is a complex analytic space equipped with a finite flat holomorphic map $\pi : \Y \to \X$ of degree $n$ such that $\cal{A}$ is the pushforward of its structure sheaf:
\begin{eqntag}
\label{1789153157048}
	\pi : \Y \coleq \Spec_\X (\cal{A}) \too \X
\qqtext{such that}
	\pi_\ast \cal{O}_\Y = \cal{A}
\fullstop
\end{eqntag}
More generally, for any finite holomorphic map $\pi : \Y \to \X$, the pushforward $\cal{A} \coleq \pi_\ast \cal{O}_\Y$ is a commutative unital $\cal{O}_\X$-algebra, coherent as an $\cal{O}_\X$-module by the Finite Mapping Theorem \cite[Theorem~2.7.1]{GuiCAG}; see also \cite[Theorem~I.1.67]{GreuelLossenShustin2007}.
The assignments 
\begin{eqn}
	\pi \mapstoo \pi_\ast \cal{O}_\Y
\qqtext{and}
	\cal{A} \mapstoo \Spec_\X (\cal{A})
\end{eqn}
give an anti-equivalence between complex analytic spaces finite over $\X$ and coherent commutative unital $\cal{O}_\X$-algebras \cite[Theorem~2.9.3]{GuiCAG}; see also \cite[Ch.~III, \S1(b), p.~94]{BanicaStanasila1976}.
Under this correspondence, finite flat covers of degree $n$ correspond precisely to algebras locally free of rank $n$ over $\cal{O}_\X$.

For a given finite flat map $\pi : \Y \to \X$ with $\cal{A} = \pi_\ast \cal{O}_\Y$, there is a corresponding equivalence between coherent sheaves on $\Y$ and coherent $\cal{A}$-modules on $\X$.
A coherent sheaf $\cal{L}$ on $\Y$ gives the coherent $\cal{A}$-module $\pi_\ast \cal{L}$ on $\X$, with the algebra action induced by multiplication of sections of $\cal{L}$ by holomorphic functions on $\Y$.
Conversely, a coherent $\cal{A}$-module $\cal{E}$ on $\X$ determines a coherent sheaf $\cal{L}$ on $\Y$ together with an $\cal{A}$-linear isomorphism $\iota : \pi_\ast \cal{L} \iso \cal{E}$.
Thus, the algebra action is precisely the extra information that allows a sheaf on $\X$ to be interpreted as a sheaf on the covering space $\Y$.
This equivalence of categories is one of the most important tools in our study of spectral structures, so let us recall how it works in suitable detail.

Write $\mathbf{Coh} (\Y)$ for the category whose objects are coherent $\cal{O}_\Y$-modules $\cal{L}$, and whose morphisms are $\cal{O}_\Y$-linear sheaf morphisms.
Similarly, write $\mathbf{Coh} (\X, \cal{A})$ for the category whose objects are coherent $\cal{O}_\X$-modules $\cal{E}$ equipped with an $\cal{A}$-action extending their $\cal{O}_\X$-module structure, and whose morphisms are $\cal{A}$-linear sheaf morphisms.

\begin{proposition}[finite sheaf-module correspondence]
\label{1793046131649}
Let $\pi : \Y \to \X$ be a finite holomorphic map of complex analytic spaces, and put $\cal{A} \coleq \pi_\ast \cal{O}_\Y$.
Then the pushforward functor defines an exact equivalence of categories
\begin{eqntag}
\label{1789155196091}
	\pi_\ast : \mathbf{Coh} (\Y) \iso \mathbf{Coh} (\X, \cal{A})
\fullstop
\end{eqntag}
A quasi-inverse $\pi^\star$ is the $\cal{A}$-module inverse image functor:
\begin{eqn}
	\pi^\star \cal{E}
	\coleq
	\cal{O}_\Y
	\otimes_{\pi^{-1} \cal{A}}
	\pi^{-1} \cal{E}
\fullstop{,}
\end{eqn}
where $\pi^{-1} \cal{A} \to \cal{O}_\Y$ is the canonical evaluation homomorphism.
There are natural isomorphisms
\begin{eqn}
	\cal{E}
	\iso
	\pi_\ast \pi^\star \cal{E}
\qqtext{and}
	\pi^\star \pi_\ast \cal{L}
	\iso
	\cal{L}
\fullstop
\end{eqn}
Moreover, for every coherent sheaf $\cal{L}$ on $\Y$ and every $i \geq 0$, there is a canonical isomorphism
\begin{eqn}
	\H^i (\Y, \cal{L})
	\cong
	\H^i (\X, \pi_\ast \cal{L})
\fullstop
\end{eqn}
\end{proposition}

Observe that the functor $\pi^\star$ is the \textit{$\cal{A}$-module} inverse image functor; i.e., the tensor product in its definition is taken over $\pi^{-1} \cal{A}$.
This is different from the usual \textit{$\cal{O}_\X$-module} inverse image functor $\pi^\ast$; i.e., where the tensor product is taken over $\pi^{-1} \cal{O}_\X$.
This subtle difference is crucial and it is what incorporates the prescribed algebra action into the reconstructed sheaf on $\Y$.
We will see in the proof that the stalk $\cal{A}_x$ at any $x \in \X$ decomposes according to the fibre of $\Y$ above $x$.
At any point above $x$, the corresponding factor of $\cal{A}_x$ then recovers the stalk of the $\cal{A}$-module inverse image sheaf, with its full local-ring structure, including any nilpotents.

\begin{proof}[Proof of \autoref{1793046131649}]
The $\cal{A}$-action on $\pi_\ast \cal{L}$ is induced by multiplication of sections of $\cal{L}$ by holomorphic functions on $\Y$.
The finite mapping theorem ensures that $\pi_\ast \cal{L}$ is coherent over $\cal{O}_\X$.

We first record the local description of finite direct image.
Fix any $x \in \X$, and write
\begin{eqn}
	\pi^{-1} (x) = \set{y_1, \ldots, y_k}
\fullstop
\end{eqn}
Choose disjoint open neighbourhoods $\V_i$ of the preimages $y_i$.
Since $\pi$ is proper, every sufficiently small neighbourhood $\U$ of $x$ has its full inverse image contained in the union of these neighbourhoods; i.e., $\pi^{-1} (\U) \subset \V_1 \cup \cdots \cup \V_k$.
Passing to germs therefore gives canonical identifications
\begin{eqn}
	\cal{A}_x
	\cong
	\prod_{i=1}^k \cal{O}_{\Y,y_i}
\qqtext{and}
	(\pi_\ast \cal{L})_x
	\cong
	\bigoplus_{i=1}^k \cal{L}_{y_i}
\fullstop
\end{eqn}
The first identification is an algebra isomorphism, and the second isomorphism is compatible with the componentwise action of this product algebra presentation of $\cal{A}_x$.
The same description applies to arbitrary sheaves of modules, so finite direct image is exact.

For completeness, these identifications can be seen directly.
A finite collection of germs $f_i$ at the points $y_i$ can be represented by sections $\tilde{f}_i$ of $\cal{L}$ on the disjoint neighbourhoods $\V_i$.
Hence, after shrinking $\U$, the tuple $(\tilde{f}_1, \ldots, \tilde{f}_k)$ gives a section of $\cal{L}$ over $\pi^{-1} (\U)$, and therefore a section of $\pi_\ast \cal{L}$ over $\U$.
Conversely, a section whose germs vanish at all the $y_i$ vanishes over $\pi^{-1} (\U)$ after possibly shrinking $\U$ one more time.
This proves both surjectivity and injectivity of the germ maps.

We next check that $\pi^\star \cal{E}$ is coherent.
Locally on $\X$, the module $\cal{E}$ admits a finite presentation
\begin{eqn}
	\cal{A}^{\oplus a}
	\too
	\cal{A}^{\oplus b}
	\too
	\cal{E}
	\too 0
\fullstop
\end{eqn}
Indeed, finitely many $\cal{O}_\X$-generators also generate $\cal{E}$ over $\cal{A}$.
The kernel of the resulting surjection from $\cal{A}^{\oplus b}$ is coherent over $\cal{O}_\X$, and finitely many generators of this kernel give the required presentation.
Applying the definition of the inverse image functor $\pi^\star$ shows that $\pi^\star \cal{E}$ is coherent because it acquires the finite presentation
\begin{eqn}
	\cal{O}_\Y^{\oplus a}
	\too
	\cal{O}_\Y^{\oplus b}
	\too
	\pi^\star \cal{E}
	\too 0
\fullstop
\end{eqn}

Let $\epsilon_i \in \cal{A}_x$ be the idempotent equal to $1$ in the factor $\cal{O}_{\Y,y_i}$ and $0$ in every other factor.
The evaluation homomorphism $\cal{A}_x \to \cal{O}_{\Y,y_i}$ is the projection onto that factor, so
\begin{eqn}
	(\pi^\star \cal{E})_{y_i}
	=
	\cal{O}_{\Y,y_i}
	\otimes_{\cal{A}_x}
	\cal{E}_x
	\cong
	\epsilon_i \cal{E}_x
\fullstop
\end{eqn}
Consequently, the natural map $\cal{E} \to \pi_\ast \pi^\star \cal{E}$ is, on the stalk at $x$, nothing but the isomorphism 
\begin{eqn}
	\cal{E}_x
	\too
	\Oplus_{i=1}^k \epsilon_i \cal{E}_x
\qtext{given by}
	e \mapstoo (\epsilon_1 e, \ldots, \epsilon_k e)
\fullstop
\end{eqn}
Its inverse sends $(e_1, \ldots, e_k)$ to $e_1 + \cdots + e_k$.

Likewise, the natural evaluation map $\pi^\star \pi_\ast \cal{L} \to \cal{L}$ is an isomorphism because its stalk at $y_i$ is the canonical identification
\begin{eqn}
	\cal{O}_{\Y,y_i}
	\otimes_{\cal{A}_x}
	\Oplus_{j=1}^k \cal{L}_{y_j}
	\cong
	\cal{L}_{y_i}
\fullstop
\end{eqn}
These maps are natural in $\cal{E}$ and $\cal{L}$, and therefore exhibit $\pi^\star$ and $\pi_\ast$ as quasi-inverse functors.
Since finite direct image is exact, this is an exact equivalence.

Finally, exactness of finite direct image gives $R^i \pi_\ast \cal{L} = 0$ for every $i > 0$.
The asserted cohomology isomorphisms follow from the Leray spectral sequence.
\end{proof}

\begin{corollary}[local freeness under the correspondence]
\label{1793046145022}
A coherent sheaf $\cal{L}$ on $\Y$ is locally free of rank $r$ if and only if $\pi_\ast \cal{L}$ is locally free of rank $r$ over $\cal{A}$.
Furthermore, if $\pi$ is finite flat of degree $n$, then $\pi_\ast \cal{L}$ is a vector bundle of rank $nr$ over $\X$ whenever $\cal{L}$ is locally free of rank $r$ over $\Y$.
In particular, line bundles on $\Y$ correspond precisely to $\cal{A}$-modules locally free of rank one, and their pushforwards carry the natural spectral structure.
\end{corollary}

\begin{proof}
Suppose $\cal{L}$ is locally free of rank $r$ on $\Y$.
Fix $x \in \X$ and trivialise $\cal{L}$ on disjoint neighbourhoods of the finitely many points of $\pi^{-1} (x)$.
Properness of $\pi$ allows us to choose a neighbourhood $\U$ of $x$ whose inverse image is contained in their union.
The chosen trivialisations then give
\begin{eqn}
	\cal{L}|_{\pi^{-1} (\U)}
	\cong
	\cal{O}_{\pi^{-1} (\U)}^{\oplus r}
\fullstop
\end{eqn}
Pushing forward gives an $\cal{A}_\U$-linear isomorphism $(\pi_\ast \cal{L})_\U \cong \cal{A}_\U^{\oplus r}$.
Conversely, applying $\pi^\star$ to a local $\cal{A}$-module trivialisation gives a local $\cal{O}_\Y$-module trivialisation of the same rank.
Finally, if $\pi$ is finite flat of degree $n$, then $\cal{A}$ is locally free of rank $n$ over $\cal{O}_\X$.
\end{proof}

\subsection{Spectral Realisations}
\label{1790944507299}

The Finite Sheaf-Module Correspondence (\autoref{1793046131649}) from the previous subsection gives spectral structures a direct geometric interpretation.
The algebra $\cal{A}$ determines a finite flat cover of $\X$, and the vector bundle $\cal{E}$ which is a locally free rank-one $\cal{A}$-module determines a line bundle on that cover.
Together, these give a presentation of the underlying vector bundle as the pushforward of a line bundle.

\begin{definition}[spectral realisation]
\label{1793121565784}
A \dfn{spectral realisation} of a vector bundle $\cal{E}$ of rank $n$ on $\X$ is a triple $(\Y \tinyoverset{\pi}{\to} \X, \cal{L}, \iota)$ consisting of a finite flat map $\pi : \Y \to \X$ of degree $n$, a line bundle $\cal{L}$ on $\Y$, and a vector bundle isomorphism $\iota : \pi_\ast \cal{L} \iso \cal{E}$.
We call $\pi : \Y \to \X$ the \dfn{spectral cover}, $\Y$ the \dfn{spectral curve}, and $\cal{L}$ the \dfn{spectral line bundle} of the realisation.
\end{definition}

Given a spectral structure $(\cal{A}, \rho)$ on $\cal{E}$, take
\begin{eqn}
	\pi : \Y \coleq \Spec_\X (\cal{A}) \too \X
\qqtext{and}
	\cal{L} \coleq \pi^\star \cal{E}
	= \cal{O}_\Y \otimes_{\pi^{-1} \cal{A}} \pi^{-1} \cal{E}
\fullstop
\end{eqn}
Since $\cal{E}$ is locally free of rank one over $\cal{A}$, \autoref{1793046145022} shows that $\cal{L}$ is a line bundle on $\Y$.
The Finite Sheaf-Module Correspondence (\autoref{1793046131649}) then furnishes this data with a canonical $\cal{A}$-linear isomorphism $\iota : \pi_\ast \cal{L} \iso \cal{E}$.
We call the resulting triple $(\Y \tinyoverset{\pi}{\to} \X, \cal{L}, \iota)$ the \dfn{spectral realisation} of $(\cal{A}, \rho)$, and also refer to the cover $\pi : \Y \to \X$ itself as the \dfn{spectral realisation} of the spectral algebra $\cal{A}$.
Since the pushforward functor $\pi_\ast$ of \eqref{1789155196091} is fully faithful, the pair $(\cal{L}, \iota)$ is unique up to a unique isomorphism compatible with $\iota$.

Thus, the Finite Sheaf-Module Correspondence assigns to every spectral structure on $\cal{E}$ its spectral realisation:
\begin{eqntag}
\label{1789155313329}
	(\cal{A}, \rho)
	\mapstoo
	(\Y \tinyoverset{\pi}{\to} \X, \cal{L}, \iota)
\fullstop
\end{eqntag}
This correspondence can be extended to an equivalence of groupoids.
Namely, consider the groupoid whose objects are spectral realisations of $\cal{E}$ in the sense of \autoref{1793121565784}.
An isomorphism $(\Y \tinyoverset{\pi}{\to} \X, \cal{L}, \iota) \iso (\Y' \tinyoverset{\pi'}{\to} \X, \cal{L}', \iota')$ is a pair $(f, \lambda)$ consisting of an isomorphism $f : \Y \iso \Y'$ over $\X$ (i.e., satisfying $\pi' \circ f = \pi$) and an $\cal{O}_\Y$-linear isomorphism $\lambda : \cal{L} \iso f^\ast \cal{L}'$ such that the induced vector bundle isomorphism $\pi_\ast \lambda : \pi_\ast \cal{L} \iso \pi_\ast f^\ast \cal{L}' \cong \pi'_\ast \cal{L}'$ intertwines the identifications $\iota$ and $\iota'$:
\begin{eqn}
\begin{tikzcd}
	\pi_\ast \cal{L}
		\ar[r, "\pi_\ast \lambda"]
		\ar[d, "\iota"']
& 	\pi'_\ast \cal{L}'
		\ar[d, "\iota'"]	
\\
	\cal{E}
		\ar[r, equal]
&
	\cal{E}
\fullstop
\end{tikzcd}
\end{eqn}
Let us denote this groupoid of spectral realisations over $\cal{E}$ by
\begin{eqn}
	\mathbf{SpReal} (\cal{E})
	\quad:\quad
	\set{ (\Y \tinyoverset{\pi}{\to} \X, \cal{L}, \iota) 
		~\left|~ 
		\substack{ \DISPLAY{ \pi \quad \text{finite flat of degree $n$}} 
			\\ \DISPLAY{\cal{L} \in \mathbf{Pic} (\Y) \qtext{and} \iota : \pi_\ast \cal{L} \iso \cal{E} } 
			} \right. }
\fullstop
\end{eqn}
Here, $\mathbf{Pic} (\Y)$ is the groupoid of invertible $\cal{O}_\Y$-modules, called the \dfn{Picard groupoid} of $\Y$.
Unfixing the vector bundle $\cal{E}$ yields the groupoid of all degree-$n$ spectral realisations over $\X$:
\begin{eqn}
	\mathbf{SpReal}_n (\X)
	\quad:\quad
	\set{ (\Y \tinyoverset{\pi}{\to} \X, \cal{L}) 
		~\left|~ 
		\substack{ \DISPLAY{ \pi ~ \text{finite flat of degree $n$}} 
			\\ \DISPLAY{\cal{L} \in \mathbf{Pic} (\Y) } 
			} \right. }
\fullstop
\end{eqn}
Isomorphisms in this groupoid are similarly pairs $(f, \lambda) : (\Y \tinyoverset{\pi}{\to} \X, \cal{L}) \to (\Y' \tinyoverset{\pi'}{\to} \X, \cal{L}')$ such that $\pi' \circ f = \pi$ and $\pi_\ast \lambda : \pi_\ast \cal{L} \iso \pi_\ast f^\ast \cal{L}' \cong \pi'_\ast \cal{L}'$ but without the compatibility condition involving $\iota$ and $\iota'$.

\begin{proposition}[equivalence of spectral structures and realisations]
\label{1789153885363}
For every vector bundle $\cal{E}$, there is a canonical equivalence of groupoids
\begin{eqn}
	\mathbf{SpStr} (\cal{E}) 
	\cong
	\mathbf{SpReal} (\cal{E})
\qtext{sending}
	(\cal{A}, \rho) \mapstoo (\Y \tinyoverset{\pi}{\to} \X, \cal{L}, \iota)
\fullstop
\end{eqn}
More generally, there is a canonical equivalence of groupoids
\begin{eqn}
	\mathbf{SpStr}_n (\X)
	\cong
	\mathbf{SpReal}_n (\X)
\qtext{sending}
	(\cal{E}, \cal{A}, \rho) \mapstoo (\Y \tinyoverset{\pi}{\to} \X, \cal{L})
\fullstop
\end{eqn}
\end{proposition}

\begin{proof}
This equivalence is induced by the correspondence \eqref{1789155313329}.
On objects, the inverse correspondence sends $(\Y \tinyoverset{\pi}{\to} \X, \cal{L}, \iota) \mapsto (\cal{A}, \rho)$ by setting $\cal{A} \coleq \pi_\ast \cal{O}_\Y$ and transporting the natural multiplication action $\rho : \pi_\ast \cal{O}_\Y \to \cal{End} (\pi_\ast \cal{L}) \cong \cal{End} (\cal{E})$ through the isomorphism $\iota$.
See \autoref{1789283249898} for more details.

On morphisms, suppose we are given an isomorphism $\alpha : \cal{A} \iso \cal{A}'$ compatible with the two actions on $\cal{E}$; i.e., $\rho = \rho' \circ \alpha$.
Contravariance of the relative Spec gives an isomorphism $f \coleq \Spec_\X (\alpha^{-1}) : \Y \iso \Y'$.
Under the equivalence \eqref{1789155196091}, the pullback bundle $f^\ast \cal{L}'$ corresponds to $\cal{E}$ endowed with the $\cal{A}$-module structure $\rho = \rho' \circ \alpha$.
Hence, the identity of $\cal{E}$ is $\cal{A}$-linear and determines a unique isomorphism $\lambda : \cal{L} \iso f^\ast \cal{L}'$.
This construction is functorial and gives the required correspondence on morphisms.
\end{proof}

\begin{example}[split spectral structures: realisation]
\label{1791460489043}
The spectral realisation of the split spectral structure of \autoref{1791434155987} is given by the disjoint union of $n$ copies of $\X$ itself,
\begin{eqn}
	\Y = \Spec_\X (\cal{A}) = \X_1 \sqcup \cdots \sqcup \X_n
\qtext{where}
	\X_i \cong \X
\fullstop
\end{eqn}
The corresponding line bundle $\cal{L}$ on $\Y$ is determined by its restriction to each copy $\X_i$, this restriction being the line bundle $\cal{L}_i$.
\end{example}

\begin{example}[cyclic spectral structures: realisation]
\label{1791460744565}
The cyclic spectral structures of rank $2$ from \autoref{1791437029550} are realised by double covers of $\X$.
Namely, $\pi : \Y \to \X$ is the double cover locally given by the equation $w^2 = s$, branched over the zeros of $s$.
If we take $\cal{L} = \cal{O}_\Y$, then we obtain $\cal{E} = \pi_\ast \cal{O}_\Y = \cal{O}_\X \oplus \cal{M}^{-1}$ with its canonical $\cal{A}$-module structure.
\end{example}

To complete the picture, let us introduce the groupoid of finite flat degree-$n$ covers $\pi : \Y \to \X$ where an isomorphism $(\pi : \Y \to \X) \iso (\pi' : \Y' \to \X)$ is a holomorphic isomorphism $f : \Y \iso \Y'$ satisfying $\pi' \circ f = \pi$:
\begin{eqn}
	\mathbf{Cov}_n (\X)
	\quad:\quad
	\set{ \pi : \Y \to \X \quad \text{finite flat of degree $n$} \Big. }
\fullstop
\end{eqn}
There is a forgetful functor that drops the invertible sheaf $\cal{L}$:
\begin{eqn}
	G : \mathbf{SpReal}_n (\X) \too \mathbf{Cov}_n (\X)
\qtext{sending}
	(\Y \tinyoverset{\pi}{\to} \X, \cal{L}) \mapstoo (\Y \tinyoverset{\pi}{\to} \X)
\fullstop
\end{eqn}
Its strict fibre $G^\inv (\pi)$ above any $\pi \in \mathbf{Cov}_n (\X)$ is precisely the Picard groupoid $\mathbf{Pic} (\Y)$ and this is canonically equivalent to its essential fibre $G^{[\inv]} (\pi)$:
\begin{eqntag}
\label{1791390038522}
	G^{[\inv]} (\pi) \cong G^\inv (\pi) = \mathbf{Pic} (\Y)
\fullstop
\end{eqntag}
Combined with the canonical equivalence in the above proposition, we obtain the diagram
\begin{eqntag}
\label{1790940189747}
\begin{tikzcd}
&	\mathbf{Pic} (\Y)
		\ar[r, hook]
&	\mathbf{SpStr}_n (\X)
		\ar[dl, "G"']
		\ar[dr, "F"]
&	\mathbf{SpStr} (\cal{E})
		\ar[l, hook']
\\
	(\Y \tinyoverset{\pi}{\to} \X)
		\ar[r, hook]
&	\mathbf{Cov}_n (\X)
&&	\mathbf{Bun}_n (\X)
&	\cal{E}
		\ar[l, hook']
\end{tikzcd}
\end{eqntag}
Thus, the groupoid of spectral structures $\mathbf{SpStr}_n (\X)$ interpolates between the geometry of finite flat covers of $\X$ of degree $n$ and the geometry of vector bundles over $\X$ of rank $n$.
In a subsequent paper, we will fix a smooth connected cover $\pi: \Y \to \X$ and study the pushforward functor $\pi_\ast : \mathbf{Pic} (\Y) \to \mathbf{Bun}_n (\X)$, including the geometry of its image.

The functors from $\mathbf{Pic} (\Y)$ and $\mathbf{SpStr} (\cal{E})$ into $\mathbf{SpStr}_n (\X)$ in the top row of diagram \eqref{1790940189747} are faithful but in general not full.
The reason is that, although these functors include the objects of each fibre into the ambient groupoid $\mathbf{SpStr}_n (\X)$, the latter may have more morphisms between those objects.
In the fibre over $\cal{E}$, morphisms must induce the identity on $\cal{E}$; in the fibre over $\pi$, they must induce the identity on $\pi$.
For essential fibres, these conditions are interpreted through the chosen identifications with the fixed objects.
But the ambient groupoid $\mathbf{SpStr}_n (\X)$ allows additional morphisms inducing nonidentity automorphisms of the underlying bundle or cover.
Forgetting the fixed identification preserves distinct morphisms, which makes the functor faithful, but can introduce additional morphisms between its image objects, obstructing fullness.

\section{Smooth Spectral Structures}
\label{1791480675169}

In general, the analytic space $\Y$ realising a spectral structure is not a Riemann surface because it may be nonreduced, singular, or disconnected.

\begin{definition}[smooth, connected spectral structure]
\label{1789157644037}
A spectral structure $\cal{A}$ on $\cal{E}$ is called \dfn{smooth} if its spectral realisation $\Y = \Spec_\X (\cal{A})$ is smooth; i.e., a possibly disconnected compact Riemann surface.
It is called \dfn{connected} if $\Y$ is connected.
\end{definition}

The full subgroupoids of $\mathbf{SpStr} (\cal{E})$ and $\mathbf{SpStr}_n (\X)$ formed by smooth and/or connected spectral structures will be denoted by adding the superscripts `$\rm{sm}$' and/or `$\rm{con}$'; e.g.,
\begin{eqn}
	\mathbf{SpStr}^{\rm{sm},\rm{con}} (\cal{E}) \subset \mathbf{SpStr}^\rm{sm} (\cal{E}) \subset \mathbf{SpStr} (\cal{E})
\fullstop
\end{eqn}
For example, the spectral realisation of the split spectral structure on $\cal{E} = \cal{L}_1 \oplus \cdots \oplus \cal{L}_n$ (cf. \autoref{1791434155987}) is the split cover $\Y = \X_1 \sqcup \cdots \sqcup \X_n$ with $\X_i \cong \X$ (cf. \autoref{1791460489043}).
It is always smooth because $\X$ is smooth, but always disconnected unless $n = 1$.

The next result says that every vector bundle admits a smooth and connected spectral structure.
In particular, this demonstrates that the groupoids above are never empty.
In fact, any vector bundle admits a spectral structure with a particularly nice geometric realisation as follows.

\begin{definition}[simply branched cover]
\label{1791408863410}
A smooth cover $\pi : \Y \to \X$ is called \dfn{simply branched} if all of its ramification points are simple, and no two distinct ramification points are contained in the same fibre.
\end{definition}

\begin{proposition}[existence of smooth connected spectral structures]
\label{1790777298686}
Every vector bundle $\cal{E}$ on $\X$ admits a smooth and connected spectral structure; in particular, $\mathbf{SpStr}^{\rm{sm},\rm{con}} (\cal{E})$ is nonempty for every $\cal{E}$.
That is to say, if $\cal{E}$ is any vector bundle on $\X$ of rank $n$, then there exists a smooth and connected cover $\pi : \Y \to \X$ of degree $n$, a line bundle $\cal{L}$ on $\Y$, and an isomorphism $\iota : \pi_\ast \cal{L} \iso \cal{E}$.
Furthermore, this cover $\pi$ can always be chosen with simple branching.
\end{proposition}

We give two proofs of this proposition.
The first, presented immediately below, is based on a projective-geometric construction suggested to the author by Anand Patel.
The second, given in \ref{1791465731919}, uses the Higgs bundle methods.
The two constructions are genuinely different: as we show in \autoref{1793629098096} below, in rank at least three the first proof method produces spectral structures that cannot be induced by \textit{any} Higgs field.

\begin{proof}[Proof (based on a suggestion by Anand Patel)]
\ResetClaims
For $n = 1$, take $\pi = \id_\X$ and $\cal{L} = \cal{E}$.
Henceforth, suppose $n \geq 2$.
Since $\X$ is a compact Riemann surface, we may regard it as a smooth projective algebraic curve and, by GAGA, regard $\cal{E}$ as an algebraic vector bundle.
We carry out the construction algebraically and then pass to the associated analytic spaces.

Choose a line bundle $\cal{N}$ on $\X$ of sufficiently large degree that $\cal{F} \coleq \cal{E} \otimes \cal{N}$ is globally generated and has positive degree.
Write
\begin{eqn}
	\cal{M} \coleq \det \cal{F}
\qqtext{and}
	\V \coleq \CC^{n+1}
\fullstop
\end{eqn}
In particular, $\deg \cal{M} = \deg \cal{F} > 0$.
We first construct the required cover and a line bundle whose pushforward is $\cal{F}$.
The twist by $\cal{N}$ will be removed at the end.

\begin{claim}
\label{1793632382349}
The bundle $\cal{F}$ is generated by $n+1$ global sections.
\end{claim}

To see this, put
\begin{eqn}
	\W \coleq \H^0 (\X, \cal{F})^{\oplus (n+1)}
\fullstop
\end{eqn}
A tuple $\boldsymbol{s} = (s_0,\ldots,s_n) \in \W$ determines a homomorphism
\begin{eqn}
	\V \otimes \cal{O}_\X \too \cal{F}
\qqtext{sending}
	e_j \mapsto s_j
\fullstop{,}
\end{eqn}
where $e_0, \ldots, e_n$ is the standard basis of $\V$.
This homomorphism is surjective precisely when the values $s_0(x),\ldots,s_n(x)$ span $\cal{F}(x)$ for every $x \in \X$.

Fix $x \in \X$.
Since $\cal{F}$ is globally generated, evaluation at $x$ gives a surjective linear map
\begin{eqn}
	\op{ev}_x : \W \too \Hom \big(\V, \cal{F}(x) \big)
\fullstop
\end{eqn}
The nonsurjective maps in the target form a closed locus of codimension two.
Indeed, every such map has image contained in a hyperplane of the $n$-dimensional vector space $\cal{F}(x)$.
The hyperplanes form a projective space of dimension $n-1$, and the maps from $\V$ into a fixed hyperplane form a vector space of dimension $(n+1)(n-1)$.
A map of rank $n-1$ determines its image hyperplane uniquely, so the nonsurjective locus has dimension $(n-1) + (n+1)(n-1) = n(n+1) - 2$.
Its inverse image under $\op{ev}_x$ therefore has codimension two in $\W$.

Consider the closed incidence locus
\begin{eqn}
	\J \coleq
	\set{
		(x,\boldsymbol{s}) \in \X \times \W
		~\Big|~
		s_0(x), \ldots, s_n(x)
		\text{ do not span } \cal{F}(x)
	}
\fullstop
\end{eqn}
Every fibre of $\J \to \X$ has dimension $\dim \W - 2$, so $\dim \J = \dim \W - 1$.
The projection $\X \times \W \to \W$ is proper because $\X$ is projective.
Consequently, the image of $\J$ is closed, and its dimension is at most $\dim \W - 1$.
It is therefore a proper subset of $\W$.
Any tuple outside this image generates $\cal{F}$ everywhere.
This proves \autoref{1793632382349}.

Choose such a tuple.
The resulting surjection $\V \otimes \cal{O}_\X \to \cal{F}$ has a line bundle as its kernel.
Taking determinants, and using the standard basis of $\V$ to trivialise $\det \V$, identifies this kernel with $\cal{M}^{-1}$.
Thus we obtain an exact sequence
\begin{eqn}
	0
	\too
	\cal{M}^{-1}
	\tinyoverset{\alpha}{\too}
	\V \otimes \cal{O}_\X
	\too
	\cal{F}
	\too
	0
\fullstop
\end{eqn}

Let us use the quotient convention for projective bundles: a point of $\PP (\V)$ represents a one-dimensional quotient $\V \to \ell$.
Such a quotient factors through $\cal{F}(x)$ precisely when it vanishes on the image of $\alpha_x$.
This gives a closed embedding $\PP (\cal{F}) \inj \X \times \PP (\V)$ whose fibre over $x \in \X$ is the hyperplane
\begin{eqn}
	\H_x \coleq \PP (\cal{F}(x)) \subset \PP (\V) \simeq \PP^n
\fullstop
\end{eqn}

More precisely, let $p$ and $q$ be the projections from $\X \times \PP (\V)$.
Pulling back $\alpha$ and composing with the tautological quotient gives a homomorphism
\begin{eqn}
	p^\ast \cal{M}^{-1}
	\too \V \otimes \cal{O}_{\X \times \PP (\V)}
	\too q^\ast \cal{O}_{\PP (\V)} (1)
\fullstop
\end{eqn}
Equivalently, this is a section of $p^\ast \cal{M} \otimes q^\ast \cal{O}_{\PP (\V)} (1)$ and its zero divisor is $\PP (\cal{F})$.
For example, if $\cal{M}^{-1}$ is trivialised over an open subset $\U \subset \X$ and $\alpha(1) = a_0 e_0 + \cdots + a_n e_n$, then $\PP (\cal{F})$ is given over $\U$ by the equation $a_0(x) Z_0 + \cdots + a_n(x) Z_n = 0$.
The functions $a_0, \ldots, a_n$ have no common zero because $\alpha$ is the inclusion of a line subbundle.

Choose homogeneous coordinates on $\PP (\V)$ and fix the rational normal embedding
\begin{eqn}
	\nu : \PP^1 \inj \PP (\V)
\qqtext{sending}
	[s : t] \mapsto [s^n : s^{n-1}t : \cdots : t^n]
\fullstop
\end{eqn}
Write $\G \coleq \PGL (\V)$.
For each $g \in \G$, put $\nu_g \coleq g \circ \nu$ and form the fibre product
\begin{eqn}
	\Y_g \coleq
	\PP (\cal{F}) \times_{\PP (\V)} \PP^1
\fullstop
\end{eqn}
Here the maps to $\PP (\V)$ are the second projection and $\nu_g$, respectively.
Thus $\Y_g$ is naturally a closed subscheme of $\S \coleq \X \times \PP^1$.
Let $p_\X$ and $p_{\PP^1}$ denote the two projections from $\S$.

Since $\nu_g^\ast \cal{O}_{\PP (\V)} (1) \simeq \cal{O}_{\PP^1} (n)$, pulling back the equation of $\PP (\cal{F})$ shows that $\Y_g$ is the zero divisor of a section of $p_\X^\ast \cal{M} \otimes p_{\PP^1}^\ast \cal{O}_{\PP^1} (n)$.
Let $\pi_g : \Y_g \to \X$ be the restriction of $p_\X$.
For each $x \in \X$, the scheme-theoretic fibre is the divisor $\nu_g^{-1} (\H_x)$ on $\PP^1$.
A rational normal curve is not contained in any hyperplane, so the equation defining this divisor is a nonzero section of $\cal{O}_{\PP^1} (n)$.
Consequently, every fibre has length exactly $n$, counting multiplicities.

This also gives a direct local verification that $\pi_g$ is finite and flat of degree $n$.
Fix $x_0 \in \X$ and choose a point of $\PP^1$ outside the finite divisor $\nu_g^{-1} (\H_{x_0})$.
Use this point as the point at infinity of an affine coordinate $t$.
After shrinking to an open neighbourhood $\U$ of $x_0$ and trivialising $\cal{M}$ over $\U$, the equation of $\Y_g$ has a nonvanishing leading coefficient.
Dividing by this coefficient writes the equation as
\begin{eqn}
	P(x,t) = t^n + b_1(x)t^{n-1} + \cdots + b_n(x) = 0
\fullstop
\end{eqn}
Therefore,
\begin{eqn}
	(\pi_g)_\ast \cal{O}_{\Y_g}|_\U
	\cong \cal{O}_\U[t] \big/ \inner{P(x,t)}
\end{eqn}
is locally free with basis $1, t, \ldots, t^{n-1}$.
Thus $\pi_g$ is finite and flat of degree $n$, and $\Y_g$ is a relative effective Cartier divisor in $\S$.

We next choose $g$ so that $\Y_g$ is smooth.
The projective bundle $\PP (\cal{F})$ is smooth of dimension $n$, and $\PP^1$ is smooth of dimension one.
The connected algebraic group $\G$ acts transitively on $\PP (\V)$.
Kleiman's Transversality Theorem \cite[Theorem~2]{Kle74}, applied to the two maps to $\PP (\V)$, therefore shows that, for $g$ in a nonempty Zariski-open subset of $\G$, the fibre product $\Y_g$ is smooth of dimension
\begin{eqn}
	\dim \PP (\cal{F}) + \dim \PP^1 - \dim \PP (\V)
	= n + 1 - n
	= 1
\fullstop
\end{eqn}
The fibre product is nonempty because its projection to $\X$ has degree $n$.

We must also arrange simple branching.
For each $d \geq 0$, let
\begin{eqn}
	\P_d \coleq |\cal{O}_{\PP^1} (d)| \simeq \PP^d
\end{eqn}
be the complete linear system of effective divisors of degree $d$ on $\PP^1$.
Let $\B_n \subset \P_n$ be the locus of divisors having either a point of multiplicity at least three or two distinct points of multiplicity at least two.

\begin{claim}
\label{1793632343816}
$\B_n$ is closed and has codimension at least two.	
\end{claim}

If $n \geq 3$, the divisors containing a triple point form the image of the proper map
\begin{eqn}
	\PP^1 \times \P_{n-3}
	\too \P_n
\qqtext{sending}
	(p, \D) \mapsto 3p + \D
\fullstop
\end{eqn}
Its image is closed and has dimension at most $1 + (n-3) = n-2$.
If $n \geq 4$, consider also the proper map
\begin{eqn}
	\PP^1 \times \PP^1 \times \P_{n-4}
	\too \P_n
\qqtext{sending}
	(p,q,\D) \mapsto 2p + 2q + \D
\fullstop
\end{eqn}
Its image is closed and has dimension at most $2 + (n-4) = n-2$.
Allowing $p = q$ adds only divisors containing a quadruple point, which already belong to the triple-point locus.
Thus the union of these two images is precisely $\B_n$.
The corresponding image is omitted when $n < 3$ or $n < 4$, respectively.
In particular, $\B_2$ is empty.
This proves \autoref{1793632343816}.

Restriction along $\nu$ identifies linear forms on $\PP (\V)$ with homogeneous polynomials of degree $n$ on $\PP^1$.
Indeed, the coordinate sections $Z_0,\ldots,Z_n$ restrict to the basis $s^n, s^{n-1}t, \ldots, t^n$ of $\H^0 (\PP^1, \cal{O}_{\PP^1} (n))$.
Consequently, restriction gives a projective-linear isomorphism from the space of hyperplanes in $\PP (\V)$ to $\P_n$.

For fixed $x \in \X$, consider the map
\begin{eqn}
	\psi_x : \G \too \P_n
\qqtext{sending}
	g \mapsto [\nu_g^{-1} (\H_x)]
\fullstop
\end{eqn}
Since $\nu_g^{-1} (\H_x) = \nu^{-1} (g^{-1} \H_x)$, this is the map $g \mapsto g^{-1} \H_x$ followed by the projective-linear isomorphism just described.
The action of $\G$ on the space of hyperplanes is transitive, and its orbit maps are smooth.
Therefore, $\psi_x$ is smooth and surjective.
It follows that $\psi_x^{-1} (\B_n)$ has codimension at least two in $\G$ whenever it is nonempty.

Now consider the closed incidence locus
\begin{eqn}
	\I \coleq
	\set{
		(x,g) \in \X \times \G
		~\Big|~
		\nu_g^{-1} (\H_x) \in \B_n
	}
\fullstop
\end{eqn}
Every fibre of $\I \to \X$ has dimension at most $\dim \G - 2$, so $\dim \I \leq \dim \G - 1$.
Since $\X$ is projective, the image of $\I$ in $\G$ is closed.
The dimension estimate shows that it is proper.
Thus, for $g$ in a nonempty Zariski-open subset of $\G$, every fibre of $\pi_g$ consists either of $n$ distinct simple points or of one double point and $n-2$ distinct simple points.

The group $\G$ is irreducible, so this open subset meets the nonempty open subset on which $\Y_g$ is smooth.
Choose $g$ satisfying both conditions.
We check that $\pi_g$ is then simply branched.

At a simple point of a fibre, the derivative of a local defining equation with respect to the fibre coordinate is nonzero.
The implicit function theorem therefore expresses the curve locally as a graph over $\X$, so $\pi_g$ is unramified there.
Now let $y$ be the double point of a fibre over $x \in \X$.
Choose local coordinates $z$ on $\X$ centred at $x$ and $t$ on $\PP^1$ centred at the second coordinate of $y$.
If $\Y_g$ is locally given by $f(z,t)=0$, then the multiplicity-two condition says that
\begin{eqn}
	f(0,0) = f_t(0,0) = 0
\qqtext{and}
	f_{tt} (0,0) \neq 0
\fullstop
\end{eqn}
Since $\Y_g$ is smooth, its defining equation has nonzero differential.
Thus $f_z(0,0) \neq 0$.
The implicit function theorem writes the curve locally as $z = \varphi(t)$, and differentiating $f(\varphi(t),t)=0$ gives
\begin{eqn}
	\varphi'(0) = 0
\qqtext{and}
	\varphi''(0)
	= - \frac{f_{tt}(0,0)}{f_z(0,0)}
	\neq 0
\fullstop
\end{eqn}
Hence $\varphi(t) = t^2 u(t)$ for a holomorphic function $u$ with $u(0) \neq 0$.
After choosing a local square root of $u$ and replacing $t$ by $t \sqrt{u(t)}$, the map has the local form $z = t^2$.
Thus every ramification point is simple.
Since every fibre contains at most one double point, distinct ramification points have distinct images.
Therefore, $\pi_g$ is simply branched.

Fix this choice of $g$, and write $\Y \coleq \Y_g$ and $\pi \coleq \pi_g$.
We next prove that $\Y$ is connected.
Let $i : \Y \inj \S$ be the inclusion.
The equation of $\Y$ gives
\begin{eqn}
	\cal{I}_\Y
	\simeq p_\X^\ast \cal{M}^{-1} \otimes p_{\PP^1}^\ast \cal{O}_{\PP^1} (-n)
\end{eqn}
and an exact sequence
\begin{eqntag}
\label{1793633087093}
	0
	\too
	\cal{I}_\Y
	\too
	\cal{O}_\S
	\too
	i_\ast \cal{O}_\Y
	\too
	0
\fullstop
\end{eqntag}
Since $\deg \cal{M} > 0$, we have $\H^0 (\X, \cal{M}^{-1}) = 0$.
Also, $\H^0 (\PP^1, \cal{O}_{\PP^1} (-n))=0$.
The Künneth formula therefore gives $\H^0 (\S,\cal{I}_\Y)=0$ and
\begin{eqn}
	\H^1 (\S,\cal{I}_\Y) \simeq
	\big( \H^1 (\X,\cal{M}^{-1}) \otimes \H^0 (\PP^1,\cal{O}_{\PP^1} (-n)) \big)
	\oplus
	\big( \H^0 (\X,\cal{M}^{-1}) \otimes \H^1 (\PP^1,\cal{O}_{\PP^1} (-n)) \big)
	= 0
\fullstop
\end{eqn}
Taking cohomology of the ideal sequence \eqref{1793633087093}, we obtain
\begin{eqn}
	\H^0 (\Y, \cal{O}_\Y)
	\cong \H^0 (\S, \cal{O}_\S)
	\cong \CC
\fullstop
\end{eqn}
But a holomorphic function on each connected component of a smooth compact curve is constant.
Thus the dimension of $\H^0 (\Y, \cal{O}_\Y)$ equals the number of connected components of $\Y$.
Consequently, $\Y$ is connected, and hence is a compact Riemann surface.

It remains to construct a line bundle on $\Y$ whose pushforward is $\cal{E}$.
Choose an isomorphism
\begin{eqn}
	\theta_g :
	\nu_g^\ast \cal{O}_{\PP (\V)} (1)
	\iso \cal{O}_{\PP^1} (n),
\end{eqn}
and use this choice to write the equation and ideal sequence of $\Y$.
Restriction along $\nu_g$, followed by $\theta_g$, gives an isomorphism
\begin{eqn}
	r_g :
	\V
	\cong \H^0 (\PP (\V),\cal{O}_{\PP (\V)} (1))
	\too \H^0 (\PP^1,\cal{O}_{\PP^1} (n))
\fullstop
\end{eqn}
This is an isomorphism because $\nu_g$ is a rational normal embedding.

Put $\cal{K} \coleq p_{\PP^1}^\ast \cal{O}_{\PP^1} (n)|_\Y$.
Tensoring the ideal sequence \eqref{1793633087093} by $p_{\PP^1}^\ast \cal{O}_{\PP^1} (n)$ gives
\begin{eqntag}
\label{1793633349618}
	0 \too p_\X^\ast \cal{M}^{-1}
	\too p_{\PP^1}^\ast \cal{O}_{\PP^1} (n)
	\too i_\ast \cal{K}
	\too 0
\fullstop
\end{eqntag}
We have
\begin{eqn}
	(p_\X)_\ast p_\X^\ast \cal{M}^{-1}
	\cong \cal{M}^{-1}
\qtext{and}
	R^1 (p_\X)_\ast p_\X^\ast \cal{M}^{-1}
	\cong
	\cal{M}^{-1} \otimes
	\H^1 (\PP^1,\cal{O}_{\PP^1})
	=0
\fullstop{,}
\end{eqn}
and
\begin{eqn}
	(p_\X)_\ast
	p_{\PP^1}^\ast \cal{O}_{\PP^1} (n)
	\cong
	\H^0 (\PP^1,\cal{O}_{\PP^1} (n))
	\otimes \cal{O}_\X
\fullstop
\end{eqn}
Since $p_\X \circ i=\pi$, pushing forward the sequence \eqref{1793633349618} therefore gives an exact sequence
\begin{eqn}
	0 \too \cal{M}^{-1}
	\too
	\H^0 (\PP^1,\cal{O}_{\PP^1} (n))
	\otimes \cal{O}_\X
	\too \pi_\ast \cal{K}
	\too 0
\fullstop
\end{eqn}

Under the identification $r_g$, the first map in this sequence is precisely $\alpha$.
Indeed, on an open subset where $\alpha (1) = a_0 e_0 + \cdots + a_n e_n$, the equation of $\Y$ is the section $a_0 (x) r_g (e_0) + \cdots + a_n (x) r_g (e_n)$ of $\cal{O}_{\PP^1} (n)$.
The map from $p_\X^\ast \cal{M}^{-1}$ in the sequence \eqref{1793633349618} sends the local generator to this section.
After pushing forward, it therefore sends
\begin{eqn}
	1 \mapstoo
	\sum_{j=0}^n a_j\,r_g(e_j)
	= r_g \big(\alpha(1)\big)
\fullstop
\end{eqn}
Thus, after identifying the middle bundles, the pushforward sequence is the original presentation
\begin{eqn}
	0 \too \cal{M}^{-1}
	\xrightarrow{\alpha} \V \otimes \cal{O}_\X
	\too \cal{F}
	\too 0
\fullstop
\end{eqn}
Comparing cokernels gives $\pi_\ast \cal{K} \simeq \cal{F} = \cal{E} \otimes \cal{N}$.

Finally, define $\cal{L} \coleq \cal{K} \otimes \pi^\ast \cal{N}^{-1}$.
The projection formula gives
\begin{eqn}
	\pi_\ast \cal{L}
	\cong
	(\pi_\ast \cal{K}) \otimes \cal{N}^{-1}
	\simeq
	(\cal{E} \otimes \cal{N}) \otimes \cal{N}^{-1}
	\cong \cal{E}
\fullstop
\end{eqn}
Choosing the resulting isomorphism $\iota : \pi_\ast \cal{L} \iso \cal{E}$ gives the required spectral realisation.
The cover $\pi : \Y \to \X$ is smooth, connected, and simply branched, as required.
\end{proof}

\begin{remark}[allowing the cover to vary]
\label{1789399711381}
It is important to distinguish what \autoref{1790777298686} \textit{says} from what it does \textit{not} say.
\autoref{1790777298686} says that every vector bundle $\cal{E}$ admits \textit{some} smooth connected spectral structure, but the associated cover $\Y$ may depend on $\cal{E}$.
However, it does \textit{not} say that every vector bundle can be obtained by pushing forward a line bundle along a given cover $\pi : \Y \to \X$.
This is typically \textit{not} true.
In fact, this locus is --- except in completely trivial situations --- always a proper subspace in every bundle degree component which may or may not be dense.
Also, the connectedness requirement matters here: for example, thanks to the Birkhoff-Grothendieck Theorem, the fully split cover of $\PP^1$ captures every vector bundle of the correct rank as the pushforward of a line bundle; see \autoref{1791434155987} and \autoref{1791460489043}.
\end{remark}

\subsection{Recognition}
\label{1792146117234}

We now discuss some practical ways to characterise smoothness of spectral structures.
First, there is a simple and intrinsic algebraic criterion on the defining algebra $\cal{A}$ to guarantee smoothness.
Recall that $\cal{A}$ is a normal $\cal{O}_\X$-algebra if, for every point $x \in \X$, the stalk $\cal{A}_x$ is a normal ring.
This means $\cal{A}_x$ is reduced and integrally closed in its total ring of fractions $\op{Frac} (\cal{A}_x)$.
Equivalently, it means that every element of $\op{Frac} (\cal{A}_x)$ that is integral over $\cal{A}_x$ already lies in $\cal{A}_x$.

\begin{lemma}[smooth as normal]
\label{1789157786307}
$\Y = \Spec_\X (\cal{A})$ is a possibly disconnected smooth compact complex curve if and only if the algebra $\cal{A}$ is normal.
\end{lemma}

\begin{proof}
The first fact is that normality of $\cal{A}$ is equivalent to normality of $\Y$.
This is because, for any $x \in \X$, there is canonical isomorphism of stalks at $x$:
\begin{eqntag}
\label{1791566069221}
	\cal{A}_x = (\pi_\ast \cal{O}_\Y)_x \cong \prod_{y \in \pi^{-1} (x)} \cal{O}_{\Y, y}
\fullstop
\end{eqntag}
A finite product of rings is normal whenever each factor is normal.
Consequently, $\cal{A}_x$ is normal if and only if $\cal{O}_{\Y, y}$ is normal for every $y \in \pi^{-1} (x)$.

The second fact is the special property in dimension one that normal complex curves are necessarily smooth.
To see this, take $y \in \Y$ and suppose $\Y$ is normal.
Since $\Y$ is one-dimensional, the local ring $\cal{O}_{\Y, y}$ is a one-dimensional Noetherian normal local domain, hence a discrete valuation ring, and hence a regular local ring of dimension one \cite[Proposition~9.2]{AtiyahMacdonald1969}.
Therefore, every point of $\Y$ is regular.
But for a reduced complex analytic curve, regularity is equivalent to smoothness.
So $\Y$ is a complex manifold of dimension one; i.e., a possibly disconnected Riemann surface.

The third and final fact concerns compactness.
Since $\pi : \Y \to \X$ is finite, it is a proper map.
But since $\X$ itself is compact, then so is $\Y = \pi^\inv (\X)$.
\end{proof}

In \autoref{1789283249898} we saw that a line bundle on a finite cover gives its pushforward a canonical action of the cover algebra by multiplication.
When the cover is smooth, this construction has a useful converse: a faithful action of $\pi_\ast \cal{O}_\Y$ on a vector bundle whose rank equals the degree of the cover automatically makes that bundle locally free of rank one over the cover algebra.
In addition, if the cover is connected, any such action is automatically faithful.
For related recognition criteria in terms of Cartan subalgebra bundles, including a parabolic version for ramified covers, see \cite[Theorems~2.4 and~3.1]{AuffarthBiswas2018}.

\begin{proposition}[recognition from algebra actions]
\label{1791408088098}
Let $\pi : \Y \to \X$ be a smooth compact finite cover of degree $n$, possibly disconnected.
Put $\cal{A} = \pi_\ast \cal{O}_\Y$ and let $\cal{E}$ be a vector bundle of rank $n$ on $\X$.
Then every injective $\cal{O}_\X$-algebra homomorphism $\rho : \cal{A} \to \cal{End} (\cal{E})$ equips $\cal{E}$ with the spectral structure $(\cal{A}, \rho)$; i.e., $\cal{E}$ is automatically locally free of rank one over $\cal{A}$.
In addition, if $\Y$ is connected, then any $\cal{O}_\X$-algebra homomorphism $\rho : \cal{A} \to \cal{End} (\cal{E})$ is automatically injective.

Consequently, any injective $\cal{O}_\X$-algebra homomorphism $\rho : \cal{A} \to \cal{End} (\cal{E})$ induces a line bundle $\cal{L}$ on $\Y$ together with an $\cal{A}$-linear isomorphism $\pi_\ast \cal{L} \iso \cal{E}$, and its degree is
\begin{eqntag}
\label{1791410001614}
	\deg \cal{L} = \deg \cal{E} + \tfrac{1}{2} \deg \R
\fullstop{,}
\end{eqntag}
where $\R = \Ram (\pi)$ is the ramification divisor of $\pi$.
If $\Y$ is disconnected, then $\deg \cal{L}$ is the sum of the degrees of $\cal{L}$ on each connected component of $\Y$.
\end{proposition}

\begin{proof}
Let $\rho : \cal{A} \to \cal{End} (\cal{E})$ be any $\cal{O}_\X$-algebra homomorphism.
The equivalence of categories between coherent sheaves on $\Y$ and coherent $\cal{A}$-modules on $\X$ associates to the $\cal{A}$-module $\cal{E}$ a coherent sheaf $\cal{L}$ on $\Y$, together with an $\cal{A}$-linear isomorphism $\iota : \pi_\ast \cal{L} \iso \cal{E}$.
Therefore, $\cal{E}$ is locally free of rank one over $\cal{A}$ if and only if $\cal{L}$ is a line bundle.

First, suppose that $\Y$ is connected.
Let $\KK \coleq \cal{M} (\X)$ and $\LL \coleq \cal{M} (\Y)$ be the fields of global meromorphic functions on the two curves.
The pullback by $\pi$ allows us to regard this as a finite field extension $\KK \subset \LL$ of degree $n$.
Let $\EE \coleq \cal{M} (\X, \cal{E})$ denote the vector space of global meromorphic sections of $\cal{E}$.
Since $\cal{E}$ has rank $n$, this is a $\KK$-vector space of dimension $n$.

At the same time, a global meromorphic section of $\cal{A} = \pi_\ast \cal{O}_\Y$ is nothing but a global meromorphic function on $\Y$, so $\cal{M} (\X, \cal{A}) = \LL$.
The algebra action $\rho$ makes $\EE$ into a vector space over $\LL$.
Incidentally, this turns the algebra action $\rho$ into a homomorphism from the field $\LL$ to the algebra $\End_\KK (\EE)$.
Such a homomorphism is necessarily injective because a nonzero kernel would be a nonzero proper ideal of which a field has none.
This forces the map $\rho : \cal{A} \to \cal{End} (\cal{E})$ to be injective.
For otherwise, if it has a kernel, then this kernel must be a torsion subsheaf of $\cal{A}$ because it vanishes after localisation.
But $\cal{A}$ is torsion-free, so this proves the second stated assertion.

The dimension of $\EE$ over $\LL$ satisfies
\begin{eqn}
	n
	= \dim_\KK \EE
	= [\LL : \KK] \dim_\LL \EE
	= n \dim_\LL \EE
\fullstop
\end{eqn}
Thus, $\dim_\LL \EE = 1$.
But by the finite module correspondence, the dimension of $\EE$ over $\LL$ is precisely the rank of $\cal{L}$ on $\Y$.
Thus, $\cal{L}$ has rank one on the connected curve $\Y$.

Now, let us allow $\Y$ to be disconnected, and assume that $\rho : \cal{A} \to \cal{End} (\cal{E})$ is injective.
Write $\Y = \Y_1 \sqcup \cdots \sqcup \Y_m$, let $n_i$ be the degree of the restriction $\pi : \Y_i \to \X$, and put $\LL_i \coleq \cal{M} (\Y_i)$.
Then $[\LL_i : \KK] = n_i$ and $n_1 + \cdots + n_m = n$.

In this case, a global meromorphic function on $\Y$ is an $m$-tuple of global meromorphic functions, one on each connected component.
Therefore, the algebra of global meromorphic sections of $\cal{A}$ decomposes into the product $\cal{M} (\X, \cal{A}) = \LL_1 \times \cdots \times \LL_m$.
Its action on $\EE$ is still faithful because $\rho$ is injective.
Let $\epsilon_i$ be the idempotent equal to $1$ on the $i$-th factor $\LL_i$ and $0$ on the others.
These idempotents determine a $\KK$-vector space decomposition $\EE = \EE_1 \oplus \cdots \oplus \EE_m$ where $\EE_i \coleq \epsilon_i \EE$.
Additionally, each $\EE_i$ is a vector space over $\LL_i$ of positive dimension $d_i$ because otherwise the idempotent $\epsilon_i$ would be acting as zero on all of $\EE$, contrary to the faithfulness of the action.
And we have the identity
\begin{eqn}
	n 
	= \dim_\KK \EE
	= \sum_{i=1}^m [\LL_i : \KK] \dim_{\LL_i} \EE_i
	= \sum_{i=1}^m n_i d_i
\fullstop
\end{eqn}
But since $n_1 + \cdots + n_m = n$, this implies
\begin{eqn}
	\sum_{i=1}^m n_i (d_i - 1) = 0
\fullstop
\end{eqn}
Each summand is nonnegative and every $n_i$ is positive, so this forces $d_i = 1$ for each $i$.
But, under the finite module correspondence, the dimension of $\EE_i$ over $\LL_i$ is precisely the rank of $\cal{L}$ on $\Y_i$.
Thus, $\cal{L}$ has rank one on each connected component of the curve $\Y$.

We claim now that $\cal{L}$ is torsion-free.
Let $\cal{F} \subset \cal{L}$ be the torsion subsheaf.
It is necessarily supported at finitely many points of $\Y$.
Since the finite pushforward $\pi_\ast$ is exact, $\iota (\pi_\ast \cal{F})$ is a subsheaf of $\cal{E}$ supported at finitely many points.
But a vector bundle on a smooth curve $\X$ is torsion-free, so $\pi_\ast \cal{F} = 0$.
Since the finite pushforward $\pi_\ast$ is also faithful, it follows that $\cal{F} = 0$.
Thus, $\cal{L}$ is a torsion-free rank-one sheaf, hence a line bundle.

Finally, identity \eqref{1791410001614} is a standard application of the Riemann-Roch and Riemann-Hurwitz formulas \cite[Ch.~IV, \S\S1--2]{Hartshorne1977}.
\end{proof}

\enlargethispage{20pt}
When $\Y$ is disconnected, the explicit requirement in \autoref{1791408088098} that $\rho$ be injective is indispensable.
For example, if $\Y = \X \sqcup \X$, then $\cal{A} = \cal{O}_\X \oplus \cal{O}_\X$ acts on any rank-two bundle by $(a,b) \mapsto a \id$.
This action is not faithful and does not make the bundle free of rank one over $\cal{A}$, hence does not determine a spectral structure.

\subsection{Local But Not Global Generators}
\label{1792146177359}

A spectral structure is an action by an algebra $\cal{A}$ that makes the vector bundle $\cal{E}$ locally free of rank one, which means $\cal{E}$ is always locally generated by a single element as an $\cal{A}$-module.
But this does not mean that $\cal{A}$ itself can in general be generated by a single element, even locally.
However, for \textit{smooth} spectral structures, this turns out to be true: the next result shows that every smooth spectral structure is locally generated by a single regular endomorphism of the vector bundle, and this description applies even near branch points.
We establish this local description and explain why such generators need not exist globally.

\begin{proposition}[local regular generators]
\label{1791453151472}
Every smooth spectral structure is locally generated by a regular endomorphism.
More precisely, let $(\cal{A}, \rho)$ be a smooth spectral structure on a vector bundle $\cal{E}$ of rank $n$.
Then every point $x \in \X$ has an open neighbourhood $\U$ such that there is a section $a \in \H^0 (\U, \cal{A})$ with the property that
\begin{eqntag}
\label{1791569347339}
	\cal{A}_\U = \cal{O}_\U [a]
\qqtext{and}
	\rho (\cal{A}_\U) = \cal{O}_\U [\Phi]
\qtext{where}
	\Phi \coleq \rho (a)
\fullstop
\end{eqntag}
The local endomorphism $\Phi$ is regular and it induces an $\cal{O}_\U$-algebra isomorphism
\begin{eqn}
	\cal{A}_\U \cong \cal{O}_\U [t] \big/ \inner{ P_a (t) }
\fullstop
\end{eqn}
where $P_a (t)$ is the characteristic polynomial of $\mu_a : \cal{A}_\U \to \cal{A}_\U$, the multiplication by $a$.
\end{proposition}

To be clear, $\cal{O}_\U [a]$ in \eqref{1791569347339} is, by definition, the $\cal{O}_\U$-subalgebra of $\cal{A}_\U$ generated by the element $a$.
Similarly, $\cal{O}_\U [\Phi]$ is the $\cal{O}_\U$-subalgebra of $\cal{End} (\cal{E})$ generated by $\Phi$.

\begin{proof}
Let $\pi : \Y \to \X$ be the smooth spectral realisation of $\cal{A}$ so that $\cal{A} = \pi_\ast \cal{O}_\Y$.
Let $\pi^{-1} (x) = \set{ y_1, \ldots, y_k }$ and let $r_i$ be the ramification index of $\pi$ at $y_i$, so that $r_1 + \cdots + r_k = n$.
Choose a coordinate $w$ around $x$ and coordinate $z_i$ around each $y_i$ such that $w = z_i^{r_i}$ is the local equation of $\pi$ near $y_i$.
Consider the fibre algebra
\begin{eqn}
	\cal{A} (x)
	=
	\cal{A}_x / \frak{m}_x \cal{A}_x
	\cong \prod_{y \in \pi^{-1} (x)} \cal{O}_{\Y, y} / \frak{m}_x \cal{O}_{\Y, y}	
	\cong \prod_{i = 1}^k \CC [z_i] \big/ \inner{ z_i^{r_i} }
\fullstop
\end{eqn}
Choose any $k$ pairwise distinct complex numbers $c_1, \ldots, c_k \in \CC$ and consider the element $h \coleq \big( z_1 + c_1, \ldots, z_k + c_k \big) \in \cal{A} (x)$ in this product algebra.
We claim that $1, h, h^2, \ldots, h^{n-1}$ forms a basis of the vector space $\cal{A}(x)$.
Indeed, the homomorphism $\CC [t] \to \cal{A} (x)$ sending $t \mapsto h$ has a kernel:
\begin{eqn}
	\K 
	\coleq
	\prod_{i=1}^k \inner[\big]{ (t - c_i)^{r_i} }
	=
	\inner[\Big]{ (t-c_1)^{r_1} \cdots (t-c_k)^{r_k} }
\fullstop
\end{eqn}
This equality follows from the Chinese Remainder Theorem because the ideals $\inner{ (t - c_i)^{r_i} }$ are pairwise comaximal.
So we get an isomorphism $\CC [t] / \K \to \cal{A}(x)$ given by sending $t \mapsto h$.
Meanwhile, the polynomial $(t-c_1)^{r_1} \cdots (t-c_k)^{r_k}$ is monic of degree $n$, so the quotient vector space has a basis $1, t, t^2, \ldots, t^{n-1}$.
Then the isomorphism induces the basis $1, h, h^2, \ldots, h^{n-1}$ on $\cal{A}(x)$.

Next, we lift $h$ to a section $a$ of $\cal{A}$ on an open neighbourhood $\U \subset \X$ of $x$.
Then the induced holomorphic bundle map
\begin{eqn}
	\cal{O}^{\oplus n}_\U \too \cal{A}_\U
\qtext{sending}
	(f_0, \ldots, f_{n-1}) \mapstoo f_0 + f_1 a + f_2 a^2 + \cdots + f_{n-1} a^{n-1}
\end{eqn}
is an isomorphism on the fibre at $x$, and hence in a neighbourhood of $x$.
After possibly shrinking $\U$, we therefore obtain an isomorphism $\cal{A}_\U \cong \cal{O}_\U [a]$.

Next, we choose a local $\cal{A}$-module generator $e$ of $\cal{E}$ over $\U$; we may have to shrink $\U$ to be able to do so.
Then the sections $e, \Phi (e), \Phi^2 (e), \ldots, \Phi^{n-1} (e)$ where $\Phi \coleq \rho (a)$ are the images of the basis $1, a, \ldots, a^{n-1}$ under the map $b \mapsto \rho(b)e$.
Consequently, they form a basis of generators of $\cal{E}_\U$.
Thus, the section $e$ determines a cyclic vector for $\Phi$ in every fibre of $\cal{E}_\U$, so the endomorphism $\Phi$ is regular.

Let $\mu_a : \cal{A}_\U \to \cal{A}_\U$ denote the multiplication by $a$, and let $P_a (t)$ be its characteristic polynomial.
Then, $P_a (a) = 0$.
Indeed, $P_a (\mu_a) = 0$ by the Cayley-Hamilton Theorem.
But $P_a (\mu_a)$ is an endomorphism of $\cal{A}_\U$, which happens to be zero, so we can apply it to $1 \in \cal{A}_\U$, getting $P_a (\mu_a) 1 = 0$.
Expanding the lefthand side gives precisely the identity $P_a (a) = 0$.

Thus, evaluation at $a$ defines an algebra homomorphism
\begin{eqn}
    \cal{O}_\U [t] \big/ \inner{P_a(t)} \too \cal{A}_\U
\qqtext{sending}
    [f(t)] \mapstoo f(a).
\end{eqn}
Since $\cal{A}_\U$ has rank $n$, the polynomial $P_a (t)$ is monic of degree $n$.
So the lefthand side $\cal{O}_\U [t] \big/ \inner{P_a(t)}$ is locally free with basis $1, t, \ldots, t^{n-1}$, and the homomorphism sends it to the basis $1, a, \ldots, a^{n-1}$ of $\cal{A}_\U$.
Therefore, it is actually an isomorphism of $\cal{O}_\U$-algebras.
\end{proof}

Let us emphasise that in general the local algebra generators constructed in \autoref{1791453151472} cannot be chosen globally.
Indeed, if $\Y$ is smooth, compact, and connected, then the only global sections of $\cal{A}$ are scalars:
\begin{eqn}
	\H^0 (\X, \cal{A}) \cong \H^0 (\Y, \cal{O}_\Y) \simeq \CC
\fullstop
\end{eqn}
A scalar section can only generate the scalar subalgebra $\cal{O}_\X \cdot 1 \subset \cal{A}$ which is strictly smaller unless $n = 1$.
Equivalently, every global endomorphism in the image $\rho (\cal{A}) \subset \cal{End} (\cal{E})$ is a scalar multiple of the identity.
Incidentally, this obstruction is exactly why the theory of Higgs bundles involves the twisting line bundle $\Lambda$ --- its purpose is to enable $\cal{End} (\cal{E}) \otimes \Lambda$ to admit global sections.

\subsection{Cyclic Covers}
\label{1792146206910}

In this subsection, we give several explicit nontrivial examples of smooth spectral structures.
Cyclic covers provide a particularly simple but illuminating set of examples.
We examine their smoothness, connectedness, and ramification, and illustrate how the choice of a line bundle on the cover affects its pushforward.

\begin{example}[double covers]
\label{1789286127875}
The cyclic spectral structures of rank $2$ from \autoref{1791437029550} are realised by double covers of $\X$ which can be arranged to be smooth and connected.
Recall that we have a line bundle $\cal{M}$ on $\X$ of degree $k$, that $s \in \H^0 (\X, \cal{M}^{\otimes 2})$, and that $\pi : \Y \to \X$ is the double cover locally given by the equation $w^2 = s$, branched over the zeros of $s$.
See \autoref{1791460744565}.
If all zeros of $s$ are simple, then $\Y$ is a smooth double cover of $\X$ with simple branch points at the $2k$ zeros of $s$; this is easy to verify in local coordinates.

Even if $\Y$ is smooth, it is in general disconnected.
However, if $k > 0$, or if $k = 0$ but the line bundle $\cal{M}$ is nontrivial, then $\Y$ is connected with genus $g_\Y = 2g_\X - 1 + k$, where $g_\X$ is the genus of $\X$.
Indeed, if $k > 0$, then $\H^0 (\X, \cal{M}^{-1}) = 0$, and so
\begin{eqn}
	\H^0 (\Y, \cal{O}_\Y)
	\cong \H^0 (\X, \pi_\ast \cal{O}_\Y)
	\cong \H^0 (\X, \cal{O}_\X) \oplus \H^0 (\X, \cal{M}^{-1})
	\cong \CC
\fullstop
\end{eqn}
On the other hand, if $k = 0$, then the section $s$ is nonvanishing and determines a trivialisation of $\cal{M}^{\otimes 2}$.
If $\cal{M}$ is nontrivial, then $\cal{M}^{-1}$ is also nontrivial of degree zero and has no nonzero sections, so the same calculation proves connectedness.
But if $\cal{M}$ is trivial, then $s$ is a nonzero constant; its two distinct square roots split the cover $\Y$ into two copies of $\X$, so $\Y$ in this case is smooth but disconnected.
\end{example}

\begin{example}[elliptic double cover of $\PP^1$]
\label{1789296853723}
Let us make \autoref{1789286127875} more concrete.
If $\X = \PP^1$, take $\cal{M} = \cal{O} (k)$ with $k \geq 1$, together with a square-free section $s \in \H^0 (\PP^1, \cal{O} (2k))$.
Then $\pi : \Y \to \PP^1$ is a smooth and connected double cover branched at $2k$ points, and $\Y$ has genus $k - 1$.

For example, take the line bundle $\cal{M} = \cal{O} (2)$ together with the square-free quartic $s(z) = z (z-1) (z-2) (z-3)$, where $z$ is an affine coordinate on $\PP^1$.
Then 
\begin{eqn}
	\cal{A} = \cal{O} \oplus \cal{O} (-2)
\end{eqn}
with multiplication determined by the relation $w^2 = s(z)$.
The corresponding double cover $\Y$ is an elliptic curve given explicitly in the affine chart $z$ by 
\begin{eqn}
	\Y ~:\qquad w^2 = z (z-1) (z-2) (z-3)
\fullstop
\end{eqn}
The projection $\pi : \Y \to \PP^1$ is $(z,w) \mapsto z$, which is branched at $z = 0, 1, 2, 3$.

If we take $\cal{L} = \cal{O}_\Y$, then 
\begin{eqn}
	\cal{E} = \pi_\ast \cal{O}_\Y \simeq \cal{O} \oplus \cal{O} (-2)
\fullstop
\end{eqn}
On the other hand, if we take any other line bundle $\cal{L} \not\cong \cal{O}_\Y$ of degree $0$, then
\begin{eqn}
	\cal{E} = \pi_\ast \cal{L} \simeq \cal{O} (-1) \oplus \cal{O} (-1)
\fullstop
\end{eqn}
To see this, write $\cal{E} = \cal{O} (a) \oplus \cal{O} (b)$ with $a \geq b$ by the Birkhoff-Grothendieck Theorem.
Then 
\begin{eqn}
	h^0 (\cal{O}(a)) + h^0 (\cal{O}(b)) = h^0 (\cal{E}) = h^0 (\cal{L})
\fullstop
\end{eqn}
At the same time, we must have $a + b = -2$ because $\deg \cal{E} = \deg \cal{L} - 2 = -2$.
If $\cal{L} = \cal{O}_\Y$, then $h^0 (\cal{L}) = 1$.
But if $a \geq 1$ then $h^0 (\cal{O}(a)) \geq 2$ whilst if $0 > a \geq b$ then both summands have no nonzero global sections.
This leaves the only possibility $a = 0$, $b = -2$.
If $\cal{L}$ is nontrivial, then $h^0 (\cal{L}) = 0$, and this forces $a = b = -1$.

Take any line bundle $\cal{L}$ of degree $1$.
For example, take the line bundle $\cal{L} = \cal{O}_\Y (p)$ for any $p \in \Y$ such as $p = (z,w) = (0,0)$.
Then 
\begin{eqn}
	\cal{E} = \pi_\ast (\cal{O}_\Y (p)) \simeq \cal{O} \oplus \cal{O} (-1)
\fullstop
\end{eqn}
This follows because $\deg \cal{E} = \deg \cal{L} - 2 = -1$ so $a + b = -1$.
At the same time, $h^0 (\cal{L}) = 1$ because $\Y$ has genus $1$ and $\deg \cal{L} = 1$.
Then the only possibility is $a = 0$, $b = -1$.
\end{example}

\begin{example}[double cover of an elliptic curve]
\label{1789296880368}
Here is another concrete realisation of \autoref{1789286127875}.
If $\X$ is an elliptic curve, take $\cal{M}$ to be a line bundle of degree $k > 0$.
Then $\pi : \Y \to \X$ is a smooth double cover of genus $k+1$ branched over the $2k$ zeros of $s$.

For example, if $\X$ is given by the equation 
\begin{eqn}
	y^2 = x^3 - x
\end{eqn}
over an affine chart $x$, let $\infty \in \X$ denote the point at infinity (i.e., the unique preimage of the point at infinity $x = \infty$) and take $\cal{M} = \cal{O}_\X (\infty)$.
Now $x$ is a meromorphic function $x : \X \to \PP^1$ with a double pole at $\infty$, hence a global section of $\cal{M}^{\otimes 2} = \cal{O}_\X (2\infty)$.

Take, for example, the section 
\begin{eqn}
	s = x - 2 \in \H^0_\X (\cal{O}_\X (2\infty))
\fullstop
\end{eqn}
Its zeros are the two points $(x,y) = (2, \pm \sqrt{6})$, and they are simple because $2$ is not a branch value of the map $x : \X \to \PP^1$.
Then we take $\cal{A} = \cal{O}_\X \oplus \cal{O}_\X (- \infty)$ with multiplication determined by the equation $w^2 = s = x - 2$.
It follows that the double cover $\Y \to \X$ is obtained by adjoining the square root of $x - 2$.
In an affine chart, $\Y$ is given by
\begin{eqn}
	y^2 = x^3 - x
\qqtext{and}
	w^2 = x - 2
\fullstop
\end{eqn}
The projection $\pi : \Y \to \X$ sends $(x,y,w) \mapsto (x,y)$ and has degree $2$.
By Riemann-Hurwitz, $\Y$ has genus $2$.
\end{example} 

\begin{example}[cyclic triple covers]
\label{1789293744134}
It is easy to generalise \autoref{1789286127875} in order to produce spectral structures of higher rank.
For example, let $\cal{M}$ be a line bundle of degree $k \geq 0$ on $\X$ and suppose $s \in \H^0 (\X, \cal{M}^{\otimes 3})$ is a section with simple zeros.
Define the algebra 
\begin{eqn}
	\cal{A} = \cal{O}_\X \oplus \cal{M}^{-1} \oplus \cal{M}^{-2}
\end{eqn}
with multiplication determined by the map $s : \cal{M}^{-3} \to \cal{O}_\X$.
It defines a smooth compact triple cover $\pi : \Y \to \X$ locally given by the equation $w^3 = s$.
There are $3k$ zeros of $s$, each with just one point in the fibre above it with ramification index $3$.
If $k > 0$, or $k = 0$ but $\cal{M}$ is nontrivial, then $\Y$ is connected and has genus $3g_\X - 2 + 3k$, where $g_\X$ is the genus of $\X$.

Smoothness of $\Y$ follows from the assumption that the zeros of $s$ are simple and can be checked in local coordinates.
Connectedness of $\Y$ when $k > 0$ follows from the fact that $\cal{M}$ has positive degree, so $\H^0 (\X, \cal{M}^{-i}) = 0$ for all $i \geq 1$, and so
\begin{eqn}
	\H^0 (\Y, \cal{O}_\Y)
	\cong \H^0 (\X, \pi_\ast \cal{O}_\Y)
	\cong \H^0 (\X, \cal{O}_\X) \oplus \H^0 (\X, \cal{M}^{-1}) \oplus \H^0 (\X, \cal{M}^{-2})
	\cong \CC
\fullstop
\end{eqn}

If $k = 0$, then the section $s$ is nonvanishing and determines a trivialisation of $\cal{M}^{\otimes3}$.
If $\cal{M}$ is nontrivial, both $\cal{M}^{-1}$ and $\cal{M}^{-2}$ are nontrivial of degree zero and have no nonzero sections, so the same calculation proves connectedness.
If $\cal{M}$ is trivial, $s$ is a nonzero constant; its three distinct cube roots split the covering into three copies of $\X$.
In this case, $\Y$ is smooth but disconnected.
\end{example}

\section{The Off-Diagonal Bundle}
\label{1792146312894}

Associated with every spectral structure $(\cal{A}, \rho)$ on a vector bundle $\cal{E}$ is a natural quotient of the endomorphism bundle $\cal{End} (\cal{E})$ by the spectral algebra.
Its cohomology provides the vanishing condition defining clean spectral structures and will be one of the main inputs into the smoothness results of \autoref{1792401370132}.
It will also play a crucial role in the deformation theory of bundles obtained as pushforwards of line bundles along a fixed cover which we will study in a subsequent paper.

\begin{definition}[off-diagonal bundle]
\label{1790711090032}
The \dfn{off-diagonal bundle} associated with the spectral structure $(\cal{A}, \rho)$ on $\cal{E}$ is the quotient sheaf
\begin{eqntag}
\label{1790711147900}
	\cal{Q} \coleq \cal{End} (\cal{E}) / \rho(\cal{A})
\fullstop
\end{eqntag}
\end{definition}

The off-diagonal bundle $\cal{Q}$ is an $\cal{O}_\X$-module, but it is not in general a quotient algebra because $\rho(\cal{A}) \subset \cal{End} (\cal{E})$ need not be an ideal.
However, it \textit{is} always a vector bundle.

\begin{lemma}[the off-diagonal exact sequence]
\label{1789332843734}
Let $(\cal{A}, \rho)$ be any spectral structure on a vector bundle $\cal{E}$ of rank $n$.
Then the associated off-diagonal bundle $\cal{Q} = \cal{End} (\cal{E}) / \rho (\cal{A})$ is a vector bundle of rank $n (n-1)$.
In particular, there is a canonical short exact sequence of vector bundles on $\X$:
\begin{eqntag}
\label{1789348440662}
	0 \too \cal{A} \overset{\rho}{\too} \cal{End} (\cal{E}) \too \cal{Q} \too 0
\fullstop
\end{eqntag}
\end{lemma}

\begin{proof}
Take an open subset $\U \subset \X$ over which $\cal{E}$ is free of rank $1$ as an $\cal{A}$-module, and choose an $\cal{A}$-module trivialisation $\phi : \cal{E}_\U \iso \cal{A}_\U$.
Then, under this identification, the spectral structure $\rho$ becomes an inclusion
\begin{eqn}
	i \coleq \Ad_\phi \circ \rho : \cal{A}_\U \inj \cal{End} (\cal{A}_\U)
\qtext{sending}
	a \mapsto m_a
\fullstop{,}
\end{eqn}
where $m_a : \cal{A} \to \cal{A}$ is the multiplication-by-$a$ map in $\cal{A}$; i.e., $m_a (b) = ab$, and where $\Ad_\phi$ is conjugation by $\phi$; i.e., $\Ad_\phi (T) = \phi \circ T \circ \phi^\inv$.
This inclusion map has a left inverse 
\begin{eqn}
	r : \cal{End} (\cal{A}_\U) \to \cal{A}_\U
	\qtext{given by}
	r (T) \coleq T(1)
\fullstop{,}
\end{eqn}
where $1 \in \cal{A}$ is the algebra unit.
Indeed, $(r \circ i) (a) = r(m_a) = m_a (1) = a$.
Therefore, the algebra $\cal{A}_\U$ sits inside $\cal{End} (\cal{A}_\U)$ as a direct summand:
\begin{eqntag}
\label{1790712637915}
	\cal{End} (\cal{A}_\U) \cong \cal{im} (i) \oplus \cal{ker} (r)
\fullstop
\end{eqntag}
Indeed, any endomorphism $T \in \cal{End} (\cal{A}_\U)$ decomposes uniquely as
\begin{eqn}
	T = m_{T(1)} + (T - m_{T(1)})
\end{eqn}
where $m_{T(1)}$ is in the image of $i$ and $T - m_{T(1)}$ is in the kernel of $r$.
Consequently,
\begin{eqntag}
\label{1790712698863}
	\cal{Q}_\U 
		\cong \cal{End} (\cal{A}_\U) \big/ i (\cal{A}_\U)
		\cong \cal{ker} (r)
\fullstop
\end{eqntag}
But this means that $\cal{Q}_\U$ is a direct summand of a locally free $\cal{O}_\U$-module, hence must itself be locally free over $\cal{O}_\U$.
This shows that $\cal{Q}$ is a locally free $\cal{O}_\X$-module.
Finally, the ranks of $\cal{End} (\cal{E})$ and $\cal{A}$ as $\cal{O}_\X$-modules are $n^2$ and $n$, respectively, which yields the number $n^2 - n = n (n-1)$ for the rank of $\cal{Q}$.
\end{proof}

It is worth stressing that the local splittings \eqref{1790712637915} and \eqref{1790712698863} used in this proof are \textit{not} canonical: they depend on a choice of a local $\cal{A}$-module generator of $\cal{E}$.
In general, the short exact sequence \eqref{1789348440662} has no canonical splitting.

\begin{example}[split spectral structures: the off-diagonal bundle]
\label{1791452789065}
For the split spectral structure of $\cal{E} = \cal{L}_1 \oplus \cdots \oplus \cal{L}_n$ (cf. \autoref{1791434155987}), the off-diagonal bundle $\cal{Q} = \cal{End} (\cal{E}) / \rho (\cal{A})$ is canonically identified with the subbundle of $\cal{End} (\cal{E})$ consisting of strictly off-diagonal endomorphisms of $\cal{E}$:
\begin{eqn}
	\cal{Q} = \Oplus_{\substack{i,j=1 \\ i \neq j}}^n \cal{Hom} (\cal{L}_j, \cal{L}_i)
\fullstop
\qedhere
\end{eqn}
\end{example}

For smooth spectral structures, we can give even more detailed information about this vector bundle as follows.

\begin{proposition}[degree and Euler characteristic]
\label{1790716932644}
For any spectral structure $(\cal{A}, \rho)$ on a vector bundle $\cal{E}$ of rank $n$, the degree of the off-diagonal bundle $\cal{Q}$ satisfies
\begin{eqn}
	\deg \cal{Q} = - \deg \cal{A}
\fullstop
\end{eqn}
If the spectral realisation $\pi : \Y \to \X$ is smooth and connected, then the degree of $\cal{Q}$ and hence its Euler characteristic are
\begin{eqntag}
\label{1791394092196}
	\deg \cal{Q} = \tfrac{1}{2} r
\qtext{and}
	\chi (\cal{Q}) = \tfrac{1}{2} r - n (n-1) (g_\X - 1) = (g_\Y - 1) - n^2 (g_\X - 1)
\end{eqntag}
where $g_\Y$ and $g_\X$ are the genera of $\Y$ and $\X$, and $r = \deg \R$ is the degree of the ramification divisor $\R = \sfop{Ram} (\pi)$.
More generally, if $\Y$ is smooth but possibly disconnected with $m$ connected components $\Y_i$ of genus $g_{\Y_i}$, then these formulas remain valid if we interpret the genus of $\Y$ as
\begin{eqntag}
\label{1792494116976}
	g_\Y = 1 + \sum_{i=1}^m (g_{\Y_i} - 1)
\fullstop
\end{eqntag}
\end{proposition}

These numerical invariants of $\cal{Q}$ do not depend on the degree of the vector bundle $\cal{E}$.

\begin{proof}
To calculate the degree of $\cal{Q}$, we appeal to the short exact sequence \eqref{1789348440662}.
First, we find that $\deg \cal{Q} = \deg \cal{End} (\cal{E}) - \deg \cal{A} = - \deg \cal{A}$ because $\deg \cal{End} (\cal{E}) = 0$ for any $\cal{E}$.
In the smooth case, $\cal{A} = \pi_\ast \cal{O}_\Y$.
Using the Riemann-Roch formula together with the identity $\chi (\X, \pi_\ast \cal{L}) = \chi (\Y, \cal{L})$, we can calculate the degree of the pushforward vector bundle $\pi_\ast \cal{O}_\Y$ to be $1 - g_\Y + n(g_\X-1)$.
Using the Riemann-Hurwitz formula, $\deg \cal{A} = -r/2$.
Finally, the Euler characteristic of a vector bundle on a Riemann surface is $\chi (\cal{Q}) = \deg \cal{Q} - \op{rk} (\cal{Q}) (g_\X - 1)$.
From \autoref{1789332843734}, the rank of $\cal{Q}$ is $n(n-1)$.
\end{proof}

\subsection{Fibrewise Description}
\label{1791375997197}

Let us now describe the off-diagonal bundle $\cal{Q}$ even more explicitly fibre by fibre.
Thus, suppose $(\cal{A}, \rho)$ is a smooth spectral structure on a vector bundle $\cal{E}$ of rank $n$, let $\pi : \Y \to \X$ be its spectral realisation, let $\cal{L}$ be such that $\cal{E} \cong \pi_\ast \cal{L}$, and put
\begin{eqn}
	\cal{Q} = \cal{End} (\cal{E}) / \rho (\cal{A})
\fullstop
\end{eqn}


\subsubsection*{Unramified fibres}

First, let $x \in \X$ be a point \textit{not} on the branch locus of $\pi$, so that the fibre above it has $n$ distinct points: $\pi^\inv (x) = \set{y_1, \ldots, y_n} \subset \Y$.
Denote the fibres of $\cal{E}, \cal{L}, \cal{A}, \cal{Q}$ at these points by
\begin{eqn}
	\E \coleq \cal{E}(x)
\fullstop{,}
\qquad
	\L_i \coleq \cal{L} (y_i)
\fullstop{,}
\qquad
	\A \coleq \cal{A} (x)
\fullstop{,}
\qquad
	\Q \coleq \cal{Q} (x)
\fullstop{,}
\end{eqn}
for $i = 1, \ldots, n$.
Then the $n$-dimensional vector space $\E$ has a canonical direct sum decomposition labelled by the points in the fibre above $x$, and correspondingly the vector space $\End (\E)$ of endomorphisms gets a canonical direct sum decomposition labelled by ordered pairs of points in the fibre above $x$:
\begin{eqntag}
\label{1789345317274}
	\E \cong \Oplus_{i=1}^n \L_i
\qtext{and so}
	\End (\E) \cong \Oplus_{i,j = 1}^n \Hom \big(\L_j, \L_i \big)
\fullstop
\end{eqntag}
With respect to this decomposition, the $n$-dimensional algebra $\A$ acts diagonally on $\E$.
This means that the image $\rho (\A)$ of $\rho : \A \to \End (\E)$ consists of the diagonal components, and the quotient $\Q = \End (\E) / \rho (\A)$ is canonically identified with the complement of $\rho (\A)$ which consists of all the off-diagonal components:
\begin{eqntag}
\label{1789345747755}
	\rho (\A) \cong \Oplus_{i=1}^n \End (\L_i)
\qtext{and}
	\Q \cong \Oplus_{\substack{i,j = 1 \\ i \neq j}}^n \Hom \big(\L_j, \L_i \big)
\fullstop
\end{eqntag}
In words, away from branch points, the spectral algebra $\cal{A}$ describes the diagonal part of the endomorphism bundle with respect to the spectral decomposition, whilst $\cal{Q}$ describes its off-diagonal part.

\subsubsection*{Ramified fibres}

If $x \in \X$ is a branch point, assume first that the fibre above it consists of just a single point $y \in \Y$ with ramification index $n \geq 2$.
We can choose local coordinates $z$ centred at $x$ and $w$ centred at $y$ such that $z = w^n$ is the local equation for $\Y$.
This determines an isomorphism
\begin{eqn}
	\A \iso \CC [w] / \inner{w^n}
\end{eqn}
as $\CC$-algebras, and the fibre $\E$ is naturally a module over $\A$.
Upon choosing a local trivialisation of the line bundle $\cal{L}$ near $y$, we get an isomorphism $\E \cong \A$ as $\A$-modules.
Multiplication by $w$ is therefore represented by the nilpotent Jordan block matrix
\begin{eqn}
	N = \mtx{ 0 & 1 \\ & 0 & 1 \\ & & \ddots & \ddots \\ & & & 0 & 1 \\ & & & & 0 }
\end{eqn}
in the basis $w^{n-1}, w^{n-2}, \ldots, w, 1$.
Consequently, the ramified analogues of the diagonal and off-diagonal descriptions \eqref{1789345747755} are the following vector spaces:
\begin{eqntag}
\label{1789353182353}
	\rho (\A) \simeq \CC [N] \coleq \Span \set{ I, N, \ldots, N^{n-1} }
\qtext{and}
	\Q \simeq \End ( \E ) \big/ \CC[N]
\fullstop
\end{eqntag}

More generally, suppose $\pi^\inv (x) = \set{y_1, \ldots, y_k}$ with ramification indices $r_1, \ldots, r_k$ so that $r_1 + \ldots + r_k = n$.
Then the fibre $\A$ decomposes accordingly:
\begin{eqntag}
\label{1790707173313}
	\A \simeq \bigoplus_{i = 1}^k \CC [w_i] \big/ \inner{w_i^{r_i}}
\fullstop
\end{eqntag}
Correspondingly, $\E$ decomposes canonically into the generalised eigenspaces indexed by the points in the fibre above $x$, and the vector space of endomorphisms has the resulting block form:
\vspace{-10pt}
\begin{eqntag}
\label{1790707168072}
	\E = \bigoplus_{i=1}^{k} \E_i
\qtext{with}
	\dim (\E_i)  = r_i
\qqtext{and}
	\End (\E) \cong \Oplus_{i,j = 1}^k \Hom (\E_j, \E_i)
\fullstop
\vspace{-5pt}
\end{eqntag}
Meanwhile, $\rho (\A) \subset \End (\E)$ is block-diagonal with polynomial blocks:
\begin{eqntag}
\label{1790707289213}
	\rho (\A) \simeq \Oplus_{i=1}^k \CC [N_i] = \CC [N_1] \oplus \cdots \oplus \CC [N_k]
\fullstop{,}
\end{eqntag}
with the $i$-th block $\CC [N_i] \subset \End (\E_i)$ generated by the local nilpotent action of $w_i$ on $\E_i$.
Therefore, $\Q$ has contributions of two kinds:
\begin{eqntag}
\label{1789354199774}
	\Q
	\simeq \left( \Oplus_{i \neq j} \Hom (\E_j, \E_i) \right)
	\oplus
	\left( \Oplus_i \End (\E_i) \big/ \CC [N_i] \right)
\fullstop
\end{eqntag}
The terms inside the first big bracket are pure off-block-diagonal terms, the same as in the unramified case \eqref{1789345747755}.
The other terms are the ramified analogues within each Jordan block as in \eqref{1789353182353}.

Most of these decompositions can actually be made canonical without choosing any coordinates or trivialisations.
So we summarise the findings of this subsection in the following statement.

\enlargethispage{20pt}

\begin{proposition}[fibrewise decomposition of the spectral action]
\label{1790706598412}
Let $\pi : \Y \to \X$ be the smooth realisation of a spectral structure $(\cal{A}, \rho)$ on $\cal{E}$, and let $\cal{L}$ be its spectral line bundle.
Fix $x \in \X$, write $\pi^\inv (x) = \set{ y_1, \ldots, y_k }$, and let $r_i$ be the ramification index of $\pi$ at $y_i$ so that $r_1 + \cdots + r_k = n$.
Let $\frak{m}_x$ be the maximal ideal of $\cal{O}_{\X,x}$, and let $\frak{m}_x \cal{O}_{\Y, y_i}$ be the ideal in $\cal{O}_{\Y, y_i}$ generated by the image of $\frak{m}_x$ under the pullback map $\pi^\ast : \cal{O}_{\X,x} \to \cal{O}_{\Y, y_i}$.
Then the fibre $\A \coleq \cal{A} (x)$ has a canonical $\CC$-algebra decomposition into local Artinian algebras $\A_i$, one for each point of $\Y$ above $x$:
\begin{eqntag}
\label{1791397292013}
	\A \cong \Oplus_{i = 1}^k \A_i
\qqtext{where}
	\A_i \coleq \cal{O}_{\Y, y_i} \big/ \frak{m}_x \cal{O}_{\Y, y_i}	
\fullstop
\end{eqntag}
Let $\epsilon_i \in \A$ be the idempotent equal to $1$ on the factor $\A_i$ and $0$ on all other factors.
Then the $r_i$-dimensional vector subspace $\E_i \coleq \rho (\epsilon_i) \E$ of the fibre $\E \coleq \cal{E} (x)$ is a free $\A_i$-module of rank one, and the spectral action $\rho$ determines the canonical vector space decompositions
\begin{eqntag}
\label{1791397660619}
	\E = \Oplus_{i=1}^k \E_i
\qtext{and}
	\End (\E) = \Oplus_{i,j = 1}^k \Hom (\E_j, \E_i)
\fullstop
\end{eqntag}
If $\rho_i = \rho |_{\A_i} : \A_i \to \End (\E_i)$ denotes the induced spectral action, then there is a canonical vector space decomposition of the fibre $\Q \coleq \cal{Q} (x)$:
\begin{eqn}
\label{1791398338256}
	\Q \cong
	\Bigg( \Oplus_{\substack{i,j = 1 \\ i \neq j}}^k \Hom (\E_j, \E_i) \Bigg)
	\oplus
	\Bigg( \Oplus_{i = 1}^k \End (\E_i) / \rho_i (\A_i) \Bigg)
\fullstop
\end{eqn}
If $x$ is not a branch point, then $k = n$, all $r_i = 1$, and every quotient $\End (\E_i) / \rho_i (\A_i)$ is zero.
\end{proposition}

\begin{corollary}[Jordan block description of the fibres]
\label{1791400768742}
Upon choosing local coordinates $z$ near $x$ and $w_i$ near $y_i$, as well as a local generator of $\cal{L}$ near $y_i$, the algebra $\A_i$ and its action can be written as
\begin{eqntag}
\label{1791401354603}
	\A_i \cong \CC [w_i] / \inner{ w_i^{r_i} }
\qqtext{and}
	\rho_i (\A_i) \cong \CC [N_i]
\fullstop{,}
\end{eqntag}
where $N_i$ is a nilpotent Jordan block of size $r_i$.
\end{corollary}

\begin{proof}
These statements follow from the preceding calculations by observing that, in local coordinates, the ideal of $\frak{m}_x \cal{O}_{\Y, y_i}$ is simply $\inner{w_i^{r_i}}$.
The idempotent $\epsilon_i$ is nothing but the identity element $1$ of $\A_i$.
\end{proof}

\vspace{-5pt}
\enlargethispage{20pt}

Our terminology \textit{off-diagonal bundle} is literally accurate over the unramified locus: its fibre consists precisely of the off-diagonal matrix entries with respect to the spectral decomposition of $\E$.
At a branch point, however, the quotient also contains internal directions $\End(\E_i)/\CC [N_i]$ within each ramified block.

\subsection{Clean Spectral Structures}
\label{1792426526441}

The off-diagonal exact sequence \eqref{1789348440662} leads to an effective way to compare the cohomology of the spectral algebra with the cohomology of the endomorphism bundle.
That is, since $\pi$ is finite, we have $\H^1(\Y,\cal{O}_\Y) \cong \H^1(\X,\cal{A})$, and the long exact sequence in cohomology becomes:
\begin{eqn}
    \H^1 (\Y, \cal{O}_\Y)
    \cong
    \H^1 (\X, \cal{A})
    \too
    \H^1 \big( \X, \cal{End} (\cal{E}) \big)
    \too
    \H^1 (\X, \cal{Q})
    \too 0
\fullstop
\end{eqn}
We can observe that $\H^1 (\X, \cal{Q}) = 0$ precisely when every class in $\H^1 \big( \X, \cal{End} (\cal{E}) \big)$ comes from the spectral algebra; i.e., when the following natural map is surjective:
\begin{eqntag}
\label{1792428063597}
	\H^1 (\Y, \cal{O}_\Y)
	\cong
	\H^1 (\X, \cal{A})
	\too
	\H^1 \big(\X, \cal{End} (\cal{E}) \big)
\fullstop
\end{eqntag}
The surjectivity of this map has a useful geometric interpretation via deformation theory.
Indeed, first-order deformations of $\cal{L}$ on the fixed curve $\Y$ are parametrised by the cohomology group $\H^1 (\Y, \cal{O}_\Y)$, whilst first-order deformations of $\cal{E}$ on the fixed curve $\X$ are parametrised by $\H^1 (\X, \cal{End} (\cal{E}))$ \cite[Proposition~2.6 and Theorem~2.7]{Hartshorne2010}.
Under these identifications, the map above is the map on first-order deformations induced by the pushforward functor $\pi_\ast$: an infinitesimal change in the transition functions of $\cal{L}$ induces an infinitesimal change in the transition functions of $\pi_\ast \cal{L}$ which can be regarded as action on $\cal{E}$ through the spectral action $\rho : \cal{A} \to \cal{End} (\cal{E})$.
In this interpretation, the vanishing of $\H^1 (\X, \cal{Q})$ means that every first-order deformation of $\cal{E}$ is induced by a deformation of $\cal{L}$ on the same cover.
We shall pick up this point and study the fixed-cover deformation problem in much more detail in a subsequent paper.

This motivates the following definition, in which we combine the cohomological vanishing condition with smoothness, connectedness, and simple branching.

\begin{definition}[clean spectral structure]
\label{1792404625757}
A spectral structure $(\cal{A}, \rho)$ on a vector bundle $\cal{E}$ over $\X$ is called \dfn{clean} if its spectral realisation $\pi : \Y \to \X$ is smooth, connected, and simply branched, and its off-diagonal bundle $\cal{Q}$ satisfies
\begin{eqn}
    \H^1 (\X, \cal{Q}) = 0
\fullstop
\end{eqn}
\end{definition}

The existence of smooth connected covers $\pi : \Y \to \X$ for which the natural map $\H^1 (\Y, \cal{O}_\Y) \to \H^1 \big(\X, \cal{End} (\cal{E}) \big)$ is surjective is already contained in the theorem of Hirschowitz-Narasimhan \cite[Theorem~6.4]{HirschowitzNarasimhan1994}.
The following proposition follows the same general strategy, but establishes the required cohomological surjectivity by a direct interpolation argument and also explains that the spectral cover can always be arranged to be simply branched.

\begin{proposition}[existence of clean spectral structures]
\label{1792167000673}
Every vector bundle $\cal{E}$ on $\X$ admits a clean spectral structure.
Thus, the groupoid $\mathbf{SpStr}^\rm{clean} (\cal{E})$ is nonempty for all $\cal{E}$.
\end{proposition}

We will prove it in \smartref{1791465731919} using the methods of Higgs bundles.

\subsection{Trace Duality}
\label{1791466758188}

There is a useful intrinsic description of the off-diagonal bundle in terms of the natural \textit{trace pairing} on $\cal{End} (\cal{E})$:
\vspace{-5pt}
\enlargethispage{20pt}
\begin{eqn}
	\tr : \cal{End} (\cal{E}) \otimes \cal{End} (\cal{E}) \too \cal{O}_\X
\qtext{given by}
	\Phi \otimes \Psi \mapstoo \tr (\Phi \circ \Psi)
\fullstop
\end{eqn}
This is a perfect pairing, meaning that it is nondegenerate on every fibre, or equivalently it induces an isomorphism
\vspace{-5pt}
\begin{eqntag}
\label{1791572230994}
	\cal{End} (\cal{E}) \iso \cal{End} (\cal{E})^\vee
\qtext{sending}
	\Phi \mapstoo \tr (\Phi \circ \dummy)
\fullstop
\end{eqntag}
Given a spectral structure $(\cal{A}, \rho)$ on $\cal{E}$, let us define the following \dfn{trace-annihilator} in $\cal{End} (\cal{E})$ of the spectral algebra $\cal{A}$ determined by the action $\rho$:
\begin{eqntag}
\label{1791457710150}
	\cal{A}^\perp 
	\coleq
	\rho(\cal{A})^\perp
	=
	\set{ \Phi \in \cal{End} (\cal{E}) 
		~\Big|~ \text{$\tr \big( \Phi \circ \rho (a) \big) = 0$ for all $a \in \cal{A}$} }
	\subset \cal{End} (\cal{E})
\fullstop
\end{eqntag}

\begin{proposition}[trace-annihilator]
\label{1791453955298}
Let $(\cal{A}, \rho)$ be a spectral structure on a vector bundle $\cal{E}$ and let $\cal{Q} = \cal{End} (\cal{E}) /\rho(\cal{A})$.
Then there is a canonical isomorphism of vector bundles
\begin{eqntag}
\label{1791458038531}
	\cal{A}^\perp \iso \cal{Q}^\vee
\qtext{sending}
	\xi \mapstoo \Big( [\Phi] \mapsto \tr ( \xi \circ \Phi ) \Big)
\fullstop
\end{eqntag}
Consequently, $\cal{A}^\perp$ fits into a canonical short exact sequence of vector bundles on $\X$,
\begin{eqntag}
\label{1791458083502}
	0
	\too
	\cal{A}^\perp
	\too
	\cal{End} (\cal{E})
	\overset{\Tr}{\too}
	\cal{A}^\vee
	\too 0
\qtext{where}
	\Tr : \Phi \mapstoo \Big( a \mapsto \tr \big( \Phi \circ \rho (a) \big) \Big)
\fullstop
\end{eqntag}
Furthermore, if the algebra $\cal{A}$ is locally generated by a section $a$ over $\U$, then $\Phi \coleq \rho (a)$ induces an isomorphism
\begin{eqntag}
\label{1793101070963}
	\cal{Q}_\U \iso \cal{A}^\perp_\U
\qqtext{sending}
	[\Psi] \mapstoo [\Phi, \Psi]
\fullstop
\end{eqntag}
\end{proposition}

\begin{proof}
Recall that $\cal{Q}$ fits into the short exact sequence
\begin{eqntag}
\label{1791572591116}
	0 \too \cal{A} \overset{\rho}{\too} \cal{End} (\cal{E}) \too \cal{Q} \too 0
\fullstop
\end{eqntag}
Dualising, we get the short exact sequence
\begin{eqntag}
\label{1791574384938}
	0
	\too
	\cal{Q}^\vee
	\too
	\cal{End} (\cal{E})^\vee
	\overset{\rho^\vee}{\too}
	\cal{A}^\vee
	\too
	0
\fullstop
\end{eqntag}
Exactness follows from the local splitting of \eqref{1791572591116}.
The map $\rho^\vee$ sends a functional to its restriction along the image of $\rho$; i.e., $\rho^\vee$ sends $\xi \in \cal{End} (\cal{E})^\vee$ to the functional $a \mapsto \xi \big( \rho (a) \big)$.
Under the trace pairing isomorphism \eqref{1791572230994}, this map becomes precisely the map $\Tr : \cal{End} (\cal{E}) \to \cal{A}^\vee$ defined in \eqref{1791458083502}.
Its kernel consists exactly of the endomorphisms $\Phi$ such that $\tr \big( \Phi \circ \rho (a) \big) = 0$ for every local section $a \in \cal{A}$.
But this is exactly the definition of the trace-annihilator $\cal{A}^\perp$ from \eqref{1791457710150}.
Consequently, this short exact sequence gives a canonical isomorphism $\cal{Q}^\vee \cong \cal{A}^\perp$.

For the final assertion, observe that every section of $\rho(\cal{A}_\U)$ commutes with $\Phi$.
Since trace is cyclic, we get:
\begin{eqn}
	\tr \big( [\Phi, \Psi] \rho(b) \big) 
	= \tr \big( \Psi [ \rho(b), \Phi ] \big)
	= 0
\end{eqn}
for any $\Psi \in \cal{End} (\cal{E})$ and any $b \in \cal{A}$.
So the commutator induces the map $\cal{Q}_\U \to \cal{A}^\perp_\U$ as in \eqref{1793101070963}.
By \autoref{1791381265048}, it is an isomorphism on every fibre because the element $a(x)$ generates $\cal{A}(x)$ for every $x \in \U$.
\end{proof}

We remark that although the trace pairing is perfect, the orthogonal complement $\rho(\cal{A})^\perp$ need not be complementary to $\rho (\cal{A})$, because the restriction of the trace pairing to $\rho(\cal{A})$ can degenerate at a branch point.
For example, if $\cal{A} (x) = \CC [t] / \inner{t^2}$, then multiplication by $t$ is represented by the nonzero $2 \times 2$ nilpotent matrix $N$ which satisfies $\tr (N) = \tr (N^2) = 0$.
Thus, $N$ belongs to both $\rho (\cal{A} (x))$ and to its trace annihilator $\rho(\cal{A}(x))^\perp$.

\section{Higgs Bundles and Spectral Structures}
\label{1792146501847}

\subsection{Regular Higgs Bundles}
\label{1792146546398}

Higgs bundles provide an important source of spectral structures.
Indeed, the classical spectral correspondence for regular Higgs bundles naturally equips the underlying vector bundle with a spectral structure \cite{Hitchin1987, MR998478, Markman1994}.
In fact, one of the principal motivations behind the present project was to isolate and study this structure on its own; i.e., retaining the spectral algebra and its action without having to choose a particular Higgs field.

Fix a line bundle $\Lambda$ on $\X$, and let $(\cal{E}, \Phi) = (\cal{E}, \Phi : \cal{E} \to \cal{E} \otimes \Lambda)$ be a $\Lambda$-twisted Higgs bundle of rank $n$.
Recall that the zero locus of its characteristic polynomial $\det (\lambda - \Phi)$ defines a curve $\Y$ inside the total space $\sfop{tot} (\Lambda)$ called the \textit{spectral curve}, and with it a finite flat cover $\pi : \Y \to \X$ of degree $n$.
Take the spectral algebra $\cal{A} \coleq \pi_\ast \cal{O}_\Y$.
As an $\cal{O}_\X$-module,
\begin{eqn}
	\cal{A}
	=
	\pi_\ast \cal{O}_\Y 
	\cong
	\cal{O}_\X \oplus \Lambda^{-1} \oplus \Lambda^{-2} \oplus \cdots \oplus \Lambda^{-(n-1)}
\fullstop
\end{eqn}
The Higgs field $\Phi$ canonically defines an $\cal{O}_\X$-algebra morphism $\rho : \cal{A} \to \cal{End} (\cal{E})$.
To see this, take a local trivialisation of $\Lambda$ so that $\Phi$ is locally represented by an endomorphism $T \in \cal{End} (\cal{E})$.
Then the spectral curve $\Y$ is locally cut out by the characteristic polynomial $\det(tI - T)$, and consequently the spectral algebra has the local presentation 
\begin{eqntag}
\label{1792409375201}
	\cal{A} \simeq \cal{O}_\X [t] / \inner[\big]{\det(tI - T)}
\fullstop
\end{eqntag}
Then $\rho$ is induced by the map $t \mapsto T$, and this is well-defined by the Cayley-Hamilton Theorem.
Thus every Higgs bundle carries a canonical action of its spectral algebra on its underlying vector bundle.

Recall that $\Phi$ is called \textit{regular} if its characteristic and minimal polynomials coincide.
Equivalently and in more elementary terms, at every $x \in \X$ and with respect to any local trivialisation of $\Lambda$, the fibrewise endomorphism $T_x \in \End (\cal{E} (x))$ representing $\Phi_x$ has exactly one Jordan block for each eigenvalue.
Equivalently again, $T_x$ admits a cyclic vector: there exists some $v \in \cal{E} (x)$ such that the vectors $v, T_x v, T^2_x v, \ldots, T^{n-1}_x v$ span the vector space $\cal{E} (x)$.

\begin{proposition}[spectral structures from regular Higgs bundles]
\label{1790763564459}
Let $(\cal{E}, \Phi)$ be a $\Lambda$-twisted Higgs bundle of rank $n$, and let $\pi : \Y \to \X$ be its spectral cover.
Then $\cal{A} \coleq \pi_\ast \cal{O}_\Y$ is a spectral algebra and $\Phi$ induces an action $\rho : \cal{A} \to \cal{End} (\cal{E})$. 
The pair $(\cal{A}, \rho)$ is a spectral structure on $\cal{E}$ if and only if $\Phi$ is regular.

In this case, the locally free rank-one $\cal{A}$-module $\cal{E}$ is the pushforward $\cal{E} \cong \pi_\ast \cal{L}$ of the spectral line bundle, $\rho$ is the natural action of $\pi_\ast \cal{O}_\Y$ on $\pi_\ast \cal{L}$ induced by the $\cal{O}_\Y$-multiplication on $\cal{L}$, and the spectral algebra admits a canonical global $\Lambda$-twisted generator $a \in \H^0(\X, \cal{A} \otimes \Lambda)$ satisfying 
\begin{eqn}
	(\rho \otimes \id_\Lambda) (a) = \Phi
\fullstop
\end{eqn}
\end{proposition}

To clarify, if $\Phi = T \otimes e$ in a local trivialisation $e$ of $\Lambda$, then $a = [t] \otimes e$, where $[t]$ is the local generator in \eqref{1792409375201}.

\begin{proof}
Suppose $\Phi$ is regular, so that every $T_x$ as above admits a cyclic vector $v_x \in \cal{E} (x)$.
Then the linear map $\cal{A} (x) \to \cal{E} (x)$ sending $f(t) \mapsto f(T_x) v_x$ is an isomorphism.
Choose a local holomorphic section $v$ of $\cal{E}$ extending $v_x$.
The induced $\cal{O}_\X$-linear map $\cal{A} \to \cal{E}$ is an isomorphism at $x$ and therefore remains an isomorphism after shrinking the neighbourhood because both bundles have rank $n$.
Therefore, $\cal{E}$ is locally free of rank one over $\cal{A}$, implying by definition that $(\cal{A}, \rho)$ is a spectral structure on $\cal{E}$.
Conversely, if $\cal{E}$ is locally free of rank one over $\cal{A}$, then a local $\cal{A}$-module generator of $\cal{E}$ determines a cyclic vector for every $T_x$.
Finally, under the correspondence explained in \ref{1790944507299}, this rank-one $\cal{A}$-module corresponds to a line bundle $\cal{L}$ on $\Y$ with $\cal{E} \cong \pi_\ast \cal{L}$.
\end{proof}

By passing from regular Higgs bundles to spectral structures we forget the choice of a twisted generator of the spectral algebra.
The following proposition describes exactly what we must supply in order to reverse this passage for a given spectral structure; compare \cite[\S17]{DonagiGaitsgory2002}.

\begin{proposition}[global twisted generators]
\label{1792167000670}
Let $(\cal{A}, \rho)$ be a spectral structure on a vector bundle $\cal{E}$ of rank $n$, let $\pi : \Y \to \X$ be its spectral cover, and fix a line bundle $\Lambda$ on $\X$.
The following conditions are equivalent:
\begin{enumerate}
\item The spectral structure $(\cal{A}, \rho)$ is isomorphic to the spectral structure induced by some regular $\Lambda$-twisted Higgs field $\Phi$ on $\cal{E}$.
\label{1792410292789}
\item There exists a section $a \in \H^0 (\X, \cal{A} \otimes \Lambda)$ whose representative in every local trivialisation of $\Lambda$ generates $\cal{A}$ as an $\cal{O}_\X$-algebra.
\label{1792410296505}
\item There exists a closed holomorphic embedding $i : \Y \hookrightarrow \tot (\Lambda)$ over $\X$.
\label{1792410299907}
\end{enumerate}
In fact, the sets of Higgs fields $\Phi$, global sections $a$, and holomorphic embeddings $i$ in these three conditions are in canonical bijection:
\begin{eqn}
	\set[\Big]{\textup{Higgs fields $\Phi$ as in \autoref*{1792410292789}}}
	\cong
	\set[\Big]{\textup{sections $a$ as in \autoref*{1792410296505}}}
	\cong
	\set[\Big]{\textup{embeddings $i$ as in \autoref*{1792410299907}}}
\fullstop
\end{eqn}
Precisely, the Higgs field associated with a section $a$ is $\Phi = (\rho \otimes \id_\Lambda) (a)$.
\end{proposition}

\begin{proof}
First, suppose $a \in \H^0 (\X, \cal{A} \otimes \Lambda)$ satisfies condition \autoref*{1792410296505}, and put $\Phi \coleq (\rho \otimes \id_\Lambda) (a)$.
Take a trivialising open neighbourhood $\U \subset \X$ for $\Lambda$, choose a generator $e \in \Lambda (\U)$, and write $a_\U = b \otimes e$.
Then $\Phi_\U = T \otimes e$ where $T \coleq \rho (b)$, and $\cal{A}_\U = \cal{O}_\U [b]$.
Since $\cal{E}$ is locally free of rank one over $\cal{A}$, a local $\cal{A}$-module generator of $\cal{E}$ is a cyclic vector for $T$ in every fibre.
So $\Phi$ is regular.
Moreover, the induced spectral action has image $\cal{O}_\U [T] = \rho (\cal{A}_\U)$, so its spectral structure is isomorphic to $(\cal{A}, \rho)$.

Conversely, suppose a regular $\Lambda$-twisted Higgs field $\Phi$ induces a spectral structure isomorphic to $(\cal{A}, \rho)$.
In the same local notation, the image of its spectral action is $\cal{O}_\U [T] = \rho (\cal{A}_\U)$.
Since $\rho$ is injective, there is a unique section $b \in \cal{A} (\U)$ such that $T = \rho (b)$, and this section generates $\cal{A}_\U$ as an $\cal{O}_\U$-algebra.
If $e' = he$ is another generator of $\Lambda_\U$, then
\begin{eqn}
	T' = h^{-1} T
\qqtext{and hence}
	b' = h^{-1} b
\fullstop
\end{eqn}
Therefore, the local sections $b \otimes e$ glue into a global section $a$ of $\cal{A} \otimes \Lambda$.
These constructions are inverse to one another.

Next, by the projection formula, a section $a \in \H^0 (\X, \cal{A} \otimes \Lambda)$ is the same thing as a section $\widetilde{a} \in \H^0 (\Y, \pi^\ast \Lambda)$.
The latter determines a holomorphic map over $\X$,
\begin{eqn}
	i_a : \Y \too \tot (\Lambda)
\qtext{given by}
	y \mapstoo \widetilde{a} (y) \in \Lambda \big(\pi (y) \big)
\fullstop
\end{eqn}
Suppose $a$ satisfies condition \autoref*{1792410296505} and let us use the local notation $a_\U = b \otimes e$ again.
Since $b$ generates the rank-$n$ algebra $\cal{A}_\U$, its powers $1, b, \ldots, b^{n-1}$ form an $\cal{O}_\U$-basis.
Indeed, they span each fibre by the Cayley-Hamilton Theorem and hence form a basis locally.
If $P_b (t)$ is the characteristic polynomial of multiplication by $b$, the Cayley-Hamilton Theorem gives an algebra isomorphism
\begin{eqn}
	\cal{O}_\U [t] \big/ \inner{P_b (t)} \iso \cal{A}_\U
\qtext{by}
	t \mapstoo b
\fullstop
\end{eqn}
So the holomorphic map $i_a$ is a closed embedding that, in this trivialisation, identifies $\Y_\U$ with the closed analytic subspace of $\U \times \CC$ cut out by the monic equation $P_b (t) = 0$:
\begin{eqn}
	i_a (\Y_\U) \iso \set[\big]{ P_b (t) = 0 } \subset \U \times \CC
\fullstop
\end{eqn}

Conversely, suppose $i : \Y \inj \tot (\Lambda)$ is a closed embedding over $\X$.
Let $a \in \H^0 (\X, \cal{A} \otimes \Lambda)$ be the corresponding section.
Fix $x \in \X$, choose a trivialising neighbourhood $\U$ for $\Lambda$ with a generator $e$, and write $a_\U = b \otimes e$.
The trivialisation $e(x)$ identifies $\Lambda(x)$ with $\CC$; write $t$ for its coordinate.
By the correspondence defining $a$, the pullback of $t$ along $i_x : \Y_x \inj \Lambda(x)$ is the element $b(x)$ of
$\cal{A}(x) \cong \H^0(\Y_x, \cal{O}_{\Y_x})$.
Polynomial interpolation at the finitely many points of $\Y_x$, including the derivatives prescribed by their multiplicities, shows that every function on $\Y_x$ is represented by a unique polynomial of degree less than $n$.
Thus, the elements $1, b(x), \ldots, b(x)^{n-1}$ form a basis of the fibre algebra $\cal{A} (x)$.

Consequently, the bundle map
\begin{eqn}
    \cal{O}_\U^{\oplus n} \too \cal{A}_\U
    \qtext{sending}
    (f_0, \ldots, f_{n-1})
    \mapstoo
    f_0 + f_1 b + \cdots + f_{n-1} b^{n-1}
\end{eqn}
is an isomorphism on the fibre at $x$.
Its determinant is therefore nonzero at $x$, and hence remains nonzero after shrinking $\U$.
It follows that $1, b, \ldots, b^{n-1}$ form a local frame of $\cal{A}_\U$, so $\cal{A}_\U = \cal{O}_\U [b]$.
The canonical correspondence between sections and maps over $\X$ therefore identifies precisely the twisted algebra generators with the closed embeddings.
\end{proof}

\begin{example}[regular and non-regular Higgs fields]
\label{1790757898859}
The following two rank-two examples have the same local spectral algebra $\cal{O}_\X [t] / \inner{t^2}$, but only one of these actions defines a spectral structure.

The simplest example of a non-regular Higgs field is $\Phi = 0$ on a rank-$2$ bundle $\cal{E}$.
Notably, $\Phi$ has two distinct Jordan blocks for the eigenvalue $0$.
Locally, its characteristic polynomial is $t^2$, so the spectral algebra is locally $\cal{A} \simeq \cal{O}_\X [t] / \inner{t^2}$.
But $t$ acts trivially on $\cal{E}$ because $T = 0$ in any local trivialisation of $\Lambda$, so any $a = f + gt \in \cal{A}$ acts on any local section $e \in \cal{E}$ by $a \cdot e = fe$.
In other words, $\cal{A} \cdot e = \cal{O}_\X \cdot e$ for every local section $e \in \cal{E}$, so no single section of $\cal{E}$ can generate the rank-$2$ bundle $\cal{E}$ as an $\cal{A}$-module.
Thus $\cal{E}$ cannot be locally free of rank one over $\cal{A}$, and so the canonical algebra action $\rho : \cal{A} \to \cal{End} (\cal{E})$ does not define a spectral structure.

In contrast, consider a rank-$2$ Higgs bundle $(\cal{E}, \Phi)$ locally represented by the endomorphism consisting of the single Jordan block
\begin{eqn}
	\Phi \loceq \mtx{ 0 & 1 \\ 0 & 0 }
\fullstop
\end{eqn}
This has the same characteristic polynomial $t^2$ as the zero Higgs field, and hence the same spectral algebra $\cal{A} \simeq \cal{O}_\X [t] / \inner{t^2}$.
However, the Higgs field $\Phi$ is regular, so the induced algebra action $\rho : \cal{A} \to \cal{End} (\cal{E})$ defines a spectral structure on $\cal{E}$.
Indeed, if $e_1, e_2$ is the local basis of $\cal{E}$ with respect to which $\Phi$ is the Jordan block above, then $e_2$ determines a cyclic vector in every fibre because $\Phi e_2 = e_1$.
Consequently, $\cal{E} \simeq \cal{A}$ as an $\cal{A}$-module.
\end{example}

\subsection{Groupoids of Higgs Bundles and Spectral Structures}
\label{1793135088636}

The correspondence in \autoref{1792167000670} also has a natural formulation in terms of groupoids.
Fix a line bundle $\Lambda$ on $\X$ and introduce the groupoid of regular $\Lambda$-twisted Higgs bundles on $\X$ of rank $n$:
\begin{eqn}
	\mathbf{Higgs}_{n,\Lambda}^\rm{reg} (\X)
	\quad:\quad
	\set[\Big]{ (\cal{E}, \Phi) \text{ regular $\Lambda$-twisted Higgs bundles of rank-$n$} }
\fullstop
\end{eqn}
Recall that an isomorphism $(\cal{E}, \Phi) \iso (\cal{E}', \Phi')$ is a vector bundle isomorphism $\phi : \cal{E} \iso \cal{E}'$ satisfying
\begin{eqn}
	(\phi \otimes \id_\Lambda) \circ \Phi
	=
	\Phi' \circ \phi
\fullstop
\end{eqn}
For a fixed vector bundle $\cal{E}$, we similarly consider the groupoid of regular $\Lambda$-twisted Higgs fields on $\cal{E}$ where isomorphisms are required to induce the identity on $\cal{E}$:
\begin{eqn}
	\mathbf{Higgs}_{\Lambda}^\rm{reg} (\cal{E})
	\quad:\quad
	\set[\Big]{ \Phi : \cal{E} \to \cal{E} \otimes \Lambda \text{ regular Higgs field} }
	\subset \H^0 \big( \X, \cal{End} (\cal{E}) \otimes \Lambda \big)
\fullstop
\end{eqn}
This is a discrete groupoid: its objects are regular global sections of $\cal{End} (\cal{E}) \otimes \Lambda$, and its only morphisms are identities.

By \autoref{1790763564459}, there is a natural forgetful functor
\begin{eqntag}
\label{1792167464000}
	H : \mathbf{Higgs}_{n,\Lambda}^\rm{reg} (\X) \too \mathbf{SpStr}_n (\X)
\qtext{sending}
	(\cal{E}, \Phi) \mapstoo (\cal{E}, \cal{A}, \rho)
\fullstop
\end{eqntag}
where $(\cal{A}, \rho)$ is the canonical spectral structure in $\cal{E}$ given by pushing forward the spectral line bundle from the spectral curve of $\Phi$.
On morphisms, $H$ sends a morphism $\phi : (\cal{E}, \Phi) \to (\cal{E}', \Phi')$ to a morphism $(\phi, \alpha)$, where $\alpha$ is the unique algebra isomorphism satisfying 
\begin{eqn}
	\rho' \circ \alpha = \Ad_\phi \circ \rho
\fullstop
\end{eqn}
Indeed, if $T$ and $T'$ locally represent $\Phi$ and $\Phi'$, then $T' = \phi \circ T \circ \phi^{-1}$, so conjugation identifies their generated algebras.
Fixing $\cal{E}$ gives the corresponding forgetful functor
\begin{eqn}
	H : \mathbf{Higgs}_\Lambda^\rm{reg} (\cal{E}) \too \mathbf{SpStr} (\cal{E})
\fullstop
\end{eqn}
These functors fit into a commutative square
\begin{eqn}
\begin{tikzcd}
	\mathbf{Higgs}_\Lambda^\rm{reg} (\cal{E})
		\ar[r, hook]
		\ar[d, "H"']
&	\mathbf{Higgs}_{n,\Lambda}^\rm{reg} (\X)
		\ar[d, "H"]
\\
	\mathbf{SpStr} (\cal{E})
		\ar[r, hook]
&	\mathbf{SpStr}_n (\X)
\end{tikzcd}
\end{eqn}

The forgetful functor $H$ obviously commutes with forgetting to $\mathbf{Bun}_n (\X)$, and the fibre of this forgetful functor over any bundle $\cal{E}$ is
\begin{eqn}
	(F \circ H)^{-1} (\cal{E})
	= \mathbf{Higgs}_{\Lambda}^\rm{reg} (\cal{E})
	\cong (F \circ H)^{[-1]} (\cal{E})
\fullstop
\end{eqn}
At the same time, by \autoref{1790763564459} the spectral curve $\Y$ of a Higgs field $\Phi$ is the same as the spectral realisation of the canonical spectral algebra $\cal{A}$ acting on $\cal{E}$ induced by $\Phi$.
Combined with the diagram \eqref{1790940189747}, we therefore have the following picture:
\begin{eqntag}
\label{1792401591304}
\begin{tikzcd}
&
&	\mathbf{Bun}_n (\X)
\\
	\mathbf{Higgs}^\rm{reg}_{n, \Lambda} (\X)
		\ar[r, "H"]
		\ar[rru, bend left]
		\ar[rrd, bend right]
&	\mathbf{SpStr}_n (\X)
		\ar[ur, "F"]
		\ar[dr, "G"']
\\
&
&	\mathbf{Cov}_n (\X)
\end{tikzcd}
\end{eqntag}
This diagram shows that, like spectral structures, regular Higgs bundles also interpolate between the geometry of finite covers of $\X$ and the geometry of vector bundles on $\X$.
But, from this point of view, spectral structures provide a more intrinsic intermediary.
Indeed, the groupoid of Higgs bundles requires a choice of a twisting line bundle $\Lambda$, but the groupoid of spectral structures and the two forgetful functors from spectral structures to covers and vector bundles require no additional choices.

As \autoref{1792167000670} explains, to build a Higgs bundle from a spectral structure, we must additionally choose a line bundle $\Lambda$ and a global $\Lambda$-twisted generator of the spectral algebra.
For any line bundle $\Lambda$ on $\X$ and any spectral algebra $\cal{A}$ of rank $n$, introduce the set of global $\Lambda$-twisted generators of $\cal{A}$:
\begin{eqn}
	\sfop{Gen}_\Lambda (\cal{A})
	\coleq
	\set{ a \in \H^0 (\X, \cal{A} \otimes \Lambda) ~\Big|~ \text{local representatives of $a$ generate $\cal{A}$} }
\fullstop
\end{eqn}
Let $\mathbf{Gen}_\Lambda (\cal{A})$ denote the discrete groupoid over $\sfop{Gen}_\Lambda (\cal{A})$; i.e., its objects are $\sfop{Gen}_\Lambda (\cal{A})$ and its only morphisms are the identity morphisms.
Next, introduce the groupoid
\begin{eqn}
	\mathbf{SpStrGen}_\Lambda (\X)
	\quad:\quad
	\set[\Big]{ (\cal{E}, \cal{A}, \rho, a) ~\Big|~ a \in \sfop{Gen}_\Lambda (\cal{A}) }
\fullstop
\end{eqn}
of spectral structures on $\X$ equipped with a global $\Lambda$-twisted generator.
Its objects are tuples $(\cal{E}, \cal{A}, \rho, a)$ consisting of a spectral structure and a section $a \in \sfop{Gen}_\Lambda (\cal{A})$.
An isomorphism between two such tuples is an isomorphism $(\phi, \alpha)$ of the underlying spectral structures satisfying $(\alpha \otimes \id_\Lambda) (a) = a'$.

\begin{proposition}[forgetting the twisted generator]
\label{1792167464001}
The groupoid of regular $\Lambda$-twisted Higgs bundles is canonically equivalent to the groupoid of spectral structures equipped with a global $\Lambda$-twisted generator:
\begin{eqn}
	\mathbf{Higgs}_{n,\Lambda}^\rm{reg} (\X)
	\cong
	\mathbf{SpStrGen}_\Lambda (\X)
\fullstop
\end{eqn}
Under this equivalence, the forgetful functor $H$ simply forgets the global twisted generator $a$:
\begin{eqn}
	H : (\cal{E}, \cal{A}, \rho, a) \mapstoo (\cal{E}, \cal{A}, \rho)
\fullstop
\end{eqn}
In particular, $H$ is faithful.
Furthermore, for every spectral structure $(\cal{A}, \rho)$ on $\cal{E}$, there are canonical equivalences of groupoids
\begin{eqn}
	H^{[-1]} (\cal{E}, \cal{A}, \rho)
	\cong
	H^{[-1]} (\cal{A}, \rho)
	\cong
	\mathbf{Gen}_\Lambda (\cal{A})
\fullstop
\end{eqn}
If nonempty, the subset $\sfop{Gen}_\Lambda (\cal{A}) \subset \H^0 (\X, \cal{A} \otimes \Lambda)$ is Zariski-open, hence a smooth complex manifold of dimension $h^0 (\X, \cal{A} \otimes \Lambda)$.
\end{proposition}

Thus, by \autoref{1792167000670}, the essential image of $H$ consists precisely of the spectral structures whose spectral covers admit an embedding over $\X$ into $\tot (\Lambda)$.
In particular, the functor $H$ need not be essentially surjective.
It need not be full either: an isomorphism of spectral structures lifts to an isomorphism of Higgs bundles exactly when it also carries one chosen twisted generator to the other.

\begin{proof}
By \autoref{1790763564459}, a regular Higgs field $\Phi$ determines a spectral structure $(\cal{A}, \rho)$ together with a canonical twisted generator $a \in \H^0 (\X, \cal{A} \otimes \Lambda)$.
Locally, $a$ is represented by the class of $t$ in the presentation \eqref{1792409375201} of $\cal{A}$.
Conversely, a tuple $(\cal{E}, \cal{A}, \rho, a)$ determines the regular Higgs field $\Phi = (\rho \otimes \id_\Lambda) (a)$ by \autoref{1792167000670}.
The two constructions are inverse up to the canonical identification of the spectral algebras.

An isomorphism $(\phi, \alpha)$ of spectral structures intertwines the resulting Higgs fields precisely when $(\alpha \otimes \id_\Lambda) (a) = a'$.
Indeed, this equality gives
\begin{eqn}
	\Phi'
	=
	(\rho' \otimes \id_\Lambda) (a')
	=
	(\Ad_\phi \otimes \id_\Lambda) \big( (\rho \otimes \id_\Lambda) (a) \big)
	=
	(\phi \otimes \id_\Lambda) \circ \Phi \circ \phi^{-1}
\fullstop
\end{eqn}
Conversely, this identity and the injectivity of $\rho'$ force $(\alpha \otimes \id_\Lambda) (a) = a'$.
This proves the equivalence on morphisms as well as on objects.
Faithfulness follows because a morphism of the enriched tuples $(\cal{E}, \cal{A}, \rho, a)$ is already determined by the underlying morphism of spectral structures $(\cal{E}, \cal{A}, \rho)$.

Next, an object of the essential fibre over $(\cal{E}, \cal{A}, \rho)$ is a regular Higgs bundle $(\cal{E}', \Phi')$ together with an isomorphism $(\phi, \alpha) : H (\cal{E}', \Phi') \iso (\cal{E}, \cal{A}, \rho)$.
Transporting its canonical generator $a'$ through the algebra isomorphism $\alpha$ gives an element $a \in \sfop{Gen}_\Lambda (\cal{A})$.
Conversely, every such $a$ gives a Higgs field on $\cal{E}$ together with its canonical spectral identification.
Two objects of the essential fibre are isomorphic if and only if their transported generators agree.
When they agree, the required Higgs-bundle isomorphism is uniquely determined by their identifications with $\cal{E}$.
In particular, all automorphism groups in the essential fibre are trivial.
The same argument with the bundle $\cal{E}$ fixed proves the other asserted equivalence.

Finally, for any $a \in \H^0 (\X, \cal{A} \otimes \Lambda)$, its powers define a bundle morphism
\begin{eqn}
	\epsilon_a : \cal{O}_\X \oplus \Lambda^{-1} \oplus \cdots \oplus \Lambda^{-(n-1)} \too \cal{A}
\qtext{by}
	(u_0, \ldots, u_{n-1}) \mapstoo \sum_{i=0}^{n-1} u_i a^i
\fullstop
\end{eqn}
The local representatives of $a$ generate $\cal{A}$ if and only if $\epsilon_a$ is an isomorphism.
Indeed, the fact that local representatives of $a$ generate $\cal{A}$ implies that the sections $1, a, \ldots, a^{n-1}$ span every fibre by the Cayley-Hamilton Theorem.
The converse is immediate.
Taking determinants gives a section
\begin{eqn}
	\det \epsilon_a \in \H^0 \big( \X, \det \cal{A} \otimes \Lambda^{n(n-1)/2} \big)
\fullstop
\end{eqn}
If $\sfop{Gen}_\Lambda (\cal{A})$ is nonempty, choose one generator $a_0$.
Then $\det \epsilon_{a_0}$ is nowhere vanishing and trivialises this line bundle.
Since $\X$ is compact and connected, every other determinant section has the form $\det \epsilon_a = p(a) \det \epsilon_{a_0}$ for some scalar $p(a) \in \CC$.
The entries of $\epsilon_a$ depend polynomially on $a$, so $p$ is a polynomial function on the finite-dimensional vector space $\H^0 (\X, \cal{A} \otimes \Lambda)$.
Consequently, we have the equality of sets
\begin{eqn}
	\sfop{Gen}_\Lambda (\cal{A})
	=
	\set{ a \in \H^0 (\X, \cal{A} \otimes \Lambda) ~\big|~ p(a) \neq 0 }
\fullstop
\end{eqn}
Since $p(a_0) = 1$, this is a nonempty Zariski-open subset.
\end{proof}

\subsubsection*{Forgetting the identification}

As in the discussion of the forgetful functor \eqref{1791407291197}, it matters whether an identification with the fixed spectral structure is part of the data.
The essential fibre $H^{[-1]} (\cal{E}, \cal{A}, \rho)$ includes such an identification, which is why its automorphism group is trivial.
If, instead, we take the full subgroupoid of $\mathbf{Higgs}_{n,\Lambda}^\rm{reg} (\X)$ consisting of Higgs bundles whose spectral structures are merely isomorphic to $(\cal{E}, \cal{A}, \rho)$, this groupoid is equivalent to the action groupoid
\begin{eqn}
	\Aut (\cal{E}, \cal{A}, \rho)
	\ltimes
	\sfop{Gen}_\Lambda (\cal{A})
	\rightrightarrows
	\sfop{Gen}_\Lambda (\cal{A})
\fullstop{,}
\end{eqn}
with action given by $(\phi, \alpha) : a \mapsto (\alpha \otimes \id_\Lambda) (a)$, where $\Aut (\cal{E}, \cal{A}, \rho)$ is the automorphism group of $(\cal{E}, \cal{A}, \rho)$ in $\mathbf{SpStr}_n (\X)$.
Indeed, choosing an identification with $(\cal{E}, \cal{A}, \rho)$ produces a twisted generator $a$, and changing that identification acts on the generator by precisely the rule $a \mapsto (\alpha \otimes \id_\Lambda) (a)$.
The morphisms in the resulting action groupoid are exactly the isomorphisms of spectral structures which carry one generator to the other.

The stabiliser of any $a \in \sfop{Gen}_\Lambda (\cal{A})$ is canonically
\begin{eqn}
	\Aut (\cal{E}, \Phi_a)
	\cong
	\H^0 (\X, \cal{A})^\times
	\cong
	\H^0 (\Y, \cal{O}_\Y)^\times
\fullstop{,}
\end{eqn}
where $\Phi_a = (\rho \otimes \id_\Lambda) (a)$.
To see this, if $(\alpha \otimes \id_\Lambda) (a) = a$, then the fact that $a$ is a generator forces $\alpha = \id_\cal{A}$.
So the stabiliser consists precisely of the $\cal{A}$-linear automorphisms of $\cal{E}$.
But since $\cal{E}$ is locally free of rank one over the commutative algebra $\cal{A}$, these are exactly multiplication by global units of $\cal{A}$.
In particular, when $\Y$ is smooth and connected, the automorphism group of every corresponding Higgs bundle is $\CC^\times$.

\begin{example}[forgetting the eigenvalues]
\label{1792167464002}
Suppose $\cal{E} = \cal{L}_1 \oplus \cdots \oplus \cal{L}_n$ has its split spectral structure with algebra $\cal{A} = \cal{O}_\X^{\oplus n}$, and take $\Lambda = \cal{O}_\X$.
Since $\X$ is compact and connected, a global section of $\cal{A}$ is a tuple of constants $(c_1, \ldots, c_n) \in \CC^n$.
Such a tuple generates the algebra in every fibre if and only if its entries are pairwise distinct.
Consequently,
\begin{eqn}
	\sfop{Gen}_{\cal{O}_\X} (\cal{A})
	=
	\set{ (c_1, \ldots, c_n) \in \CC^n ~\big|~ c_i \neq c_j \text{ whenever } i \neq j }
\fullstop
\end{eqn}
The corresponding regular Higgs field $\Phi_a$ acts on $\cal{L}_i$ as $c_i \id_{\cal{L}_i}$; i.e., the global generator $a = (c_1, \ldots, c_n) \in \CC^n$ corresponds to the diagonal Higgs field
\begin{eqn}
	\Phi_a = \diag \big( c_1, \ldots, c_n \big)
\fullstop
\end{eqn}
Thus the forgetful functor $H$ retains the decomposition into the eigenline bundles but forgets the eigenvalues.
An identification with the fixed spectral structure remembers the individual factors of $\cal{A}$, so its essential fibre is the space of all possible $n$-tuples of distinct eigenvalues, which is the ordered configuration space displayed above.
\end{example}

\subsection{Existence of Spectral Structures via Higgs Bundles}
\label{1791465731919}

In this subsection, we give a second proof of \autoref{1790777298686} and establish \autoref{1792167000673} using the methods of Higgs bundles.
The construction relies on choosing a sufficiently positive twist and a suitable regular Higgs field whose associated spectral structure has the desired properties.
In particular, every vector bundle admits a clean spectral structure induced by a regular Higgs field.
However, a given clean spectral structure need not itself be induced by a Higgs field; see \autoref{1792434634960}.

The construction has two main steps.
The first step is to show that a general suitably twisted Higgs field has a smooth, connected, and simply branched spectral cover.

\begin{lemma}[smooth connected simply branched spectral curves]
\label{1791390497007}
Let $\cal{E}$ be a vector bundle of rank $n \geq 2$ on $\X$, and let $\Lambda$ be a line bundle of positive degree such that $\cal{End} (\cal{E}) \otimes \Lambda$ is globally generated.
Then a general section
\begin{eqn}
	\Phi \in \H^0 \big( \X, \cal{End} (\cal{E}) \otimes \Lambda \big)
\end{eqn}
defines a regular $\Lambda$-twisted Higgs field whose spectral cover $\pi : \Y \to \X$ is smooth, connected, and simply branched.
In particular, the spectral correspondence provides a line bundle $\cal{L}$ on $\Y$ and an isomorphism $\iota : \pi_\ast \cal{L} \iso \cal{E}$.
\end{lemma}

A ``general section'' means any section belonging to a suitable nonempty Zariski-open subset of the relevant space of global sections, in this case $\H^0 \big( \X, \cal{End} (\cal{E}) \otimes \Lambda \big)$.
The hypothesis that $\cal{End} (\cal{E}) \otimes \Lambda$ is globally generated can always be fulfilled by taking a sufficiently high power of any positive-degree line bundle.
This lemma provides an alternative proof of \autoref{1790777298686}, independent of the projective-geometric construction given in \ref{1791480675169}.
Indeed, if $\cal{E}$ has rank $n \geq 2$, choose a positive-degree line bundle $\Lambda_0$ and take $N$ sufficiently large that $\cal{End} (\cal{E}) \otimes \Lambda_0^N$ is globally generated~\cite[Ch.~II, \S7; Ch.~IV, Corollary~3.3]{Hartshorne1977}.
\autoref{1791390497007} then gives a smooth connected spectral structure on $\cal{E}$, with simple branching, induced
by a regular twisted Higgs field.
For rank one, the same assertion follows by taking the identity cover $\pi = \id_\X$ and $\cal{L} = \cal{E}$.

The second main step is to show that, in fact, this power can be increased so high that a general Higgs field also gives the required vanishing of off-diagonal cohomology.

\begin{lemma}[clean spectral structures from regular Higgs fields]
\label{1792430412943}
Let $\cal{E}$ be a vector bundle of rank $n$ on $\X$, let $\Lambda_0$ be any line bundle on $\X$ of positive degree, and set $\Lambda \coleq \Lambda_0^N$ for an integer $N$.
Then, for every sufficiently large $N$, a general section
\begin{eqn}
	\Phi \in \H^0 \big( \X, \cal{End} (\cal{E}) \otimes \Lambda \big)
\end{eqn}
defines a regular $\Lambda$-twisted Higgs field whose induced spectral structure is clean.
That is, its spectral cover is smooth, connected, and simply branched, and its off-diagonal bundle $\cal{Q}$ satisfies
\begin{eqn}
	\H^1 (\X, \cal{Q}) = 0
\fullstop
\end{eqn}
\end{lemma}

This lemma establishes the stronger existence result of \autoref{1792167000673} for clean spectral structures.
We do not know if the alternative proof of \autoref{1790777298686} presented in \ref{1791480675169} can be generalised to the situation of clean spectral structures.

\begin{proof}[Proof of \autoref{1791390497007}]
Put $\V \coleq \H^0 \big(\X, \cal{End} (\cal{E}) \otimes \Lambda \big)$.
Since $\cal{End} (\cal{E}) \otimes \Lambda$ is assumed to be globally generated, the evaluation map
\begin{eqn}
	\op{ev} : \X \times \V \too \tot \big(\cal{End} (\cal{E}) \otimes \Lambda \big)
\qqtext{given by}
	(x, \Phi) \mapstoo \Phi(x)
\end{eqn}
is a holomorphic submersion.
We use this to show that a general section avoids the singular locus of the matrix discriminant and meets its smooth locus transversely.

First consider the space of complex $n {\times} n$ matrices, $\sfop{Mat}_n(\CC)$.
Let $\D$ be the discriminant hypersurface consisting of matrices with a repeated eigenvalue; equivalently,
\begin{eqn}
	\D = \set{ T \in \sfop{Mat}_n(\CC) ~\Big|~ \Delta(T) = 0 }
\qqtext{where}
	\Delta(T) \coleq \op{disc} \big( \det(tI-T) \big)
\fullstop
\end{eqn}
Its smooth locus consists precisely of the matrices having one eigenvalue of multiplicity two, represented by a single Jordan block of size two, with all other eigenvalues simple.
Here is a verification, including the distinction between a Jordan block and a semisimple double eigenvalue.

Write $p_T(t) \coleq \det(tI-T)$ for the characteristic polynomial of $T$.
Suppose first that $p_T$ has exactly one double root $\lambda$ and all its other roots are simple.
We examine when the differential of $\Delta$ is nonzero at such a matrix $T$.

Near $p_T$ in the space of monic polynomials, the two roots near $\lambda$ determine a holomorphically varying monic quadratic factor.
The remaining roots are simple, so they vary holomorphically and remain distinct from one another and from the roots of this quadratic factor.
Consequently, the discriminant of the full polynomial is a nowhere-vanishing holomorphic factor times the discriminant of the quadratic factor.

To make this explicit, shift the variable by $\lambda$ and write the quadratic factor as
\begin{eqn}
	u^2 + au + b
\qqtext{where}
	u = t - \lambda
\fullstop
\end{eqn}
At the original polynomial $p_T$, we have $a = b = 0$.
The discriminant of this quadratic factor is $a^2 - 4b$, whose differential at $a = b = 0$ sends a variation $(\dot{a}, \dot{b})$ to $-4\dot{b}$.
On the other hand, if
\begin{eqn}
	p_T(t) = (t-\lambda)^2 r(t)
\qqtext{with}
	r(\lambda) \neq 0
\fullstop{,}
\end{eqn}
then a variation $\dot{p}$ of the full polynomial satisfies $\dot{p}(\lambda) = r(\lambda) \dot{b}$.
Thus, the differential of the polynomial discriminant at $p_T$ is a nonzero scalar multiple of the functional $\dot{p} \mapsto \dot{p}(\lambda)$.

Now let $\dot{T}$ be a variation of the matrix $T$.
Differentiating the determinant gives the corresponding variation of its characteristic polynomial:
\begin{eqn}
	\dot{p}_T (t)
	\coleq
	\left.
		\frac{\del}{\del \epsilon}
	\right|_{\epsilon=0}
	\det \big( tI - T - \epsilon \dot{T} \big)
	=
	-\tr \big( \op{adj}(tI-T) \dot{T} \big)
\fullstop{,}
\end{eqn}
where $\op{adj}$ denotes the adjugate matrix.
Consequently, the differential of $\Delta$ at $T$ is a nonzero scalar multiple of the linear functional
\begin{eqn}
	\dot{T} \mapstoo -\tr \big( \op{adj}(\lambda I-T) \dot{T} \big)
\fullstop
\end{eqn}
The trace pairing on matrices is nondegenerate, so this functional is nonzero if and only if $\op{adj}(\lambda I-T)$ is nonzero.
Since $\lambda$ is an eigenvalue, this is equivalent to $\lambda I-T$ having rank $n-1$, or equivalently to the $\lambda$-eigenspace being one-dimensional.

There are just two possible Jordan types for an eigenvalue of multiplicity two.
If it is represented by a single Jordan block of size two, its eigenspace is one-dimensional and the differential of $\Delta$ is nonzero.
If it is represented by two Jordan blocks of size one, its eigenspace is two-dimensional, the adjugate vanishes, and so does the differential of $\Delta$.
This explains the distinction between the two cases.

Next, suppose that $p_T$ has fewer than $n-1$ distinct roots.
Then $\deg \gcd (p_T, p_T') \geq 2$.
The Sylvester matrix defining the resultant $\op{Res}(p_T,p_T')$ therefore has corank at least two.
All its cofactors vanish, so the differential of its determinant vanishes.
Since the discriminant of a monic polynomial is, up to a fixed sign, this resultant, its differential vanishes at $p_T$.
It follows that the differential of $\Delta$ also vanishes at $T$.

We now show that the matrices where this differential vanishes form a subset of codimension at least two in $\sfop{Mat}_n(\CC)$.
For this, we count the dimensions of the relevant Jordan strata.
Fix a Jordan type with $k$ distinct eigenvalues, and let $c$ be the dimension of the centraliser of a matrix of that type.
Its conjugacy orbit has dimension $n^2-c$, and its $k$ distinct eigenvalues supply $k$ further parameters.
Thus, the corresponding Jordan stratum has dimension
\begin{eqn}
	n^2 - c + k
\fullstop
\end{eqn}
We always have $c \geq n$.
Indeed, the block-diagonal matrices that act polynomially on each individual Jordan block already form a commuting space whose dimension is the sum of the block sizes, namely $n$.
Therefore, if $k \leq n-2$, the stratum has dimension at most $n^2-2$.

The only remaining case where the differential vanishes has $k=n-1$ and one semisimple double eigenvalue.
Its centraliser consists of an arbitrary endomorphism of the two-dimensional eigenspace, together with scalars on the remaining $n-2$ eigenlines.
Hence,
\begin{eqn}
	c = 4 + (n-2) = n+2
\fullstop{,}
\end{eqn}
and the corresponding stratum has dimension $n^2-3$.
Since there are only finitely many Jordan types, the union of all these strata has dimension at most $n^2-2$.

This also verifies that $\Delta$ defines a reduced hypersurface.
Indeed, if $\Delta$ had a repeated irreducible factor, its differential would vanish along the entire corresponding hypersurface component, contradicting the dimension bound just obtained.
The usual differential criterion for smoothness therefore applies to $\D$.
We have proved both the asserted description of $\D^{\rm{sm}}$ and the codimension estimate
\begin{eqn}
	\op{codim}_{\sfop{Mat}_n(\CC)}
	\D^{\rm{sing}}
	\geq 2
\fullstop
\end{eqn}

We now apply these matrix calculations to the bundle
$\cal{End}(\cal{E}) \otimes \Lambda$.
The subsets $\D$, $\D^{\rm{sm}}$, and $\D^{\rm{sing}}$ are invariant under conjugation and multiplication by a nonzero scalar.
These are precisely the changes arising from different local trivialisations of $\cal{E}$ and $\Lambda$.
Thus, abusing notation, they define corresponding subsets $\D, \D^\rm{sm}, \D^\rm{sing}$ in the total space $\tot \big( \cal{End}(\cal{E}) \otimes \Lambda \big)$ independently of these choices.

Since the evaluation map $\op{ev}$ is a submersion, the preimage $\op{ev}^{-1}(\D^{\rm{sing}})$ has codimension at least two in $\X \times \V$.
In particular,
\begin{eqn}
	\dim \op{ev}^{-1}(\D^{\rm{sing}}) \leq \dim \V - 1
\fullstop
\end{eqn}
The projection $\op{pr}_\V : \X \times \V \to \V$ is proper because $\X$ is compact.
The Proper Mapping Theorem \cite[Ch.~III, Theorem~2.11]{BanicaStanasila1976} therefore shows that the image of this preimage is a closed analytic subset of $\V$.
Its dimension is at most $\dim \V - 1$, so it cannot be all of $\V$.
Consequently, its complement
\begin{eqn}
	\V_0
	\coleq
	\V \setminus
	\op{pr}_\V
	\big( \op{ev}^{-1}(\D^{\rm{sing}}) \big)
\end{eqn}
is nonempty and open.
Every section $\Phi \in \V_0$ avoids $\D^{\rm{sing}}$ at every point of $\X$.

It remains to arrange transversality to $\D^{\rm{sm}}$.
Consider the incidence space
\begin{eqn}
	\W
	\coleq
	\set{ (x,\Phi) \in \X \times \V_0 ~\Big|~ \Phi(x) \in \D^{\rm{sm}} }
\fullstop
\end{eqn}
Since $\op{ev}$ is a submersion and $\D^{\rm{sm}}$ is a smooth hypersurface, $\W$ is smooth of dimension $\dim \V$.
Moreover, no section in $\V_0$ meets $\D^{\rm{sing}}$, so $\W$ is also the preimage of the closed subset $\D$ inside $\X \times \V_0$.
It is therefore closed in $\X \times \V_0$, and the projection $q : \W \to \V_0$ is proper.

We claim that $q$ is a submersion at $(x,\Phi)$ if and only if $\Phi$ meets $\D^{\rm{sm}}$ transversely at $x$.
To see this, choose a local coordinate $z$ on $\X$ and local trivialisations of $\cal{E}$ and $\Lambda$.
Write $T_\Phi(z)$ for the matrix representing $\Phi$ in these trivialisations.
Then a local defining equation for $\W$ is
\begin{eqn}
	\delta(z,\Phi)
	\coleq
	\Delta \big( T_\Phi(z) \big)
	=
	0
\fullstop
\end{eqn}
The derivative of $\delta$ in the $\V$-direction is surjective onto $\CC$.
Indeed, evaluation at $x$ is surjective onto the matrix space, and the differential of $\Delta$ is nonzero at $T_\Phi(x)$.
The tangent space to $\W$ is given by the equation
\begin{eqn}
	\frac{\del \delta}{\del z} \dot{z} + \d_\Phi \delta (\dot{\Phi}) = 0
\fullstop
\end{eqn}
Thus, every variation $\dot{\Phi}$ can be lifted to a tangent vector of $\W$ precisely when $\partial \delta / \partial z$ is nonzero.
This is exactly the condition that $\Phi$ meet $\D^{\rm{sm}}$ transversely, proving the claim.

The critical locus of $q$ is a closed analytic subset of $\W$.
Since $q$ is proper, its image is a closed analytic subset of $\V_0$.
By Sard's Theorem \cite[Theorem~4.1]{Sard1942}, this image cannot fill $\V_0$.
Its complement is therefore a nonempty open subset $\V_1 \subset \V_0$.
Every $\Phi \in \V_1$ avoids $\D^{\rm{sing}}$ and meets $\D^{\rm{sm}}$ transversely.

We also need that these good choices form a Zariski-open subset of $\V$.
The discriminant of a Higgs field is a global section
\begin{eqn}
	\Delta_\Phi \in \H^0 \big( \X, \Lambda^{n(n-1)} \big)
\fullstop
\end{eqn}
The two conditions just established are equivalent to this section having only simple zeros.
Indeed, a point where $\Phi$ meets $\D^{\rm{sing}}$ gives a zero at which the derivative of the discriminant also vanishes.
At a point where $\Phi$ meets $\D^{\rm{sm}}$, the zero is simple precisely when the intersection is transverse.

Now, $\X$ is a smooth projective algebraic curve, and by GAGA \cite[Théorème~3]{Serre1956GAGA} we may regard $\cal{E}$ and $\Lambda$ as algebraic bundles.
The discriminant depends polynomially on $\Phi$.
Consequently, the condition that $\Delta_\Phi$ vanish to order at least two at $x$ defines a closed algebraic subset of $\X \times \V$.
Locally, it is given by the equations
\begin{eqn}
	\delta(z,\Phi) = 0
\qqtext{and}
	\frac{\partial \delta}{\partial z}(z,\Phi) = 0
\fullstop
\end{eqn}
Since $\X$ is proper, the projection of this subset to $\V$ is closed.
Its complement is precisely $\V_1$, which is therefore Zariski open as well as nonempty.

Choose any $\Phi \in \V_1$, and let $\Y \subset \tot(\Lambda)$ be its spectral curve.
In a local trivialisation of $\Lambda$, its equation is
\begin{eqn}
	p_\Phi (z,t) \coleq \det \big( tI-T_\Phi(z) \big) = 	0
\fullstop
\end{eqn}
This polynomial is monic of degree $n$ in $t$.
Thus, the projection $\pi : \Y \to \X$ is finite flat of degree $n$, and its algebra is locally free with basis $1,t,\ldots,t^{n-1}$.
Globally, this gives
\begin{eqn}
	\cal{A}
	\coleq
	\pi_\ast \cal{O}_\Y
	\cong
	\cal{O}_\X
	\oplus \Lambda^{-1}
	\oplus \cdots
	\oplus \Lambda^{-(n-1)}
\fullstop
\end{eqn}
Since $\pi$ is finite and $\X$ is compact, $\Y$ is compact.

We now check that $\Y$ is smooth and that $\pi$ is simply branched.
If $T_\Phi(x)$ has distinct eigenvalues, every root of $p_\Phi(x,t)$ is simple.
The Holomorphic Implicit Function Theorem then expresses $\Y$ near this fibre as $n$ disjoint smooth sheets, each mapping isomorphically to a neighbourhood of $x$.

Otherwise, $T_\Phi(x)$ belongs to $\D^{\rm{sm}}$.
It has exactly one double eigenvalue, represented by a single Jordan block of size two, and all its other eigenvalues are simple.
The simple eigenvalues again give smooth unramified sheets.
The two roots near the double eigenvalue determine a holomorphically varying quadratic factor.
After completing the square and choosing $z=0$ at $x$, its equation becomes $u^2 - b(z) = 0$ with $b(0) = 0$.

The discriminant of the full characteristic polynomial is a nowhere-vanishing holomorphic factor times $b(z)$.
Since $\Delta_\Phi$ has a simple zero at $x$, we have $b'(0) \neq 0$.
We may therefore use $b(z)$ as a new coordinate on the base, in which the ramified part of the cover has the familiar form $u^2 = z$.
It follows that $\Y$ is smooth here and that the ramification point is simple.
There is exactly one ramification point over this branch value because all the other sheets are unramified.
Thus, $\pi$ is simply branched.

The positivity of $\Lambda$ now gives connectedness.
For every $i \geq 1$, the line bundle $\Lambda^{-i}$ has negative degree and hence no nonzero global sections.
The decomposition of $\cal{A}$ therefore gives
\begin{eqn}
	\H^0 (\Y, \cal{O}_\Y)
	\cong
	\H^0 (\X, \cal{A})
	\cong
	\H^0 (\X, \cal{O}_\X)
	\cong
	\CC
\fullstop
\end{eqn}
A smooth compact complex curve with $k$ connected components has global holomorphic functions $\CC^k$.
Consequently, $\Y$ has exactly one connected component.

Finally, every matrix $T_\Phi(x)$ is regular.
Indeed, it either has distinct eigenvalues, or has a single Jordan block at its double eigenvalue and one-dimensional blocks at its remaining eigenvalues.
Thus, $\Phi$ is a regular Higgs field.

To see the corresponding spectral line bundle explicitly, fix a point $x \in \X$, choose a local trivialisation of $\Lambda$, and write $T$ for the local endomorphism representing $\Phi$.
Choose a cyclic vector $v_x \in \cal{E}(x)$ for $T(x)$ and extend it to a local holomorphic section $v$ of $\cal{E}$.
The sections $v, Tv, T^2v, \ldots, T^{n-1}v$ form a basis on the fibre at $x$.
Their determinant is therefore nonzero at $x$ and remains nonzero after shrinking the neighbourhood.
They consequently form a local frame of $\cal{E}$.
The Cayley-Hamilton action then gives a local $\cal{A}$-linear isomorphism
\begin{eqn}
	\cal{A} \iso \cal{E}
\qqtext{sending}
	[f(t)] \mapstoo f(T)v
\fullstop
\end{eqn}
Thus, $\cal{E}$ is locally free of rank one over $\cal{A}$.
The Finite Sheaf-Module Correspondence now supplies a line bundle $\cal{L}$ on $\Y$ together with an $\cal{A}$-linear isomorphism $\iota : \pi_\ast \cal{L} \iso \cal{E}$.
\end{proof}

\begin{proof}[Proof of \autoref{1792430412943}]
When $n = 1$, there is nothing to prove: every section $\Phi \in \H^0 (\X, \Lambda_0^N)$ is a regular Higgs field for every $N$, and its spectral cover is the graph of $\Phi$ in $\tot (\Lambda_0^N)$.
The projection $\pi : \Y \to \X$ is therefore an isomorphism.
Under this identification, $\cal{L} = \cal{E}$,~ $\cal{A} = \cal{O}_\X$, and $\cal{Q} = 0$ hence it has no cohomology and so the spectral structure is clean.

More generally, suppose now that $\cal{E}$ has rank $n \geq 2$.
Choose a line bundle $\Lambda_0$ of positive degree, and pick an integer $N$ sufficiently large that $\cal{End} (\cal{E}) \otimes \Lambda$ with $\Lambda \coleq \Lambda_0^N$ is globally generated.
Recall that such a choice is always possible by taking a sufficiently high power of any positive-degree line bundle and applying the usual global-generation theorem for coherent sheaves on compact Riemann surfaces.

By \autoref{1791390497007}, there is a nonempty Zariski-open subset
\begin{eqn}
	\V \subset \H^0 \big( \X, \cal{End} (\cal{E}) \otimes \Lambda \big)
\end{eqn}
whose elements are regular Higgs fields inducing smooth, connected, and simply branched spectral structures.
It remains to arrange the vanishing of the off-diagonal cohomology.
For this, introduce the vector space
\begin{eqn}
	\W \coleq \H^0 \big( \X, \cal{End} (\cal{E}) \otimes \omega_\X \big)
\fullstop
\end{eqn}
By \autoref{1791453955298} and Serre duality \cite[Corollary~III.7.7]{Hartshorne1977}, the vanishing condition $\H^1 (\X, \cal{Q}) = 0$ is equivalent to
\begin{eqn}
	\H^0 \big( \X, \cal{A}^\perp \otimes \omega_\X \big) = 0
\fullstop
\end{eqn}
If $\W = 0$, then every $\Phi \in \V$ already induces a clean spectral structure, and the assertion follows.
We may therefore assume that $\W \neq 0$.

Choose finitely many distinct points $x_1, \ldots, x_m \in \X$ and endomorphisms $T_k \in \End (\cal{E} (x_k))$ such that the linear map
\begin{eqn}
	\W \too \bigoplus_{k=1}^m \omega_\X (x_k)
\qtext{given by}
	S \mapstoo \big( \tr (S(x_k) \circ T_k) \big)_{k=1}^m
\end{eqn}
is injective.
Here and below, $\omega_\X (x_k)$ denotes the fibre of the canonical bundle at $x_k$.
To make these choices, take any nonzero $S \in \W$ and choose a point $x_1$ at which $S(x_1) \neq 0$.
The matrix trace pairing is perfect, so there is some $T_1$ such that $\tr (S(x_1) \circ T_1) \neq 0$.
Thus this first evaluation has a kernel of smaller dimension than $\W$.
If that kernel is nonzero, repeat the construction with a nonzero section in the kernel and a new point $x_2$.
A nonzero holomorphic section is nonzero away from finitely many points, so the new point can always be chosen distinct from those already used.
Each step strictly decreases the dimension of the common kernel, so the procedure terminates after finitely many steps.

Put $\D \coleq x_1 + \cdots + x_m$ and make sure the power $N$ in $\Lambda = \Lambda_0^N$ is large enough that $\cal{End} (\cal{E}) \otimes \Lambda$ is globally generated and
\begin{eqn}
	\H^1 \big( \X, \cal{End} (\cal{E}) \otimes \Lambda (-\D) \big) = 0
\fullstop
\end{eqn}
Both conditions hold for every sufficiently large $N$ \cite[Proposition~III.5.3]{Hartshorne1977}.
For each such $N$, apply \autoref{1791390497007} with $\Lambda = \Lambda_0^N$, and let $\V$ denote the resulting nonempty Zariski-open subset of $\H^0 \big( \X, \cal{End} (\cal{E}) \otimes \Lambda \big)$.
Restriction to the divisor $\D$ induces an exact sequence which shows that the evaluation map
\begin{eqn}
	\H^0 \big( \X, \cal{End} (\cal{E}) \otimes \Lambda \big)
	\too
	\bigoplus_{k=1}^m \End (\cal{E} (x_k)) \otimes \Lambda (x_k)
\end{eqn}
is surjective.
Choose nonzero vectors $\lambda_k \in \Lambda (x_k)$.
There is therefore a Higgs field $\Phi_0$ such that $\Phi_0 (x_k) = T_k \otimes \lambda_k$ for every $k$.

For any $\Lambda$-twisted Higgs field $\Phi$, consider the linear map
\begin{eqn}
	B_\Phi : \W \too \H^0 (\X, \omega_\X \otimes \Lambda)
\qtext{given by}
	S \mapstoo \tr (S \circ \Phi)
\fullstop
\end{eqn}
If $B_{\Phi_0} (S) = 0$, evaluation at every $x_k$ gives $\tr (S(x_k) \circ T_k) = 0$.
The choice of the $T_k$ therefore implies $S = 0$, so $B_{\Phi_0}$ is injective.
Since the matrix entries of $B_\Phi$ depend linearly on $\Phi$, injectivity is a nonempty Zariski-open condition on the vector space $\H^0 (\X, \cal{End} (\cal{E}) \otimes \Lambda)$.
Its intersection with $\V$ is again nonempty and Zariski-open because the ambient vector space $\H^0 \big( \X, \cal{End} (\cal{E}) \otimes \Lambda \big)$ is irreducible.

For any $\Phi$ in this intersection, let $(\cal{A}, \rho)$ be the induced spectral structure on $\cal{E}$.
A section $S \in \H^0 (\X, \cal{A}^\perp \otimes \omega_\X)$ pairs to zero with every local section of $\rho (\cal{A})$ under the trace pairing.
In particular, it pairs to zero with every local endomorphism representing $\Phi$, because these local endomorphisms belong to $\rho (\cal{A})$.
Thus $B_\Phi (S) = \tr (S \circ \Phi) = 0$.
The injectivity of $B_\Phi$ implies $S = 0$, proving the required vanishing.
\end{proof}

The construction in the proof of \autoref{1791390497007} gives some helpful additional numerical information.
For any positive-degree line bundle $\Lambda$ such that $\cal{End} (\cal{E}) \otimes \Lambda$ is globally generated, the sections $\Phi$ which lead to a smooth and simply branched spectral cover $\pi : \Y \to \X$ form a nonempty dense open subset of $\H^0 (\X, \cal{End}(\cal{E}) \otimes \Lambda)$.
The discriminant of $\Phi$ is a global section of the line bundle $\Lambda^{n(n-1)}$ and has only simple zeros, one for each simple ramification point.
Consequently, the covers $\Y$ obtained this way satisfy the following numerical constraint on their genus $g_\Y$ and number $r$ of ramification points:
\begin{eqntag}
\label{1791391782181}
	g_\Y = 1 + n(g_\X - 1) + r/2
\qqtext{and}
	r = n (n-1) \deg \Lambda
\fullstop{,}
\end{eqntag}
where $g_\X$ is the genus of $\X$.
Therefore, an increasingly positive twist $\Lambda$ leads to an increasingly ramified smooth spectral realisation of the same vector bundle.

\subsection{Non-Higgs Spectral Structures}
\label{1793634279019}

The constructions in \ref{1791465731919} give spectral structures with a global twisted generator, namely the generating Higgs field.
We now show that every vector bundle of rank at least three also admits smooth, connected, and simply branched spectral structures without a global twisted generator for any choice of twisting line bundle.
The projective-geometric construction in the proof of \autoref{1790777298686} presented in \ref{1791480675169} provides these examples.
We also give a determinant obstruction to the existence of a global twisted generator and illustrate it with two explicit examples.
The second example shows that a non-Higgs spectral structure can also be clean.

\begin{proposition}[existence of non-Higgs spectral structures]
\label{1793629098096}
Let $\cal{E}$ be a vector bundle of rank $n \geq 3$ on $\X$.
Then $\cal{E}$ admits a smooth and connected spectral structure with simply branched spectral cover which cannot be induced by any regular twisted Higgs field, for any choice of twisting line bundle.
\end{proposition}

\begin{proof}
Choose a sufficiently positive line bundle $\cal{N}$ on $\X$, and let
\begin{eqn}
	\pi : \Y \to \X
	\quad,
	\qquad
	\cal{L} \in \Pic (\Y)
	\quad,
	\qquad
	\pi_\ast \cal{L} \simeq \cal{E}
\end{eqn}
be a spectral realisation obtained by the projective-geometric construction in the proof of \autoref{1790777298686} presented in \ref{1791480675169}.
Thus $\Y$ is smooth and connected, and $\pi$ is simply branched.
Write
\begin{eqn}
	\cal{M} \coleq \det (\cal{E} \otimes \cal{N})
	\quad,
	\qquad
	m \coleq \deg \cal{M} > 0
	\quad,
	\qquad
	\cal{A} \coleq \pi_\ast \cal{O}_\Y
\fullstop
\end{eqn}

Recall that $\Y$ is a divisor in $\X \times \PP^1$ with
\begin{eqn}
	\cal{O}_{\X \times \PP^1} (\Y)
	\simeq p_\X^\ast \cal{M} \otimes
	p_{\PP^1}^\ast \cal{O}_{\PP^1} (n)
\fullstop
\end{eqn}
So its ideal sequence is
\begin{eqn}
	0 \too
	p_\X^\ast \cal{M}^{-1} \otimes
	p_{\PP^1}^\ast \cal{O}_{\PP^1} (-n)
	\too \cal{O}_{\X \times \PP^1}
	\too i_\ast \cal{O}_\Y
	\too 0
\fullstop
\end{eqn}
Here $i : \Y \inj \X \times \PP^1$ is the inclusion.
Since
\begin{eqn}
	\H^0 (\PP^1, \cal{O}_{\PP^1} (-n)) = 0
\qqtext{and}
	\H^1 (\PP^1, \cal{O}_{\PP^1}) = 0
\fullstop{,}
\end{eqn}
pushing forward along $p_\X$ gives an exact sequence
\begin{eqn}
	0 \too \cal{O}_\X
	\too \cal{A}
	\too \cal{M}^{-1} \otimes
		\H^1 (\PP^1, \cal{O}_{\PP^1} (-n))
	\too 0
\fullstop
\end{eqn}
But $h^1 (\PP^1, \cal{O}_{\PP^1} (-n)) = n-1$, so
\begin{eqn}
	\cal{A} / \cal{O}_\X
	\simeq (\cal{M}^{-1})^{\oplus (n-1)}
\fullstop
\end{eqn}

Suppose, towards a contradiction, that this spectral structure is induced by a regular $\Lambda$-twisted Higgs field for some line bundle $\Lambda$ on $\X$.
Then, by \autoref{1792167000670}, the algebra $\cal{A}$ admits a global $\Lambda$-twisted generator.
The powers of this generator give an isomorphism
\begin{eqn}
	\cal{O}_\X \oplus \Lambda^{-1} \oplus \cdots
	\oplus \Lambda^{-(n-1)}
	\iso \cal{A}
\end{eqn}
whose restriction to $\cal{O}_\X$ is the unit map.
Passing to the quotient, we obtain
\begin{eqn}
	\cal{A} / \cal{O}_\X
	\cong \Lambda^{-1} \oplus \Lambda^{-2}
	\oplus \cdots \oplus \Lambda^{-(n-1)}
\fullstop
\end{eqn}
Comparing degrees in the two descriptions of $\cal{A} / \cal{O}_\X$ gives
\begin{eqn}
	-(n-1)m
	= - \frac{n(n-1)}{2} \deg \Lambda
	\quad,
	\qquad
	\deg \Lambda = \frac{2m}{n} > 0
\fullstop
\end{eqn}
Since $n \geq 3$, it follows that
\begin{eqn}
	\deg (\Lambda \otimes \cal{M}^{-1})
	= \frac{2-n}{n} m < 0
\fullstop
\end{eqn}
Consequently, $\Hom (\Lambda^{-1}, \cal{M}^{-1}) \cong \H^0 (\X, \Lambda \otimes \cal{M}^{-1}) = 0$.
So there is no nonzero homomorphism
\begin{eqn}
	\Lambda^{-1} \too
	(\cal{M}^{-1})^{\oplus (n-1)}
\fullstop
\end{eqn}
On the other hand, the description of $\cal{A} / \cal{O}_\X$ coming from the twisted generator exhibits $\Lambda^{-1}$ as a direct summand, giving precisely such a homomorphism.
This is a contradiction.

Therefore, $\cal{A}$ admits no global twisted generator for any line bundle $\Lambda$.
By \autoref{1792167000670}, the spectral structure cannot be induced by any regular twisted Higgs field.
\end{proof}

\begin{remark}
\label{1793629534346}
The restriction $n \geq 3$ in \autoref{1793629098096} is essential.
For a degree-two cover, the trace map splits the inclusion $\cal{O}_\X \inj \cal{A}$, giving $\cal{A} \cong \cal{O}_\X \oplus (\cal{A} / \cal{O}_\X)$.
The trace-free summand, viewed as a twisted section of $\cal{A}$, is a global algebra generator.
Thus every rank-two spectral structure is induced by a regular twisted Higgs field.

For $n \geq 3$, the proposition shows that the distinction between spectral structures and Higgs-induced spectral structures persists over every vector bundle, even when the spectral cover is required to be smooth, connected, and simply branched.
\end{remark}

The existence of local regular generators does not imply the existence of a global twisted generator, even if the twisting line bundle is allowed to vary.
There is already a necessary condition at the level of determinants.

\begin{corollary}[a determinant obstruction]
\label{1792167000671}
Suppose a spectral structure $(\cal{A}, \rho)$ on a rank-$n$ bundle is induced by a regular $\Lambda$-twisted Higgs field.
Then
\begin{eqn}
	\det \cal{A} \cong \Lambda^{-n(n-1)/2}
\fullstop
\end{eqn}
If its spectral cover $\pi : \Y \to \X$ is smooth with ramification divisor $\R = \Ram (\pi)$, then
\begin{eqn}
	\deg \R = n (n-1) \deg \Lambda
\fullstop
\end{eqn}
\end{corollary}

\begin{proof}
As an $\cal{O}_\X$-module, the spectral algebra of a $\Lambda$-twisted Higgs field is
\begin{eqn}
	\cal{A} \cong \cal{O}_\X \oplus \Lambda^{-1} \oplus \cdots \oplus \Lambda^{-(n-1)}
\fullstop
\end{eqn}
Taking determinants gives the first assertion because $1 + \cdots + (n-1) = n (n-1)/2$.
For a smooth spectral cover, \autoref{1790716932644} gives $\deg \cal{A} = - \deg \R / 2$, and the second assertion follows.
\end{proof}

\begin{example}[smooth spectral structure without a global twisted generator]
\label{1792167000672}
Let us take $\X = \Y = \PP^1$ and the degree-three cover
\begin{eqn}
	\pi : \PP^1 \too \PP^1
\qtext{given by}
	z \mapstoo z^3
\fullstop
\end{eqn}
Compare with \cite[\S3]{BoulterRayan2025}, where the same cover is used to construct Higgs fields whose spectral curves are singular.

Put $\cal{A} \coleq \pi_\ast \cal{O}_\Y$ and equip $\cal{E} \coleq \cal{A}$ with its natural multiplication action $\rho : \cal{A} \to \cal{End} (\cal{A})$.
Then $(\cal{A}, \rho)$ is a smooth connected spectral structure of rank three.
Its ramification divisor is $\R = 2[0] + 2[\infty]$, so $\deg \R = 4$.
If this spectral structure were induced by a regular $\Lambda$-twisted Higgs field for some line bundle $\Lambda$ on $\PP^1$, then \autoref{1792167000671} would imply $4 = 3 \cdot 2 \cdot \deg \Lambda$.
This is impossible because $\deg \Lambda$ is an integer.

Thus, the smooth spectral structure $(\cal{A}, \rho)$ admits local regular generators by \autoref{1791453151472}, but no global twisted generator for any line bundle $\Lambda$.
Equivalently, its cover cannot be embedded in the total space of any line bundle on $\X$.
Of course, this does not prevent the same vector bundle $\cal{E}$ from carrying other spectral structures induced by regular twisted Higgs fields, as guaranteed by \autoref{1790777298686}.
\end{example}

The construction in \autoref{1793629098096} does not, by itself, establish the cohomological vanishing required for clean spectral structures.
The next example shows that non-Higgs spectral structures can nevertheless be clean.

\enlargethispage{20pt}
\begin{example}[clean spectral structure without a global twisted generator]
\label{1792434634960}
This example shows that a clean spectral structure need not admit a global twisted generator for any choice of twisting line bundle.

Take $\X = \Y = \PP^1$ and the degree-three cover
\begin{eqn}
	\pi : \PP^1 \too \PP^1
\qtext{given by}
	z \mapstoo z^2 + \frac{1}{z}
\fullstop
\end{eqn}
Put $\cal{A} \coleq \pi_\ast \cal{O}_\Y$ and $\cal{L} \coleq \cal{O}_\Y (2)$, and equip $\cal{E} \coleq \pi_\ast \cal{L}$ with the canonical spectral action $\rho : \cal{A} \to \cal{End} (\cal{E})$.
We will show that this spectral structure is clean, but admits no global twisted generator for any line bundle on $\X$.

The spectral realisation is smooth and connected because $\Y = \PP^1$.
To examine its branching, observe that
\begin{eqn}
	\pi' (z) = \frac{2z^3 - 1}{z^2}
\fullstop
\end{eqn}
Thus the three roots $z_1, z_2, z_3$ of $z^3 = 1/2$ are simple ramification points.
Their images are $\pi (z_k) = 3z_k^2$, which are pairwise distinct.
There is one further ramification point at $\infty$, with
ramification index two, because $\pi$ has a pole of order two there.
The pole at $0$ is simple, so $\pi$ is unramified there.
Consequently, $\pi$ is simply branched and its ramification divisor is
\begin{eqn}
	\R = [z_1] + [z_2] + [z_3] + [\infty]
\qqtext{and hence}
	\deg \R = 4
\fullstop
\end{eqn}

We next show that $\cal{E} \simeq \cal{O}_\X^{\oplus 3}$.
Since $\pi$ is finite, Riemann-Roch gives
\begin{eqn}
	\deg \cal{E} + 3 = \chi (\X, \cal{E}) = 	\chi (\Y, \cal{L}) = 3
\fullstop
\end{eqn}
Thus $\deg \cal{E} = 0$.
On the other hand, the projection formula and the fact that $\pi^\ast \cal{O}_\X (1) \simeq \cal{O}_\Y (3)$ give
\begin{eqn}
	\H^0 \big( \X, \cal{E} (-1) \big)
	\simeq
	\H^0 \big( \Y, \cal{L} \otimes \pi^\ast \cal{O}_\X (-1) \big)
	\simeq
	\H^0 \big( \Y, \cal{O}_\Y (-1) \big)
	= 0
\fullstop
\end{eqn}
By the Birkhoff-Grothendieck Theorem, we can write 
\begin{eqn}
	\cal{E} \simeq \cal{O}_\X (k_1) \oplus \cal{O}_\X (k_2) \oplus \cal{O}_\X (k_3)
\fullstop
\end{eqn}
The vanishing above implies $k_i \leq 0$ for every $i$, whilst $\deg \cal{E} = 0$ implies $k_1 + k_2 + k_3 = 0$.
Therefore, all three integers are zero.
It follows that
\begin{eqn}
    \H^1 \big( \X, \cal{End} (\cal{E}) \big)
    \simeq
    \H^1 \big( \PP^1, \cal{O}_{\PP^1}^{\oplus 9} \big)
    = 0
\fullstop
\end{eqn}
The off-diagonal exact sequence then gives $\H^1 (\X, \cal{Q}) = 0$.
Thus $(\cal{A}, \rho)$ is clean.

Nevertheless, if this spectral structure were induced by a regular $\Lambda$-twisted Higgs field for some line bundle $\Lambda$ on $\X$, then \autoref{1792167000671} would imply $4 = \deg \R = 3 \cdot 2 \cdot \deg \Lambda$.
This is impossible because $\deg \Lambda$ is an integer.
Consequently, this clean spectral structure admits no global twisted generator for any line bundle $\Lambda$.
Equivalently, its cover cannot be embedded over $\X$ in the total space of any line bundle on $\X$.

Of course, the same vector bundle $\cal{E}$ admits other clean spectral structures induced by regular twisted Higgs fields, as constructed in the proof of \autoref{1792430412943}.
\end{example}

\subsection{The Off-Diagonal Form}
\label{1792146586506}

By \autoref{1791453955298}, the trace pairing identifies $\cal{A}^\perp$ canonically with $\cal{Q}^\vee$.
On the other hand, $\cal{Q}$ itself is in general not isomorphic to $\cal{A}^\perp$.
A local choice of a regular generator of $\cal{A}$ provides an identification between $\cal{Q}$ and $\cal{A}^\perp$ through the commutator map, but this identification is only local and depends on this choice of a local generator.
A regular $\Lambda$-twisted Higgs field supplies additional data: its commutator induces a global isomorphism \begin{eqn}
	\cal{Q} \iso \cal{A}^\perp \otimes \Lambda
\fullstop
\end{eqn}
The following proposition describes the resulting perfect alternating $\Lambda$-valued form on $\cal{Q}$.

\begin{proposition}[off-diagonal form of a regular Higgs field]
\label{1791458921327}
Suppose $(\cal{E}, \Phi)$ is a $\Lambda$-twisted regular Higgs bundle of rank $n$, and let $(\cal{A}, \rho)$ be the induced spectral structure.
Then the off-diagonal bundle $\cal{Q} = \cal{End} (\cal{E}) / \rho (\cal{A})$ acquires a perfect $\Lambda$-valued pairing
\begin{eqntag}
\label{1791459765927}
	\eta : \cal{Q} \otimes \cal{Q} \too \Lambda
\qtext{given by}
	\eta \big( [T_1], [T_2] \big) \coleq \tr \big( \Phi \circ [T_1, T_2] \big)
\fullstop
\end{eqntag}
This pairing and its determinant determine bundle isomorphisms
\begin{eqntag}
\label{1791577928214}
	\cal{Q} \iso \cal{Q}^\vee \otimes \Lambda
\qqtext{and}
	\det \cal{Q} \iso \Lambda^{n(n-1)/2}
\fullstop
\end{eqntag}
Consequently, the degree and Euler characteristic of $\cal{Q}$ can be expressed in terms of the twist $\Lambda$:
\begin{eqntag}
\label{1791459662959}
	\deg \cal{Q} = \tfrac{1}{2} n (n-1) \deg \Lambda
\qqtext{and}
	\chi (\cal{Q}) = \tfrac{1}{2} n (n-1) ( \deg \Lambda - 2 g_\X + 2 )
\fullstop
\end{eqntag}
\end{proposition}

Here, $T_1$ and $T_2$ are local endomorphisms representing the classes $[T_1]$ and $[T_2]$ in $\cal{Q}$.
Also, to be explicit, the isomorphism $\cal{Q} \iso \cal{Q}^\vee \otimes \Lambda$ is given by 
\begin{eqn}
	[T] \mapstoo \eta \big( [T], \dummy \big) = \Big( S \mapsto \eta ([T], [S]) \Big)
\fullstop
\end{eqn}

\begin{proof}
First, take a trivialising neighbourhood $\U \subset \X$ for $\Lambda$ and pick a generator $e \in \Lambda(\U)$.
Let $T$ be the endomorphism representing the Higgs field $\Phi$ in this trivialisation; i.e., $\Phi_\U = T \otimes e$.
Then the spectral action has image $\cal{O}_\U [T] \subset \cal{End} (\cal{E})$, and every section of $\rho (\cal{A}_\U)$ commutes with $T$.

For local endomorphisms $S, R$ of $\cal{E}_\U$, define
\begin{eqn}
	\bar{\eta} (S,R)
	\coleq
	\tr \big( T \circ [S,R] \big) e
\fullstop
\end{eqn}
Bilinearity and skewsymmetry of the commutator and linearity of the trace make $\bar{\eta}$ bilinear and skew-symmetric:
\begin{eqn}
	\bar{\eta} (S,R) = - \bar{\eta} (R,S)
\qtext{and}
	\bar{\eta} (S + S', R) = \bar{\eta}(S,R) + \bar{\eta} (S',R)
\fullstop
\end{eqn}
We claim that $\bar{\eta}$ depends only on the equivalence classes of $S$ and $R$ in $\cal{Q}_\U$.
Indeed, cyclicity of the trace gives 
\begin{eqn}
	\tr \big( T \circ [S, R] \big)
	=
	\tr \big( [T,S] \circ R \big)
	=
	- \tr \big( [T,R] \circ S \big)
\fullstop
\end{eqn}
So if either $S$ or $R$ belongs to $\rho (\cal{A}_\U)$, then either $[T,S] = 0$ or $[T,R] = 0$, so this expression vanishes.
Consequently, $\bar{\eta}$ descends to the quotient where it defines a skew-symmetric bilinear form $\eta : \cal{Q}_\U \otimes \cal{Q}_\U \to \Lambda$.

This definition is independent of the chosen local generator $e$ of $\Lambda$.
Indeed, if $e' = he$ is another generator for some nonvanishing function $h \in \cal{O}_\U^\times$, then the endomorphism $T'$ representing $\Phi$ in this trivialisation is $T' = h^{-1} T$.
Consequently, 
\begin{eqn}
	\tr \big( T \circ [S,R] \big) e = \tr \big( T' \circ [S,R] \big) e'
\fullstop{,}
\end{eqn}
so $\bar{\eta}$ and hence $\eta$ are well-defined without choosing a local generator of $\Lambda$.
This also means that these local forms glue together into a global form $\eta$ exactly as in \eqref{1791459765927}.

We claim that this form is nondegenerate, which is enough to check on fibres.
Fix $x \in \X$ and use $e(x)$ to identify the fibre $\Lambda (x)$ with $\CC$.
Suppose an endomorphism $S \in \End (\cal{E}(x))$ pairs to zero with every endomorphism $R \in \End (\cal{E} (x))$.
Then we must have 
\begin{eqn}
	\tr \big( [T(x), S] \circ R \big) = 0
\qquad
	\forall R \in \End (\cal{E}(x))
\fullstop
\end{eqn}
But the trace pairing is perfect, so $[T(x), S] = 0$.
Since the centraliser of the regular endomorphism $T(x)$ is the algebra 
\begin{eqn}
	\CC [T(x)] = \rho (\cal{A}(x))
\fullstop{,}
\end{eqn}
the endomorphism $S$ belongs to the spectral algebra $\rho (\cal{A}(x))$.
So $[S] = 0$ in $\cal{Q} (x)$.
So $\eta$ is nondegenerate on every fibre, and this shows that it is a perfect pairing and hence defines the isomorphism \eqref{1791577928214}.

Finally, $\cal{Q}$ has rank $n (n-1)$.
Since the skew-symmetric form $\eta$ is nondegenerate, its top nonzero exterior power defines a nowhere-vanishing global section
\begin{eqn}
	\eta^{\wedge \rm{top}} 
	= \eta^{\wedge n (n-1)/2} \in \H^0 \big(\X, (\det \cal{Q})^\vee \otimes \Lambda^{n(n-1)/2} \big)
\fullstop
\end{eqn}
So it defines an isomorphism $\det \cal{Q} \cong \Lambda^{n(n-1)/2}$, which implies the first formula in \eqref{1791459662959}, and the second formula follows from the Riemann-Roch formula.
\end{proof}

By Serre duality, this proposition yields an interesting consequence for spectral structures induced by \textit{classical} Higgs fields; i.e., those valued in the canonical bundle $\Lambda = \omega_\X$.

\begin{corollary}[balanced off-diagonal cohomology]
\label{1791459820103}
If a spectral structure on the vector bundle $\cal{E}$ is induced by a regular classical Higgs field 
\begin{eqn}
	\Phi \in \H^0 \big(\X, \cal{End} (\cal{E}) \otimes \omega_\X \big)
\fullstop{,}
\end{eqn}
then the off-diagonal bundle $\cal{Q}$ has balanced cohomology:
\begin{eqntag}
\label{1791459900368}
	\H^1 (\X, \cal{Q})^\vee \cong \H^0 (\X, \cal{Q})
\qtext{and so}
	h^0 (\X, \cal{Q}) = h^1 (\X, \cal{Q})
\fullstop
\end{eqntag}
\end{corollary}

\begin{proof}
Serre duality and the isomorphism \eqref{1791577928214} of \autoref{1791458921327} with $\Lambda = \omega_\X$ give:
\begin{eqn}
	\H^1 (\X, \cal{Q})^\vee \cong \H^0 (\X, \cal{Q}^\vee \otimes \omega_\X) \cong \H^0 (\X, \cal{Q})
\fullstop
\qedhere
\end{eqn}
\end{proof}

\section{Spectral Structures in Families}
\label{1791914002217}

We now extend the definitions and constructions of the preceding sections to holomorphic families over arbitrary analytic spaces, allowing the base curve itself to vary and the parameter space to be nonreduced.
We begin with the Relative Finite Sheaf-Module Correspondence, which explains how the correspondence of \ref{1792143161338} (namely, \autoref{1793046131649}) behaves under arbitrary changes of parameter space.
We then apply it to spectral structures and their realisations, discuss the fibrewise interpretation of smoothness, connectedness, and cleanliness, and prove local existence of clean spectral structures.

Let $\S$ be any complex analytic space, possibly nonreduced, and let
\begin{eqn}
    p : \X \too \S
\end{eqn}
be a proper smooth surjective morphism of relative dimension one with connected fibres.
Thus, in this section, $\X$ denotes the total space of a family of compact Riemann surfaces
\begin{eqn}
    \X_s \coleq \X \times_\S \set{s}
\qqtext{for every}
    s \in \S
\fullstop
\end{eqn}
For a sheaf $\cal{F}$ on $\X$, denote by $\cal{F}_s$ its restriction to $\X_s$.
For an open subset $\U \subset \S$, put $\X_\U \coleq p^{-1} (\U)$ and write $\cal{F}_\U$ for the restriction of $\cal{F}$ to $\X_\U$.

An important nuance is that identities between morphisms over $\X$ must hold as equalities of sheaf morphisms rather than merely when restricted to every fibre $\X_s$.
The reason is that restricting to fibres is problematic because such restrictions may annihilate nonzero terms involving nilpotent functions on $\S$.
On the other hand, a morphism between vector bundles on $\X$ is an isomorphism whenever it is a vector space isomorphism on the fibre at every point $x \in \X$: its determinant is then a unit in every local ring.

\subsection{Relative Finite Sheaf-Module Correspondence}
\label{1793124287763}

The Finite Sheaf-Module Correspondence (\autoref{1793046131649}) applies to the total spaces of analytic families, including singular and nonreduced total spaces.
We now show that it is compatible with arbitrary changes of parameter space.
This compatibility is what allows spectral realisations of a family of vector bundles to pull back to spectral realisations of the pulled-back family.
The following result holds for any holomorphic map $p : \X \to \S$, without the properness or smoothness assumptions imposed above.

\begin{proposition}[relative finite sheaf-module correspondence]
\label{1793046154727}
Let $p : \X \to \S$ be any holomorphic map of complex analytic spaces, let $\pi : \Y \to \X$ be finite, and put $\cal{A} \coleq \pi_\ast \cal{O}_\Y$.
Given any holomorphic map $g : \S' \to \S$, form $\X' \coleq \X \times_\S \S'$ and $\Y' \coleq \Y \times_\S \S' \cong \Y \times_\X \X'$.
Write the induced Cartesian square as
\begin{eqn}
\begin{tikzcd}
	\Y'
		\ar[r, "\tilde{f}"]
		\ar[d, "\pi'"']
&
	\Y
		\ar[d, "\pi"]
\\
	\X'
		\ar[r, "f"']
&
	\X
\fullstop
\end{tikzcd}
\end{eqn}
Then $\pi'$ is finite, and there are canonical isomorphisms
\begin{eqn}
	\cal{A}' \coleq f^\ast \cal{A}
	\iso
	\pi'_\ast \cal{O}_{\Y'}
\qqtext{and}
	\Y'
	\cong
	\Spec_{\X'} (\cal{A}')
\fullstop
\end{eqn}
For every coherent sheaf $\cal{L}$ on $\Y$, the natural base-change morphism is an isomorphism
\begin{eqn}
	b :
	f^\ast \pi_\ast \cal{L}
	\iso
	\pi'_\ast \tilde{f}^\ast \cal{L}
\fullstop
\end{eqn}
It is $\cal{A}'$-linear under the preceding identification of algebras.
Consequently, the Sheaf-Module Correspondence commutes with arbitrary base change: for every $\cal{E} \in \mathbf{Coh} (\X, \cal{A})$, there is a canonical isomorphism
\begin{eqn}
	\tilde{f}^\ast \pi^\star \cal{E}
	\cong
	\pi'^\star f^\ast \cal{E}
\fullstop
\end{eqn}
These identifications are natural and compatible with composition of base changes.
Moreover, a coherent sheaf $\cal{L}$ on $\Y$ is flat over $\S$ if and only if $\pi_\ast \cal{L}$ is flat over $\S$.
\end{proposition}

\begin{proof}
Finiteness is preserved by base change.
The finite analytic base change in the case of the structure sheaf gives the canonical algebra isomorphism
\begin{eqn}
	f^\ast \cal{A}
	=
	f^\ast \pi_\ast \cal{O}_\Y
	\iso
	\pi'_\ast \cal{O}_{\Y'}
\fullstop
\end{eqn}
See \cite[Corollary~2.8.3]{GuiCAG} and \cite[Ch.~III, \S1(b), p.~94]{BanicaStanasila1976}.
The identification of $\Y'$ with $\Spec_{\X'} (\cal{A}')$ then follows from the finite analytic-spectrum correspondence.
We explain how this algebra isomorphism gives the base-change assertion for every coherent sheaf.
The morphism $b$ is obtained by pulling back local sections of $\cal{L}$.
Its compatibility with multiplication by local functions on $\Y$ shows that it is $\cal{A}'$-linear.

Locally on $\X$, the module $\pi_\ast \cal{L}$ admits a finite presentation over $\cal{A}$.
Applying the exact equivalence of \autoref{1793046131649} gives, over the preimage of that neighbourhood, a finite presentation
\begin{eqn}
	\cal{O}_\Y^{\oplus a}
	\too
	\cal{O}_\Y^{\oplus b}
	\too
	\cal{L}
	\too 0
\fullstop
\end{eqn}
Both functors $\cal{L} \mapsto f^\ast \pi_\ast \cal{L}$ and $\cal{L} \mapsto \pi'_\ast \tilde{f}^\ast \cal{L}$ are right exact: pullback is right exact, and finite direct image is exact.
The natural comparison between them is an isomorphism on the finite free sheaves in this presentation, by the structure-sheaf case.
It is therefore an isomorphism on their cokernels.
This proves the assertion for $\cal{L}$ without a flatness assumption.

Apply this isomorphism to $\cal{L} = \pi^\star \cal{E}$.
Using $\pi_\ast \pi^\star \cal{E} \cong \cal{E}$ gives an $\cal{A}'$-linear isomorphism
\begin{eqn}
	\pi'_\ast \tilde{f}^\ast \pi^\star \cal{E}
	\cong
	f^\ast \cal{E}
\fullstop
\end{eqn}
The equivalence for $\pi'$ therefore identifies $\tilde{f}^\ast \pi^\star \cal{E}$ with $\pi'^\star f^\ast \cal{E}$.
Naturality and compatibility with composition follow from the canonical pullback maps and associativity of tensor products.
Finally, if $x \in \X$ and $s = p(x)$, the stalk decomposition in the proof of \autoref{1793046131649} is an isomorphism of $\cal{O}_{\S,s}$-modules:
\begin{eqn}
	(\pi_\ast \cal{L})_x
	\cong
	\bigoplus_{y \in \pi^{-1} (x)}
	\cal{L}_y
\fullstop
\end{eqn}
A finite direct sum is flat if and only if every summand is flat.
\end{proof}

\subsection{Spectral Structures, Realisations, and Base Change}
\label{1792146709610}

We now apply the relative correspondence to spectral structures and their realisations.
The definitions and constructions of \smartref{1790715973901} and \smartref{1790944507299} apply to vector bundles on the total space $\X$, and \autoref{1793046154727} provides their compatibility with arbitrary base change.

Thus, for example, a \dfn{spectral structure} on a rank-$n$ vector bundle $\cal{E}$ over the analytic family $\X$ is a pair $(\cal{A}, \rho)$ consisting of a locally free commutative unital $\cal{O}_\X$-algebra $\cal{A}$ and an algebra homomorphism $\rho : \cal{A} \to \cal{End} (\cal{E})$ that makes $\cal{E}$ a locally free rank-one $\cal{A}$-module.
For each $s \in \S$, it restricts to a spectral structure $(\cal{A}_s, \rho_s)$ on the vector bundle $\cal{E}_s$ over $\X_s$.

Likewise, the \dfn{spectral realisation} of $(\cal{A}, \rho)$ consists of a finite flat degree-$n$ map
\begin{eqntag}
\label{1791476831216}
    \pi:\Y\coleq\Spec_\X(\cal{A})\too\X
\qqtext{such that}
    \pi_\ast\cal{O}_\Y\cong\cal{A}
\fullstop{,}
\end{eqntag}
together with a line bundle $\cal{L}$ on $\Y$ and an $\cal{A}$-linear isomorphism $\iota : \pi_\ast \cal{L} \iso \cal{E}$.
The composite projection $q \coleq p \circ \pi : \Y \to \S$ is proper and flat because $\pi$ is finite flat and $p$ is proper and smooth, hence flat.
Each fibre $\Y_s \coleq \Y \times_\S \set{s}$ is a compact curve, possibly singular, nonreduced, or disconnected, with a finite flat degree-$n$ map $\pi_s : \Y_s \to \X_s$.
Also, neither family need be holomorphically locally trivial over $\S$.
These constructions are compatible with arbitrary base change in $\S$, as we prove in \autoref{1791551111381} below.

However, numerical assertions (such as genus, Euler characteristic, or degree) as well as assertions about cohomology, smoothness, connectedness, or cleanliness, must be interpreted fibrewise on $\S$; see \smartref{1792146750001} for more information.
In particular, the explicit descriptions of ramified fibres in \smartref{1791375997197} must also be understood in this sense.
The corresponding existence result for families is local on $\S$ and will be proved in \autoref{1791544458956}.
The final exception concerns the observation following \autoref{1791434155987} about split spectral structures on vector bundles over $\PP^1$.
That observation employs the Birkhoff-Grothendieck Theorem which guarantees a decomposition of each vector bundle $\cal{E}_s$ over $\PP^1 \times \set{s} \cong \PP^1$, but does not guarantee a decomposition of the vector bundle $\cal{E}$ even over the trivial family $\PP^1 \times \S$.

\clearpage
\begin{corollary}[compatibility with base change]
\label{1791551111381}
Given a holomorphic map $g : \S' \to \S$, form the base change $\X'$ along $g$ with projections
\begin{eqn}
	p' : \X' \coleq \X \times_\S \S' \too \S'
\qqtext{and}
	f : \X' \too \X
\fullstop
\end{eqn}
If $(\cal{A},\rho)$ is a spectral structure on $\cal{E}$ over $\X$, then $(\cal{A}', \rho') \coleq (f^\ast \cal{A}, f^\ast \rho)$ is a spectral structure on $\cal{E}'\coleq f^\ast\cal{E}$ over $\X'$.
Its spectral realisation is canonically identified with
\begin{eqn}
	\pi' : \Y' \coleq \Y \times_\S \S' \too \X'
\qtext{and}
	\cal{L}' \coleq \tilde{f}^\ast \cal{L}
\fullstop{,}
\end{eqn}
where $\tilde{f} : \Y' \to \Y$ is the induced projection.
Under this identification, $\iota'$ is the composite
\begin{eqn}
	\Big( \iota' : \pi'_\ast \cal{L}' \iso \cal{E}' \Big)
	=
	\Big(
    \pi'_\ast \tilde{f}^\ast \cal{L}
    \iso
    f^\ast \pi_\ast \cal{L}
	\tinyunderset{f^\ast \iota}{\iso}
    f^\ast \cal{E}
    \Big)
\fullstop{,}
\end{eqn}
where the first map is the finite base-change isomorphism.
Finally, the off-diagonal bundle $\cal{Q}'$ of $(\cal{A}', \rho')$ is canonically identified with the pullback of the off-diagonal bundle $\cal{Q}$ of $(\cal{A}, \rho)$:
\begin{eqn}
	\cal{Q}' \cong f^\ast \cal{Q}
\fullstop
\end{eqn}
\end{corollary}

In particular, taking $\S' = \set{s}$ for any $s \in \S$ gives canonical isomorphisms
\begin{eqn}
	\Y_s \cong \Spec_{\X_s} (\cal{A}_s)
\fullstop{,}
	\qqqquad
	\cal{A}_s \cong (\pi_s)_\ast \cal{O}_{\Y_s}
\fullstop{,}
	\qqqquad
	\cal{E}_s \cong (\pi_s)_\ast \cal{L}_s
\fullstop
\end{eqn}
Likewise, $\cal{Q}_s$ is the off-diagonal bundle of $(\cal{A}_s, \rho_s)$ over $\X_s$.

\begin{proof}
Pulling back a local $\cal{A}$-module trivialisation $\cal{E} \cong \cal{A}$ gives a local $\cal{A}'$-module trivialisation $\cal{E}' \cong \cal{A}'$.
Thus, $(\cal{A}', \rho')$ is a spectral structure on $\cal{E}'$.

By \autoref{1793046154727}, the relative spectrum of $\cal{A}'$ is canonically $\Y' = \Y \times_\S \S'$, and the pulled-back spectral line bundle is $\cal{L}' = \tilde{f}^\ast \cal{L}$.
Writing $b : f^\ast \pi_\ast \cal{L} \iso \pi'_\ast \tilde{f}^\ast \cal{L}$ for the canonical finite base-change isomorphism, the induced identification with $\cal{E}'$ is $\iota' = f^\ast \iota \circ b^{-1}$.
This is $\cal{A}'$-linear, and since $\pi'_\ast$ is fully faithful it gives the asserted canonical identification of spectral realisations.

Finally, the defining sequence \eqref{1789348440662} for $\cal{Q}$ is locally split by the proof of \autoref{1789332843734}.
Pullback therefore preserves its exactness, and its cokernel is both $f^\ast \cal{Q}$ and $\cal{Q}'$.
\end{proof}

There is also a useful compatibility with tensoring by line bundles from the parameter space.

\begin{proposition}[twisting by a line bundle from parameter space]
\label{1791564860452}
Let $\cal{M}$ be a line bundle on $\S$, and put $\cal{N} \coleq p^\ast \cal{M}$.
If $(\cal{A}, \rho)$ is a spectral structure on a vector bundle $\cal{E}$ over $\X$ then the canonical identification $\cal{End} (\cal{E}) \cong \cal{End} (\cal{E} \otimes \cal{N})$ transports the spectral structure $(\cal{A}, \rho)$ to the vector bundle $\cal{E} \otimes \cal{N}$.
If $\cal{L}$ is the spectral line bundle over $\Y$ for $\cal{E}$, then $\cal{L} \otimes \pi^\ast \cal{N}$ is the spectral line bundle for $\cal{E} \otimes \cal{N}$.
\end{proposition}

\begin{proof}
The canonical identification of endomorphism bundles sends a local endomorphism $T$ of $\cal{E}$ to $T \otimes \id_\cal{N}$.
Thus, the transported action is $\rho (a) \otimes \id_\cal{N}$.
Locally trivialising $\cal{N}$ identifies $\cal{E} \otimes \cal{N}$ with $\cal{E}$ as $\cal{A}$-modules, so the former is still locally free of rank one over $\cal{A}$.

For the second assertion, the projection formula gives a chain of $\cal{A}$-linear isomorphisms
\begin{eqn}
	\pi_\ast \big( \cal{L} \otimes \pi^\ast \cal{N} \big)
	\cong
	\pi_\ast \cal{L} \otimes \cal{N}
	\iso
	\cal{E} \otimes \cal{N}
\fullstop{,}
\end{eqn}
where the second arrow is $\iota \otimes \id_\cal{N}$.
The action on the last term is precisely the transported action above.
Consequently, the finite sheaf-module correspondence identifies its spectral line bundle with $\cal{L} \otimes \pi^\ast \cal{N}$.
\end{proof}

\subsection{Smoothness, Connectedness, and Cleanliness}
\label{1792146750001}

The smoothness and connectedness conditions of \ref{1791480675169} as well as the cleanliness condition of \ref{1792426526441} require a relative interpretation.

\begin{definition}[fibrewise smoothness, connectedness, and cleanliness]
\label{1791482794969}
A spectral structure $(\cal{A}, \rho)$ on $\cal{E}$ over $\X$ is called \dfn{smooth} if every fibre $\Y_s$ of its spectral realisation $\Y$ is smooth; \dfn{connected} if every fibre $\Y_s$ is connected; \dfn{clean} if $\H^1 (\X_s, \cal{Q}_s) = 0$ and $\Y_s$ is smooth, connected, and simply branched for every $s \in \S$.
We then say that $\Y$ enjoys these properties \dfn{fibrewise} over $\S$.
\end{definition}

We stress that these conditions concern only the fibres of $q : \Y \to \S$ rather than the total space $\Y$ itself.
In particular, the total space of a fibrewise smooth and/or connected family may be singular, nonreduced, or disconnected: this may happen when $\S$ is singular, nonreduced, or disconnected; see \autoref{1791541569294}.
Nor does the reverse implication hold: the total space $\Y$ may be smooth and connected despite having singular or disconnected fibres; see \autoref{1791542183202}.

\begin{lemma}[fibrewise and relative smoothness]
\label{1791493100668}
The spectral realisation $\Y$ is fibrewise smooth if and only if the map $q : \Y \to \S$ is smooth of relative dimension one.
\end{lemma}

\begin{proof}
Concretely, $q$ being smooth of relative dimension one means that it is locally a projection $\DD \times \U \to \U$ where $\DD$ is an open disc and $\U$ is open in $\S$.
Its fibres are therefore smooth curves.

Conversely, the map $q = p \circ \pi$ is flat because both $\pi$ and $p$ are flat.
By assumption, every fibre of $q$ is smooth of pure dimension one.
The general analytic smoothness criterion then says that $q$ is smooth of relative dimension one; see \cite[Théorème 3.1, (i) $\Rightarrow$ (iii)]{1791539610207}.
\end{proof}

It follows that, if $\S$ is smooth and $\Y$ is fibrewise smooth, then $\Y$ is smooth.
Indeed, the local product description in the proof is then a product of two complex manifolds.
Likewise, if $\S$ is connected and $\Y$ is fibrewise connected, then $\Y$ is connected.
For otherwise, a decomposition of $\Y$ into two nonempty disjoint open and closed subsets would give such a decomposition of $\S$: their images under $q$ are closed because $q$ is proper, disjoint because every fibre is connected, and cover $\S$ because $q$ is surjective.

The normality criterion for smoothness in \autoref{1789157786307} can be applied fibrewise: $\Y$ is fibrewise smooth if and only if each $\cal{O}_{\X_s}$-algebra $\cal{A}_s$ is normal.
Yet normality of the algebra $\cal{A}$ is neither necessary nor sufficient for the fibrewise smoothness of $\Y$, as the following two examples illustrate.

\begin{example}[smooth fibres with nonreduced total space]
\label{1791541569294}
Let $\pi_0 : \Y_0 \to \X_0$ be a finite cover of smooth compact complex curves, with $\X_0$ connected, and take the fat point 
\begin{eqn}
	\S \coleq \Spec \big( \CC [\epsilon] / \inner{\epsilon^2} \big)
\fullstop
\end{eqn}
Put $\X \coleq \X_0 \times \S$ and $\Y \coleq \Y_0 \times \S$, with the induced map $\pi \coleq \pi_0 \times \id_\S$.
The projection $q : \Y \to \S$ is smooth of relative dimension one because every point of $\Y$ has a neighbourhood $\V \times \S$ where $\V \subset \Y_0$ is an open disc neighbourhood.
Every fibre is the smooth curve $\Y_0$.
But the total space $\Y$ is nonreduced because its local rings are
\begin{eqn}
	\cal{O}_{\Y,y}
	\cong
	\cal{O}_{\Y_0,y} [\epsilon] / \inner{\epsilon^2}
\fullstop
\end{eqn}
Thus, $\Y$ is fibrewise smooth even though $\cal{A} = \pi_\ast \cal{O}_\Y$ is not normal.
\end{example}

\begin{example}[smooth total space with a singular fibre]
\label{1791542183202}
Take $\S = \CC$ and $\X = \PP^1 \times \S$, with $p$ the projection to $\S$.
Let $[z_0 : z_1]$ be homogeneous coordinates on $\PP^1$, and let $\Lambda$ be the pullback of $\cal{O}_{\PP^1} (1)$ to $\X$.
In the total space of $\Lambda$, consider the family of double covers
\begin{eqn}
	\Y : \qquad w^2 = z_1^2 + s z_0^2
\fullstop
\end{eqn}
On the chart $z_0 \neq 0$, its equation is $w^2 = z^2 + s$.
The derivative of the defining equation with respect to $s$ is nonzero, so $\Y$ is smooth on this chart.
On the chart $z_1 \neq 0$, the equation is $v^2 = 1 + s u^2$.
At the remaining points, where $u = 0$, we have $v = \pm 1$, so the derivative with respect to $v$ is nonzero.
Consequently, the total space $\Y$ is smooth everywhere, and $\cal{A} = \pi_\ast \cal{O}_\Y$ is normal.
However, the fibre at $s = 0$ has a node because its local equation near $z = w = 0$ is $(w - z) (w + z) = 0$.
Thus, this family is not fibrewise smooth.
\end{example}

For fibrewise smooth covers, the relative analogue of the recognition principle in \autoref{1791408088098} is as follows.

\begin{lemma}[recognition in families]
\label{1792143470690}
Let $\pi : \Y \to \X$ be finite flat of degree $n$, with $\Y$ fibrewise smooth over $\S$, and put $\cal{A} \coleq \pi_\ast \cal{O}_\Y$.
Suppose a vector bundle $\cal{E}$ of rank $n$ on $\X$ is equipped with an algebra homomorphism $\rho : \cal{A} \to \cal{End} (\cal{E})$.
If $\rho$ is $\S$-fibrewise injective (i.e., if every restriction $\rho_s : \cal{A}_s \to \cal{End} (\cal{E}_s)$ is injective), then $\cal{E}$ is locally free of rank one over $\cal{A}$.
Moreover, if $\Y$ is fibrewise connected, then every algebra homomorphism $\rho : \cal{A} \to \cal{End} (\cal{E})$ is automatically $\S$-fibrewise injective.

In either case, there is a line bundle $\cal{L}$ on $\Y$ with an $\cal{A}$-linear isomorphism $\iota : \pi_\ast \cal{L} \iso \cal{E}$.
The pair $(\cal{L}, \iota)$ is unique up to a unique isomorphism compatible with $\iota$, and its construction commutes with arbitrary base change in $\S$.
\end{lemma}

\begin{proof}
By finite base change, $\cal{A}_s \cong (\pi_s)_\ast \cal{O}_{\Y_s}$.
Since $\Y_s$ is smooth, \autoref{1791408088098} applies to $\pi_s : \Y_s \to \X_s$ and the action $\rho_s$.
It follows that $\cal{E}_s$ is locally free of rank one over $\cal{A}_s$.
The same proposition also shows that $\rho_s$ is automatically injective when $\Y_s$ is connected.

We now claim that these fibrewise rank-one conditions imply the rank-one condition over $\X$ itself.
Fix $x \in \X$, put $s \coleq p(x)$, and choose a local $\cal{A}_s$-module generator $e_s$ of $\cal{E}_s$ near $x$.
Lift its value at $x$ to a local section $e$ of $\cal{E}$ near $x$.
The action defines an $\cal{A}$-linear map
\begin{eqn}
	\varepsilon_e : \cal{A} \too \cal{E}
\qtext{sending}
	a \mapstoo \rho (a) e
\qquad\text{locally near $x \in \X$}
\fullstop
\end{eqn}
Both sides are locally free of rank $n$ over $\cal{O}_\X$.
Moreover, on the vector space fibre at $x$, this map is the isomorphism determined by the generator $e_s (x)$.
Thus, its determinant is nonzero at $x$, and hence is a unit in the local analytic ring $\cal{O}_{\X,x}$.
After possibly shrinking the neighbourhood of $x$, the map $\varepsilon_e$ is an isomorphism.
Since it is $\cal{A}$-linear and since $x$ was arbitrary, it gives a local $\cal{A}$-module trivialisation of $\cal{E}$.
This proves the claim.

The finite sheaf-module correspondence now gives a line bundle $\cal{L}$ on $\Y$ together with an $\cal{A}$-linear isomorphism $\iota : \pi_\ast \cal{L} \iso \cal{E}$.
If another pair $(\cal{L}', \iota')$ has this property, then the $\cal{A}$-linear isomorphism $(\iota')^{-1} \circ \iota$ lifts to a unique isomorphism $\lambda : \cal{L} \iso \cal{L}'$ satisfying $\iota = \iota' \circ \pi_\ast \lambda$.
Finally, compatibility with arbitrary base change follows from \autoref{1791551111381}.
\end{proof}

\subsection{Higgs Bundles in Families}
\label{1793045933301}

All contents of \ref{1792146546398} and \ref{1792146586506} carry over to families verbatim with the following modifications.
For classical Higgs fields, the canonical bundle is replaced by the relative canonical bundle $\Lambda = \omega_{\X/\S} = \Omega^1_{\X/\S}$, and the balanced cohomology statement of \autoref{1791459820103} remains true fibrewise only:
\begin{eqn}
	\H^1 (\X_s, \cal{Q}_s)^\vee
	\cong
	\H^0 (\X_s, \cal{Q}_s)
\qqtext{so}
	h^0 (\X_s, \cal{Q}_s) = h^1 (\X_s, \cal{Q}_s)
\qquad
	\text{for all $s \in \S$}
\fullstop
\end{eqn}

Thus, for example, we can fix a line bundle $\Lambda$ on $\X$ and consider a family of rank-$n$ Higgs bundles $(\cal{E}, \Phi) = (\cal{E}, \Phi : \cal{E} \to \cal{E} \otimes \Lambda)$.
We call this family \dfn{regular} if every restriction $(\cal{E}_s, \Phi_s)$ is regular.
Equivalently, after locally trivialising $\Lambda$, the endomorphism representing $\Phi$ is regular on the vector space fibre at every point $x \in \X$.

The characteristic polynomial of $\Phi$ defines a finite flat degree-$n$ spectral cover $\pi : \Y \to \X$ inside the total space of $\Lambda$.
The Cayley-Hamilton Theorem gives an action $\rho : \cal{A} \to \cal{End} (\cal{E})$ of the spectral algebra $\cal{A} \coleq \pi_\ast \cal{O}_\Y$.
The assertion of \autoref{1790763564459} remains unchanged: this action defines a spectral structure on $\cal{E}$ if and only if the family is regular.
Indeed, suppose that the family is regular.
Take a local trivialisation of $\Lambda$ and let $T$ be the endomorphism representing $\Phi$.
Given a point $x \in \X$, lift a cyclic vector for $T (x)$ to a local section $e$ of $\cal{E}$ near $x$.
The determinant of the sections $e, Te, \ldots, T^{n-1} e$ is nonzero at $x$, and hence is a unit near $x$.
These sections therefore form a local frame, and the map $a \mapsto \rho (a) e$ gives a local $\cal{A}$-module isomorphism $\cal{A} \iso \cal{E}$.
This argument works also when $\S$ is nonreduced.
Conversely, a local $\cal{A}$-module generator determines a cyclic vector in every vector space fibre, so the family is regular.

For a regular family, these constructions are compatible with arbitrary base change in $\S$.
In the notation of \autoref{1791551111381}, the pulled-back Higgs bundle is
\begin{eqn}
	(\cal{E}', \Phi') \coleq (f^\ast \cal{E}, f^\ast \Phi)
\qqtext{with twist}
	\Lambda' \coleq f^\ast \Lambda
\fullstop
\end{eqn}
The characteristic polynomial and its Cayley-Hamilton action commute with pullback, so the spectral structure induced by $\Phi'$ is precisely $(\cal{A}', \rho') = (f^\ast \cal{A}, f^\ast \rho)$.
In particular, regularity is preserved by base change, and \autoref{1791551111381} describes the resulting spectral cover and spectral line bundle.
Likewise, the off-diagonal form $\eta$ of \autoref{1791458921327} is defined over the whole family and commutes with arbitrary base change, so the isomorphism $\cal{Q} \iso \cal{Q}^\vee \otimes \Lambda$ holds over $\X$ itself.

The equivalence of \autoref{1792167464001} also extends to families and commutes with arbitrary base change.
A global $\Lambda$-twisted generator $a$ of $\cal{A}$ determines the Higgs field
\begin{eqn}
	\Phi = (\rho \otimes \id_\Lambda) (a)
\fullstop
\end{eqn}
An isomorphism $(\phi,\alpha)$ of spectral structures intertwines the corresponding Higgs fields precisely when
\begin{eqn}
	(\alpha \otimes \id_\Lambda) (a) = a'
\fullstop
\end{eqn}
These formulas are defined over $\X$ itself and are preserved by pullback.
Thus, after fixing an identification with a given family of spectral structures $(\cal{E}, \cal{A}, \rho)$, the corresponding groupoid of regular Higgs bundles is equivalent to the discrete groupoid of global $\Lambda$-twisted generators of $\cal{A}$.

\subsection{Local Existence in Families}
\label{1792146750003}

The existence results of \autoref{1790777298686} and \autoref{1792167000673} also extend to families, locally over the parameter space.
The construction uses \autoref{1792430412943} to choose a suitable regular Higgs field on one fibre, and then extends it to all nearby fibres.
After shrinking the neighbourhood of this one fibre in the parameter space, the resulting spectral structures remain clean.

\begin{theorem}[local existence of clean spectral structures]
\label{1791544458956}
Let $\cal{E}$ be any vector bundle of rank $n$ on $\X$.
Every point $s_0 \in \S$ has a neighbourhood $\U \subset \S$ over which $\cal{E}_\U$ admits a fibrewise clean spectral structure.
Moreover, this spectral structure can be induced by a regular $\Lambda$-twisted Higgs field on $\cal{E}_\U$ for a suitable line bundle $\Lambda$ on $\X_\U$.

More explicitly, there is a finite flat degree-$n$ cover $\pi_\U : \Y_\U \to \X_\U$, a line bundle $\cal{L}_\U$ on $\Y_\U$, and an isomorphism $\iota_\U : (\pi_\U)_\ast \cal{L}_\U \iso \cal{E}_\U$.
For every $s \in \U$, the curve $\Y_s$ is smooth and connected, the cover $\pi_s : \Y_s \to \X_s$ is simply branched, and the off-diagonal bundle satisfies $\H^1 (\X_s, \cal{Q}_s) = 0$.
The resulting spectral structure and its spectral realisation are compatible with arbitrary base change in $\U$.
\end{theorem}

\begin{proof}
When $n = 1$, we can take $\U = \S$,~ $\Y = \X$,~ $\pi = \id_{\X}$, and $\cal{L} = \cal{E}$.
The spectral algebra is $\cal{O}_{\X}$ and the off-diagonal bundle is zero, so every fibre is clean.
This spectral structure is induced by the regular Higgs field $\Phi = \id_{\cal{E}}$ with twist $\Lambda = \cal{O}_{\X}$.

Suppose now that $n \geq 2$.
Fix $s_0 \in \S$ and choose a point $x_0 \in \X_{s_0}$.
Since $p$ is smooth, its local product description near $x_0$ gives a holomorphic section $\sigma : \U \to \X_\U$ through $x_0$, after restricting to a neighbourhood $\U$ of $s_0$.
Its image $\sfSigma \coleq \sigma (\U)$ is a closed effective Cartier divisor on $\X_\U$.
Indeed, the section is closed because $p$ is separated, and near its image we may choose a relative disc coordinate $z$ in which $\sigma$ is given by $z = 0$.
Its ideal is then generated by $z$, which is a non-zero-divisor even when $\U$ is nonreduced.

Let $\Lambda_0 \coleq \cal{O}_{\X_\U} (\sfSigma)$.
For any integer $N$, take 
\begin{eqn}
	\Lambda \coleq \Lambda_0^N = \cal{O}_{\X_\U} (N \sfSigma)
\qtext{so that}
	\Lambda_s \cong \cal{O}_{\X_s} \big(N \sigma(s) \big)
\fullstop
\end{eqn}
Then $\deg \Lambda_s = N$ for every $s \in \U$.
Let $N$ be sufficiently large that $\cal{End} (\cal{E}_{s_0}) \otimes \Lambda_{s_0}$ is globally generated and
\begin{eqntag}
\label{1791545577180}
	\H^1 \big(\X_{s_0}, \cal{End} (\cal{E}_{s_0}) \otimes \Lambda_{s_0} \big) = 0
\fullstop
\end{eqntag}
By \autoref{1792430412943}, after increasing $N$ if necessary, there is also a regular Higgs field
\begin{eqn}
	\Phi_{s_0}
	\in
	\H^0 \big(\X_{s_0}, \cal{End} (\cal{E}_{s_0}) \otimes \Lambda_{s_0} \big)
\end{eqn}
whose induced spectral structure is clean.
In particular, its spectral cover is smooth, connected, and simply branched, and its off-diagonal bundle satisfies
\begin{eqn}
	\H^1 (\X_{s_0}, \cal{Q}_{s_0}) = 0
\fullstop
\end{eqn}

The point is that the vanishing of cohomology \eqref{1791545577180} allows us to extend $\Phi_{s_0}$ to a family of Higgs fields over a neighbourhood of $s_0$.
That is, after possibly shrinking the neighbourhood $\U \subset \S$ of $s_0$, there is some Higgs field
\begin{eqn}
	\Phi_\U \in \H^0 (\X_\U, \cal{End} (\cal{E}_\U) \otimes \Lambda)
\end{eqn}
whose restriction to the central fibre $\X_{s_0}$ is the Higgs field $\Phi_{s_0}$.
By upper semicontinuity, after shrinking $\U$ around $s_0$, we may assume that
\begin{eqn}
	\H^1 (\X_s, \cal{End} (\cal{E}_s) \otimes \Lambda_s) = 0
\qquad
	\text{for all $s \in \U$}
\fullstop
\end{eqn}
Then the Analytic Cohomology and Base Change Theorem implies that 
\begin{eqn}
	\cal{H} \coleq p_\ast \big( \cal{End}(\cal{E}_\U) \otimes \Lambda \big)
\end{eqn}
is locally free on $\U$ and that the canonical maps
\begin{eqn}
	\cal{H} (s) \iso \H^0 \big( \X_s, \cal{End} (\cal{E}_s) \otimes \Lambda_s \big)
\end{eqn}
are isomorphisms for every $s \in \U$; see \cite[Corollary~6.6.11 and Theorem~6.7.4]{GuiCAG} and \cite[Ch.~III, Corollary~3.9 and Theorem~4.12(i)]{BanicaStanasila1976}.
After shrinking $\U$ further, we may assume that $\cal{H}$ is free.
The vector $\Phi_{s_0} \in \cal{H} (s_0)$ then extends to a section $\Phi_\U$ of $\cal{H}$ over $\U$: in a local frame, simply extend its coordinates as constant functions.
This extension is not unique because it depends on the trivialisation of $\cal{H}$, but that is immaterial for us.
This proves the claim.

We claim that the discriminant of $\Phi_s$ has only simple zeros for every $s$ sufficiently close to $s_0$.
The discriminant of $\Phi_\U$ is a section $\Delta_\U$ of $\Lambda^{n(n-1)}$.
In a relative local coordinate $z$ and a local trivialisation of this line bundle, write $g$ for the function representing $\Delta_\U$ and consider the locus $g = 0$ and $\del g / \del z = 0$.

These local loci define a closed analytic subset $\B \subset \X_\U$, independently of the choices.
Indeed, changing the relative coordinate multiplies the derivative by a unit, whilst replacing $g$ by $u g$ for a unit $u$ replaces its derivative by $u \del g / \del z + g \del u / \del z$.
Thus, the simultaneous vanishing condition is unchanged.
By construction, $\B$ is disjoint from $\X_{s_0}$ because $\Delta_{s_0}$ has only simple zeros.
But $p_\U$ is proper, so $p_\U (\B)$ is closed in $\U$.
Shrinking $\U$ to avoid this image proves the claim.

\enlargethispage{20pt}
The characteristic polynomial of $\Phi_\U$ now defines a finite flat degree-$n$ cover $\pi_\U : \Y_\U \to \X_\U$.
For every $s \in \U$, the simplicity of the zeros of $\Delta_s$ says that $\Phi_s$ avoids the singular matrix discriminant and meets its smooth locus transversely.
As in the discriminant argument in the proof of \autoref{1791390497007}, this implies that $\Y_s$ is smooth and $\pi_s$ is simply branched.
Concretely, near a branch point, the ramified part has local equation $w^2 = z$, with all remaining sheets unramified.

As an $\cal{O}_{\X_\U}$-module, the spectral algebra is
\begin{eqn}
	\cal{A}_\U
	\coleq
	(\pi_\U)_\ast \cal{O}_{\Y_\U}
	\cong
	\cal{O}_{\X_\U} \oplus \Lambda^{-1} \oplus \cdots \oplus \Lambda^{-(n-1)}
\fullstop
\end{eqn}
Since every $\Lambda_s$ has positive degree, all of its negative powers have no nonzero global sections.
Consequently, finite base change gives
\begin{eqn}
	\H^0 (\Y_s, \cal{O}_{\Y_s})
	\cong
	\H^0 (\X_s, \cal{A}_s)
	\cong
	\H^0 (\X_s, \cal{O}_{\X_s})
	\cong
	\CC
\fullstop
\end{eqn}
The smooth compact curve $\Y_s$ is therefore connected.

The Cayley-Hamilton Theorem gives an action $\rho_\U : \cal{A}_\U \to \cal{End} (\cal{E}_\U)$.
Since $\Y_\U$ is fibrewise smooth and connected over $\U$, \autoref{1792143470690} shows that this action makes
$\cal{E}_\U$ locally free of rank one over $\cal{A}_\U$.
The finite sheaf-module correspondence then gives a line bundle $\cal{L}_\U$ on $\Y_\U$ and an isomorphism $\iota_\U : (\pi_\U)_\ast \cal{L}_\U \iso \cal{E}_\U$.

It remains to arrange the vanishing of off-diagonal cohomology on every fibre over $\U$.
The off-diagonal bundle $\cal{Q}_\U = \cal{End} (\cal{E}_\U) / \rho_\U (\cal{A}_\U)$ is a vector bundle on $\X_\U$, and by \autoref{1791551111381} its restriction to each fibre $\X_s$ is the off-diagonal bundle $\cal{Q}_s = \cal{End} (\cal{E}_s) / \rho_s (\cal{A}_s)$.
In particular, its restriction to the central fibre is the off-diagonal bundle of the clean structure chosen above.
Since $\cal{Q}_\U$ is flat over $\U$, upper semicontinuity and the vanishing of $\H^1 (\X_{s_0}, \cal{Q}_{s_0})$ allow us to shrink $\U$ one more time to ensure that $\H^1 (\X_s, \cal{Q}_s) = 0$ for every $s \in \U$.
Thus the induced spectral structure is fibrewise clean.

Finally, compatibility of the spectral structure, its spectral realisation, and the off-diagonal bundle with arbitrary base change follows from \autoref{1791551111381}.
After any such base change, each fibre is identified with a fibre of the family just constructed, so it remains clean.
\end{proof}

\section{The Moduli Stack of Spectral Structures}
\label{1792401370132}

We now collect the family constructions of the previous section into a moduli stack of spectral structures.
This allows us to regard the groupoids of spectral structures and spectral realisations as the fibres of a single geometric object, and to describe how this object relates in a geometric sense to the moduli stacks of covers and vector bundles.
We work on the site $\mathbf{An}$ of complex analytic spaces, possibly nonreduced, locally of finite dimension with open coverings.
An \dfn{analytic stack} is a stack in groupoids whose diagonal is representable by analytic spaces and which admits a smooth surjective atlas by an analytic space (some authors would call this an \dfn{analytic Artin stack}) \cite[\S\S2.3 and~3.2]{PortaYu2016}.

Our first theorem (\autoref{1792494933366}) establishes analyticity of the moduli stack and representability of its forgetful morphism to vector bundles.
Our second theorem (\autoref{1792494388516}) describes distinguished smooth open substacks and shows that the forgetful morphism is smooth and surjective on the clean substack.

Throughout, let $\S$ be any complex analytic space, possibly nonreduced, and consider a proper smooth family $p : \X \to \S$ of connected compact curves.
Fix an integer $n \geq 1$.
For any holomorphic map $\T \to \S$, put
\begin{eqn}
	\X_\T \coleq \X \times_\S \T
\fullstop
\end{eqn}

\subsection{The Relative Analytic Stack of Vector Bundles}
\label{1793143183111}

Before constructing the moduli stack of spectral structures, we construct the relative analytic stack of vector bundles, allowing the base curve to vary and the parameter space to be nonreduced.

\begin{proposition}[relative analytic stack of vector bundles]
\label{1793142766507}
The moduli stack $\frak{Bun}_n (\X / \S)$ of rank-$n$ holomorphic vector bundles over the family $\X \to \S$ is an analytic stack, locally of finite presentation and smooth over $\S$.
Its formation commutes with arbitrary analytic base change: for every holomorphic map $\S' \to \S$, there is a canonical equivalence
\begin{eqn}
	\frak{Bun}_n (\X / \S) \times_\S \S'
	\cong
	\frak{Bun}_n (\X \times_\S \S' / \S')
\fullstop
\end{eqn}
\end{proposition}

To prove this proposition, we begin with a general construction of relative spaces of sections which will be used both here and in the subsequent moduli constructions.
It says that over the given proper smooth family of compact curves, the relative functor of global sections of any separated holomorphic map is representable by an analytic space, the relative space of sections.

\begin{lemma}[relative space of sections]
\label{1792526939305}
Let $\T \to \S$ be an analytic space, and let $u : \Z \to \X_\T$ be a separated holomorphic map.
For every analytic space $\T' \to \T$, put $\Z_{\T'} \coleq \Z \times_\T \T'$ and $u_{\T'} : \Z_{\T'} \to \X_{\T'}$.
Then there exists an analytic space $\M$, locally of finite presentation over $\T$, together with a universal section $\sigma_\M : \X_\M \too \Z_\M$ of the base-change projection $u_\M : \Z_\M \to \X_\M$ such that the operation of pullback of $\sigma_\M$ gives a natural bijection, for every analytic space $\T' \to \T$,
\begin{eqn}
	\Hom_\T (\T', \M)
	\cong
	\set{
		\sigma : \X_{\T'} \to \Z_{\T'}
		~\Big|~
		u_{\T'} \circ \sigma = \id_{\X_{\T'}}
	}
\fullstop
\end{eqn}
Moreover, this construction commutes with arbitrary base change: for every $\T' \to \T$, the analytic space $\M'$ determined by $u_{\T'}$ comes with a canonical isomorphism 
\begin{eqn}
	\M' \cong \M \times_\T \T'
\end{eqn}
which intertwines the universal sections $\sigma_{\M}$ and $\sigma_{\M'}$ in the sense that there is a commutative diagram
\begin{eqn}
\begin{tikzcd}
	\Z_\M
		\ar[d]
&	\Z_\M \times_\T \T'
		\ar[r, "\cong"]
		\ar[d]
		\ar[l]
&	\Z_{\M'}
		\ar[d]
\\
	\X_\M
		\ar[u, bend left, "\sigma_\M"]
&	\X_\M \times_\T \T'
		\ar[r, "\cong"']
		\ar[u, bend left, "\sigma'_\M"]
		\ar[l]
&	\X_{\M'}
		\ar[u, bend right, "\sigma_{\M'}"']
\fullstop
\end{tikzcd}
\end{eqn}
\end{lemma}

To be clear, the section $\sigma'_\M$ in the above commutative diagram is the pullback of $\sigma_\M$: i.e., $\sigma'_\M = \sigma_\M \times_\T \id_{\T'}$.
That is to say, if $\sigma_\M$ takes a point $x \in \X_\M$ to $z \in \Z_\M$, then $\sigma'_\M$ takes any point $(x,t')$ to $(z,t')$.

\begin{proof}[Proof of \autoref{1792526939305}]
We construct $\M$ by regarding a section of $u : \Z \to \X_\T$ as a closed analytic subspace of $\Z$ whose projection to $\X_\T$ via $u$ is an isomorphism.

By the Relative Douady Theorem \cite[Théorème~2]{Pourcin1969}, there is an analytic space $\R \to \T$ and a universal closed analytic subspace $	\W \subset \Z \times_\T \R$ which is proper and flat over $\R$.
Its universal property is the following.
For every analytic space $\T' \to \T$, morphisms $\T' \to \R$ over $\T$ correspond naturally to closed analytic subspaces $\W' \subset \Z \times_\T \T'$ which are proper and flat over $\T'$.
Under this correspondence, $\W'$ is the pullback $\W \times_\R \T'$, regarded as a closed analytic subspace of
$\Z \times_\T \T'$.

The morphism $u : \Z \to \X_\T$ induces the morphism $u_\R : \Z \times_\T \R \to \X_\T \times_\T \R$ which then induces an $\R$-morphism $h : \W \to \X_\T \times_\T \R$ by restriction to the subspace $\W$.
Both its source and its target are proper and flat over $\R$.
For the source, this is part of the universal property of the relative Douady space.
For the target, it follows by base change from the properness and flatness of $\X_\T \to \T$.
If $r \in \R$ lies over $t \in \T$, the induced morphism on fibres is the restriction of $u$; i.e., $h_r : \W_r \too \X_t$.

By the openness of the isomorphism locus \cite[\S10.1, Proposition~1]{Douady1966}, the subset
\begin{eqn}
	\M
	\coleq
	\set[\Big]{r \in \R ~\big|~ \text{$h_r$ is an isomorphism} }
	\subset \R
\end{eqn}
is an open analytic subspace.
Consequently, the restriction of $h$ to $\M$ is an isomorphism of analytic spaces
\begin{eqn}
	h_{\M} : \W \times_\R \M \iso \X_\T \times_\T \M
\fullstop
\end{eqn}

Let $j_\M : \W \times_\R \M \inj \Z \times_\T \M$ be the restriction to $\M$ of the inclusion map $\W \inj \Z \times_\T \R$.
Then the composite map
\begin{eqn}
	\sigma_\M
	\coleq
	j_\M
	\circ
	h_\M^{-1}
	:
	\X_\T \times_\T \M
	\too
	\Z \times_\T \M
\end{eqn}
is a section of $u_\M : \Z \times_\T \M \to \X_\T \times_\T \M$ because $j_\M$ is an inclusion and $h_\M$ is the restriction of $u_\M$.
Thus, under the identification given by $h_\M$, the universal closed subspace $\W \times_\R \M \subset \Z \times_\T \M$ is precisely the image of this section $\sigma_\M$.

We verify that $\sigma_\M$ has the required universal property.
Let $\T' \to \T$ be any analytic space, and let $\sigma : \X_{\T'} \to \Z_{\T'} \coleq \Z \times_\T \T'$ be a section of $u' : \Z_{\T'} \to \X_{\T'}$.
Since $u : \Z \to \X_\T$ is separated, its base change is separated, so $\sigma$ is a closed immersion.
Its image $\sfSigma \subset \Z_{\T'}$ is therefore a closed analytic subspace.
Moreover, the projection $u'$ identifies $\sfSigma$ with $\X_{\T'}$, so $\sfSigma$ is proper and flat over $\T'$.

The universal property of $\R$ now gives a unique morphism $f_\sigma : \T' \to \R$ over $\T$ whose pullback of $\W$ is $\sfSigma$.
For every $t' \in \T'$, the projection $\sfSigma_{t'} \to \X_{t'}$ is an isomorphism.
Consequently, $f_\sigma (t') \in \M$ for every $t' \in \T'$.
But since $\M$ is open in $\R$, the morphism $f_\sigma$ factors uniquely through $\M$.

Conversely, any $\T$-morphism $f : \T' \to \M$ gives a section of $\Z_{\T'} \to \X_{\T'}$ by pulling back $\sigma_\M$, and these two constructions are inverse to each other.
Indeed, the section recovered from $\sfSigma$ is $\sigma$, whilst the closed subspace associated with the pullback of $\sigma_\M$ is the pullback of $\W$; uniqueness in the relative Douady universal property therefore recovers $f$.
The constructions also commute with arbitrary base change in $\T'$.
Therefore, the open analytic space $\M$ represents the desired functor of sections with universal section $\sigma_\M$.

The assertion about base change also follows from this universal property.
Indeed, for any analytic space $\U \to \T'$, morphisms $\U \to \M \times_\T \T'$ over $\T'$ correspond to sections of $u_\U : \Z_\U \to \X_\U$.
Thus, $\M \times_\T \T'$ represents the same functor as the analytic space $\M'$ corresponding to $u_{\T'}$.
Uniqueness of the representing space gives the canonical isomorphism and identifies the universal sections.

Finally, the holomorphic map $\M \to \T$ locally factors as a closed analytic embedding into a finite-dimensional relative polydisc over $\T$, followed by the projection to $\T$.
The defining ideal of this embedding is coherent and hence locally finitely generated.
Therefore, $\M \to \T$ is locally of finite presentation.
This completes the proof of \autoref{1792526939305}.
\end{proof}

\begin{proof}[Proof of \autoref{1793142766507}]
For every analytic space $\T \to \S$, the groupoid $\frak{Bun}_n (\X / \S) (\T)$ consists of rank-$n$ vector bundles on $\X_\T$ and their isomorphisms.
Vector bundles and their isomorphisms glue over open coverings, so these groupoids form a stack.
We prove representability of its diagonal and then construct a smooth surjective atlas.

\caps{Step 1: Representability of the diagonal.}
Let $\cal{E}$ and $\cal{F}$ be rank-$n$ vector bundles on $\X_\T$.
Consider the open analytic subspace
\begin{eqn}
	\Z \subset \tot \big( \cal{Hom} (\cal{E}, \cal{F}) \big)
\end{eqn}
consisting of invertible linear maps between their fibres.
Invertibility is the nonvanishing of the determinant, so this is indeed an open subspace.
Its projection $u : \Z \to \X_\T$ is separated because it is the restriction of a vector-bundle projection.

For every analytic space $\T' \to \T$, write $\cal{E}_{\T'}$ and $\cal{F}_{\T'}$ for the pulled-back bundles on $\X_{\T'}$.
Sections of $u_{\T'} : \Z_{\T'} \to \X_{\T'}$ are precisely isomorphisms $\cal{E}_{\T'} \iso \cal{F}_{\T'}$.
By \autoref{1792526939305}, this functor is represented by an analytic space $\M \to \T$, locally of finite presentation.
Therefore, the relative isomorphism functor between $\cal{E}$ and $\cal{F}$ is representable by an analytic space.
This proves representability of the diagonal.

\caps{Step 2: Spaces of generating sections.}
The construction of an atlas is local on $\S$.
After restricting to a neighbourhood of any chosen point of $\S$, smoothness of $p : \X \to \S$ gives a holomorphic section $\sigma : \S \to \X$.
Its image $\sfSigma \coleq \sigma (\S)$ is an effective Cartier divisor: locally in a relative coordinate $z$, it is given by $z = 0$, and $z$ is a non-zero-divisor even over a nonreduced base.
Put $\Lambda \coleq \cal{O}_\X (\sfSigma)$, so that $\Lambda$ has degree one on every fibre.

For integers $m \geq 1$ and $N \geq n$, put $\cal{P} = \cal{P}_{m,N} \coleq (\Lambda^{-m})^{\oplus N}$, and consider the relative Grassmannian of rank-$(N-n)$ vector subbundles of $\cal{P}$:
\begin{eqn}
	\G = \G_{m,N}
	\coleq
	\sfop{Gr}_\X \big( N-n, \cal{P} \big)
\fullstop
\end{eqn}
Suppressing the pullbacks to $\G$ from the notation, the Grassmannian $\G$ carries the universal exact sequence
\begin{eqn}
	0
	\too
	\cal{K}
	\too
	\cal{P}
	\too
	\cal{E}
	\too 0
\fullstop{,}
\end{eqn}
where $\cal{E}$ is a vector bundle of rank $n$.

The projection $\G \to \X$ is separated.
Therefore, \autoref{1792526939305}, applied with $\T = \S$, gives an analytic space $\M \to \S$ parametrising its sections.
Pulling back the universal exact sequence by the universal section gives a universal rank-$n$ vector-bundle quotient on $\X_{\M}$:
\begin{eqn}
	(\Lambda_{\M}^{-m})^{\oplus N}
	\too
	\cal{E}_{\M}
\fullstop
\end{eqn}
Conversely, every rank-$n$ vector-bundle quotient of $(\Lambda^{-m})^{\oplus N}$ defines a section of $\G \to \X$ by taking its kernel.
These correspondences commute with arbitrary base change.

Let $\V = \V_{m,N} \subset \M$ be the open locus where 
\begin{eqn}
	\H^1 (\X_s, \cal{E}_v \otimes \Lambda_s^m) = 0
\end{eqn}
for a parameter $v$ lying over $s \in \S$.
This condition is open by upper semicontinuity \cite[Ch.~III, Theorem~4.12(i)]{BanicaStanasila1976}.
The universal quotient gives a morphism
\begin{eqn}
	a : \V \too \frak{Bun}_n (\X / \S)
\fullstop
\end{eqn}

We claim that $a$ is representable and smooth.
To see this, let $\T \to \S$ be any analytic space and let $\cal{E}$ be a vector bundle on $\X_\T$.
Let $\U \subset \T$ be the open subset where
\begin{eqn}
	\H^1 (\X_t, \cal{E}_t \otimes \Lambda_t^m) = 0
\qquad
	\forall t \in \U
\fullstop
\end{eqn}
By analytic cohomology and base change, the sheaf
\begin{eqn}
	\cal{H}
	\coleq
	(p_{\U})_\ast
	\big( \cal{E}|_{\X_{\U}} \otimes \Lambda_{\U}^m \big)
\end{eqn}
is locally free on $\U$, and its formation commutes with arbitrary analytic base change \cite[Ch.~III, Theorem~3.4 and Corollaries~3.5 and~3.9]{BanicaStanasila1976}.

An object of the fibre product $\T \times_{\frak{Bun}_n (\X / \S)} \V$ is a quotient together with an identification of its target with the prescribed bundle $\cal{E}$.
It is therefore exactly an $N$-tuple of sections of $\cal{E}|_{\X_\U} \otimes \Lambda_\U^m$ which generates that bundle.
The specified identification removes all automorphisms.
Consequently, this fibre product is the open subspace of $\tot (\cal{H}^{\oplus N}) \to \U$ consisting of the generating $N$-tuples.

\clearpage

To check openness, consider the universal evaluation map over this total space.
Its cokernel is coherent, and the image of its support under the proper curve projection is closed.
The complement of this image is precisely the generating locus.
Thus, the fibre product is an analytic space, smooth over the open subspace $\U$, and hence smooth over $\T$.
This proves the claim.

\caps{Step 3: Surjectivity of the atlas.}
Fix a vector bundle $\cal{E}$ on a fibre $\X_s$.
For sufficiently large $m$, the positive twist $\cal{E} \otimes \Lambda_s^m$ is globally generated and has vanishing first cohomology.
Choosing finitely many generating sections gives a quotient $(\Lambda_s^{-m})^{\oplus N} \to \cal{E}$, and hence a point of $\V$ over $\cal{E}$.

The same construction works locally for any family $\cal{E}$ on $\X_\T$.
Indeed, near a chosen $t \in \T$, cohomology and base change allow the selected generators on $\X_t$ to extend to sections over a neighbourhood of $t$.
Their evaluation remains surjective after shrinking, by properness of the curve projection.
Thus, every such family lifts locally to some $\V$.

Taking the disjoint union of these spaces $\V_{m,N}$ over all $m,N$ and over an open cover of the original base $\S$ gives a smooth surjective atlas.
Together with Step 1, this proves that $\frak{Bun}_n (\X / \S)$ is an analytic stack.
The chart spaces are locally of finite presentation over $\S$ because they are complex analytic spaces over $\S$.
Therefore, the stack is locally of finite presentation over $\S$.

\caps{Step 4: Smoothness over $\S$.}
We prove that each chart $\V = \V_{m,N} \to \S$ is smooth by the analytic infinitesimal lifting criterion \cite[Théorème~3.1]{1791539610207}.
It suffices to solve lifting problems for small extensions $\R' \to \R$ of local Artinian $\CC$-algebras, with kernel $\I$ annihilated by the maximal ideal of $\R'$.
Indeed, every Artinian surjection factors into such extensions.

Put $\T \coleq \Spec (\R)$ and $\T' \coleq \Spec (\R')$.
Suppose we are given a map $\T' \to \S$ and a point of $\V$ over $\T$, represented by a quotient
\begin{eqn}
	(\Lambda_\T^{-m})^{\oplus N} \too \cal{E}_\T
\fullstop
\end{eqn}
Write $\X_0$, $\cal{E}_0$, and $\Lambda_0$ for the restrictions to the closed point.
Flatness of $\X_{\T'}$ over $\T'$ identifies the ideal of $\X_\T$ in $\X_{\T'}$ with $\cal{O}_{\X_0} \otimes_\CC \I$.
The obstruction to lifting $\cal{E}_\T$ to a vector bundle on $\X_{\T'}$, with its prescribed identification over $\T$, lies in
\begin{eqn}
	\H^2 \big( \X_0, \cal{End} (\cal{E}_0) \big) \otimes_\CC \I
	= 0
\fullstop
\end{eqn}
Here the vanishing holds because $\X_0$ is a compact Riemann surface.
Concretely, on a Stein open cover, local lifts of the transition matrices fail to satisfy the cocycle identity by a \v{C}ech $2$-cocycle with values in $\cal{End} (\cal{E}_0) \otimes_\CC \I$.
Its vanishing cohomology class allows us to adjust the lifts so that they glue to a vector bundle $\cal{E}_{\T'}$; see also \cite[\href{https://stacks.math.columbia.edu/tag/08VR}{Tag 08VR}]{stacks-project}.

The original quotient is equivalent to $N$ generating sections of $\cal{E}_\T \otimes \Lambda_\T^m$.
The exact sequence
\begin{eqn}
	0
	\too
	\cal{E}_0 \otimes \Lambda_0^m \otimes_\CC \I
	\too
	\cal{E}_{\T'} \otimes \Lambda_{\T'}^m
	\too
	\cal{E}_\T \otimes \Lambda_\T^m
	\too 0
\end{eqn}
and the defining vanishing $\H^1 (\X_0, \cal{E}_0 \otimes \Lambda_0^m) = 0$ show that all these sections lift.
Their lifted evaluation remains surjective by Nakayama's Lemma.
It therefore defines a quotient $(\Lambda_{\T'}^{-m})^{\oplus N} \to \cal{E}_{\T'}$ restricting to the original one.
The closed fibre is unchanged, so this quotient still lies in $\V_{m,N}$.
This solves the lifting problem and proves that $\V \to \S$ is smooth.
Smoothness of $\frak{Bun}_n (\X / \S)$ over $\S$ follows from its smooth surjective atlas.

Finally, the base-change equivalence follows directly from the definition: both sides classify vector bundles on $\X \times_\S \T$ for analytic spaces $\T \to \S'$.
\end{proof}

\subsection{The Analytic Stack of Spectral Structures}
\label{1793143227123}

We now use the relative analytic stack of vector bundles from \autoref{1793142766507} and the construction of the relative space of sections of \autoref{1792526939305} to construct the moduli stack of spectral structures.
The essential step is to represent spectral structures on a prescribed family of vector bundles by an analytic space.

\begin{definition}
\label{1792401376079}
The \dfn{moduli stack of spectral structures} of rank $n$ over the family $\X \to \S$ is the stack $\frak{SpStr}_n (\X / \S)$ over the analytic site $\mathbf{An}$ whose groupoid over an analytic space $\T \to \S$ is
\begin{eqn}
	\frak{SpStr}_n (\X / \S) (\T)
	\coleq
	\mathbf{SpStr}_n (\X_\T)
\fullstop
\end{eqn}
\end{definition}

Thus, an object is a triple $(\cal{E}, \cal{A}, \rho)$ consisting of a vector bundle $\cal{E}$ of rank $n$ on $\X_\T$ and a spectral structure $(\cal{A}, \rho)$ on $\cal{E}$ in the sense of \smartref{1792146709610}.
In particular, $\cal{A}$ is a commutative unital $\cal{O}_{\X_\T}$-algebra, locally free of rank $n$, and $\rho : \cal{A} \to \cal{End} (\cal{E})$ makes $\cal{E}$ locally free of rank $1$ over $\cal{A}$.
Isomorphisms are the same as in \autoref{1789150206981}: they are pairs $(\phi, \alpha)$ consisting of a vector bundle isomorphism $\phi : \cal{E} \iso \cal{E}'$ and an $\cal{O}_{\X_\T}$-algebra isomorphism $\alpha : \cal{A} \iso \cal{A}'$ such that $\Ad_\phi \circ \rho = \rho' \circ \alpha$.
The pullback functors are given by \autoref{1791551111381}.
The fact that this is indeed a stack over $\mathbf{An}$ follows immediately from effective descent for vector bundles and their morphisms because the algebra identities and the rank-one module condition are local.

For a fixed compact Riemann surface $\X$, we take $\S = \ast = \Spec (\CC)$ and write simply $\frak{SpStr}_n (\X)$.
Its groupoid over a point is precisely the groupoid $\mathbf{SpStr}_n (\X)$ introduced in \autoref{1789150206981}.
More generally, base change to any point $s \in \S$ gives the moduli stack $\frak{SpStr}_n (\X_s)$.

\begin{theorem}[analyticity and representability]
\label{1792494933366}
The moduli stack of spectral structures $\frak{SpStr}_n (\X / \S)$ is an analytic stack locally of finite presentation over $\S$, and the forgetful morphism
\begin{eqn}
	F : \frak{SpStr}_n (\X / \S) \too \frak{Bun}_n (\X / \S)
\end{eqn}
is representable by analytic spaces and locally of finite presentation.
Arbitrary analytic base change commutes with the formation of this moduli stack and the forgetful morphism $F$: for every holomorphic map $\S' \to \S$, there is a canonical equivalence
\begin{eqn}
	\frak{SpStr}_n (\X / \S) \times_\S \S'
	\cong
	\frak{SpStr}_n (\X \times_\S \S' / \S')
\fullstop
\end{eqn}
\end{theorem}

\begin{proof}[Proof of \autoref{1792494933366}]
Let $\T \to \S$ be an analytic space, and let $\cal{E}$ be a vector bundle of rank $n$ on $\X_\T$ representing a morphism $\T \to \frak{Bun}_n (\X / \S)$.
We construct an analytic space $\M = \M (\cal{E})$ together with a natural equivalence of stacks
\begin{eqntag}
\label{1793028752430}
	\T \times_{\frak{Bun}_n (\X / \S)} \frak{SpStr}_n (\X / \S)
	\cong
	\M
\fullstop
\end{eqntag}
An object over a test space $\T' \to \T$ of the fibre product on the left consists of a spectral structure $(\cal{E}', \cal{A}', \rho')$ together with an isomorphism $\phi : \cal{E}' \iso \cal{E}_{\T'}$, where $\cal{E}_{\T'}$ is the pullback of $\cal{E}$ to $\X_{\T'}$.
Conjugating the action by $\phi$ gives the spectral structure on $\cal{E}_{\T'}$ via
\begin{eqn}
	\cal{A}' \too \cal{End} (\cal{E}_{\T'})
	\qtext{sending}
	a \mapstoo \phi \circ \rho'(a) \circ \phi^{-1}
\fullstop
\end{eqn}
Compatibility with the specified identification $\phi$ means that isomorphisms between the resulting spectral structures must induce the identity on the vector bundle $\cal{E}_{\T'}$.

\enlargethispage{30pt}
Thus, a spectral structure $(\cal{A}, \rho)$ on the given vector bundle $\cal{E}$ is determined, up to the unique isomorphism inducing the identity on $\cal{E}$, by its image
\begin{eqn}
	\cal{B} \coleq \rho (\cal{A}) \subset \cal{End} (\cal{E})
\fullstop
\end{eqn}
This is a vector subbundle of rank $n$.
Indeed, if $e$ is a local $\cal{A}$-module generator of $\cal{E}$, then the composite map
\begin{eqn}
	\cal{A}
	\overset{\rho}{\too}
	\cal{End} (\cal{E})
	\too
	\cal{E}
\qqtext{given locally by}
	a \mapstoo \rho (a) e
\end{eqn}
\clearpage
is an isomorphism.
Evaluation at $e$ therefore gives a local splitting of $\rho$.
Conversely, a commutative unital subalgebra $\cal{B} \subset \cal{End} (\cal{E})$ which is a vector subbundle of rank $n$ defines a spectral structure precisely when $\cal{E}$ is locally free of rank one over $\cal{B}$.

It also follows that the fibre product on the lefthand side of \eqref{1793028752430} has no automorphisms: indeed, an automorphism inducing the identity on $\cal{E}$ must also induce the identity on the subalgebra $\cal{B} \subset \cal{End} (\cal{E})$.

Consider the relative Grassmannian $\G \coleq \sfop{Gr}_{\X_\T} \big( n, \cal{End} (\cal{E}) \big)$ of rank-$n$ vector subbundles of $\cal{End} (\cal{E})$.
After suppressing from the notation all pullbacks from $\X_\T$ to $\G$, we write the universal off-diagonal exact sequence on $\G$ as follows:
\begin{eqn}
	0
	\too
	\cal{B}
	\too
	\cal{End} (\cal{E})
	\overset{q}{\too}
	\cal{Q}
	\too 0
\fullstop
\end{eqn}
Consider the closed analytic subspace $\C \subset \G$ parametrising the subbundles $\cal{B} \subset \cal{End} (\cal{E})$ that are commutative unital subalgebras.
Concretely, it can be defined as the common vanishing locus of three maps as follows.
First, the identity endomorphism on $\cal{E}$ is a global section of $\cal{End} (\cal{E})$, thus defining a map $i : \cal{O}_\G \to \cal{End} (\cal{E})$.
Second, the commutator on $\cal{End} (\cal{E})$ gives a map $c : \cal{B} \otimes \cal{B} \to \cal{End} (\cal{E})$ sending $b_1 \otimes b_2 \mapsto [b_1, b_2]$.
Third, the composition on $\cal{End} (\cal{E})$ gives a map $m : \cal{B} \otimes \cal{B} \to \cal{End} (\cal{E})$ sending $b_1 \otimes b_2 \mapsto b_1 \circ b_2$.
Let $\C$ be the common vanishing locus of $c$ and the composite maps $q \circ i : \cal{O}_\G \to \cal{Q}$ and $q \circ m : \cal{B} \otimes \cal{B} \to \cal{Q}$:
\begin{eqn}
	\C \coleq \sfop{Zero} (q \circ i) \cap \sfop{Zero} (q \circ m) \cap \sfop{Zero} (c) \subset \G
\fullstop
\end{eqn}

The first vanishing condition ensures that $\cal{B}$ contains the identity, the second ensures that it is closed under multiplication, and the third ensures that it is commutative.
These are vanishing conditions for morphisms of vector bundles on $\G$, so $\C$ is defined by finitely many local holomorphic equations.

Next, we impose the condition that $\cal{E}$ be locally free of rank one over $\cal{B}$.
This defines an open analytic subspace $\Z \subset \C$ equipped with a separated holomorphic map $u : \Z \to \X_\T$.
To see this, consider the total space of $\cal{E}$ over $\C$.
At a vector $e \in \tot (\cal{E})$, the evaluation map $\cal{B} \to \cal{E}$ sending $b \mapsto be$ is an isomorphism precisely when its determinant is nonzero, which therefore defines an open subset of $\tot (\cal{E})$.
Its image in $\C$ is open because the bundle projection $\tot (\cal{E}|_\C) \to \C$ is an open map; define $\Z$ to be this image.
A local choice of such a vector gives a $\cal{B}$-linear local isomorphism $\cal{B} \iso \cal{E}$.
The restriction of the relative Grassmannian's projection $\G \to \X_\T$ to a projection $u : \Z \to \X_\T$ is a separated holomorphic map because $\Z$ is an open subspace of a closed analytic subspace of $\G$.

Consequently, for every analytic space $\T' \to \T$, spectral structures on the pulled-back bundle $\cal{E}_{\T'}$ correspond naturally to sections of the separated holomorphic map $u_{\T'} : \Z_{\T'} \to \X_{\T'}$.
By \autoref{1792526939305}, this relative functor of sections is represented by an analytic space $\M$, locally of finite presentation over $\T$.
The preceding correspondence therefore gives a natural equivalence
\begin{eqntag}
\label{1793101181533}
    \T \times_{\frak{Bun}_n (\X / \S)} \frak{SpStr}_n (\X / \S)
    \cong
    \M
\fullstop
\end{eqntag}
Since $\T$ and $\cal{E}$ were arbitrary, the forgetful morphism $F$ is representable by analytic spaces and locally of finite presentation.

Finally, by \autoref{1793142766507}, choose a smooth surjective atlas $\V \too \frak{Bun}_n (\X / \S)$ and form the base change
\begin{eqn}
	\U \coleq \V \times_{\frak{Bun}_n (\X / \S)} \frak{SpStr}_n (\X / \S)
\fullstop
\end{eqn}
This is an analytic space by the representability of $F$.
The projection $\U \to \frak{SpStr}_n (\X / \S)$ is smooth and surjective because it is the base change of the chosen atlas.
The diagonal of $\frak{SpStr}_n (\X / \S)$ is representable because $F$ and the diagonal of $\frak{Bun}_n (\X / \S)$ are representable.
Thus, $\frak{SpStr}_n (\X / \S)$ is an analytic stack.

Moreover, $\U \to \V$ is locally of finite presentation, and $\V \to \S$ is smooth by \autoref{1793142766507} and the choice of atlas.
Their composite $\U \to \S$ is therefore locally of finite presentation.
This proves that $\frak{SpStr}_n (\X / \S)$ is locally of finite presentation over $\S$.

For the base-change assertion, let $\S' \to \S$ be any holomorphic map.
For every analytic space $\T \to \S'$, the canonical identification
\begin{eqn}
	(\X \times_\S \S') \times_{\S'} \T
	\cong
	\X \times_\S \T
\end{eqn}
identifies the groupoids of spectral structures on the two sides.
These identifications are natural in $\T$ and compatible with the pullback functors of \autoref{1791551111381}, so they give the asserted equivalence of stacks.
Since they preserve the underlying vector bundles, they also identify the corresponding forgetful morphism with the base change of $F$.
This completes the proof of \autoref{1792494933366}.
\end{proof}

For a fixed vector bundle $\cal{E}$ on a fixed compact Riemann surface $\X$, the essential fibre of $F$ is an analytic space:
\begin{eqntag}
\label{1793046205260}
	\frak{SpStr} (\cal{E})
	\coleq
	\frak{SpStr}_n (\X)
	\times_{\frak{Bun}_n (\X)}
	\ast
\fullstop{,}
\end{eqntag}
where $\ast \to \frak{Bun}_n (\X)$ is the point classifying $\cal{E}$.
An object of this fibre consists of a spectral structure together with an identification of its underlying vector bundle with $\cal{E}$.
As in \eqref{1791390052260}, this identification is part of the data.
Consequently, an automorphism in the fibre must induce the identity on $\cal{E}$, and hence also on its spectral algebra because the spectral action is injective.

On the other hand, in the moduli stack of spectral structures, scalar automorphisms of the underlying vector bundle are allowed.
Every spectral structure therefore has a central subgroup $\CC^\times$ in its automorphism group.
Thus, $\frak{SpStr}_n (\X)$ is not an analytic Deligne-Mumford stack (some authors would call it an analytic \textit{Artin stack}), even though the essential fibres of $F$ are analytic spaces.

\begin{example}[obstructed spectral structure]
\label{1792516769027}
This example shows that deformations of spectral structures can be obstructed even if the underlying vector bundle is allowed to vary.

Take $\S = \ast$ and let $\X$ be a compact Riemann surface of positive genus.
We will construct a first-order deformation of a rank-two spectral structure which does not extend to second order.
Consider the algebra $\cal{A}_0 \coleq \cal{O}_\X [u] / \inner{u^2}$ and equip $\cal{E}_0 \coleq \cal{A}_0$ with its multiplication action $\rho_0 : \cal{A}_0 \to \cal{End} (\cal{E}_0)$.
This defines a spectral structure whose spectral realisation is the doubled curve $\Y_0 \coleq \Spec_\X (\cal{A}_0)$ with spectral line bundle $\cal{L}_0 = \cal{O}_{\Y_0}$.

Since $\X$ has positive genus, we can choose an open cover $\X = \bigcup_i \U_i$ and a holomorphic Čech cocycle $(c_{ij})$ representing a nonzero class in $\H^1 (\X, \cal{O}_\X)$.
Consider the order-1 and order-2 fat points
\begin{eqn}
	\DD_1 \coleq \Spec \big( \CC [t] / \inner{t^2} \big)
\qqtext{and}
	\DD_2 \coleq \Spec \big( \CC [t] / \inner{t^3} \big)
\fullstop
\end{eqn}
Over each $\U_i \times \DD_1$, consider the algebra
\begin{eqn}
	\cal{A}_{1,i}
	\coleq
	\cal{O}_{\U_i \times \DD_1} [u_i] / \inner{u_i^2 - t}
\fullstop
\end{eqn}
Then glue these algebras on double overlaps by $u_i = (1 + t c_{ij}) u_j$.
These transition maps preserve the defining relation because
\begin{eqn}
	\big( (1 + t c_{ij}) u_j \big)^2
	=
	(1 + 2t c_{ij}) t
	=
	t
	\mod
	t^2
\fullstop
\end{eqn}
Their cocycle condition follows from the additive cocycle condition for $(c_{ij})$.
Thus, they define a commutative unital algebra $\cal{A}_1$ on $\X \times \DD_1$, locally free of rank two with local basis $1, u_i$.
Its restriction to $t = 0$ is $\cal{A}_0$.
Taking $\cal{E}_1 \coleq \cal{A}_1$ with its multiplication action therefore gives a first-order deformation of $(\cal{E}_0, \cal{A}_0, \rho_0)$.

We claim that this spectral structure does not extend over $\DD_2$.
Suppose otherwise.
Then its spectral algebra would give a commutative unital algebra $\cal{A}_2$ on $\X \times \DD_2$, locally free of rank two, together with an identification of its restriction to $\DD_1$ with $\cal{A}_1$.
The trace map gives a canonical decomposition
\begin{eqn}
	\cal{A}_2
	=
	\cal{O}_{\X \times \DD_2}
	\oplus
	\cal{N}_2
\qqtext{where}
	\cal{N}_2 \coleq \cal{ker} \big( \tr : \cal{A}_2 \to \cal{O}_{\X \times \DD_2} \big)
\fullstop
\end{eqn}
Notice that $\cal{N}_2$ is a line bundle because the trace of the identity equals $2$.
Its restriction to $\DD_1$ is the trace-zero line subbundle locally generated by $u_i$.

After refining the cover if necessary, choose local generators $v_i$ of $\cal{N}_2$ which restrict to $u_i$ on $\DD_1$.
The square of a trace-zero element of a rank-two algebra is scalar by the Cayley-Hamilton Theorem.
Consequently, the defining relations and transition maps must have the form
\begin{eqn}
	v_i^2 = t + t^2 b_i
\qqtext{and}
	v_i = (1 + t c_{ij} + t^2 d_{ij}) v_j
\fullstop{,}
\end{eqn}
where $b_i$ is holomorphic on $\U_i$ and $d_{ij}$ is holomorphic on $\U_i \cap \U_j$.
Comparing the squares on an overlap gives
\begin{eqn}
	t + t^2 b_i
	=
	(1 + t c_{ij} + t^2 d_{ij})^2 (t + t^2 b_j)
	=
	t + t^2 (b_j + 2c_{ij})
	\mod
	t^3
\fullstop
\end{eqn}
Therefore, $b_i - b_j = 2c_{ij}$, which means $(c_{ij})$ is a coboundary, contrary to its choice.
So $\cal{A}_1$ cannot extend over $\DD_2$.

It follows that $\frak{SpStr}_2 (\X)$ is not formally smooth --- and hence is not smooth --- at $(\cal{E}_0, \cal{A}_0, \rho_0)$.
Notice that the obstruction persists even though the underlying vector bundle is allowed to vary.
\end{example}

\begin{example}[obstructions in higher rank]
\label{1792530283401}
The preceding \autoref{1792516769027} extends to every rank $n > 2$.
Consider the product algebra
\begin{eqn}
	\cal{B}_1
	\coleq
	\cal{A}_1
	\times
	\cal{O}_{\X \times \DD_1}
	\times \cdots \times
	\cal{O}_{\X \times \DD_1}
\fullstop{,}
\end{eqn}
with $n-2$ factors of $\cal{O}_{\X \times \DD_1}$, and equip $\cal{E}_1 \coleq \cal{B}_1$ with its multiplication action.
Any lift over $\DD_2$ would lift the idempotents defining this product decomposition, since idempotents lift uniquely across nilpotent ideals \cite[\href{https://stacks.math.columbia.edu/tag/00J9}{Tag 00J9}]{stacks-project}.
Its rank-two direct factor would then give a lift of $\cal{A}_1$, which we have shown to be impossible.
Thus, $\frak{SpStr}_n (\X)$ is not smooth for any $n \geq 2$ when $\X$ has positive genus.

More generally, if a family $p : \X \to \S$ has a fibre $\X_s$ of positive genus, then the structure morphism $\frak{SpStr}_n (\X / \S) \to \S$ is not smooth at the corresponding points of its fibre over $s$.
Indeed, its base change to $s$ is $\frak{SpStr}_n (\X_s)$, and smoothness is preserved by base change.
\end{example}

\subsection{Spectral Realisations and Distinguished Smooth Open Substacks}
\label{1793136202900}

The full substacks of smooth, connected, and/or clean spectral structures are determined by \autoref{1791482794969} and denoted by adding a superscript `sm', `con', and `clean`, respectively.

\begin{theorem}[distinguished smooth open substacks]
\label{1792494388516}
The moduli stack of spectral structures $\frak{SpStr}_n (\X / \S)$ has the following distinguished open analytic substacks, each of which is smooth over $\S$:
\begin{eqntag}
\label{1793009278980}
	\frak{SpStr}_n^\rm{clean} (\X / \S)
	\subset \frak{SpStr}_n^{\rm{sm},\rm{con}} (\X / \S)
	\subset \frak{SpStr}_n^{\rm{sm}} (\X / \S)
	\subset \frak{SpStr}_n (\X / \S)
\fullstop
\end{eqntag}
Arbitrary analytic base change commutes with the formation of these open substacks.
In general, none of these inclusions is dense; in fact, the second inclusion is open-and-closed.
Furthermore, the forgetful morphism $F$ restricted to clean spectral structures,
\begin{eqntag}
\label{1792514157603}
	F :
	\frak{SpStr}_n^\rm{clean} (\X / \S)
	\too
	\frak{Bun}_n (\X / \S)
\fullstop{,}
\end{eqntag}
is smooth and surjective, and is representable by analytic spaces locally of finite presentation.
\end{theorem}

The following immediate corollary is the stack-theoretic formulation of \autoref{1791544458956}:

\begin{corollary}[local lifting]
\label{1793219634470}
Every holomorphic map $\T \to \frak{Bun}_n (\X / \S)$ lifts locally on $\T$ through clean spectral structures via the forgetful functor $F$:
\begin{eqntag}
\label{1792580234202}
\begin{tikzcd}
&	\frak{SpStr}_n^\rm{clean} (\X / \S)
		\ar[d, "F"]
\\
	\T
		\ar[r]
		\ar[ur, dashed]
&	\frak{Bun}_n (\X / \S)
\fullstop
\end{tikzcd}
\end{eqntag}
\end{corollary}

The proof of this theorem goes through the geometry of spectral realisations.
Namely, the family construction in \ref{1792146709610} yields the stack of spectral realisations over the analytic site $\mathbf{An}$:
\begin{eqn}
	\frak{SpReal}_n (\X / \S)
\end{eqn}
Concretely, an object over $\T$ is a pair $(\Y \tinyoverset{\pi}{\to} \X_\T, \cal{L})$ where $\pi$ is finite flat of degree $n$ and $\cal{L}$ is a line bundle on $\Y$.
The corresponding vector bundle is $\cal{E} = \pi_\ast \cal{L}$ over the family $\X_\T$.
Isomorphisms are isomorphisms of covers together with compatible isomorphisms of their line bundles.
Pullback is defined by base change of the cover and pullback of the line bundle.
The equivalence of groupoids from \autoref{1789153885363} is compatible with base change, and therefore gives an equivalence of stacks as follows.

\begin{lemma}[description as spectral realisations]
\label{1792502248471}
There is a canonical equivalence of stacks over the analytic site, compatible with arbitrary analytic base change:
\begin{eqn}
	\frak{SpStr}_n (\X / \S)
	\cong
	\frak{SpReal}_n (\X / \S)
\fullstop
\end{eqn}
In particular, $\frak{SpReal}_n (\X / \S)$ is an analytic stack locally of finite presentation over $\S$, and its formation commutes with arbitrary analytic base change.
\end{lemma}

\begin{proof}
For every analytic space $\T \to \S$, the finite sheaf-module correspondence gives the equivalence between spectral structures on $\X_\T$ and their spectral realisations, as in \autoref{1789153885363}.
The inverse sends $(\Y \tinyoverset{\pi}{\to} \X_\T, \cal{L})$ to $\cal{E} \coleq \pi_\ast \cal{L}$ and $\cal{A} \coleq \pi_\ast \cal{O}_\Y$ with the multiplication action $\rho : \cal{A} \to \cal{End} (\cal{E})$.
By \autoref{1791551111381}, these constructions commute canonically with arbitrary base change.
They therefore give an equivalence of categories fibred in groupoids over the analytic site.
So $\frak{SpReal}_n (\X / \S)$ is a stack because $\frak{SpStr}_n (\X / \S)$ is.
A smooth surjective atlas for $\frak{SpStr}_n (\X / \S)$ gives a smooth surjective atlas for $\frak{SpReal}_n (\X / \S)$ through the equivalence, and representability of the diagonal is preserved.
Analyticity and local finite presentation over $\S$ now follow from \autoref{1792494933366}.
The base-change assertion follows from the compatibility of the equivalence with base change and \autoref{1792494933366}.
\end{proof}

\enlargethispage{30pt}

\begin{proof}[Proof of \autoref{1792494388516}]
We prove this by working with spectral realisations through the equivalence of \autoref{1792502248471}.

\caps{Step 1: The substack of smooth spectral structures.}
Consider a family of spectral structures parametrised by an analytic space $\T \to \S$, with spectral realisation $\pi : \Y \to \X_\T$.
The map $q \coleq p_\T \circ \pi : \Y \to \T$ is proper and flat.
Its non-smooth locus is closed in $\Y$, and its image in $\T$ is closed because $q$ is proper.
The complement of this image is precisely the set of parameters $t \in \T$ for which $\Y_t$ is smooth.
Thus, $\frak{SpStr}_n^{\rm{sm}} (\X / \S)$ is an open analytic substack.

To prove smoothness over $\S$, we use the infinitesimal lifting criterion \cite[Théorème~3.1]{1791539610207}, applied on smooth atlases.
The substack is locally of finite presentation by \autoref{1792494933366}.
As in the proof of \autoref{1793142766507}, it suffices to solve lifting problems across small extensions, with the identification over the original parameter space prescribed and the map to $\S$ fixed.

So let $\R' \to \R$ be a small extension of local Artinian $\CC$-algebras, with kernel $\I$ annihilated by the maximal ideal of $\R'$.
Put $\T \coleq \Spec (\R)$ and $\T' \coleq \Spec (\R')$.
Suppose we are given a map $\T' \to \S$ and a smooth spectral structure over $\X_\T$.
Write its spectral realisation as $\pi_\T : \Y_\T \to \X_\T$ with spectral line bundle $\cal{L}_\T$.
The goal is to lift the cover and the line bundle over the target $\X_{\T'}$.

Let $\pi_0 : \Y_0 \to \X_0$ denote the cover on the closed fibre.
Choose disjoint small coordinate discs around its branch values.
After shrinking these discs, their inverse images are disjoint unions of discs, and on each component the closed-fibre map has the form $x = y^r$ where $r \geq 1$ is the ramification index of that component.
Since $\X_{\T'} \to \T'$ is smooth, we may choose these discs inside relative coordinate neighbourhoods of $\X_{\T'}$.

The coordinate $y$ lifts to a holomorphic function on the corresponding component of $\Y_\T$ because the underlying disc is Stein and the thickening is nilpotent.
The sections $1, y, \ldots, y^{r-1}$ form a basis of the direct-image algebra on the closed fibre, and hence also a basis over $\T$ by Nakayama's Lemma.
Consequently, this component has a monic presentation
\begin{eqn}
	P (x,y) = y^r + a_{r-1} (x) y^{r-1} + \cdots + a_0 (x) = 0
\fullstop{,}
\end{eqn}
whose reduction to the closed fibre is $y^r - x = 0$.
Lift the coefficients of $P$ from $\R$ to $\R'$.
The resulting monic equation defines a finite flat cover of the same degree over the corresponding relative disc in $\X_{\T'}$.
Its source is smooth over $\T'$ because the derivative of its defining equation with respect to $x$ reduces to $-1$.

Over the complement of the branch values, the given cover $\pi_\T$ is finite étale.
Finite étale analytic covers lift uniquely, up to unique isomorphism, across nilpotent thickenings.
Indeed, they are locally finite disjoint unions of copies of the target, with locally constant permutations as transition functions, and a nilpotent thickening does not change the underlying topological space.
The same observation shows that isomorphisms between finite étale covers lift uniquely.

The covers constructed over the relative discs are finite étale on their overlaps with this complement.
The original overlap identifications therefore lift uniquely, and their lifts satisfy the cocycle condition.
They glue the local constructions to a finite flat cover $\pi_{\T'} : \Y_{\T'} \too \X_{\T'}$ which restricts to $\pi_\T$ and whose source is smooth over $\T'$.
If $\pi_0$ is unramified, the finite étale lifting argument applies directly to the entire cover.

It remains to lift $\cal{L}_\T$ to $\Y_{\T'}$.
Flatness over $\R'$ identifies the ideal of $\Y_\T$ in $\Y_{\T'}$ with $\cal{O}_{\Y_0} \otimes_\CC \I$.
The exact sequence of sheaves of units for this square-zero thickening places the obstruction to lifting $\cal{L}_\T$ in $\H^2 (\Y_0, \cal{O}_{\Y_0}) \otimes_\CC \I = 0$.
The vanishing holds because $\Y_0$ is a compact curve.
Thus, there is a line bundle $\cal{L}_{\T'}$ lifting $\cal{L}_\T$.

Finally, put $\cal{E}_{\T'} \coleq (\pi_{\T'})_\ast \cal{L}_{\T'}$ and $\cal{A}_{\T'} \coleq (\pi_{\T'})_\ast \cal{O}_{\Y_{\T'}}$.
Since $\pi_{\T'}$ is finite flat of degree $n$, the sheaf $\cal{E}_{\T'}$ is a vector bundle of rank $n$.
The multiplication action of $\cal{A}_{\T'}$ makes it a spectral structure, and compatibility with finite base change identifies its restriction with the given structure over $\T$.
This solves the lifting problem and proves that $\frak{SpStr}_n^\rm{sm} (\X / \S)$ is smooth over $\S$.

\caps{Step 2: The substack of smooth and connected spectral structures.}
For a family in $\frak{SpStr}_n^\rm{sm} (\X / \S)$, the map $q : \Y \to \T$ is proper and smooth.
The number of connected components of its fibres is locally constant on $\T$ \cite[Ch.~III, Proposition~3.12]{BanicaStanasila1976}.
Consequently, the set of parameters for which $\Y_t$ is connected is both open and closed.
This gives the asserted open-and-closed analytic substack, which is smooth over $\S$ because $\frak{SpStr}_n^\rm{sm} (\X / \S)$ is smooth over $\S$.

\caps{Step 3: The substack of clean spectral structures.}
Consider a family in the substack $\frak{SpStr}_n^{\rm{sm},\rm{con}} (\X / \S)$.
The relative differential 
\begin{eqn}
	\d{\pi} : \pi^\ast \Omega^1_{\X_\T/\T} \to \Omega^1_{\Y/\T}
\end{eqn}
is a map of line bundles whose restriction to every fibre is nonzero, and its zero locus is the relative ramification divisor $\R \subset \Y$.
Since its restriction to each fibre is an effective Cartier divisor, $\R$ is flat over $\T$.
It is also proper with finite fibres, and hence finite flat over $\T$.

\clearpage

Suppose that $\pi_t$ is simply branched at some parameter $t \in \T$.
Then $\R_t$ is reduced, so the map $\R \to \T$ is finite étale after shrinking to a neighbourhood around $t$.
Locally on $\T$, it is therefore a disjoint union of sections.
Their images under $\pi$ are pairwise distinct at $t$, and remain pairwise distinct after a further shrinking.
Thus, simple branching is an open condition on $\T$.

The off-diagonal sheaf $\cal{Q}$ is a vector bundle on $\X_\T$, and its definition commutes with arbitrary base change by \autoref{1791551111381}.
In particular, it is flat over $\T$, and its restriction to $\X_t$ is the off-diagonal bundle of the spectral structure on that fibre.
Upper semicontinuity of cohomology for the proper map $p_\T : \X_\T \to \T$ shows that the subset of $t \in \T$ satisfying $\H^1 (\X_t, \cal{Q}_t) = 0$ is open; see \cite[Theorem~6.7.4]{GuiCAG} and \cite[Ch.~III, Theorem~4.12(i)]{BanicaStanasila1976}.
Intersecting this open set with the simply branched locus gives precisely the clean locus.
Therefore, it defines an open analytic substack of $\frak{SpStr}_n^{\rm{sm},\rm{con}} (\X / \S)$, which is smooth over $\S$.

\caps{Step 4: The forgetful morphism \eqref{1792514157603}.}
Since the substack $\frak{SpStr}_n^{\rm{clean}} (\X / \S)$ is open, representability and local finite presentation of \eqref{1792514157603} follow from \autoref{1792494933366}.
We prove smoothness by the infinitesimal lifting criterion.

We use the same notation $\R' \to \R$, $\I$, $\T$, and $\T'$ as in Step 1.
Suppose we are given a vector bundle $\cal{E}_{\T'}$ on $\X_{\T'}$ together with a clean spectral structure on its restriction $\cal{E}_\T$.
Write the corresponding spectral realisation as $\pi_\T : \Y_\T \to \X_\T$ with spectral line bundle $\cal{L}_\T$ and identification $\iota_\T : (\pi_\T)_\ast \cal{L}_\T \iso \cal{E}_\T$.

The same argument as in Step 1 allows us to lift this cover to a finite flat cover $\pi_{\T'} : \Y_{\T'} \to \X_{\T'}$ extending $\pi_\T$, and the vanishing of $\H^2 (\Y_0, \cal{O}_{\Y_0})$ over the closed fibre gives a line bundle $\tilde{\cal{L}}_{\T'}$ extending $\cal{L}_\T$.
Put $\tilde{\cal{E}}_{\T'} \coleq (\pi_{\T'})_\ast \tilde{\cal{L}}_{\T'}$.
Through $\iota_\T$, both $\tilde{\cal{E}}_{\T'}$ and $\cal{E}_{\T'}$ are lifts of the same vector bundle $\cal{E}_\T$, with identified restrictions.

The difference between these two lifts is a class 
\begin{eqn}
	\delta \in \H^1 \big( \X_0, \cal{End} (\cal{E}_0) \big) \otimes_\CC \I
\end{eqn}
over the closed fibre.
On the other hand, the isomorphism classes of lifts of $\cal{L}_\T$ to the chosen cover $\Y_{\T'}$, with identified restriction to $\Y_\T$, form a torsor under $\H^1 (\Y_0, \cal{O}_{\Y_0}) \otimes_\CC \I$ \cite[\href{https://stacks.math.columbia.edu/tag/08VR}{Tag 08VR}]{stacks-project}.
Changing the lift of $\cal{L}_\T$ changes its pushforward through the natural map
\begin{eqn}
	\H^1 (\Y_0, \cal{O}_{\Y_0})
	\cong
	\H^1 (\X_0, \cal{A}_0)
	\too
	\H^1 \big( \X_0, \cal{End} (\cal{E}_0) \big)
\fullstop
\end{eqn}
Indeed, multiplying the transition functions of the lifted line bundle by a cocycle of units of the form $1+\alpha_{ij}$ changes the transition functions of its pushforward by the corresponding multiplication endomorphisms under the spectral action $\rho_0 : \cal{A}_0 \to \cal{End} (\cal{E}_0)$ over the closed fibre.

Since the spectral structure $(\cal{A}_0, \rho_0)$ is clean, we have $\H^1 (\X_0, \cal{Q}_0) = 0$.
Then the long exact sequence in cohomology induced by the off-diagonal exact sequence shows that the displayed map is surjective.
We may consequently modify $\tilde{\cal{L}}_{\T'}$ to a lift $\cal{L}_{\T'}$ whose pushforward has the same deformation class as $\cal{E}_{\T'}$.
Thus, there is an isomorphism $\iota_{\T'} : (\pi_{\T'})_\ast \cal{L}_{\T'} \iso \cal{E}_{\T'}$ which restricts to the prescribed isomorphism $\iota_\T$.
Transporting the multiplication action through $\iota_{\T'}$ gives the required lift of the spectral structure.

The lifted structure is fibrewise clean because the spectral structure over the closed fibre is unchanged.
This solves every infinitesimal lifting problem for \eqref{1792514157603}, proving that the restricted forgetful morphism is smooth.
Apply \autoref{1791544458956} to the vector bundle on $\X_\T$ corresponding to any holomorphic map $\T \to \frak{Bun}_n (\X / \S)$.
It gives a clean spectral structure locally around every point of $\T$.
This proves surjectivity of $F$.

Finally, the base-change equivalence of \autoref{1792494933366} preserves the fibrewise conditions of \autoref{1791482794969} and therefore restricts to each of the distinguished open substacks.
\end{proof}

\begin{proof}[Proof of \autoref{1793219634470}]
By \autoref{1792494388516}, the base change of $F$ along $\T \to \frak{Bun}_n (\X / \S)$ is a smooth surjective holomorphic map of analytic spaces.
Such a map admits local holomorphic sections, which give the required lifts.
\end{proof}

\begin{corollary}[sufficient criterion for smoothness of $F$]
\label{1793046223705}
In fact, the forgetful morphism
\begin{eqn}
	F :
	\frak{SpStr}_n (\X / \S)
	\too
	\frak{Bun}_n (\X / \S)
\end{eqn}
is smooth at every point represented by a smooth spectral structure $(\cal{A}, \rho)$ on a vector bundle $\cal{E}$ over a fibre $\X_s$ whose off-diagonal bundle satisfies $\H^1 (\X_s, \cal{Q}) = 0$.
No connectedness or simple branching assumption is required.
\end{corollary}

\begin{proof}
Tracing through the proof of \autoref{1792494388516}, we can observe that the infinitesimal lifting argument for $F$ first lifts the smooth cover and its spectral line bundle, and then modifies the line-bundle lift to correct its pushforward.
The first step requires only smoothness of the spectral curve, whilst the correction uses only the vanishing of $\H^1 (\X_s, \cal{Q})$.
The same argument therefore proves that $F$ is smooth at every such point.
\end{proof}

The following proposition gives more information about the substack of smooth spectral structures.
It says that this substack is \textit{not} the same as the smooth locus of the moduli stack of spectral structures (except in trivial scenarios). 

\begin{proposition}[smooth spectral curves and smooth moduli points]
\label{1792514390194}
The moduli substack of smooth spectral structures is contained in the smooth locus over $\S$ of the moduli stack of all spectral structures.
In fact, if $n \geq 2$ and $\S$ is nonempty, this containment is strict:
\begin{eqntag}
\label{1793011620151}
	\frak{SpStr}_n^{\rm{sm}} (\X / \S)
	\subsetneq
	\big( \frak{SpStr}_n (\X / \S) \big)^{\rm{sm}}
\fullstop
\end{eqntag}
On the other hand, for $n = 1$, all the inclusions in \eqref{1793009278980} are equalities, and all these substacks coincide with the relative Picard stack:
\begin{eqntag}
\label{1793009445991}
	\frak{SpStr}_1^{\rm{clean}} (\X / \S)
	= \frak{SpStr}_1^{\rm{sm,con}} (\X / \S)
	= \frak{SpStr}_1^{\rm{sm}} (\X / \S)
	= \frak{SpStr}_1 (\X / \S)
	= \frak{Pic} (\X/\S)
\fullstop
\end{eqntag}
\end{proposition}

Recall that the relative Picard stack $\frak{Pic} (\X/\S)$ is the smooth analytic stack whose objects over $\T \to \S$ are line bundles on $\X_\T$.

\begin{proof}
By \autoref{1792494388516}, $\frak{SpStr}_n^\rm{sm} (\X / \S)$ is open in $\frak{SpStr}_n (\X / \S)$ and smooth over $\S$, so it is contained in the smooth locus of the structure morphism $\frak{SpStr}_n (\X / \S) \to \S$.
That this inclusion is always strict when $n \geq 2$ and $\S \neq \emptyset$ is demonstrated in \autoref{1793011572262} below which gives a general but explicit construction of a singular spectral structure which is nevertheless a smooth point of the moduli stack of spectral structures.
Finally, for $n = 1$, every spectral algebra $\cal{A}$ is canonically identified with $\cal{O}_{\X_\T}$ via the unit map $\cal{O}_{\X_\T} \to \cal{A}$, and therefore the off-diagonal bundle in this case is always zero: $\cal{Q} = \cal{End} (\cal{E}) / \cal{A} \cong \cal{O}_{\X_\T} / \cal{O}_{\X_\T}$.
Consequently, its spectral cover is the identity cover, and this leads directly to the identity \eqref{1793009445991}.
\end{proof}

\begin{example}[unobstructed singular spectral structure]
\label{1793011572262}
To demonstrate that the inclusion \eqref{1793011620151} of \autoref{1792514390194} is always strict, we give a general explicit construction of a spectral structure that is singular but nevertheless a smooth point of the moduli stack of spectral structures.
This construction is made even more concrete in \autoref{1792517279249} below.

First, take $n = 2$.
Fix $s_0 \in \S$, put $\X_0 \coleq \X_{s_0}$, and adapt the cyclic double cover construction from \autoref{1791437029550}.
Namely, choose a sufficiently positive line bundle $\cal{M}$ on $\X_0$ such that $\H^1 (\X_0, \cal{M}^{\otimes 2}) = 0$ and such that there is a section $\sigma \in \H^0 (\X_0, \cal{M}^{\otimes 2})$ with simple zeros except for exactly one double zero at some point $x_0 \in \X_0$.
Such a section can be obtained by taking a general section of $\cal{M}^{\otimes 2} (-2x_0)$ which does not vanish at $x_0$.
Define the algebra 
\begin{eqn}
	\cal{A}_0 \coleq \cal{O}_{\X_0} \oplus \cal{M}^{-1}
\end{eqn}
with multiplication determined by $\sigma : \cal{M}^{-2} \to \cal{O}_{\X_0}$.
Take $\cal{E}_0 \coleq \cal{A}_0$ and equip it with its canonical multiplication action $\rho_0 : \cal{A}_0 \to \cal{End} (\cal{E}_0)$.
The spectral cover $\pi_0 : \Y_0 \to \X_0$ of $\cal{A}_0$ has an ordinary node above the double zero of $\sigma$, so this spectral structure $(\cal{E}_0, \cal{A}_0, \rho_0)$ does not belong to $\frak{SpStr}_2^\rm{sm} (\X / \S)$.

Nevertheless, we argue that $(\cal{E}_0, \cal{A}_0, \rho_0)$ belongs to the smooth locus of $\frak{SpStr}_2 (\X / \S)$ over $\S$.
To see this, consider a small extension $\R' \to \R$ of local Artinian $\CC$-algebras, with kernel $\I$ annihilated by the maximal ideal $\frak{m}$ of $\R'$:
\begin{eqntag}
\label{1793048015827}
	0 \to \I \to \R' \to \R \to 0
\qqtext{and}
	\frak{m} \I = 0
\fullstop
\end{eqntag}
Fix a morphism $\Spec (\R') \to \S$ whose closed point maps to $s_0$, and put
\begin{eqn}
	\X_{\R'} \coleq \X \times_\S \Spec (\R')
\qqtext{and}
	\X_\R \coleq \X \times_\S \Spec (\R)
\fullstop
\end{eqn}
Suppose that a deformation $(\cal{E}_\R, \cal{A}_\R, \rho_\R)$ of the spectral structure $(\cal{E}_0, \cal{A}_0, \rho_0)$ over $\X_\R$ is given, with a specified identification of its restriction to $\X_0$ with $(\cal{E}_0, \cal{A}_0, \rho_0)$.
For a local section $a$ of $\cal{A}_\R$, let $\tr(a)$ denote the trace of multiplication by $a$ on the rank-two vector bundle $\cal{A}_\R$.
Since $\tr(1) = 2$, the map $\frac{1}{2} \tr$ is a retraction of the unit inclusion $\cal{O}_{\X_\R} \inj \cal{A}_\R$.
Its kernel is therefore a line bundle, giving the decomposition
\begin{eqn}
	\cal{A}_\R
	=
	\cal{O}_{\X_\R} \oplus \cal{N}_\R
\qqtext{where}
	\cal{N}_\R \coleq \cal{ker} (\tr)
\fullstop
\end{eqn}
The line bundle $\cal{M}_\R \coleq \cal{N}_\R^{-1}$ deforms $\cal{M}$ because the trace-zero summand $\cal{N}_\R$ on the closed fibre is $\cal{M}^{-1}$.
By the Cayley-Hamilton Theorem, the square of a trace-zero element is scalar.
Thus, multiplication on $\cal{M}_\R^{-1}$ is determined by a section $\sigma_\R \in \H^0 (\X_\R, \cal{M}_\R^{\otimes 2})$.

The goal is to lift the spectral structure $(\cal{E}_\R, \cal{A}_\R, \rho_\R)$ to a spectral structure $(\cal{E}_{\R'}, \cal{A}_{\R'}, \rho_{\R'})$ over $\X_{\R'}$.
To do so, we lift the line bundle $\cal{M}_\R$ to a line bundle $\cal{M}_{\R'}$ over $\X_{\R'}$, and then lift the section $\sigma_\R$ to a section $\sigma_{\R'} \in \H^0 (\X_{\R'}, \cal{M}_{\R'}^{\otimes 2})$.
This will define the candidate spectral algebra $\cal{A}_{\R'}$ lifting the spectral algebra $\cal{A}_\R$.
Then we will lift the spectral line bundle $\cal{L}_\R$ over the realisation $\Y_\R$ of $\cal{A}_\R$ to a line bundle $\cal{L}_{\R'}$ over $\Y_{\R'}$, and push it forward to define the vector bundle $\cal{E}_{\R'}$.

Consider the ideal $\cal{I}$ of $\X_\R$ in $\X_{\R'}$, defined by the short exact sequence
\begin{eqntag}
\label{1793048887757}
	0 \to \cal{I} \to \cal{O}_{\X_{\R'}} \to \cal{O}_{\X_\R} \to 0
\fullstop{,}
\end{eqntag}
We claim that there is a canonical isomorphism $\cal{I} \cong \cal{O}_{\X_0} \otimes \I$.
Indeed, since $\X_{\R'}$ is flat over $\R'$, tensoring the exact sequence \eqref{1793048015827} with the structure sheaf of $\X_{\R'}$ is exact.
So we get $\cal{I} \cong \cal{O}_{\X_{\R'}} \otimes_{\R'} \I \cong \cal{O}_{\X_0} \otimes_\CC \I$.
The second identification holds because $\frak{m} \I = 0$, so the $\R'$-action on $\I$ factors through the residue field $\CC$.

This ideal satisfies $\cal{I}^2 = 0$ because $\I^2 = 0$ since $\frak{m} \I = 0$.
As a result, composing the inclusion $\cal{I} \inj \cal{O}_{\X_{\R'}}$ with the exponential map $\exp : \cal{O}_{\X_{\R'}} \to \cal{O}_{\X_{\R'}}^\times$ is simply the map $a \mapsto 1 + a$.
Consider the short exact sequence
\begin{eqn}
	0 \to \cal{I} \to \cal{O}_{\X_{\R'}}^\times \to \cal{O}_{\X_\R}^\times \to 1
\fullstop
\end{eqn}
Taking cohomology and recalling that $\H^1 (\X_\R, \cal{O}_{\X_\R}^\times) = \Pic (\X_\R)$ and $\H^1 (\X_{\R'}, \cal{O}_{\X_{\R'}}^\times) = \Pic (\X_{\R'})$, we get the exact sequence 
\begin{eqn}
	\Pic (\X_{\R'}) \too \Pic (\X_\R) \too \H^2 (\X_0, \cal{O}_{\X_0}) \otimes \I
\fullstop
\end{eqn}

Therefore, the line bundle $\cal{M}_\R$ on $\X_\R$ can always be lifted to a line bundle $\cal{M}_{\R'}$ on $\X_{\R'}$ because the obstruction to this lifting lies in $\H^2 (\X_0, \cal{O}_{\X_0}) \otimes_\CC \I$ which vanishes since $\X_0$ is a curve.
After choosing a lift $\cal{M}_{\R'}$, the obstruction to lifting the section $\sigma_\R$ lies in $\H^1 (\X_0, \cal{M}^{\otimes 2}) \otimes_\CC \I$ which also vanishes by assumption.
Thus, the spectral algebra $\cal{A}_\R$ and its finite flat cover $\Y_\R$ admit lifts over the prescribed target $\X_{\R'}$.
The spectral line bundle $\cal{L}_\R$ over $\Y_\R$ also lifts to a line bundle $\cal{L}_{\R'}$ over $\Y_{\R'}$ because its obstruction lies in $\H^2 (\Y_0, \cal{O}_{\Y_0}) \otimes_\CC \I = 0$.
Put $\cal{E}_{\R'} \coleq (\pi_{\R'})_\ast \cal{L}_{\R'}$, and let $\rho_{\R'}$ be the multiplication action of $\cal{A}_{\R'}$ on $\cal{E}_{\R'}$.
Then $(\cal{E}_{\R'}, \cal{A}_{\R'}, \rho_{\R'})$ lifts the given spectral structure $(\cal{E}_\R, \cal{A}_\R, \rho_\R)$.
This proves smoothness of $\frak{SpStr}_2 (\X / \S)$ over $\S$ at the chosen point $(\cal{E}_0, \cal{A}_0, \rho_0)$.

For $n \geq 3$, take the disjoint union of this double cover with $n-2$ copies of the identity cover $\X_0 \to \X_0$, and take the structure sheaf as the spectral line bundle.
Writing $\cal{B}_0$ for the double-cover algebra, the resulting spectral algebra is
\begin{eqn}
	\cal{A}_0
	=
	\cal{B}_0 \times \prod_{i=1}^{n-2} \cal{O}_{\X_0}
\fullstop
\end{eqn}
The projections onto these factors are multiplication by mutually orthogonal idempotents whose sum is $1$.

Now consider any infinitesimal deformation $\cal{A}_\R$ of this algebra.
Idempotents lift uniquely across nilpotent thickenings, and their
local lifts glue by uniqueness.
Thus, the displayed decomposition lifts to a decomposition of
$\cal{A}_\R$ into one rank-two factor $\cal{B}_\R$ and $n-2$
rank-one factors.
Each rank-one factor is canonically $\cal{O}_{\X_\R}$ via its unit map.
Consequently, the deformed cover remains a disjoint union of a
double cover and $n-2$ identity covers, and its spectral line bundle
restricts to a line bundle on each piece.

Over a prescribed target $\X_{\R'}$, the double-cover piece and its
spectral line bundle lift by the preceding argument.
On each identity-cover piece, only a line bundle needs to be lifted;
its obstruction lies in
$\H^2 (\X_0, \cal{O}_{\X_0}) \otimes_\CC \I = 0$.
Taking the disjoint union of these lifts solves the given lifting problem.
The original spectral cover contains the nodal double cover, so this
is a smooth point of the moduli stack lying outside
$\frak{SpStr}_n^\rm{sm} (\X / \S)$.
Together with the rank-two case, this proves the claim for every $n \geq 2$.
\end{example}

\begin{example}[unobstructed singular spectral structure over $\PP^1$]
\label{1792517279249}
This example makes the rank-two construction of \autoref{1793011572262} concrete.
Take $\S = \ast$ and $\X = \PP^1$, with homogeneous coordinates $[x_0 : x_1]$, and consider $\cal{M} = \cal{O}_{\PP^1} (2)$.
Take the global section
\begin{eqn}
	\sigma \coleq x_1^2 (x_1^2 - x_0^2) \in \H^0 \big(\X, \cal{O}_{\PP^1} (4) \big)
\qtext{and define}
	\cal{A}_0 \coleq \cal{O}_{\PP^1} \oplus \cal{O}_{\PP^1} (-2)
\end{eqn}
with multiplication determined by $\sigma : \cal{O}_{\PP^1} (-4) \to \cal{O}_{\PP^1}$, as explained in detail in \autoref{1791437029550}.
Equip $\cal{E}_0 \coleq \cal{A}_0$ with its canonical multiplication action $\rho_0$.
Let $\pi_0 : \Y_0 \to \X$ be the spectral realisation of $\cal{A}_0$.
Then $\cal{L}_0 \coleq \cal{O}_{\Y_0}$ is the spectral line bundle.

On the affine chart $x_0 \neq 0$, with coordinate $z = x_1 / x_0$, this double cover has equation
\begin{eqn}
	w^2 = z^2 (z^2 - 1)
\fullstop
\end{eqn}
The point $(z,w) = (0,0)$ is an ordinary node of $\Y_0$, but the remaining zeros of $\sigma$ are simple.
Also, $\sigma$ does not vanish at infinity, so there are no other singularities.
Moreover, $z^2 (z^2 - 1)$ is not a square in the field of meromorphic functions $\CC (z)$, so $\Y_0$ is reduced and irreducible.

Nevertheless, $\frak{SpStr}_2 (\X)$ is smooth at the point $(\cal{E}_0, \cal{A}_0, \rho_0)$.
To verify this, let $\R' \to \R$ be a small extension of local Artinian $\CC$-algebras, with kernel $\I$ annihilated by the maximal ideal of $\R'$.
Suppose that a deformation of this spectral structure over $\X_\R \coleq \PP^1 \times \Spec (\R)$ is given, with spectral algebra $\cal{A}_\R$ and spectral line bundle $\cal{L}_\R$.
Then the trace map gives a canonical decomposition
\begin{eqn}
	\cal{A}_\R
	=
	\cal{O}_{\X_\R} \oplus \cal{N}_\R
\qqtext{where}
	\cal{N}_\R \coleq \cal{ker} (\tr)
\fullstop
\end{eqn}
Here, $\cal{M}_\R \coleq \cal{N}^{-1}_\R$ is a line bundle deforming $\cal{M} = \cal{O}_{\PP^1} (2)$.
Since $\H^1 (\PP^1, \cal{O}_{\PP^1}) = 0$, every such infinitesimal deformation $\cal{M}_\R$ is isomorphic to $\cal{O}_{\X_\R} (2)$.
Choose an isomorphism compatible with the given identification on the closed fibre.

By the Cayley-Hamilton Theorem, the square of every trace-zero element of a rank-two algebra is scalar.
Thus, multiplication on $\cal{N}_\R$ is determined by a section
\begin{eqn}
	\sigma_\R
	\in
	\H^0 \big( \X_\R, \cal{O}_{\X_\R} (4) \big)
	\cong
	\H^0 \big( \PP^1, \cal{O}_{\PP^1} (4) \big) \otimes_\CC \R
\end{eqn}
which reduces to $\sigma$.
Lift all five coefficients of the given section $\sigma_\R$ from $\R$ to $\R'$ to obtain a section $\sigma_{\R'}$.
The algebra
\begin{eqn}
	\cal{A}_{\R'}
	\coleq
	\cal{O}_{\X_{\R'}} \oplus \cal{O}_{\X_{\R'}} (-2)
\end{eqn}
with multiplication determined by $\sigma_{\R'}$ then lifts $\cal{A}_\R$.
Its relative spectrum gives a finite flat cover $\pi_{\R'} : \Y_{\R'} \to \X_{\R'}$ lifting the given cover over $\R$.

The exact sequence of sheaves of units for $\Y_\R \subset \Y_{\R'}$ places the obstruction to lifting the spectral line bundle $\cal{L}_\R$ in $\H^2 (\Y_0, \cal{O}_{\Y_0}) \otimes_\CC \I$.
This cohomology group vanishes because $\Y_0$ is a compact albeit singular curve.
Consequently, $\cal{L}_\R$ lifts to a line bundle $\cal{L}_{\R'}$ on $\Y_{\R'}$.
Its pushforward $\cal{E}_{\R'} \coleq (\pi_{\R'})_\ast \cal{L}_{\R'}$, equipped with the multiplication action $\rho_{\R'}$ of $\cal{A}_{\R'}$, lifts the given spectral structure $(\cal{E}_\R, \cal{A}_\R, \rho_\R)$.

Thus, every infinitesimal lifting problem at $(\cal{E}_0, \cal{A}_0, \rho_0)$ has a solution.
Since $\frak{SpStr}_2 (\X)$ is locally of finite presentation, it is smooth at this point, despite the fact that the corresponding spectral curve is singular.
\end{example}

\begin{example}[distinguished open substacks need not be dense]
\label{1793123354876}
This example demonstrates that in general none of the distinguished open substack inclusions in \eqref{1793009278980} of \autoref{1792494388516} is dense if $n \geq 2$.
Take $\S = \ast$ and fix any compact Riemann surface $\X$; for the substack of clean spectral structures, we will additionally assume that its genus is $g_\X \geq 2$.

First, note that the degree of any given spectral algebra $\cal{A}$ is locally constant in families.
Moreover, if $\cal{A}$ is smooth, then its degree is necessarily nonpositive.
Indeed, if the spectral realisation $\pi : \Y \to \X$ is smooth, then Riemann-Hurwitz gives 
\begin{eqn}
	\deg (\cal{A}) = \deg (\pi_\ast \cal{O}_\Y) = - r/2 \leq 0
\end{eqn}
where $r = \deg \R$ is the total ramification degree of $\pi$.
This identity remains true even if $\Y$ is disconnected by applying Riemann-Hurwitz to each connected component.
This also follows directly from \autoref{1790716932644}.

On the other hand, choose any vector bundle $\cal{V}$ on $\X$ of rank $n-1$ and positive degree.
Take 
\begin{eqn}
	\cal{A} \coleq \cal{O}_\X \oplus \cal{V}
\qqtext{with algebra structure}
	(f,v) \cdot (g,u) \coleq (fg, fu + gv)
\fullstop
\end{eqn}
Notice that this satisfies $(f,v) \cdot (f,-v) = (f^2, 0)$.
Taking $\cal{E} = \cal{A}$ with this multiplication action defines a spectral structure $(\cal{E}, \cal{A}, \rho)$ over $\X$ with $\deg (\cal{A}) > 0$.
Thus, the locus in $\frak{SpStr}_n (\X)$ of spectral structures of strictly positive degree is a nonempty open-and-closed substack disjoint from $\frak{SpStr}_n^\rm{sm} (\X)$.
Consequently, $\frak{SpStr}_n^\rm{sm} (\X)$ is not dense in $\frak{SpStr}_n (\X)$.

The substack of smooth and connected spectral structures $\frak{SpStr}_n^{\rm{sm}, \rm{con}} (\X)$ is not dense in $\frak{SpStr}_n^\rm{sm} (\X)$.
Indeed, by \autoref{1792494388516}, the locus of smooth spectral structures whose spectral curves are disconnected is an open-and-closed substack of $\frak{SpStr}_n^\rm{sm} (\X)$.
This locus is nonempty: the trivial degree-$n$ cover $\pi : \Y \coleq \X \sqcup \cdots \sqcup \X \to \X$, equipped with the spectral line bundle $\cal{O}_\Y$, gives a smooth spectral structure on $\cal{E} = \cal{O}_\X^{\oplus n}$ with disconnected spectral curve.
Thus, $\frak{SpStr}_n^\rm{sm} (\X)$ contains a nonempty open substack which is disjoint from $\frak{SpStr}_n^{\rm{sm},\rm{con}} (\X)$.

Finally, assume that $g_\X \geq 2$.
The substack $\frak{SpStr}_n^{\rm{clean}} (\X)$ of clean spectral structures is in general not dense in $\frak{SpStr}_n^{\rm{sm},\rm{con}} (\X)$ because the locus of smooth and connected spectral structures of degree $0$ is a nonempty open substack in its complement.
Indeed, by \autoref{1790716932644}, for any smooth spectral structure, the off-diagonal bundle $\cal{Q}$ has Euler characteristic 
\begin{eqn}
	\chi(\cal{Q}) = - \deg (\cal{A}) - n(n-1)(g_\X - 1)
\fullstop
\end{eqn}
If this spectral structure is clean so that $\H^1 (\X, \cal{Q}) = 0$, then $\chi (\cal{Q}) = h^0 (\X, \cal{Q}) \geq 0$; i.e., 
\begin{eqn}
	\deg (\cal{A}) \leq - n (n-1) (g_\X -1)
\fullstop
\end{eqn}
Now, take a connected étale cover $\pi : \Y \to \X$ of degree $n$.
Such covers exist for every $n \geq 2$, and every spectral structure realised by a line bundle on such a cover satisfies 
\begin{eqn}
	h^1 (\X, \cal{Q}) = h^0 (\X, \cal{Q}) + n(n-1)(g_\X - 1) > 0
\fullstop
\end{eqn}
The condition $\deg (\cal{A}) = 0$ is open-and-closed within the locus of smooth covers $\pi$ because $r = -2 \deg(\cal{A})$.
It therefore defines a nonempty open substack of $\frak{SpStr}_n^{\rm{sm},\rm{con}} (\X)$ containing no clean spectral structures.
\end{example}

There is a good reason why, in \autoref{1792494388516}, we speak of `smooth and connected' spectral structures but not of `connected' spectral structures.
Connected spectral structures do form a full substack $\frak{SpStr}_n^\rm{con} (\X/\S)$, but connectedness alone need not define an open substack or even an analytic substack.
The following \autoref{1793154664668} exhibits both failures and explains why a full substack need not admit a smooth analytic atlas.

\begin{example}[connected substack need not be analytic]
\label{1793154664668}
Let $\X$ be a fixed compact Riemann surface and take $\T = \CC$ with coordinate $t$.
On $\X \times \T$, consider
\begin{eqn}
	\cal{A} \coleq \cal{O}_{\X \times \T} [z] / \inner{z^2 - t}
\qqtext{and}
	\cal{E} \coleq \cal{A}
\fullstop{,}
\end{eqn}
with the canonical multiplication action.
For $t \neq 0$, the polynomial $z^2 - t$ has two distinct roots, so the spectral realisation $\Y_t$ is the disjoint union of two copies of $\X$.
For $t=0$, the spectral realisation is the order-$1$ fat curve over $\X$:
\begin{eqn}
	\Y_0 = \Spec_\X \big( \cal{O}_\X [z] / \inner{z^2} \big)
\fullstop
\end{eqn}
Its underlying topological space is $\X$ itself, so it is connected.
Thus, the connected locus in this family $\cal{A}$ of spectral structures is the single point $t = 0$, which of course is not open in $\T$.

There is a further issue involving nilpotent functions on the parameter space.
For a holomorphic map $\U \to \T$, let $u \in \H^0 (\U, \cal{O}_\U)$ be the pullback of $t$.
Then the pulled-back family $\cal{A}_\U$ is connected precisely when $u$ vanishes at every point of $\U$; equivalently, $u$ is locally nilpotent.
Consequently, the fibre product
\begin{eqn}
	\T \times_{\frak{SpStr}_2 (\X)} \frak{SpStr}_2^\rm{con} (\X)
\end{eqn}
is the functor of locally nilpotent functions, or the formal completion of the line $\T$ at $0$, viewed on ordinary analytic parameter spaces.
In particular, it contains the maps $t \mapsto \epsilon$ from all fat points 
\begin{eqn}
	\DD_k \coleq \Spec \big( \CC [\epsilon] / \inner{\epsilon^k} \big)
	\qqtext{for}
	k \geq 1
\fullstop
\end{eqn}
No fixed nilpotence bound can capture all of these maps; i.e., these maps cannot all factor through any single closed subspace of $\T$ defined by $t^m = 0$ because, for $k > m$, the image $\epsilon$ of $t$ satisfies $\epsilon^m \neq 0$ in $\CC [\epsilon] / \inner{\epsilon^k}$.

To see directly why this functor has no smooth analytic atlas, suppose towards a contradiction that there is one, given by an analytic space $\V$.
Choose a point $v_0 \in \V$ above the object defined by $t = 0$.
Let $t_\V \in \H^0 (\V, \cal{O}_\V)$ denote the pullback of $t$.
Since $\V$ maps to this functor, the section $t_\V$ is locally nilpotent.
Consequently, after restricting to a neighbourhood $\V_0$ of $v_0$, there is an integer $m \geq 1$ such that 
\vspace{-5pt}
\enlargethispage{20pt}
\begin{eqn}
	(t_\V|_{\V_0})^m = 0
\qtext{in}
	\H^0 (\V_0, \cal{O}_{\V_0})
\fullstop
\end{eqn}

For every $k \geq 1$, the morphism $\DD_k \to \T$ defined by $t \mapsto \epsilon$ belongs to this functor.
Its restriction to the closed point is the object defined by $t = 0$, for which we have chosen the lift $v_0$.
The closed immersions $\DD_1 \inj \DD_2 \inj \cdots \inj \DD_k$ are successive small extensions.
Formal smoothness of the atlas therefore allows us to lift the morphism $\DD_k \to \T$ successively to a morphism $\ell_k : \DD_k \to \V$ whose closed point maps to $v_0$.
Since $\DD_k$ has only one underlying point, this lift factors through $\V_0$.
Pulling back the displayed equality of sections along $\ell_k$ gives
\begin{eqn}
	\epsilon^m = 0
	\quad\text{in}\quad
	\CC [\epsilon] / \inner{\epsilon^k}
\fullstop
\end{eqn}
For $k > m$, this is impossible.
Therefore, this functor has no smooth analytic atlas, and the full substack of connected spectral structures in this example is not analytic.
\end{example}

\subsection{The Universal Relative Picard Stack}
\label{1793136253999}

Now, let us consider the stack of covers $\frak{Cov}_n (\X / \S)$ over the analytic site $\mathbf{An}_\S$.
Its objects over any test space $\T \to \S$ are finite flat covers $\pi : \Y \to \X_\T$ of degree $n$.
Isomorphisms are isomorphisms of covers over $\X_\T$, and the pullback is defined by base change of the cover.
Thus, any family of covers $\pi : \Y \to \X_\T$ is classified by a morphism $c_\pi : \T \to \frak{Cov}_n (\X / \S)$.
In summary:
\begin{eqn}
	\frak{Cov}_n (\X / \S)
	\quad:\quad
	\set[\Big]{ \text{finite flat covers $\pi : \Y \to \X_\T$ of degree $n$} }_{\T \in \mathbf{An}_\S}
\fullstop
\end{eqn}
Meanwhile, thanks to the equivalence between spectral structures and their realisations (\autoref{1792502248471}), the two forgetful functors in \eqref{1790940189747} globalise to morphisms of stacks:
\begin{eqntag}
\label{1792401198207}
\begin{tikzcd}
&	\frak{SpStr}_n (\X / \S)
		\ar[dl, "G"']
		\ar[dr, "F"]
\\
	\frak{Cov}_n (\X / \S)
&&	\frak{Bun}_n (\X / \S)
\fullstop
\end{tikzcd}
\end{eqntag}

\begin{proposition}[smoothness and surjectivity of $G$]
\label{1792511507140}
The stack $\frak{Cov}_n (\X / \S)$ is an analytic stack locally of finite presentation over $\S$, and the forgetful morphism $G$ is smooth and surjective.
Arbitrary analytic base change commutes with the formation of this stack and the forgetful morphism $G$: for every holomorphic map $\S' \to \S$, there is a canonical equivalence
\begin{eqn}
	\frak{Cov}_n (\X / \S) \times_\S \S'
	\cong
	\frak{Cov}_n (\X \times_\S \S' / \S')
\fullstop
\end{eqn}
\end{proposition}

\begin{proof}
Finite flat covers of degree $n$ correspond to commutative unital algebras which are locally free of rank $n$.
Effective descent for vector bundles and their morphisms therefore gives effective descent for these covers.
This proves that $\frak{Cov}_n (\X / \S)$ is a stack over the analytic site.

To see that it is analytic, consider the forgetful morphism
\begin{eqn}
	\frak{Cov}_n (\X / \S)
	\too
	\frak{Bun}_n (\X / \S)
\qqtext{sending}
	\big( \pi : \Y \to \X_\T \big) \mapstoo \pi_\ast \cal{O}_\Y
\fullstop
\end{eqn}
For any test space $\T \to \S$, fix a vector bundle $\cal{A}$ of rank $n$ on $\X_\T$.
The base change of this forgetful morphism along $\T \to \frak{Bun}_n (\X / \S)$ parametrises all commutative unital algebra structures on $\cal{A}$.
Indeed, the identification of the underlying vector bundle with $\cal{A}$ transports the algebra structure to $\cal{A}$, and any isomorphism preserving this identification is the identity.
Such structures consist of two morphisms
\begin{eqn}
	e : \cal{O}_{\X_\T} \too \cal{A}
\qqtext{and}
	\mu : \cal{A} \otimes \cal{A} \too \cal{A}
\end{eqn}
satisfying the unit, commutativity, and associativity identities.
Specifically, introduce the vector bundle
\begin{eqn}
    \cal{H}
    \coleq
    \cal{A} \oplus \cal{Hom} (\cal{A} \otimes \cal{A}, \cal{A})
\fullstop
\end{eqn}
Then a pair $(e, \mu)$ is a global section of this vector bundle $\cal{H}$, or equivalently a global section of its total space $u : \H \coleq \tot (\cal{H}) \to \X_\T$.
Suppressing the pullbacks to $\H$ from the notation, let $(e,\mu)$ denote the tautological section of $u^\ast \cal{H}$.
Let $\Z \subset \tot (\cal{H})$ be the closed analytic subspace defined by the three identities
\begin{eqn}
	\mu \circ (e \otimes \id)
	=
	\id
	=
	\mu \circ (\id \otimes e)
\fullstop{,}\qquad
	\mu \circ \varepsilon = \mu
\fullstop{,}\qquad
	\mu \circ (\mu \otimes \id)
	=
	\mu \circ (\id \otimes \mu)
\fullstop{,}
\end{eqn}
where $\id : \cal{A} \to \cal{A}$ is the identity endomorphism and $\varepsilon : \cal{A} \otimes \cal{A} \to \cal{A} \otimes \cal{A}$ interchanges the two tensor factors.
These are finitely many local holomorphic equations in the entries of $e$ and $\mu$.

For every analytic space $\T' \to \T$, sections of $\Z_{\T'} \to \X_{\T'}$ correspond naturally to commutative unital algebra structures on $\cal{A}_{\T'}$.
Since $\Z \to \X_\T$ is separated, \autoref{1792526939305} represents this functor by an analytic space $\M$ locally of finite presentation over $\T$.
Consequently, the forgetful morphism to $\frak{Bun}_n (\X / \S)$ is representable by analytic spaces and locally of finite presentation.
Pulling back an analytic atlas for $\frak{Bun}_n (\X / \S)$ gives an analytic atlas for $\frak{Cov}_n (\X / \S)$, and its diagonal is representable for the same reason.
This proves both analyticity and local finite presentation of $\frak{Cov}_n (\X / \S)$ over $\S$.

Now we demonstrate that $G$ is smooth and surjective.
The morphism $G$ is locally of finite presentation because its source and target are analytic stacks locally of finite presentation over $\S$, thanks to \autoref{1792494933366} and \autoref{1792502248471}.
We verify its formal smoothness by lifting line bundles over prescribed families of covers.

Let $\R' \to \R$ be a small extension of local Artinian $\CC$-algebras, with kernel $\I$ annihilated by the maximal ideal of $\R'$.
Put $\T \coleq \Spec (\R)$ and $\T' \coleq \Spec (\R')$.
By the preceding equivalence, a lifting problem for $G$ consists of a prescribed family of covers $\pi_{\T'} : \Y_{\T'} \to \X_{\T'}$ together with a line bundle $\cal{L}_\T$ on its restriction $\Y_\T \coleq \Y_{\T'} \times_{\T'} \T$.
We must find a line bundle $\cal{L}_{\T'}$ on $\Y_{\T'}$ together with an isomorphism $\cal{L}_{\T'}|_{\Y_\T} \iso \cal{L}_\T$.

The composite $\Y_{\T'} \to \T'$ is proper and flat.
Let $\cal{J}$ be the defining ideal of $\Y_\T$ in $\Y_{\T'}$.
Flatness gives
\begin{eqn}
	\cal{J}
	\cong
	\cal{O}_{\Y_{\T'}} \otimes_{\R'} \I
	\cong
	\cal{O}_{\Y_0} \otimes_\CC \I
\fullstop{,}
\end{eqn}
where $\Y_0$ is the closed fibre; the second identification holds because the maximal ideal of $\R'$ annihilates $\I$.
Moreover, $\cal{J}^2 = 0$, so there is an exact sequence of sheaves of abelian groups
\begin{eqn}
	0
	\too \cal{J}
	\too \cal{O}_{\Y_{\T'}}^\times
	\too \cal{O}_{\Y_\T}^\times
	\too 1
\fullstop{,}
\end{eqn}
whose first map sends $a$ to $1+a$.
Taking cohomology gives the exact sequence
\begin{eqn}
	\Pic (\Y_{\T'})
	\too
	\Pic (\Y_\T)
	\too
	\H^2 (\Y_0, \cal{O}_{\Y_0}) \otimes_\CC \I
\fullstop
\end{eqn}
The last group vanishes because $\Y_0$ is a compact curve.
Consequently, $\cal{L}_\T$ lifts to a line bundle on $\Y_{\T'}$ with the required identification on $\Y_\T$.
This solves the lifting problem, proving that $G$ is formally smooth and hence smooth.

Every family of covers $\pi : \Y \to \X_\T$ admits the line bundle $\cal{O}_\Y$.
The resulting spectral realisation maps to $\pi$ under $G$, so $G$ is surjective.

The base-change equivalence follows directly from the definition: for every analytic space $\T \to \S'$, both sides classify finite flat degree-$n$ covers of $\X \times_\S \T$.
Under the equivalence of \autoref{1792502248471}, the morphism $G$ forgets the spectral line bundle, so its compatibility with base change follows as well.
\end{proof}

The stack of covers $\frak{Cov}_n (\X / \S)$ carries a universal family of covers 
\begin{eqn}
	q : \frak{Y}_n
	\too
	\frak{Cov}_n (\X / \S)
\end{eqn}
defined as follows.
Consider the stack $\frak{Y}_n$ over $\mathbf{An}_\S$ whose objects over $\T \to \S$ are finite flat degree-$n$ covers $\pi : \Y \to \X_\T$ together with a global section $y : \T \to \Y$ of the composite projection $\Y \to \X_\T \to \T$.
Isomorphisms are isomorphisms of covers $f : \Y \to \Y'$ intertwining the specified sections in the sense that $f \circ y = y'$.
In summary:
\begin{eqn}
	\frak{Y}_n
	\quad:\quad
	\set{ 
	(\Y \tinyoverset{\pi}{\to} \X_\T, y) 
	~\left|~
	\substack{
	\displaystyle
	\text{finite flat cover $\pi : \Y \to \X_\T$ of degree $n$}
	\\
	\displaystyle
	\text{global section 
	$\begin{tikzcd}[ampersand replacement=\&, column sep = small]
		\Y
		\ar[r,"\pi"]
	\&	\X_\T
		\ar[r]
	\&	\T
		\ar[ll, bend left, "y"]
	\end{tikzcd}$
	}
	}
	\right.
	}_{\T \in \mathbf{An}_\S}
\fullstop
\end{eqn}
Forgetting the section $y : \T \to \Y$ defines a forgetful morphism $q : \frak{Y}_n \to \frak{Cov}_n (\X/\S)$.
This morphism is representable by proper flat families of curves.
Indeed, fix $c_\pi : \T \to \frak{Cov}_n (\X / \S)$ classifying a family $\pi : \Y \to \X_\T$.
For any test space $\U \to \T$, an object over $\U$ of the fibre product $\T \times_{\frak{Cov}_n (\X / \S)} \frak{Y}_n$ consists of a cover $\Z \to \X_\U$, a section $z : \U \to \Z$, and a specified isomorphism of covers $f : \Z \iso \Y_\U$.
Transporting the section along $f$ gives a section $f \circ z : \U \to \Y_\U$.
Such a section is equivalently a morphism $\U \to \Y$ over $\T$.
Conversely, every such morphism gives an object of the fibre product by taking the cover $\Y_\U$ with its corresponding section and the identity isomorphism.

These constructions are inverse up to the specified isomorphisms.
Compatibility with $f$ leaves no nonidentity automorphisms, and the correspondence is natural in $\U$.
Thus, the fibre product is represented by $\Y$, giving us the following canonical identification over $\T$:
\begin{eqntag}
\label{1793034624749}
	c_\pi^\ast \frak{Y}_n
	=
	\T \times_{\frak{Cov}_n (\X / \S)} \frak{Y}_n
	\cong
	\Y
\fullstop
\end{eqntag}
Since $\Y \to \T$ is proper and flat, this also proves that $q$ is representable by proper flat families of curves.

This universal family of covers $q : \frak{Y}_n \to \frak{Cov}_n (\X/\S)$ even has its own relative Picard stack, denoted by $\frak{Pic} \big( \frak{Y}_n / \frak{Cov}_n (\X / \S)	\big)$.
It parametrises morphisms $\T \to \frak{Cov}_n (\X / \S)$ together with a line bundle on the corresponding family $\Y$ obtained by pullback of $\frak{Y}_n$ as in \eqref{1793034624749}:
\begin{eqn}
	\frak{Pic} \big( \frak{Y}_n / \frak{Cov}_n (\X / \S)	\big)
\quad:\quad
	\set{ (c_\pi, \cal{L})
		~\left|~
	\substack{
	\displaystyle
	c_\pi : \T \to \frak{Cov}_n (\X / \S)
	\\
	\displaystyle
	\text{$\cal{L} \in \mathbf{Pic} (\Y)$ where $\Y \cong c_\pi^\ast \frak{Y}_n$} 
	}
	\right.
	}_{\T \in \mathbf{An}_\S}
\fullstop
\end{eqn}
Equivalently, its objects over any test space $\T \to \S$ are pairs $(\pi : \Y \to \X_\T, \cal{L})$, where $\pi$ is finite flat of degree $n$ and $\cal{L}$ is a line bundle on $\Y$.
Isomorphisms are isomorphisms of covers together with compatible isomorphisms of the line bundles.
But this is precisely the stack of spectral realisations $\frak{SpReal}_n (\X / \S)$.

\begin{theorem}[description as relative Picard stack]
\label{1792502946092}
The forgetful morphism $G$ identifies $\frak{SpStr}_n (\X / \S)$ with the relative Picard stack of the universal covering curve over $\frak{Cov}_n (\X / \S)$:
\begin{eqn}
	\frak{SpStr}_n (\X / \S)
	\cong
	\frak{Pic} \big( \frak{Y}_n / \frak{Cov}_n (\X / \S)	\big)
\fullstop
\end{eqn}
That is to say, for every family of covers $\pi : \Y \to \X_\T$ classified by a morphism $\T \to \frak{Cov}_n (\X / \S)$, there is a canonical equivalence of stacks over $\T$:
\begin{eqntag}
\label{1793045330854}
	\T \times_{\frak{Cov}_n (\X / \S)}
	\frak{SpStr}_n (\X / \S)
	\cong
	\frak{Pic} (\Y / \T)
\fullstop
\end{eqntag}
Both equivalences are compatible with arbitrary analytic base change.
\end{theorem}

\enlargethispage{20pt}

\begin{proof}
Fix a family $\pi : \Y \to \X_\T$ and a morphism $\U \to \T$.
By \autoref{1792502248471}, an object over $\U$ of the fibre product in \eqref{1793045330854} can be described as a spectral realisation $(\Z \to \X_\U, \cal{L})$ together with an isomorphism of covers $f : \Z \iso \Y_\U$.
Sending this object to $f_\ast \cal{L}$ gives a line bundle on $\Y_\U$.
Conversely, a line bundle $\cal{M}$ on $\Y_\U$ gives the object $( \Y_\U \to \X_\U, \cal{M}, \id_{\Y_\U})$ of the fibre product.

These constructions are quasi-inverse.
Indeed, the specified isomorphism $f$ identifies the varying cover with $\Y_\U$, and the compatible morphisms in the fibre product become precisely isomorphisms of line bundles on $\Y_\U$.
The constructions commute canonically with pullback in $\U$, so they give an equivalence of stacks over $\T$.
Compatibility with arbitrary analytic base change in $\S$ and $\T$ follows from the same constructions and \autoref{1791551111381}.
\end{proof}

\begin{corollary}[relative dimension of $G$]
\label{1792512462793}
The relative stack dimension of the forgetful morphism $G$ at a point represented by the spectral structure with realisation $(\pi : \Y \to \X_s, \cal{L})$ is
\begin{eqn}
	h^1 (\Y, \cal{O}_\Y)
	-
	h^0 (\Y, \cal{O}_\Y)
	=
	-\chi (\cal{O}_\Y)
\fullstop
\end{eqn}
\end{corollary}

We stress that no smoothness, reducedness, or connectedness assumption on $\Y$ is required in the relative stack dimension formula.

\begin{proof}
On a fixed cover $\Y$, first-order deformations of a line bundle are parametrised by $\H^1 (\Y, \cal{O}_\Y)$ and its infinitesimal automorphisms are parametrised by $\H^0 (\Y, \cal{O}_\Y)$; see \cite[\href{https://stacks.math.columbia.edu/tag/08VR}{Tag 08VR}]{stacks-project}.
Since the relative Picard stack is smooth, the difference between the dimensions of these cohomology groups gives the asserted relative stack dimension.
\end{proof}

\begin{remark}[scalar automorphisms]
\label{1792517373880}
Suppose that the spectral curve $\Y$ is smooth and connected.
Then by \autoref{1792512462793}, the relative Picard stack contributes $-\chi(\cal{O}_\Y) = g_\Y - 1$ to the stack dimension.
The subtraction by $1$ is a record of the scalar automorphism group $\CC^\times$ of every line bundle on $\Y$.

More generally, an automorphism of a spectral realisation $(\pi : \Y \to \X, \cal{L})$ induces a deck transformation of the cover.
The kernel of the resulting homomorphism $\Aut (\pi, \cal{L}) \to \Aut_\X (\Y)$ is the group of automorphisms of $\cal{L}$; namely, $\H^0 (\Y, \cal{O}_\Y)^\times$.
For a smooth connected $\Y$, this kernel is $\CC^\times$, and the quotient is a subgroup of the finite group $\Aut_\X (\Y)$.

If $\Y$ is smooth with $N$ connected components $\Y_i$, the kernel is $(\CC^\times)^N$ and the relative stack dimension is therefore (cf. the genus convention in \eqref{1792494116976}):
\begin{eqn}
	\sum_{i=1}^N g_{\Y_i} - N
	=
	-\chi (\cal{O}_\Y)
\fullstop
\qedhere
\end{eqn}
\end{remark}

\enlargethispage{20pt}
\subsection{Comparison with the Stack of Higgs Bundles}
\label{1793045910309}

Fix a line bundle $\Lambda$ on $\X$.
The relationship between regular Higgs bundles and spectral structures was described in \autoref{1792167464001}.
Recall that the additional data supplied by a regular Higgs field is a global $\Lambda$-twisted generator of the spectral algebra, or equivalently an embedding of the spectral cover into $\tot (\Lambda)$.
In \ref{1793045933301}, we briefly explained why this description extends to families.
It therefore also gives a comparison between their moduli stacks.

Consider the analytic moduli stack of regular $\Lambda$-twisted Higgs bundles over $\X$, denoted by
\begin{eqn}
	\frak{Higgs}_{n,\Lambda}^\rm{reg} (\X / \S)
\quad:\quad
	\set{ (\cal{E}, \Phi) ~\Big|~ \text{regular $\Lambda_\T$-twisted Higgs bundles of rank-$n$ over $\X_\T$} }_{\T \in \mathbf{An}_\S}
\fullstop
\end{eqn}
Concretely, its objects over $\T \to \S$ are vector bundles $\cal{E}$ of rank $n$ over $\X_\T$ with a regular $\Lambda_\T$-twisted Higgs field $\Phi : \cal{E} \to \cal{E} \otimes \Lambda_\T$, where $\Lambda_\T$ denotes the pullback of $\Lambda$ to $\X_\T$.
Isomorphisms are vector bundle isomorphisms intertwining the Higgs fields, and pullback is defined by base change.
Spectral correspondence in families, as recalled in \ref{1793045933301}, yields a forgetful morphism
\begin{eqn}
	H :
	\frak{Higgs}_{n,\Lambda}^\rm{reg} (\X / \S)
	\too
	\frak{SpStr}_n (\X / \S)
\fullstop{,}
\end{eqn}
which fits into the following commutative diagram:
\begin{eqntag}
\label{1793045953932}
\begin{tikzcd}
&
&	\frak{Bun}_n (\X / \S)
\\
	\frak{Higgs}_{n,\Lambda}^\rm{reg} (\X / \S)
		\ar[r, "H"]
		\ar[rru, bend left]
		\ar[rrd, bend right]
&	\frak{SpStr}_n (\X / \S)
		\ar[ur, "F"]
		\ar[dr, "G"']
\\
&
&	\frak{Cov}_n (\X / \S)
\fullstop
\end{tikzcd}
\end{eqntag}

To describe this comparison more precisely, let us introduce the \textit{relative Hitchin base}.
Consider the following vector bundle on $\X$ and its bundle projection:
\begin{eqn}
	\Oplus_{i=1}^n \Lambda^{\otimes i}
\qtext{and}
	u : \tot \left( \Oplus_{i=1}^n \Lambda^{\otimes i} \right) \too \X
\fullstop
\end{eqn}
Using \autoref{1792526939305} with $\T = \S$, the relative functor of global sections of $u$ is represented by an analytic space, locally of finite presentation over $\S$, which we call the \dfn{relative Hitchin base} and denote by 
\begin{eqn}
	\frak{B} \coleq \frak{B}_{n,\Lambda} (\X / \S)
\fullstop
\end{eqn}
Concretely, for every analytic space $\T \to \S$, a morphism $\T \to \frak{B}$ over $\S$ corresponds naturally to a tuple $c = (c_1, \ldots, c_n)$ with each $c_i \in \H^0 \big( \X_\T, \Lambda_\T^{\otimes i} \big)$; in particular, for a fixed curve $\X$, this recovers the ordinary $\Lambda$-valued Hitchin base
\begin{eqn}
	\frak{B}_{n,\Lambda} (\X/\S) (\T)
	=
	\Oplus_{i=1}^n \H^0 \big( \X_\T, \Lambda^{\otimes i}_\T \big)
\qqtext{and}
	\frak{B}_{n,\Lambda} (\X)
	=
	\Oplus_{i=1}^n \H^0 \big( \X, \Lambda^{\otimes i} \big)
\fullstop
\end{eqn}
The characteristic polynomial defines the \dfn{Hitchin morphism}
\begin{eqn}
	h :
	\frak{Higgs}_{n,\Lambda}^\rm{reg} (\X / \S)
	\too
	\frak{B}_{n,\Lambda} (\X / \S)
\fullstop
\end{eqn}

Such tuples $c = (c_1, \ldots, c_n)$ are precisely sections of the base-change projection $u_\T$ as explained in \autoref{1792526939305}.
The universal section supplied by \autoref{1792526939305} gives a universal tuple of coefficients on $\X_\frak{B} = \X \times_\S \frak{B}_{n,\Lambda} (\X/\S)$ from which every tuple $c$ can be obtained by pullback.
The same lemma shows that the relative Hitchin base, together with its universal tuple, commutes with arbitrary base change on $\S$.

A tuple $c$ over $\T$ defines a finite flat degree-$n$ cover $\pi_c : \Y_c \to \X_\T$ by the monic equation
\begin{eqn}
	t^n + c_1 t^{n-1} + \cdots + c_n = 0
\end{eqn}
inside $\tot (\Lambda_\T)$.
Here, $t$ is the tautological section of the pullback of $\Lambda_\T$ to its total space.
Applying this construction to the universal tuple on $\X_\frak{B}$ gives a universal family of spectral covers that we call the \dfn{universal spectral cover}:
\begin{eqn}
	\pi_\Lambda : \Y_\Lambda \too \X_\frak{B}
\fullstop
\end{eqn}
Every tuple $c$ over $\T$ determines a classifying morphism $\T \to \frak{B}$, and the corresponding cover $\pi_c : \Y_c \to \X_\T$ is the pullback of $\pi_\Lambda$.
Forgetting their embeddings into the total space defines a morphism $i : \frak{B}_{n,\Lambda} (\X / \S) \to \frak{Cov}_n (\X / \S)$.

\enlargethispage{20pt}
\begin{theorem}[regular Higgs bundles by base change]
\label{1792546482295}
The spectral correspondence gives a Cartesian square
\begin{eqntag}
\label{1792546191429}
\begin{tikzcd}
	\frak{Higgs}_{n,\Lambda}^\rm{reg} (\X / \S)
		\ar[r, "H"]
		\ar[d, "h"']
		\ar[dr, phantom, very near start, "\ulcorner"]
&
	\frak{SpStr}_n (\X / \S)
		\ar[d, "G"]
\\
	\frak{B}_{n, \Lambda} (\X / \S)
		\ar[r, "i"']
&
	\frak{Cov}_n (\X / \S)
\fullstop
\end{tikzcd}
\end{eqntag}
Consequently, there is a canonical equivalence of stacks over the relative Hitchin base $\frak{B}$ between the stack of regular Higgs bundles and the relative Picard stack of the universal spectral cover $\Y_\Lambda$ over $\X_\frak{B}$:
\begin{eqn}
	\frak{Higgs}_{n,\Lambda}^\rm{reg} (\X / \S)
	\cong
	\frak{Pic} (\Y_\Lambda / \frak{B})
\fullstop
\end{eqn}
In particular, $\frak{Higgs}_{n,\Lambda}^\rm{reg} (\X / \S)$ is an analytic stack locally of finite presentation over $\S$.
The morphism $H$ is representable by analytic spaces and locally of finite presentation, whilst the Hitchin morphism $h$ is smooth and surjective.
All claims are compatible with arbitrary analytic base change in $\S$, with $\Lambda$ replaced by its pullback.
\end{theorem}

\begin{proof}
To deduce the Cartesian square \eqref{1792546191429}, examine directly the fibre product between $\frak{B}$ and $\frak{SpStr}_n (\X / \S)$ over the stack of covers.
An object of $\T$ consists of the following data: a tuple $c = (c_1, \ldots, c_n)$ with $c_i \in \H^0 (\X_\T, \Lambda^{\otimes i}_\T)$, a spectral realisation $(\pi : \Y \to \X_\T, \cal{L})$, and an isomorphism of covers $f : \Y \iso \Y_c$.
The isomorphism $f$ identifies this data with a line bundle $f_\ast \cal{L}$ on the embedded spectral cover $\Y_c$.
Pushing forward this line bundle and letting the tautological section $t$ act by multiplication gives a regular $\Lambda_\T$-twisted Higgs bundle with characteristic polynomial $c$.

Conversely, a regular Higgs bundle with characteristic polynomial $c$ determines a line bundle on $\Y_c$, and hence an object of the displayed fibre product.
The two constructions are quasi-inverse to each other by the family-version of spectral correspondence.
They identify the morphisms on both sides and commute canonically with arbitrary base change.
This proves that \eqref{1792546191429} is Cartesian.

To obtain the relative Picard stack description, apply \autoref{1792502946092} after base change along the morphism $i$.
By \autoref{1792511507140}, the morphism $G$ is smooth and surjective.
Its base change $h$ is therefore also smooth and surjective.
Since $\frak{B}$ is an analytic space locally of finite presentation over $\S$, the Cartesian square also shows that the regular Higgs stack $\frak{Higgs}_{n,\Lambda}^\rm{reg} (\X / \S)$ is analytic and locally of finite presentation over $\S$.

It remains to prove that $H$ is representable by analytic spaces and locally of finite presentation.
Fix a morphism $\T \to \frak{SpStr}_n (\X / \S)$ corresponding to a family of spectral structures $(\cal{E}, \cal{A}, \rho)$ on $\X_\T$.
By \autoref{1792167464001} applied in families, for every analytic space $\U \to \T$, the groupoid over $\U$ of the fibre product
\begin{eqn}
	\T
	\times_{\frak{SpStr}_n (\X / \S)}
	\frak{Higgs}_{n,\Lambda}^\rm{reg} (\X / \S)
\end{eqn}
is canonically equivalent to the discrete groupoid of sections $a \in \H^0 \big( \X_\U, \cal{A}_\U \otimes \Lambda_\U \big)$ whose local representatives generate $\cal{A}_\U$.
Here, the identification with the pulled-back spectral structure $(\cal{E}_\U, \cal{A}_\U, \rho_\U)$ is part of the data in the fibre product, which means in particular that all automorphism groups in this fibre product are trivial.

Now, apply \autoref{1792526939305} to the bundle projection $\tot (\cal{A} \otimes \Lambda_\T) \to \X_\T$.
It represents the functor of all global sections of $\cal{A} \otimes \Lambda_\T$ by an analytic space $\V$, locally of finite presentation over $\T$.
Let $a \in \H^0 \big( \X_\V, \cal{A}_\V \otimes \Lambda_\V \big)$ denote the universal section.
Its powers define the morphism
\begin{eqn}
	\epsilon_a :
	\cal{O}_{\X_\V}
	\oplus \Lambda_\V^{-1}
	\oplus \cdots
	\oplus \Lambda_\V^{-(n-1)}
	\too
	\cal{A}_\V
\fullstop
\end{eqn}
The section $a$ is a twisted generator precisely when $\epsilon_a$ is an isomorphism.

Let $\D \subset \X_\V$ be the degeneracy locus of $\epsilon_a$.
Its image in $\V$ is a closed analytic subspace because $\X_\V \to \V$ is proper.
The complement of this image represents precisely the twisted generators.
It is therefore the required base change of $H$, proving representability by analytic spaces and local finite presentation.
\end{proof}

\enlargethispage{20pt}
\begin{corollary}[fibres and essential image]
\label{1793046246686}
Let $(\cal{E}, \cal{A}, \rho)$ be a spectral structure on a fibre $\X_s$.
There is a canonical isomorphism between the essential fibre of $H$ over $(\cal{E}, \cal{A}, \rho)$ and the analytic space of global $\Lambda_s$-twisted generators of $\cal{A}$:
\begin{eqn}
	H^{[-1]} (\cal{E}, \cal{A}, \rho)
	\cong
	\sfop{Gen}_{\Lambda_s} (\cal{A})
	\subset
	\H^0 (\X_s, \cal{A} \otimes \Lambda_s)
\fullstop{,}
\end{eqn}
which is a Zariski-open subset of the vector space $\H^0 (\X_s, \cal{A} \otimes \Lambda_s)$.
Consequently, every nonempty essential fibre is smooth of dimension
\begin{eqn}
	h^0 (\X_s, \cal{A} \otimes \Lambda_s)
	=
	\sum_{i=0}^{n-1}
	h^0 \big( \X_s, \Lambda_s^{\otimes (1-i)} \big)
\fullstop
\end{eqn}
A spectral structure lies in the essential image of $H$ if and only if its spectral cover admits an embedding over $\X_s$ into the total space $\tot (\Lambda_s)$.
\end{corollary}

\begin{proof}
The description of the essential fibre, the fact that it is Zariski-open, and the embedding criterion are all simply assertions of \autoref{1792167464001} and \autoref{1792167000670} applied to $\X_s$.
If the fibre is nonempty, a twisted generator gives an isomorphism of vector bundles
\begin{eqn}
	\cal{A}
	\cong
	\cal{O}_{\X_s}
	\oplus \Lambda_s^{-1}
	\oplus \cdots
	\oplus \Lambda_s^{-(n-1)}
\fullstop
\end{eqn}
Tensoring with $\Lambda_s$ and taking global sections gives the dimension formula.
\end{proof}

This corollary shows in particular that fixing the twisting line bundle $\Lambda$ imposes restrictions on the possible spectral algebras appearing in the essential image.
Namely, every spectral structure in the essential image satisfies $\det \cal{A} \cong \Lambda_s^{-n(n-1)/2}$.
Furthermore, if its spectral cover is smooth, this also forces $r = n(n-1) \deg \Lambda_s$ by \autoref{1792167000671}.
Therefore, the example in \autoref{1792434634960} shows that $H$ need not be surjective, even onto the substack of clean spectral structures, the smallest of the distinguished substacks appearing in \autoref{1792494388516}.

\begin{example}[$H$ need not be smooth]
\label{1793046252953}
This example demonstrates that the smoothness of the nonempty fibres of $H$ does not imply that $H$ is smooth.
This is because a deformation of the spectral structure may fail to admit a deformation of its chosen twisted generator.
This failure can be seen already for smooth connected double covers as follows.

Let $\X$ be a fixed compact Riemann surface of positive genus.
Choose a sufficiently positive line bundle $\Lambda$ such that $\H^1 \big( \X, \Lambda^{\otimes 2} \big) = 0$, and choose a section $s \in \H^0 \big( \X, \Lambda^{\otimes 2} \big)$ with simple zeros.
The algebra $\cal{A} \coleq \cal{O}_\X \oplus \Lambda^{-1}$ with multiplication determined by $s : \Lambda^{-2} \to \cal{O}_\X$ defines a smooth connected double cover $\pi : \Y \to \X$.
For any line bundle $\cal{L}$ on $\Y$, the embedding of this cover into $\tot (\Lambda)$ gives a regular $\Lambda$-twisted Higgs bundle on $\cal{E} \coleq \pi_\ast \cal{L}$.

Consider the fat point $\T \coleq \Spec ( \CC [\epsilon] / \inner*{\epsilon^2} )$.
Since $\H^1 (\X, \cal{O}_\X) \neq 0$, there is a line bundle $\cal{M}_\T$ on $\X_\T$ that restricts to $\Lambda$ but which is not isomorphic to the constant lift $\Lambda_\T$.
The vanishing of $\H^1 \big( \X, \Lambda^{\otimes 2} \big)$ implies that section $s$ lifts to a section $s_\T \in \H^0 \big( \X_\T, \cal{M}_\T^{\otimes 2} \big)$.
Consequently, the algebra $\cal{A}_\T \coleq \cal{O}_{\X_\T} \oplus \cal{M}_\T^{-1}$ with multiplication determined by $s_\T$ defines a finite flat cover $\pi_\T : \Y_\T \to \X_\T$ deforming $\pi$.
The line bundle $\cal{L}$ lifts to a line bundle $\cal{L}_\T$ on $\Y_\T$, since its obstruction lies in $\H^2 (\Y, \cal{O}_\Y) = 0$.
Its pushforward therefore gives a deformation of the spectral structure.

Suppose, towards a contradiction, that this deformation lifted through $H$ with the fixed twist $\Lambda_\T$.
A global $\Lambda_\T$-twisted generator would give an isomorphism $\cal{O}_{\X_\T} \oplus \Lambda_\T^{-1} \iso \cal{A}_\T$.
Taking determinants would then imply $\Lambda_\T^{-1} \cong 	\det \cal{A}_\T = \cal{M}_\T^{-1}$, contrary to the choice of $\cal{M}_\T$.
Thus, this infinitesimal deformation of the spectral structure does not lift through $H$, proving that $H$ is not smooth at the chosen Higgs bundle.

The initial Higgs bundle may also be chosen to induce a clean spectral structure, with sufficiently positive twist, by \autoref{1792430412943}.
Every deformation over $\T$ then remains fibrewise clean because its only geometric fibre is the prescribed clean structure.
The same argument therefore shows that even the restriction of $H$ to the inverse image of the substack of clean spectral structures need not be smooth.
\end{example}

Let us conclude this subsection with two remarks.
First, for a fixed curve $\X$ and fixed twist $\Lambda$, the stack of regular Higgs bundles $\frak{Higgs}_{n,\Lambda}^\rm{reg} (\X)$ is nevertheless smooth.
Indeed, its Hitchin base is a vector space and the morphism $h$ is smooth.
This conclusion applies even at those regular Higgs bundles whose spectral curves are singular or nonreduced.

Second, the local existence theorem \autoref{1791544458956} allows the twist $\Lambda$ to be chosen to suit the given family of vector bundles.
It supplies local clean presentations induced by regular Higgs fields.
Meanwhile, the morphism $F$ and the stack of spectral structures are defined independently of any such choices.

\subsection{Dimension}
\label{1793136463825}

In this final subsection, we calculate the local stack dimension of the moduli stacks of spectral structures.
If $(\cal{E}, \cal{A}, \rho) \in \frak{SpStr}_n (\X)$ is a spectral structure, formulas for the local dimension at $(\cal{E}, \cal{A}, \rho)$ are most easily expressed in terms of the geometric invariants of the corresponding spectral realisation $(\pi, \cal{L}) \in \frak{SpReal}_n (\X)$.
Using the equivalence between spectral structures and their realisations (\autoref{1792502248471}), we allow ourselves to speak of the local stack dimension of $\frak{SpStr}_n (\X)$ at the point $(\pi, \cal{L})$:
\begin{eqn}
	\dim_{(\pi, \cal{L})} \frak{SpStr}_n (\X)
	\coleq \dim_{(\cal{E}, \cal{A}, \rho)} \frak{SpStr}_n (\X)
	= \dim_{(\pi, \cal{L})} \frak{SpReal}_n (\X)
\fullstop
\end{eqn}

\begin{theorem}[local dimension of the moduli stack]
\label{1792503824011}
Let $\X$ be a compact Riemann surface, and let $\cal{E}$ be a vector bundle on $\X$ of rank $n$ with spectral structure $(\cal{A}, \rho)$ and spectral realisation $(\pi : \Y \to \X, \cal{L})$.
Then the local dimension of the stack of spectral structures at the point $(\pi, \cal{L})$ is computable in terms of the local dimension of the stack of covers at $\pi$ and the Euler characteristic of $\Y$:
\begin{eqntag}
\label{1793105751171}
	\dim_{(\pi, \cal{L})} \frak{SpStr}_n (\X)
	= \dim_\pi \frak{Cov}_n (\X) - \chi (\cal{O}_\Y)
\fullstop
\end{eqntag}
In particular, the local dimension at $(\pi, \cal{L})$ depends only on the cover $\pi$ and not on the spectral line bundle $\cal{L}$.
More generally, for a family $p : \X \to \S$, this identity applies fibrewise over $\S$.
\end{theorem}

Explicitly, let $p : \X \to \S$ be a proper smooth family of connected compact curves.
Fix $s \in \S$ and consider a spectral realisation $(\pi_s : \Y_s \to \X_s, \cal{L}_s)$.
Base change identifies the fibres over $s$ of the two moduli stacks with $\frak{SpStr}_n (\X_s)$ and $\frak{Cov}_n (\X_s)$, respectively.
Therefore, \autoref*{1792503824011} in families means the following identity:
\begin{eqn}
	\dim_{(\pi_s, \cal{L}_s)} \frak{SpStr}_n (\X_s)
	= \dim_{\pi_s} \frak{Cov}_n (\X_s) - \chi (\cal{O}_{\Y_s})
\fullstop
\end{eqn}

\begin{proof}[Proof of \autoref{1792503824011}]
By \autoref{1792511507140}, the forgetful morphism $G : \frak{SpStr}_n (\X) \to \frak{Cov}_n (\X)$ is smooth.
By \autoref{1792512462793}, its relative stack dimension at $(\pi, \cal{L})$ is
\begin{eqn}
	h^1 (\Y, \cal{O}_\Y) - h^0 (\Y, \cal{O}_\Y) 	= - \chi (\cal{O}_\Y)
\fullstop
\end{eqn}
Therefore, the dimension formula for smooth morphisms, computed on smooth analytic atlases as in \cite[\href{https://stacks.math.columbia.edu/tag/0DRP}{Tag 0DRP}]{stacks-project}, gives precisely the identity \eqref{1793105751171}.
The fibrewise assertion follows by base change to each $s \in \S$.
\end{proof}

For smooth spectral structures, these quantities can be made much more explicit as follows.

\begin{corollary}[dimension at smooth spectral structures]
\label{1793102347934}
If, in addition, the spectral structure is smooth, let $r \coleq \deg \R$ be its total ramification degree.
Then
\begin{eqnstag}
\label{1792503826975}
	\dim_{(\pi, \cal{L})} \frak{SpStr}_n (\X)
	&= r - \chi (\cal{O}_\Y)
	= r + g_\Y - 1
\\
\label{1792503829376}
	&= \tfrac{3}{2} r + n (g_\X - 1)
	= 3 (g_\Y - 1) - 2n (g_\X - 1)
\fullstop
\end{eqnstag}
More generally, for a family $p : \X \to \S$, this identity applies fibrewise over $\S$ and computes the relative stack dimension over $\S$.
\end{corollary}

Note that $g_\Y = 1 - \chi (\cal{O}_\Y)$ with the convention of \eqref{1792494116976} when $\Y$ is disconnected.
Explicitly for a family $p : \X \to \S$, the corresponding statement concerns the structure morphism $\frak{SpStr}_n^\rm{sm} (\X / \S) \to \S$.
This morphism is smooth by \autoref{1792494388516}.
At a point represented by $(\pi_s : \Y_s \to \X_s, \cal{L}_s)$, write $r_s$ for the total ramification degree of $\pi_s$.
The relative stack dimension of the structure morphism at this point is the local dimension of its fibre:
\begin{eqn}
	\dim_{(\pi_s, \cal{L}_s)}
	\Big( \frak{SpStr}_n^\rm{sm} (\X / \S) \big/ \S \Big)
	=
	\dim_{(\pi_s, \cal{L}_s)} \frak{SpStr}_n (\X_s)
	=
	\tfrac{3}{2} r_s + n (g_{\X_s} - 1)
\fullstop
\end{eqn}

\begin{proof}[Proof of \autoref{1793102347934}]
Suppose now that $\Y$ is smooth.
Let $\cal{T}_\Y$ and $\cal{T}_\X$ denote the tangent sheaves.
Since $\pi$ is a finite map between smooth curves, its differential is injective and gives an exact sequence
\begin{eqn}
	0
	\too
	\cal{T}_\Y
	\overset{d\pi}{\too}
	\pi^\ast \cal{T}_\X
	\too
	\cal{N}_\pi
	\too
	0
\fullstop
\end{eqn}
The sheaf $\cal{N}_\pi$ is a torsion sheaf supported on the ramification points.

At a point of ramification index $k$, choose local coordinates in which the map $\pi$ is $w = z^k$.
The differential is then multiplication by $k z^{k-1}$, so the local length of $\cal{N}_\pi$ is $k-1$.
Summing over the ramification points gives $h^0 (\Y, \cal{N}_\pi) = r$ and $\H^1 (\Y, \cal{N}_\pi) = 0$.

Deformation theory of holomorphic maps with fixed target identifies first-order deformations of the cover with $\H^0 (\Y, \cal{N}_\pi)$; see \cite[\S3.1]{Naumann2018Hurwitz}.
There are no infinitesimal automorphisms of the cover because $d\pi$ is injective.
Moreover, the cover-lifting argument in the proof of \autoref{1792494388516} solves every infinitesimal lifting problem for $\pi$.
Since $\frak{Cov}_n (\X)$ is locally of finite presentation, it is therefore smooth at $\pi$.
So its local dimension is
\begin{eqn}
	\dim_\pi \frak{Cov}_n (\X)
	=
	h^0 (\Y, \cal{N}_\pi)
	=
	r
\fullstop
\end{eqn}
Substituting this into the general identity \eqref{1793105751171} proves \eqref{1792503826975}.

Finally, Riemann-Hurwitz, summed over the connected components of $\Y$ if necessary, gives
\begin{eqn}
	g_\Y - 1
	=
	n (g_\X - 1) + \tfrac{1}{2} r
\fullstop
\end{eqn}
Substituting this relation into $r + g_\Y - 1$ gives both expressions in \eqref{1792503829376}.

For a family $p : \X \to \S$, applying this calculation to each fibre $\X_s$ proves the fibrewise identity.
Since $\frak{SpStr}_n^\rm{sm} (\X / \S) \to \S$ is smooth by \autoref{1792494388516}, the same formula gives its relative stack dimension.
\end{proof}

\begin{corollary}[relative dimension of $F$ on the clean substack]
\label{1793103207826}
In addition, if the spectral structure is clean, then the relative dimension of the morphism $F : \frak{SpStr}_n^\rm{clean} (\X) \to \frak{Bun}_n (\X)$ at the point $(\pi, \cal{L})$ is
\begin{eqntag}
\label{1793046274641}
	r + h^0 (\X, \cal{Q})
	=
	r + \chi (\cal{Q})
	=
	\tfrac{3}{2} r - n(n-1) (g_\X - 1)
\fullstop
\end{eqntag}
More generally, for a family $p : \X \to \S$, this formula applies fibrewise over $\S$ and computes the relative dimension of the forgetful morphism $F : \frak{SpStr}_n^\rm{clean} (\X / \S) \to \frak{Bun}_n (\X / \S)$.
\end{corollary}

Explicitly for a family $p : \X \to \S$, the relative dimension of the smooth representable morphism $F$ is calculated as follows.
At a clean spectral realisation $(\pi_s : \Y_s \to \X_s, \cal{L}_s)$, put $\cal{E}_s \coleq (\pi_s)_\ast \cal{L}_s$, let $\cal{Q}_s$ be its off-diagonal bundle, and write $r_s$ for its total ramification degree.
The essential fibre of $F$ over $\cal{E}_s$ is the analytic space $\frak{SpStr} (\cal{E}_s)$ whose local dimension at this point is
\begin{eqn}
	\dim_{(\pi_s, \cal{L}_s)} \frak{SpStr} (\cal{E}_s)
	=
	r_s + h^0 (\X_s, \cal{Q}_s)
	=
	\tfrac{3}{2} r_s - n(n-1) (g_{\X_s} - 1)
\fullstop
\end{eqn}
Here, the identification of the underlying vector bundle with $\cal{E}_s$ is part of the fibre data.
Thus, this relative dimension measures variations of the spectral presentation with the base curve and the underlying vector bundle fixed.

\begin{proof}[Proof of \autoref{1793103207826}]
By \autoref{1792494388516}, $F$ is smooth and representable by analytic spaces.
Its relative dimension is therefore the difference between the local stack dimensions of its source and target.

The stack of vector bundles on $\X$ is smooth \cite[Proposition~1]{Heinloth2010Uniformization}, with its analytic description supplied by \cite[Corollary~6.16(1) and Remark~6.17(2)]{HolsteinPorta2025Analytification}, and Riemann-Roch gives its dimension at $\cal{E}$ as
\begin{eqn}
	h^1 \big( \X, \cal{End} (\cal{E}) \big)
	-
	h^0 \big( \X, \cal{End} (\cal{E}) \big)
	=
	-\chi \big( \cal{End} (\cal{E}) \big)
	=
	n^2 (g_\X - 1)
\fullstop
\end{eqn}
Since the substack of clean spectral structures is open, formula \eqref{1792503826975} gives its dimension at $(\pi, \cal{L})$ as $r + g_\Y - 1$.
The required relative dimension is consequently
\begin{eqn}
	r + g_\Y - 1 - n^2 (g_\X - 1)
	=
	r + \chi (\cal{Q})
\fullstop{,}
\end{eqn}
where the equality follows from the off-diagonal exact sequence.
The condition of being clean gives $\H^1 (\X, \cal{Q}) = 0$, so $\chi (\cal{Q}) = h^0 (\X, \cal{Q})$.
Finally, Riemann-Hurwitz gives
\begin{eqn}
	\chi (\cal{Q})
	=
	g_\Y - 1 - n^2 (g_\X - 1)
	=
	\tfrac{1}{2} r - n(n-1) (g_\X - 1)
\fullstop
\end{eqn}
These identities prove \eqref{1793046274641}.

For a family $p : \X \to \S$, base change to $s \in \S$ identifies $F$ with the corresponding forgetful morphism over $\X_s$.
Since the relative dimension of a smooth morphism is preserved by base change, the same calculation proves the family assertion.
\end{proof}

We remark that the dimension formulas above remain valid after fixing the degree of the underlying vector bundle.
Indeed, for a family $p : \X \to \S$, the fibrewise degree is locally constant, so prescribing this degree selects open-and-closed substacks of both $\frak{SpStr}_n (\X / \S)$ and $\frak{Bun}_n (\X / \S)$.
Restriction to these substacks changes neither the local stack dimensions nor the relative dimension of $F$ at their points.

\begin{remark}[total ramification and ramification profiles]
\label{1792517423166}
The dimension formula \eqref{1792503826975} in the case of smooth spectral structures (\autoref{1793102347934}) does not require simple branching.
There, $r$ is the total ramification degree, and a point of ramification index $k$ contributes $k-1$.
The deformations counted by $\H^0 (\Y, \cal{N}_\pi)$ allow such a point to split into several ramification points.
For example, the local deformation $w = z^k + \epsilon z$ splits a point of ramification index $k$ into $k-1$ simple ramification points for small $\epsilon \neq 0$.
Thus, a locus on which the ramification profile is prescribed can have smaller dimension than the ambient substack of smooth spectral structures.
\end{remark}

\begin{remark}[fibres over a fixed bundle]
\label{1793116276492}
The distinction between the total stack and its fibres over a fixed vector bundle $\cal{E}$ is important.
The dimension formula \eqref{1792503826975} computes the dimension of $\frak{SpStr}_n (\X)$ at a smooth spectral structure.
On the other hand, the analytic space $\frak{SpStr} (\cal{E})$ in \eqref{1793046205260} parametrises spectral structures \textit{together with} an identification of their underlying vector bundle with the prescribed bundle $\cal{E}$.
At a clean spectral structure, its dimension is the relative dimension in \eqref{1793046274641}.
\end{remark}

\enlargethispage{20pt}

\begin{remark}[unbounded dimensions]
\label{1792517441679}
For $n \geq 2$, the dimensions of the moduli stack of smooth and connected spectral structures are unbounded as the ramification degree varies, even after fixing the degree of the underlying vector bundle.

Indeed, fix an integer $d$.
Choose a line bundle $\cal{M}$ of sufficiently large degree $k$ on $\X$ and a global section $s \in \H^0 \big( \X, \cal{M}^{\otimes n} \big)$ with simple zeros.
Then the cyclic cover $\pi : \Y \to \X$ given locally by $w^n = s$ is smooth.
Over each zero of $s$, there is one point of ramification index $n$.
This also forces $\Y$ to be connected, because the connected component containing such a ramification point already has degree $n$ over $\X$.
Since $s$ has $nk$ zeros in total, its total ramification degree is $r = n(n-1)k$.

Choose a line bundle $\cal{L}$ on $\Y$ of degree $d + r/2$.
Then Riemann-Roch gives
\begin{eqn}
	\deg (\pi_\ast \cal{L})
	= \deg \cal{L} - r/2 = d
\fullstop
\end{eqn}
So the degree of the vector bundle $\cal{E} = \pi_\ast \cal{L}$ is fixed at $d$.
Nevertheless, dimension formula \eqref{1792503829376} gives the local stack dimension
\begin{eqn}
	\dim_{(\pi, \cal{L})} \frak{SpStr}_n (\X)
	=
	\tfrac{3}{2} n(n-1)k + n(g_\X - 1)
\fullstop{,}
\end{eqn}
and this tends to infinity as $k$ increases.
Thus, local finite presentation is compatible with unbounded local dimensions, and fixing $\deg \cal{E}$ alone does not give a uniform dimension bound.
\end{remark}

\clearpage
\enlargethispage{40pt}
\mbox{}
\vspace{-50pt}

\bibliographystyle{src/declare-bibliography}
\bibliography{src/bibliography}

\end{document}